\PassOptionsToPackage{unicode}{hyperref}
\PassOptionsToPackage{hyphens}{url}
\PassOptionsToPackage{dvipsnames,svgnames,x11names}{xcolor}
\documentclass[
  12pt]{article}
\usepackage{bbm}
\usepackage{enumerate}
\usepackage{xr}
\usepackage{amsmath,amssymb}
\usepackage[T1]{fontenc}
\usepackage[utf8]{inputenc}
\usepackage{textcomp}
\usepackage{lmodern}
\usepackage{placeins}
\usepackage{amsthm}

\IfFileExists{upquote.sty}{\usepackage{upquote}}{}
\IfFileExists{microtype.sty}{
  \usepackage[]{microtype}
  \UseMicrotypeSet[protrusion]{basicmath} 
}{}
\makeatletter
\@ifundefined{KOMAClassName}{
  \IfFileExists{parskip.sty}{%
    \usepackage{parskip}
  }{
    \setlength{\parindent}{0pt}
    \setlength{\parskip}{6pt plus 2pt minus 1pt}}
}{
  \KOMAoptions{parskip=half}}
\makeatother
\usepackage{xcolor}
\makeatletter
\ifx\paragraph\undefined\else
  \let\oldparagraph\paragraph
  \renewcommand{\paragraph}{
    \@ifstar
      \xxxParagraphStar
      \xxxParagraphNoStar
  }
  \newcommand{\xxxParagraphStar}[1]{\oldparagraph*{#1}\mbox{}}
  \newcommand{\xxxParagraphNoStar}[1]{\oldparagraph{#1}\mbox{}}
\fi
\ifx\subparagraph\undefined\else
  \let\oldsubparagraph\subparagraph
  \renewcommand{\subparagraph}{
    \@ifstar
      \xxxSubParagraphStar
      \xxxSubParagraphNoStar
  }
  \newcommand{\xxxSubParagraphStar}[1]{\oldsubparagraph*{#1}\mbox{}}
  \newcommand{\xxxSubParagraphNoStar}[1]{\oldsubparagraph{#1}\mbox{}}
\fi
\makeatother

\usepackage{longtable,booktabs,array}
\usepackage{calc} 
\usepackage{etoolbox}
\makeatletter
\patchcmd\longtable{\par}{\if@noskipsec\mbox{}\fi\par}{}{}
\makeatother
\IfFileExists{footnotehyper.sty}{\usepackage{footnotehyper}}{\usepackage{footnote}}
\makesavenoteenv{longtable}
\usepackage{graphicx}
\makeatletter
\def\maxwidth{\ifdim\Gin@nat@width>\linewidth\linewidth\else\Gin@nat@width\fi}
\def\maxheight{\ifdim\Gin@nat@height>\textheight\textheight\else\Gin@nat@height\fi}
\makeatother
\setkeys{Gin}{width=\maxwidth,height=\maxheight,keepaspectratio}
\makeatletter
\def\fps@figure{htbp}
\makeatother

\makeatletter
\@ifpackageloaded{caption}{}{\usepackage{caption}}
\AtBeginDocument{%
\ifdefined\contentsname
  \renewcommand*\contentsname{Table of contents}
\else
  \newcommand\contentsname{Table of contents}
\fi
\ifdefined\listfigurename
  \renewcommand*\listfigurename{List of Figures}
\else
  \newcommand\listfigurename{List of Figures}
\fi
\ifdefined\listtablename
  \renewcommand*\listtablename{List of Tables}
\else
  \newcommand\listtablename{List of Tables}
\fi
\ifdefined\figurename
  \renewcommand*\figurename{Figure}
\else
  \newcommand\figurename{Figure}
\fi
\ifdefined\tablename
  \renewcommand*\tablename{Table}
\else
  \newcommand\tablename{Table}
\fi
}
\@ifpackageloaded{float}{}{\usepackage{float}}
\floatstyle{ruled}
\@ifundefined{c@chapter}{\newfloat{codelisting}{h}{lop}}{\newfloat{codelisting}{h}{lop}[chapter]}
\floatname{codelisting}{Listing}

\makeatother
\makeatletter
\@ifpackageloaded{caption}{}{\usepackage{caption}}
\@ifpackageloaded{subcaption}{}{\usepackage{subcaption}}
\makeatother

\usepackage[]{natbib}
\usepackage{bookmark}

\IfFileExists{xurl.sty}{\usepackage{xurl}}{} 
\hypersetup{
  pdftitle={Title},
  pdfauthor={Author 1; Author 2},
  pdfkeywords={3 to 6 keywords, that do not appear in the title},
  colorlinks=true,
  linkcolor={blue},
  filecolor={Maroon},
  citecolor={Blue},
  urlcolor={Blue},
  pdfcreator={LaTeX via pandoc}}

\newcommand{\anon}{1}

\DeclareMathOperator{\Var}{Var}

\DeclareMathOperator{\Cov}{Cov}

\newtheorem{theorem}{Theorem}[section]
\newtheorem{remark}[theorem]{Remark}

\newtheorem{proposition}[theorem]{Proposition}
\newtheorem{lemma}[theorem]{Lemma}
\newtheorem{corollary}[theorem]{Corollary}

\newtheorem{example}{Example}[section]

\begin{document}

\def\spacingset#1{\renewcommand{\baselinestretch}%
{#1}\small\normalsize} \spacingset{1}


\if1\anon
{
  \title{\bf Generalized Taylor's Law for Dependent and Heterogeneous Heavy-Tailed Data}
 \author{
  Pok Him Cheng\\
  Department of Statistics, Columbia University
  \and
  Joel E. Cohen\thanks{Joel E. Cohen thanks Roseanne Benjamin for her help during this work.}\\
  Laboratory of Populations, The Rockefeller University\\
  Department of Statistics, Columbia University\\
  Climate School, Columbia University\\
  Department of Statistics, University of Chicago
  \and
  Hok Kan Ling\thanks{Hok Kan Ling gratefully acknowledges the support by NSERC Grant RGPIN/03124-2021.}\\
  Department of Mathematics and Statistics, Queen’s University
  \and
  Sheung Chi Phillip Yam\thanks{Phillip Yam acknowledges the financial supports from HKGRF-14301321 (“General Theory $\ldots$ Classical Problems”), HKGRF-14300123 (“Well-posedness $\ldots$ Applications”), and G-CUHK411/23. He also thanks The University of Texas at Dallas for the kind invitation to be a Visiting Professor in Naveen Jindal School of Management.}\\
  Department of Statistics, The Chinese University of Hong Kong
}

  \maketitle
} \fi

\if0\anon
{
  \bigskip
  \bigskip
  \bigskip
  \begin{center}
    {\LARGE\bf Generalized Taylor's Law for Dependent and Heterogeneous Heavy-Tailed Data}
\end{center}
  \medskip
} \fi

\begin{abstract}
  Taylor's law, also known as fluctuation scaling in physics and the power-law variance function in statistics, is an empirical pattern widely observed across fields including ecology, physics, finance, and epidemiology. It states that the variance of a sample scales as a power function of the mean of the sample. We study generalizations of Taylor's law in the context of heavy-tailed distributions with infinite mean and variance. We establish the probabilistic limit and analyze the associated convergence rates. Our results extend the existing literature by relaxing the i.i.d. assumption to accommodate dependence and heterogeneity among the random variables. This generalization enables application to dependent data such as time series and network-structured data. We support the theoretical developments  by extensive simulations, and the practical relevance through applications to real network data.
  \end{abstract}

\noindent%
\vfill

\newpage
\spacingset{1.8} 

\section{Introduction}
The study of heavy-tailed distributions has gained significant attention in various fields of research, including finance,  economics, operational risk management, network analysis, disaster analysis, and extreme event modeling. 
For example, \cite{rachev2003handbook} provided a comprehensive introduction to heavy-tailed distributions in finance; 
\cite{glasserman2002portfolio} explored methods for computing portfolio Value-at-Risk when underlying risk factors have heavy-tailed distributions; 
\cite{peters2015advances}  modeled heavy-tailed loss processes in operational risk; 
\cite{merz2022understanding} discussed the consequences of the heavy-tailed behavior of flood peak distributions; 
\cite{eom2015tail} studied estimation methods for the tail degree in large-scale social networks.  
When the heavy-tailed distribution satisfies $\mathbb{P}(X>x)$ = $x^{-\alpha} L(x)$ for $x \geq 0$ with tail index $\alpha \in (0,1)$ (so that the mean, variance, and all higher moments are infinite) and $L$ is slowly varying 
at infinity (grows more slowly than any power function of $x$ as $x\to\infty$), \cite{nevslehova2006infinite} considered operational risk management and \cite{pisarenko2010heavy} used  heavy-tailed distributions to  model disasters. These distributions  
exhibit rare, extreme events that occur more frequently than in normal distributions; 
they have found practical applications in real life.

Taylor's law \citep{taylor1961aggregation}, also known as fluctuation scaling in physics,
and as a power-law variance function in statistics, is an empirical pattern observed in various fields such as ecology, physics, finance, and epidemiology. 
Taylor's law states that the variance of a sample is proportional to a power of the sample mean.
For a random variable $X$ with mean $\mu_X >0$ and variance $\Var_X >0$, Taylor's law postulates that, 
as $X$ ranges over a set of two or more random variables considered to be somehow comparable,
\begin{align}\label{eq:TLpower}
    \Var_X=a\mu_X^b,
\end{align}
where $a>0$ and $b$ are real constants, the same for all random variables in the set. 
Most applications have $b>0$. 
In its logarithmic form,  Taylor's law becomes
\begin{align}\label{eq:TL}
    \log(\Var_X)=\log(a)+b\log(\mu_X)
\end{align}
where the slope $b$ in \eqref{eq:TL} is the exponent $b$ in the power-law variance-mean relationship \eqref{eq:TLpower}. 
Equivalently, 
\begin{equation*}
    \frac{\log (\Var_X)}{\log (\mu_X)} = \frac{\log (a)}{\log (\mu_X) } + b.
\end{equation*}
Replacing the population variance $\Var_X$ by the sample variance
and replacing the population mean $\mu_X$ by the sample mean
yields a sample version of Taylor's law even when the population quantities 
$\Var_X,\ \mu_X$ may be infinite.

\cite{taylor2019taylor} reviewed many applications of Taylor's law.
\cite{eisler2008fluctuation} discussed the general applicability of the scaling relationship of Taylor's law, reviewed literature, and presented new empirical data and model calculations; \cite{hanley2014fluctuation} showed that monthly crime reports in the UK followed Taylor’s law relationships approximately; \cite{giometto2015sample} extended Taylor's law to higher moments and gave conditions on the ratio of sample size per random variable and the number of random variables under estimates of the slope reflect underlying dynamics or are artifactual; \cite{cohen2013taylor} showed that  human demographic data from Norway follow Taylor's law and proposed a simple model of some features of the results.
These references focused  on light-tailed data where the population mean and variance exist. 
\cite{de2022dynamic} introduced a dynamic Taylor’s law for dependent samples using self-normalized expressions involving Bernstein blocks. They discussed variations of Taylor’s law under dependence, considering ergodic behaviors and stationary time series of  light-tailed random variables under weak dependence and strong mixing assumptions of a strictly stationary time series.

We now  review  Taylor's law for models and data following heavy-tailed 
distributions with infinite mean and therefore also infinite variance. 
To the best of our knowledge, \cite{brown2017taylor} were the first to  extend Taylor's law to distributions with infinite mean and variance. Specifically, they showed  that a modified version of Taylor's law remains valid for $\alpha$-stable distributions with tail index $\alpha \in (0, 1)$.   
\cite{cohen2020heavy} showed that heavy-tailed distributions can overwhelm the dependence among random variables; that under certain conditions, the sample correlation of regularly varying random variables converges to a random limit; 
and that the sample kurtosis  converges in distribution. 
They discussed Taylor's law for  heterogeneous and correlated heavy-tailed data, instead of focusing solely on 
independently and identically distributed (henceforth $i.i.d.$) data. 
\cite{brown2021taylor} extended the work of \cite{brown2017taylor} by demonstrating the convergence of the ratios to the sample mean of sample (local) upper and lower semivariances, higher sample moments, skewness, and kurtosis. 
They applied  their results  to various financial ratios. 

\cite{cohen2022covid} considered Taylor's law when the tail index satisfies $\alpha\in (1, 2)$, 
so that the random variable has finite mean and infinite variance. 
They argued that this setting for heavy-tailed random variables is more relevant to some real data. 
They were also the first to consider an array setting of heavy-tailed random variables and investigated the patterns of COVID-19 cases and deaths in the United States using Taylor's law. 

In this article, we focus on  nonnegative random variables with survival function $\overline{F}(x) := 1- F(x) = x^{-\alpha}l(x)$ where 
$x \geq 0,\ \alpha > 0$ and
$l(\cdot)$ is slowly varying at infinity, that is, 
$\lim_{x \to \infty}{l(\lambda x)}/{l(x)}=1$ for all $\lambda>0$.
Examples of slowly varying functions include $\log(x)$, $(\log(x))^p$ for $p \in \mathbb{R}$,  $\log(\log(x))$, and the constant function.  
\cite{seneta1976regularly}, \cite{bingham1989regular}, and \cite{foss2013introduction} give more  details about slowly varying functions.

We shall concentrate on establishing the  limit in probability $\stackrel{\mathbb{P}}{\rightarrow}$, 
rather than the weak distributional limit, of various generalizations of Taylor's law. Inspired by \cite{kolmogoroff1933grundbegriffe}, we make use of the celebrated Karamata theorem \cite[page 26, Proposition 1.5.8]{bingham1989regular} to establish the moments of the truncated heavy-tailed random variables and to approximate the heavy-tailed random variable with infinite moments. This approach allows us to impose fewer conditions than presently  found in the literature and accommodates some dependence among the random variables considered. In particular, we examine higher moments, higher central moments, and local semivariances. We also show that the covariance condition we obtained is satisfied by some weak dependence in  a time series. 
Like \cite{de2022dynamic}, who  worked on light-tailed random variables, we considered different ``mixing'' conditions of a time-series process. 
In addition to  $i.i.d.$ data
and time series, we also consider heterogeneous and correlated data. 
Since our primary condition for establishing these laws is based on the covariances of truncated random variables, our results apply to a collection of random variables, such as those on a network graph; indeed, to the best of our knowledge, this is the first article to discuss Taylor's law for a network graph structure.
Our results might be applicable in many real-life situations such as social network data.

The rest of this article is organized as follows: Section \ref{sec:main results} presents our main results on Taylor's law under weak dependence. {
We extend the results from higher (centered) moments to semivariances in Section A in the appendix.} Sections \ref{sec:heterogeneous} and \ref{sec:corr} introduce Taylor's law for heterogeneous and correlated heavy-tailed data, respectively. Section \ref{sec:network TL} generalizes Taylor's law to network data.
Section \ref{sec:application} illustrates Taylor's law for heavy-tailed network data using three real-world datasets.
Section \ref{sec:discussion} discusses and summarizes this article and proposes some possible future work. 
In the appendix, we provide simulation results illustrating the convergence discussed in the article, along with proofs of the mathematical claims in the main text.

\section{Main Results: Taylor's Law for Higher Moments with Weak Dependence}\label{sec:main results}
Unless otherwise specified, random variables $X_1,\ldots,X_n$ are nonnegative and have a common distribution $F$ with  
survival function $\overline{F}(x) := 1- F(x) = x^{-\alpha}l(x)$, where $x \geq 0,\ \alpha >0$ and $l(\cdot)$ is slowly varying at infinity. 
We do not require $X_1,\ldots,X_n$ to be independent. 
It is well-known that $\overline{F}(x)= x^{-\alpha}l(x)$ for some slowly varying function $l$ if and only if $F$ belongs to the \textit{maximum domain of attraction} of the Fr\'echet distribution, 
which has cumulative distribution function
$\Phi_\alpha (x)= \exp\{-x^{-\alpha}\}$ for $x>0$ and $\alpha>0$  \citep[Theorem 3.3.7]{embrechts2013modelling}. 
This class of distributions contains, for example, the Pareto and the $\alpha$-stable distributions. 
When $\alpha \in (0, 1)$,  the mean and variance of $F$ are infinite  \citep{nair2022fundamentals, borak2005stable}.

Let $\{a_n\}$ and $\{b_n\}$ be two real sequences. 
We write $a_n \sim b_n$ if ${a_n}/{b_n} \to 1$ as $n \to \infty$. 
We let $\sum_{i \neq j}$ denote  $\sum_{i,j=1,\ldots,n, i \neq j}$.

\subsection{Auxiliary Results}\label{subsect:auxiliary}

Two auxiliary results, Lemmas \ref{lemma:truncate_moment_asy} and \ref{lemma:existence_of_tn},  precede the main Theorems \ref{thm:moment} and \ref{thm:higher_central_moments}. 
Since the expectation of the random variables we study here 
is infinite when $\alpha \in (0, 1)$, Lemma \ref{lemma:truncate_moment_asy} considers the growth rate of the truncated moments to infinity, which we  establish using Karamata's theorem \cite[page 26, Proposition 1.5.8]{bingham1989regular}.  
\begin{lemma}[Truncated Moments]\label{lemma:truncate_moment_asy}
Suppose that $X\geq 0$ has a survival function $\overline{F}(x) = x^{-\alpha}l(x)$ for $\alpha > 0$. Let $\tilde{X} := X\mathbbm{1}(X < t_n)$, where $t_n \uparrow \infty$ as $n \rightarrow \infty$ ($t_n$ does not necessarily have to satisfy Equation (\ref{eq:choice_of_tn}) below). Then, for $p > \alpha$, as $n \rightarrow \infty$, we have
\begin{enumerate}[(a)]
    \item \begin{equation*}
    \mathbb{E}[\tilde{X}^p]  \sim \frac{\alpha }{p-\alpha} t_n^{p-\alpha}l(t_n);
\end{equation*}
\item \begin{equation*}
    \frac{\log \mathbb{E}[\tilde{X}^p]}{\log t_n} \rightarrow p- \alpha.
\end{equation*}
\end{enumerate}
\end{lemma}
The next lemma guarantees the existence of an appropriate truncation level $t_n$ 
to  use in establishing the main Theorems \ref{thm:moment} and \ref{thm:higher_central_moments}.

\begin{lemma}\label{lemma:existence_of_tn}
Let $\{c_n\}$ be a positive sequence such that $n c_n \uparrow \infty$ and $\log c_n = o(\log n)$. 
(Such a positive sequence exists, for example, $c_n:=\log(n)$ for every positive integer $n$.)
For any continuous slowly varying function $l_0 > 0$ with $\overline{F}_0(0+):=\lim_{x \downarrow 0} x^{-\alpha}l_0(x)>0$, there exists a sequence $\{t_n\}$ such that, for all sufficiently large integer $n$,
    \begin{equation}\label{eq:choice_of_tn}
        \frac{n l_0(t_n)}{t_n^\alpha} = \frac{1}{c_n}.
    \end{equation}
\end{lemma}
Hereafter we take $\{c_n\}$ to be a positive sequence such that $nc_n \uparrow \infty$ and $\log c_n = o(\log n)$ and
$\{t_n\}$ to be a sequence that satisfies (\ref{eq:choice_of_tn}).

Although the function $l(\cdot)$ in $\overline{F}(x) = x^{-\alpha}l(x)$ is not necessarily continuous due to possible
jumps of the survival function $\overline{F}$, by Proposition 1.3.4 in \cite{bingham1989regular}, there is 
a continuous slowly varying function $l_0(\cdot)$ such that $l(x) \sim l_0(x)$ and $\overline{F}_0(0+)>0$ as $x \rightarrow \infty$. 
By Lemma \ref{lemma:existence_of_tn}, there exists $\{t_n\}$ such that, for all large $n$,
\eqref{eq:choice_of_tn} holds.
Then, as $n \rightarrow \infty$,
\begin{equation}\label{eq:t_n_to_0}
    \frac{n l(t_n)}{t_n^\alpha} = \frac{n l_0(t_n)}{t_n^\alpha} \cdot \frac{l(t_n)}{l_0(t_n)} \sim \frac{1}{c_n} \rightarrow 0.
\end{equation}
Following the same argument as for \eqref{eq:choice_of_tn} and \eqref{eq:t_n_to_0}, there exists a positive sequence $\{v_n\}$ that satisfies 
\begin{equation}\label{eq:v_n_choice1}
    \frac{n l(v_n)}{v^\alpha_n} \sim c_n,
\end{equation}
and we henceforth take $\{v_n\}$ to be a positive sequence satisfying \eqref{eq:v_n_choice1}.

\subsection{Taylor's Law for Higher Moments}

Taylor’s law holds empirically when the sample variance scales approximately in direct proportion to a power of the sample mean.
In addition, sometimes  the scaling relationship between the sample mean and other higher order moments or higher centered moments is of interest. 
In particular, we may adopt higher (centered) moments in place of the variance or the mean in equation (\ref{eq:TL}). 
\cite{giometto2015sample}  applied the generalized Taylor's law for higher sample moments to tree abundance across plots in the Black Rock Forest and abundances of
carabid beetles in different sites in the Netherlands.
\cite{giometto2015sample} and \cite{brown2017taylor} give  more examples and discussion.

To establish the probability limits in our main theorems, we impose a  weak dependence condition on the truncated random variable $\breve{X}_i := X_i \mathbbm{1}(X_i < v_n)$.

\noindent
\textbf{Condition A($p$):} 
Assume the truncated random variables ${\breve{X}_i}$ satisfy
\begin{equation}\label{eq:cov_ass2}
    \sum_{i\neq j}\Cov(\breve{X}^p_i, \breve{X}^p_j) = o(v^{2p}_n c_n^2).
\end{equation}

\subsubsection{Some Sufficient Conditions for Condition A\texorpdfstring{($p$)}{(p)}}

\label{subsec:cov condition}
Condition A($p$) is relatively mild  since  it is satisfied under many common  dependence assumptions.
Here we provide some examples and sufficient conditions for Condition A($p$). 
First, if $X_1,\ldots,X_n$ are independent,  Condition A($p$) is satisfied trivially. 
Next, we consider an $m$-dependent random sequence. 
A random sequence is said to be $m$-dependent if {
there exists a nonnegative integer $m$ such that for every positive integer $k$},
$\{X_i\}_{i\leq k}$ is independent of $\{X_i\}_{i\geq k+m+1}$;
trivially, if $m=0$, $\{X_i\}_{i=1}^n$ is an independent sequence.
\begin{theorem}\label{thm:m-dependence}
Let $X_1,\ldots,X_n,\ldots$ be a sequence of 
nonnegative random variables with a common survival function $\overline{F}(x) = x^{-\alpha}l(x)$, where $l$ is slowly varying and $\alpha > 0$ (not necessarily in $(0, 1)$).
If $\{X_i,i\geq1\}$ is a $m$-dependent sequence, then Condition A($p$) holds for any $p > \alpha$.
\end{theorem}

Now we show that Condition A($p$) holds under some commonly used mixing conditions.  Let $Y:=(Y_k,k\in \mathbb{Z})$ be a sequence of random variables on a given probability space  $(\Omega, \mathcal{F}, \mathbb{P})$. For any $\sigma$-field $\mathcal{A} \subset \mathcal{F}$, let $\mathcal{F}^l_j := \sigma(Y_k, j \leq k \leq l)$ be the $\sigma$-field generated by $Y_j,\ldots,Y_l$. Let $\mathcal{L}^2_{\text{Real}}(\mathcal{A})$ denote the space of square-integrable, $\mathcal{A}$-measurable (real-valued) random variables. 
For any two $\sigma$-fields $\mathcal{A}$ and $\mathcal{B} \subset \mathcal{F}$, define the following measures of dependence as
\begin{align*}
    \alpha(\mathcal{A},\mathcal{B}) &:=\sup_{A \in \mathcal{A}, B \in \mathcal{B}} |\mathbb{P}(A\cap B)-\mathbb{P}(A)\mathbb{P}(B)|,\\
    \phi(\mathcal{A},\mathcal{B}) &:= \sup_{A \in \mathcal{A}, B \in \mathcal{B}, \mathbb{P}(A) > 0} \left| P(B|A) - P(B) \right|,\\
    \psi(\mathcal{A},\mathcal{B}) &:= \sup_{A \in \mathcal{A}, B \in \mathcal{B}, \mathbb{P}(A) > 0, \mathbb{P}(B) > 0} \left|\frac{\mathbb{P}(A \cap B)}{\mathbb{P}(A)\mathbb{P}(B)} - 1 \right|,\\
    \rho(\mathcal{A},\mathcal{B}) &:= \sup_{f \in \mathcal{L}^2_{\text{Real}}(\mathcal{A}), g \in \mathcal{L}^2_{\text{Real}}(\mathcal{B})} \left| \text{Corr}(f, g) \right|,\\
    \beta(\mathcal{A},\mathcal{B}) &:= \sup \frac{1}{2} \sum_{i=1}^{I} \sum_{j=1}^{J} \left| \mathbb{P}(A_i \cap B_j) - \mathbb{P}(A_i)\mathbb{P}(B_j) \right|,
\end{align*}
where the last supremum is taken over all pairs of (finite) partitions $\{A_1, \ldots, A_I \}$
and $\{B_1, \ldots, B_J \}$ of $\Omega$ such that $A_i \in \mathcal{A}$ for each $i$ and $B_j \in \mathcal{B}$ for each $j$. For any positive integer $n$, the strong mixing coefficient is defined as
\begin{align*}
    \alpha(n)= \sup_{j \in \mathbb{Z}} \alpha(\mathcal{F}^j_{-\infty},\mathcal{F}^{\infty}_{j+n});
\end{align*}
other types of mixing coefficients can be defined similarly by replacing $\alpha$ with their respective Greek letters.

The random sequence $Y$ is then said to be strong mixing or $\alpha$-mixing if $\alpha(n) \rightarrow 0$ as $n \rightarrow \infty$, and similarly, $\phi$-mixing if $\phi(n) \rightarrow 0$ as $n \rightarrow \infty$. 
The same applies to all the other types of mixing  above. 
These  mixing conditions are commonly considered 
\citep{zhengyan1997limit, merlevede2000functional,quintela2008hazard}. \cite{de2022dynamic}, who considered Taylor's law for a time series of 
light-tailed random variables with finite moments, also introduced strong mixing conditions and the corresponding mixing rates.

\begin{theorem}\label{thm:strong_mixing}
Suppose that a strictly stationary process $\{X_i,i\geq 1\}$ is a strong-mixing sequence of random variables with strong mixing coefficient $\alpha(k)$ such that\\ $\lim_{n \rightarrow \infty} \frac{1}{c_n^2}\sum^n_{k=1}\alpha(k) = 0 $. Then Condition A($p$) holds for any $p > 0$.
\end{theorem}
The strong mixing condition is among the weakest of all commonly used mixing conditions (e.g. compared to $\phi$-mixing) 
in the sense that if $\phi(n) \to 0$, then  $ \alpha(n) \to 0$.
\cite{Bradley2005basic}  reviews these relationships.  
Thus  if a  sequence satisfies another mixing condition  that has been introduced in this section, then Condition A($p$) holds.

\begin{corollary}\label{corollary:mixing}
Suppose that a strictly stationary process $\{X_i,i\geq 1\}$ is any one of  
\begin{enumerate}[(i)]
    \item $\psi$-mixing with $\lim_{n \to \infty} \frac{1}{c_n^2}\sum^n_{k=1} \psi(k) = 0$;
    \item $\phi$-mixing with $\lim_{n \to \infty} \frac{1}{c_n^2}\sum^n_{k=1} \phi(k) = 0$;
    \item $\rho$-mixing with $\lim_{n \to \infty} \frac{1}{c_n^2}\sum^n_{k=1} \rho(k) = 0$;
    \item $\beta$-mixing with $\lim_{n \to \infty} \frac{1}{c_n^2}\sum^n_{k=1} \beta(k) = 0$.
\end{enumerate}
Then Condition A($p$) holds for any $p > 0$.
\end{corollary}

We next demonstrate that a first-order autoregressive model known as AR(1), satisfies Condition A($p$) for any $p > 0$ under certain regularity conditions.  
\begin{theorem}\label{thm:AR_cov}
Consider a strictly stationary AR(1) process $X_t = \beta_1 X_{t-1} + \epsilon_t$, $t \in \mathbb{Z}$, where $\beta_1  \in (0, 1)$,  and $\{\epsilon_t\}$ are i.i.d. nonnegative random variables.
Suppose that $\mathbb{E}(\epsilon_1^q) < \infty$ for some $q > 0$. Furthermore, assume that $\epsilon_1$ has a density that is bounded above with respect to the Lebesgue measure that is non-strictly monotone (either non-decreasing or non-increasing) on each set in a finite partition of $\mathbb{R}$ into intervals, which may be bounded or unbounded. 
Then $\{X_t, t \geq 1\}$ satisfies Condition A($p$) 
for any $p > 0$.
\end{theorem}
In Theorem \ref{thm:AR_cov}, if $\epsilon$ has a survival function given by $\overline{F}(x) = x^{-\alpha} l(x)$ for $x > 0$, where $\alpha \in (0, 1)$ and $l$ is a slowly varying function, then $\mathbb{E}(\epsilon^q) < \infty$ for $0 < q < \alpha$. By Proposition 13.3.2 in \cite{brockwell1991time}, the corresponding AR(1) process $\{X_t\}$ is strictly stationary.

\subsubsection{Dependent Taylor's Law for Higher Moments}
For $p > 0$, denote the $p$th sample moment by $M_{n,p} := {n^{-1}}\sum^n_{i=1}X_i^p$. 
For  $h, k > \alpha$, define
\begin{equation*}\label{eq:imath}
\imath(h,k) := \frac{h-\alpha}{k-\alpha}.
\end{equation*}

\begin{lemma}\label{lemma:log_M_np/log_n}
Let $X_1,\ldots,X_n,\ldots$ be a sequence of 
nonnegative random variables with a common survival function $\overline{F}(x) = x^{-\alpha}l(x)$, 
where $l$ is slowly varying and $\alpha > 0$ (not necessarily in $(0, 1)$). 
If  Condition A($p$) holds for $p > \alpha$, then
\begin{align*}
        \frac{\log M_{n,p}}{\log n} =\frac{p-\alpha}{\alpha} +O_p\left(\frac{\log c_n}{\log n}\right) + O\left( \frac{|\log l(t_n)|}{\log n}\right),
    \end{align*}
    and thus
\begin{equation}\label{eq:estimator_moment_n}
                \frac{\log M_{n,p}}{\log n} \stackrel{\mathbb{P}}{\rightarrow}  \frac{p-\alpha}{\alpha}.
    \end{equation}
\end{lemma}

In addition to establishing the limit in probability, we also obtain the rate of convergence of $\frac{\log M_{n,h_1}}{\log M_{n,h_2}}  - \imath(h_1, h_2)$ under some conditions on the slowly varying function $l$ in Theorem \ref{thm:moment}. 
{Here, for a sequence of random variables $Y_n$, we say its rate of convergence is $r_n$ if $r_n Y_n = O_p(1)$ and $(r_n Y_n)^{-1} = O_p(1)$.}
Let $\alpha > 0$. For some $\delta \in (0, \alpha)$, 
we define  
\begin{align*}
l_{1n} &:= \left(2 \max\{l(n),l(n)^{-1},1\}  \sup_{ n^{\frac{1}{\alpha+\delta}} \leq x \leq n^{\frac{1+\delta}{\alpha-\delta}}}l(x)\right)^{1/\alpha},\\
l_{2n} &:=  \left(\frac{1}{2}\inf_{ n^{\frac{1}{\alpha+\delta}} \leq x \leq n^{\frac{1+\delta}{\alpha-\delta}}}l(x)\right)^{1/\alpha},\\
I_{n,l} &:= [n^{1/\alpha}l_{2n}, n^{1/\alpha} l_{1n}], \\
s_{1n,l} &:= \frac{|\log (\sup_{x\in I_{n,l}} l(x))| }{\log n}, \\
s_{2n,l} &:= \frac{|\log (\inf_{x\in I_{n,l}} l(x))| }{\log n}.
\end{align*}

\begin{theorem}[Taylor's Law for Higher Moments]\label{thm:moment}
Let $X_1,\ldots,X_n,\ldots$ be a sequence of 
nonnegative random variables with a common survival function $\overline{F}(x) = x^{-\alpha}l(x)$, where $l$ is slowly varying and $\alpha > 0$ (not necessarily in $(0, 1)$). Let $h_1, h_2 > \alpha$. Suppose that Condition A($p$) holds for $p = h_1$ and $p = h_2$.
\begin{enumerate}[(a)]
    \item 
    Then
    \begin{equation}\label{eq:pareto rate}
R_{n,h_1,h_2}^{\text{hm}} :=  \frac{\log M_{n,h_1}}{\log M_{n,h_2}}  - \imath(h_1, h_2) =O_p\left( \frac{\log c_n}{\log n} \right) + O \left( \frac{|\log l(t_n)|}{\log n} \right).
\end{equation}

\item  
\begin{enumerate}[(i)]
    \item  If  $\lim_{x \rightarrow \infty}l(x) = \infty$, then 
     $R_{n,h_1,h_2}^{\text{hm}} = O_p(s_{1n,l})$ and $(R_{n,h_1,h_2}^{\text{hm}})^{-1} = O_p(s_{2n,l}^{-1})$.
    \item If  $\lim_{x \rightarrow \infty}l(x) = 0$, then  $R_{n,h_1,h_2}^{\text{hm}} = O_p(s_{2n,l})$ and $(R_{n,h_1,h_2}^{\text{hm}})^{-1} = O_p(s_{1n,l}^{-1})$.

\end{enumerate}
\end{enumerate}
\end{theorem}

In view of the results in Theorem \ref{thm:moment}(b), the following remark gives sufficient conditions
on  the slowly varying function such that the rate of convergence of $R_{n,h_1,h_2}^{{hm}}$
can be established. Similar remarks also apply to the following Theorems \ref{thm:higher_central_moments}, \ref{thm:moment_hetero}, A.3, and  A.5.

\begin{remark}\label{remark:rate}
    If $\lim_{x \rightarrow \infty}l(x) = \infty$ and  there exists a constant $C>0$ independent of $n$ such that $|\log ( \sup_{x \in I_{n,l}} l(x)) | \leq C|\log \left(  \inf_{x \in I_{n,l}} l(x) \right) |$, then the rate of convergence of $R_{n,h_1,h_2}^{\text{hm}}$ is
 ${(\log n)}/{|\log \left( \sup_{x \in I_{n,l}} l(x)\right) | } $. Similarly, if $\lim_{x \rightarrow \infty} l(x) = 0$ and there exists a constant $C>0$ independent of $n$ such that $|\log ( \inf_{x \in I_{n,l}} l(x)) | \leq C|\log ( \sup_{x \in I_{n,l}} l(x)) |$, then the rate of convergence of $R_{n,h_1,h_2}^{\text{hm}}$ is ${(\log n)}/{|\log \left( \inf_{x \in I_{n,l}} l(x)\right) | } $.
\end{remark}

Several examples of slowly varying functions  illustrate the rate of convergence results in Theorem \ref{thm:moment}(b) and the discussion in Remark \ref{remark:rate}.

\begin{example}
\begin{enumerate}[(a)]
    \item Suppose that $l(x) = \exp((\log x)^\beta)$ for some $\beta \in (0,1)$, which is known to be slowly varying \citep[p. 16]{bingham1989regular}. Clearly, $l(x)$ is increasing and  $l(x) \rightarrow\infty$. Thus,
\begin{align*}
    \left|\log\left(\inf_{x\in I_{n,l}} l(x) \right) \right|=(\log (n^{1/\alpha }l_{2n}))^\beta
    = \left(\frac{1}{\alpha} \log n + \log l_{2n}\right)^\beta \leq \left( \frac{2}{\alpha} \log n \right)^\beta,
\end{align*}
as $\log n \geq  \log l_{2n}$ for all sufficiently large $n$ when $\beta < 1$. Indeed, note that 
\begin{equation*}
  \frac{\log l_{2n}}{\log n}  = \frac{ \frac{1}{\alpha}(\log\frac{1}{2} + (\log n^{1/(\alpha+\delta)} )^\beta)}{\log n} \rightarrow 0
\end{equation*}
as $n \rightarrow \infty$. On the other hand,
\begin{align*}
    \left|\log \left(\sup_{x\in I_{n,l}} l(x)\right)\right|=(\log (n^{1/\alpha }l_{1n}))^\beta \geq \left( \frac{1}{\alpha} \log n\right)^\beta.
\end{align*}
In view of Remark \ref{remark:rate},  the rate of convergence of $R_{n,h_1,h_2}^{\text{hm}}$ is $(\log n)^{1 - \beta}$.
\item Suppose that $l(x)=\exp(-(\log x)^\beta)$ for some $\beta \in (0, 1)$. It can be shown that $l$ is slowly varying, similar to the function $x \mapsto \exp((\log x)^\beta)$. Clearly, $l(x)$ is decreasing and $l(x) \rightarrow 0$. Similar to (a), we have
\begin{align*}
    \left|\log\left(\sup_{x\in I_{n,l}} l(x)\right)\right|
    = \left|-(\log (n^{1/\alpha }l_{2n}))^\beta\right| 
    \leq \left( \frac{2}{\alpha} \log n \right)^\beta
\end{align*}
and
\begin{align*}
\left|\log\left(\inf_{x\in I_{n,l}} l(x)\right)\right|
   = \left|-(\log (n^{1/\alpha }l_{1n}))^\beta\right| \geq \left( \frac{1}{\alpha} \log n\right)^\beta.
\end{align*}
In view of Remark \ref{remark:rate},  the rate of convergence of $R_{n,h_1,h_2}^{\text{hm}}$ is $(\log n)^{1 - \beta}$.

\item Suppose that $l(x)=\log(x)f(x)$ where $0 < a \leq f(x) \leq b < \infty$ such that $l$ remains slowly varying. Clearly, $l(x) \rightarrow \infty$ and 
\begin{equation*}
     \left|\log\left(\inf_{x\in I_{n,l}} l(x) \right) \right| \leq \log (b\log(n^{1/\alpha} l_{2n})) = \log b + \log \left( \frac{1}{\alpha}\log n + \log l_{2n}\right).
\end{equation*}
As
\begin{equation*}
    l_{2n} \leq \left( \frac{b}{2} \log (n^{1/(\alpha+\delta)})\right)^{1/\alpha},
\end{equation*}
we have for all sufficiently large $n$,
\begin{equation*}
    \left|\log\left(\inf_{x\in I_{n,l}} l(x) \right) \right| \leq 2 \log \log n.
\end{equation*}
On the other hand, for all sufficiently large $n$,
\begin{equation*}
         \left|\log\left(\sup_{x\in I_{n,l}} l(x) \right) \right| \geq \log (a \log (n^{1/\alpha}l_{1n})) \geq \log \log n.
\end{equation*}
In view of Remark \ref{remark:rate},  the rate of convergence of $R_{n,h_1,h_2}^{{hm}}$ is $\frac{\log n}{\log \log n}$.
\end{enumerate}
\end{example}

The next example shows that when $l(x)$ does not diverge to infinity or go to $0$, the upper bound obtained in Theorem \ref{thm:moment}(a) can be arbitrarily close to the true rate.

\begin{example}
Let $X_1,\ldots,X_n$ be $i.i.d.$ Pareto random variables, with survival function $\mathbb{P}(X_i>x)=x^{-\alpha}$ for $i=1,\ldots,n$ and $0<\alpha<1$ for $x\geq1$. In this case, the slowly varying function is $l(x)=1$.
From Theorem \ref{thm:moment}(a), for any $h_1, h_2 > \alpha$, we have
\eqref{eq:pareto rate}.
Since $l(x)=1$, $\log l(t_n) =\log 1=0$ and thus the second summand on the right side of (\ref{eq:pareto rate}) is 0. 
As a result, the upper bound of the convergence rate of 
$\frac{\log M_{n,h_1}}{\log M_{n,h_2}}  - \imath(h_1, h_2)$ to 0 is $\frac{\log c_n}{\log n}$; indeed,  $\log c_n=o(\log n)$ can be chosen to be arbitrarily slow, and the upper bound of the convergence rate can be arbitrarily close to $\frac{1}{\log n}$.\\
  \end{example}
This result coincides with that of \cite{cohen2020heavy}. 
Define $M^c_{n,2}:=n^{-1}\sum_{i=1}^n(X_i-M_{n,1})^2$ and 
$b=\frac{2-\alpha}{1-\alpha}$
and $W_n:=\frac{M^c_{n,2}}{(M_{n,1})^b}$. 
From Section 3(a) of \cite{cohen2020heavy}, we have that $W_n \overset{d}{\rightarrow }S_{\alpha/2}/S^b_\alpha$, where $S_\alpha$ is a stable random variable with index $\alpha \in (0,2)$ and shape parameter $\beta=1$, which implies that
\begin{equation}\label{eq:pareto_conv_dist}
    \log M_{n,2}^c - b \log M_{n,1} \stackrel{d}{\rightarrow} \log \left(S_{\alpha/2}/S^b_\alpha \right).
\end{equation}
From \eqref{eq:pareto_conv_dist}, we obtain
\begin{equation*}
    \log n \left(\frac{\log M_{n,1}}{\log n}\left(\frac{\log M^c_{n,2}}{\log M_{n,1}}-b \right)\right) \overset{d}{\rightarrow }\log \left(S_{\alpha/2}/S^b_\alpha \right).
\end{equation*}
By Lemma \ref{lemma:log_M_np/log_n}, $\frac{\log M_{n,1}}{\log n} \stackrel{\mathbb{P}}{\rightarrow} \frac{1-\alpha}{\alpha}$. 
Thus, the convergence rate of $\frac{\log(M^c_{n,2})}{\log(M_{n,1})}-\frac{2-\alpha}{1-\alpha}$ is $\frac{1}{\log n}$. 
Our results in this article imply that the convergence rate is arbitrarily close to $\frac{1}{\log n}$. 
Thus the two results coincide, despite the fact that we are interested in finding the probability limit with its convergence rate while \cite{cohen2020heavy} 
specified the limit of  the convergence in distribution without the rate of convergence.

\begin{remark}
    Theorem \ref{thm:moment} implies that for $h_1, h_2 > \alpha$,
\begin{equation*}
        \frac{\log M_{n,h_1}}{\log M_{n,h_2}} \stackrel{\mathbb{P}}{\rightarrow} \imath(h_1,h_2).
\end{equation*}
In particular, for $\alpha \in (0, 1)$ and $k \in \mathbb{N}$,
\begin{equation}\label{eq:estimator_TL}
        \frac{\log M_{n,k}}{\log M_{n,1}} \stackrel{\mathbb{P}}{\rightarrow} \frac{k-\alpha}{1-\alpha}.
\end{equation}
\end{remark}

\subsection{Taylor's Law for Higher Central Moments}
For integer $k \geq 2$, define the $k$th sample central moments to be
\begin{equation*}
    M^c_{n,k} := \frac{1}{n}\sum^n_{i=1}(X_i - M_{n,1})^k.
\end{equation*}

\begin{theorem}[Taylor's Law for Higher Central Moments]\label{thm:higher_central_moments}
Let $\alpha > 0$ and $k, h_1, h_2$ be integers greater than $\alpha$. 
Let $X_1,\ldots,X_n,\ldots$ be a sequence of 
nonnegative random variables with a common survival function $\overline{F}(x) = x^{-\alpha}l(x)$, where $l$ is slowly varying. 
\begin{enumerate}[(a)]
    \item Let $\alpha \in (0, 1)$. If Condition A($p$) holds for $p = 1$ and $p = k$, then
\begin{equation*}
  R_{n,k,1}^{\text{hcm1}} :=  \frac{\log |M^c_{n,k}|}{\log M_{n,1}}  - \imath(k, 1) = O_p\left( \frac{\log c_n}{\log n} \right) + O \left( \frac{|\log l(t_n)|}{\log n} \right).
\end{equation*}
\begin{enumerate}[(i)]
   \item  If $\lim_{x \rightarrow \infty}l(x) = \infty$, then 
     $R_{n,k,1}^{\text{hcm1}} = O_p(s_{1n,l})$ and $(R_{n,k,1}^{\text{hcm1}})^{-1} = O_p(s_{2n,l}^{-1})$.
    
    \item If $\lim_{x \rightarrow \infty}l(x) = 0$, then  $R_{n,k,1}^{\text{hcm1}} = O_p(s_{2n,l})$ and $(R_{n,k,1}^{\text{hcm1}})^{-1} = O_p(s_{1n,l}^{-1})$.

\end{enumerate}
   
\item 
Let $\alpha > 0$. If Condition A($p$) holds for $p = h_1$ and $p = h_2$, then
\begin{equation*}
    R_{n,h_1,h_2}^{\text{hcm2}} := \frac{\log |M^c_{n,h_1}|}{\log |M^c_{n,h_2}|}  - \imath(h_1, h_2) = O_p\left( \frac{\log c_n}{\log n} \right) + O \left( \frac{|\log l(t_n)|}{\log n} \right).
\end{equation*}
\begin{enumerate}[(i)]
   \item  If $\lim_{x \rightarrow \infty}l(x) = \infty$, then
     $R_{n,h_1,h_2}^{\text{hcm2}} = O_p(s_{1n,l})$ and $(R_{n,h_1,h_2}^{\text{hcm2}})^{-1} = O_p(s_{2n,l}^{-1})$.
    
    \item If $\lim_{x \rightarrow \infty}l(x) = 0$, then   $R_{n,h_1,h_2}^{\text{hcm2}} = O_p(s_{2n,l})$ and $(R_{n,h_1,h_2}^{\text{hcm2}})^{-1} = O_p(s_{1n,l}^{-1})$.
\end{enumerate}

\end{enumerate}

\end{theorem}

The following corollary,  an important special case of Theorem \ref{thm:higher_central_moments},  recovers a commonly known form of Taylor's law; see  \citet[Section 6]{brown2017taylor} and  \citet[Eq. (10)]{brown2021taylor}.
\begin{corollary}
Let $X_1,\ldots,X_n,\ldots$ be a sequence of 
nonnegative random variables with a common survival function $\overline{F}(x) = x^{-\alpha}l(x)$, where $l$ is slowly varying and $\alpha \in (0, 1)$. Let $V_n := M_{n,2} - M_{n,1}^2$ be the sample variance.  
If Condition A($p$) holds for $p = 1$ and $p = 2$, then
\begin{equation}\label{eq:estimator_TL_variance}
    \frac{\log V_n}{\log M_{n,1}} \stackrel{\mathbb{P}}{\rightarrow} \frac{2-\alpha}{1-\alpha}.
\end{equation}    
If $\lim_{x \to \infty}l(x)=\infty$, then
\begin{align*}
     \frac{\log V_n}{\log M_{n,1}}  -  \frac{2-\alpha}{1-\alpha} = O_p\left(s_{1n,l}\right) \quad \text{and} \quad 
     \left(\frac{\log V_n}{\log M_{n,1}}  -  \frac{2-\alpha}{1-\alpha} \right)^{-1} = O_p(s_{2n,l}^{-1}).
\end{align*}
If $\lim_{x \to \infty}l(x)=0$, then
\begin{align*}
     \frac{\log V_n}{\log M_{n,1}}  -  \frac{2-\alpha}{1-\alpha} = O_p\left(s_{2n,l}\right) \quad \text{and} \quad 
     \left(\frac{\log V_n}{\log M_{n,1}}  -  \frac{2-\alpha}{1-\alpha} \right)^{-1} = O_p(s_{1n,l}^{-1}).
\end{align*}
\end{corollary}

\section{Taylor’s Law with Heterogeneous Heavy-tailed Data}\label{sec:heterogeneous}
In this section, we consider a heterogeneous case where not all $X_i$'s 
follow the same distribution. 
\cite{cohen2020heavy} considered a similar setting for Taylor's law.
Their results focused on the scaling relationship between the sample mean and the sample variance.
We extend the results to higher moments and relax their  assumption of independence among random variables. 

Specifically, for each positive integer $n$,
we consider $u_n<n$ random variables $X_{i,U}$, $i=1,\ldots,u_n$.
Each $X_{i,U}$ has the same distribution with survival function $\overline{F}_U(x) = x^{-\alpha}l(x)$ for some $\alpha \in (0, 1)$, where $l$ is slowly varying. In addition to $X_{i,U}$'s, we also have $n -u_n$ random variables $X_{i,V}$, $i=1,\ldots,n-u_n$, where each $X_{i,V}$ has the same distribution with
\begin{equation}\label{eq:heter_heavy_tail}
    \lim_{x \rightarrow \infty} \frac{\mathbb{P}(X_{1,V} > x)}{\mathbb{P}(X_{1,U} > x)} = 0,
\end{equation}
meaning that $X_{1,U}$ has a heavier tail than $X_{1,V}$. 
Define $\overline{F}_V$ to be the survival function of $X_{1,V}$. 
We assume that $u_n/n \rightarrow p^* \in (0, 1]$ as $n\to\infty$. 
As in Section \ref{sec:main results}, let $\{t_n\}$ satisfy
\begin{equation}\label{eq:tn choice heterogeneous}
    \frac{n l(t_n)}{t_n^{\alpha}} \sim \frac{1}{c_n},
\end{equation}
and let $\{v_n\}$ satisfy
\begin{equation*}
    \frac{n l(v_n)}{v^\alpha_n} \sim c_n.
\end{equation*}
\begin{remark}
 (\ref{eq:heter_heavy_tail}) is satisfied if $\overline{F}_U(x) = x^{-\alpha_1}l(x)$ and $\overline{F}_V(x) = x^{-\alpha_2}l(x)$ with $\alpha_1 < \alpha_2$.
\end{remark}

Reusing an earlier notation for a related purpose,  we define $M_{n,p} := n^{-1}( \sum^{u_n}_{i=1}X_{i,U}^p + \sum^{n-u_n}_{i=1}X_{i,V}^p)$. 
Let $\breve{X}_{i,U} := X_{i,U} \mathbbm{1}(X_{i,U} < v_n)$ 
and $\breve{X}_{i,V} := X_{i,V} \mathbbm{1}(X_{i,V} < v_n)$. Finally, let $\breve{Y}_i := \breve{X}_{i,U}$ for $i=1,\ldots,u_n$ and $\breve{Y}_{u_n+i} := \breve{X}_{i,V}$ for $i=1,\ldots,n-u_n$.

\begin{theorem}[Taylor's Law for Higher Moments with Heterogeneous Data]\label{thm:moment_hetero}
For any positive integer $n>1$, assume that each of 
$u_n < n$ random variables $X_{i,U}$, $i=1,\ldots,u_n$, follows the same distribution with survival function $\overline{F}_U(x) = x^{-\alpha}l(x)$, where $l$ is slowly varying for some $\alpha \in (0, 1)$; 
and all of $n -u_n$ random variables $X_{i,V}$, $i=1,\ldots,n-u_n$,  have the 
same distribution such  that (\ref{eq:heter_heavy_tail}) holds. Let $h_1, h_2 > \alpha$. Assume that
\begin{equation}\label{eq:cov_cond_hetero}
    \sum_{i \neq j} \Cov(\breve{Y}_i^p, \breve{Y}_j^p) = o(v^{2p}_n c^2_n)
\end{equation}
holds for $p = h_1$ and $p = h_2$.  
\begin{enumerate}[(a)]
    \item  Then
    \begin{equation*}
   R_{n,h_1,h_2}^{\text{hetero}} := \frac{\log M_{n,h_1}}{\log M_{n,h_2}}  - \imath(h_1, h_2) =O_p\left( \frac{\log c_n}{\log n} \right) + O \left( \frac{|\log l(t_n)|}{\log n} \right).
\end{equation*}

\item  
\begin{enumerate}[(i)]
     \item  If $\lim_{x \rightarrow \infty}l(x) = \infty$, then
      $R_{n,h_1,h_2}^{\text{hetero}} = O_p(s_{1n,l})$ and $(R_{n,h_1,h_2}^{\text{hetero}})^{-1} = O_p(s_{2n,l}^{-1})$.
    
    \item If $\lim_{x \rightarrow \infty}l(x) = 0$, then   $R_{n,h_1,h_2}^{\text{hetero}} = O_p(s_{2n,l})$ and $(R_{n,h_1,h_2}^{\text{hetero}})^{-1} = O_p(s_{1n,l}^{-1})$.

\end{enumerate}

\end{enumerate}
 
\end{theorem}

\section{Taylor's Law with Correlated Heavy-tailed Data}\label{sec:corr}
In Section \ref{sec:main results}, we assumed weak dependence, Condition A($p$),
among the truncated random variables.
Here, in the spirit of \citet[Theorems 1 and 2]{cohen2022covid},
we allow  all the random variables to  be highly correlated 
provided that the dependence among them becomes weak if we condition on a $\sigma$-field.
Our results in Section \ref{sec:main results} continue to hold if
Condition A($p$)
holds when conditioning on a $\sigma$-field $\mathcal{G}$. Of course, a special case is that the random variables are conditionally independent given $\mathcal{G}$.

\begin{theorem}\label{thm:correlated_case}
 Let $X_1,\ldots,X_n,\ldots$ be a sequence of nonnegative 
 random variables with common survival function $\overline{F}(x) = x^{-\alpha}l(x)$, where $l$ is slowly varying and $\alpha > 0$ (not necessarily in $(0, 1)$).
Suppose that $\mathcal{G}$ is a $\sigma$-field and
\begin{equation}\label{eq:assumption_correlated_case}
   \mathbb{E} \left( \frac{\sum_{i \neq j}\Cov(\breve{X}_i^p, \breve{X}_j^p|\mathcal{G})}{v_n^{2p} c_n^2}\right) = o(1),
\end{equation}
for $p = h_1, h_2 > \alpha$. Then the same conclusions as in Theorem \ref{thm:moment} hold.
\end{theorem}

Two examples illustrate this framework.

\begin{example}\label{example:correlated_2020}
This example is considered in \citet[Section 3(c)]{cohen2020heavy}. 
Let $N_0$, $N_1,\ldots$ be $i.i.d.$ standard normal random variables and let $0 \leq \rho < 1$. 
Let $Z_i := \sqrt{\rho}N_0 + \sqrt{1-\rho}N_i$ for $i=1,\ldots,n$. Then each $Z_i$ is a standard normal with $\Cov(Z_i, Z_j) =\rho$ when $i\neq j$. 
For $\alpha \in (0, 1)$,  define $X_i :=1/(Z^2_i)^{1/(2\alpha)}$.
Then $X_i$ is regularly varying.
Conditional on $\mathcal{G} := \sigma(N_0)$, $Z_i$'s are independent. Thus $X_i$'s are also independent conditional on $\sigma(N_0)$. 
Condition (\ref{eq:assumption_correlated_case}) holds trivially, so our Taylor's law holds. 
\end{example}

\begin{example}\label{example:correlated_2022}
This example is considered in \citet[Theorem 2]{cohen2022covid}.
Let $\{G_i\}$ be Gaussian with mean $0$ and covariance function $\gamma(\cdot,\cdot)$, i.e., $\Cov(G_i, G_j) = \gamma(i,j)$.
Let $Z_i$'s be $i.i.d.$ Pareto-distributed with tail index $\alpha \in (0, 1)$. Let $X_i := Z_i e^{G_i}$. Then $X_i$ are identically distributed and regularly varying with index $\alpha$. Conditioning on $\mathcal{G} := \sigma(G_i, i\in\mathbb{N})$, $X_i$'s are independent. Condition (\ref{eq:assumption_correlated_case}) holds trivially and our Taylor's law holds. 
\end{example}

Our results go beyond the results 
of \cite{cohen2020heavy} and \cite{cohen2022covid} in showing that,
under our covariance condition \eqref{eq:assumption_correlated_case}
for dependent random variables, 
Taylor's law will hold for the variance and for higher moments.

\section{Taylor's Law for Heavy-tailed Network Data}\label{sec:network TL}

Although we have stated the generalizations of Taylor's law for a denumerable set
of random variables $\{X_1,X_2,\ldots,X_n, \ldots\}$, 
the results clearly also hold for a set of random variables $\{X_i: i \in \mathcal{I}_n\}$ with an arbitrary non-empty index set $\mathcal{I}_n$ when we modify the notation in the corresponding conditions accordingly. 
For example, for the Taylor's law for higher moments in Theorem \ref{thm:moment}, Condition A($p$) becomes
\begin{equation}\label{eq:network_cov_assumption}
    \sum_{i,j \in \mathcal{I}_n, i \neq j} \Cov(\breve{X}_i, \breve{X}_j) = o(v_n^{2p} c_n^2).
\end{equation}
Here is  an example where (\ref{eq:network_cov_assumption}) is satisfied in a network graph.

\begin{example}\label{Example:network}
For a natural number $n$,
let $G_n = (V_n, E_n)$ be an 
undirected or directed  
graph, where $V_n = \{v_1,\ldots,v_n\}$ is a  set of $n$  vertices and $E_n$ is a nonempty set of
edges. Define $[n] :=\{1,\ldots,n\}$ and let $U_{n,h}$ be the subset of $\mathcal{S}_n := \{(i,j): i,j \in [n]\}$ such that the unweighted distance between two nodes is $h$:
\begin{equation*}
U_{n,h} := \{ (i,j) \in \mathcal{S}_{n}:  d(v_i, v_j) = h  \},
\end{equation*}
where $d(v_i,v_j)$ is the shortest unweighted distance between $v_i$ and $v_j$ according to $E_n$. 
$U_{n,h}$ describes  a specific network structure. 
Suppose that each node $v_i$ is associated with a random variable $X_i$, 
and that each of
 $X_1,\ldots,X_n$ has a common survival function $\overline{F}(x) = x^{-\alpha}l(x)$, 
 where $l$ is slowly varying and $\alpha > 0$. Assume that 
    \begin{enumerate}
        \item $\Cov(\breve{X}_i,\breve{X}_j) = 0$ if $d(v_i,v_j) > M$, where $M$ is a nonnegative integer;
        \item for all $i \in [n]$, $|\{ j: (i,j) \in U_{n,1}\}| \leq  K$, where $K$ is a positive integer.
        (In graph-theoretic language, the degree of every vertex does not exceed $K$.)
    \end{enumerate}
    Then Condition \eqref{eq:network_cov_assumption} holds for $p > \alpha$; see Section G in the appendix.
    
\end{example}

\section{Network Data Illustration}\label{sec:application}
In this section, we illustrate Taylor's law for network data, as introduced in Section 
\ref{sec:network TL}, using three real-world datasets. To the best of our knowledge, this is 
the first study to explore the application of Taylor's law in network data. 

For each dataset, we estimate the heavy-tail index $\alpha$ using 
the stabilizing region of the alternative Hill plot \citep{resnick1997smoothing}. 
All three datasets have an estimated tail index less than~$1$. 
We also compare this estimate  with the tail index implied by Taylor's law, which is given by $(2 - \frac{\log V_n}{\log M_{n,1}}) / (1 - \frac{\log V_n}{\log M_{n,1}})$ in view of \eqref{eq:estimator_TL_variance}.
We find that the two estimates generally agree closely. 
To illustrate Taylor's law for each dataset, we construct subsamples of varying sizes from the full dataset
and, for each subsample, we plot log-variance on the vertical axis against log-mean on the  horizontal axis. 
All the plots reveal a clear linear relationship, consistent with Taylor's law for network data.
Throughout we use the natural logarithm, base $e$.

\subsection{Wikipedia Talk Dataset}
The Wikipedia Talk dataset was obtained from the Stanford Network Analysis Project\footnote{\url{https://snap.stanford.edu/data/wiki-Talk.html}}. Wikipedia is a free encyclopedia collaboratively written by volunteers worldwide. Each registered user has a talk page that serves as a space for communication and coordination. Both the user and others can edit this page to exchange messages or discuss edits to articles.

The dataset captures all user interactions via talk page edits from the inception of Wikipedia through January 2008. This network represents the communication structure within the Wikipedia community. Each node corresponds to a user, and a directed edge from user A to user B indicates that user A edited the talk page of user B. For each user, the total number of edits they made to other users' talk pages is treated as a node attribute. 
We use these values to illustrate Taylor’s law. 
The dataset contains 2,394,385 users in total. We removed users who never interacted with others, as they are not considered active participants on the platform. After removing these inactive users, the dataset consists of 147,602 users. The maximum number of edits by a single user is 100,022, while the average is 34 and the sample variance is 160,858.

\begin{figure}[ht]
\centering
\includegraphics[width=0.7\linewidth]{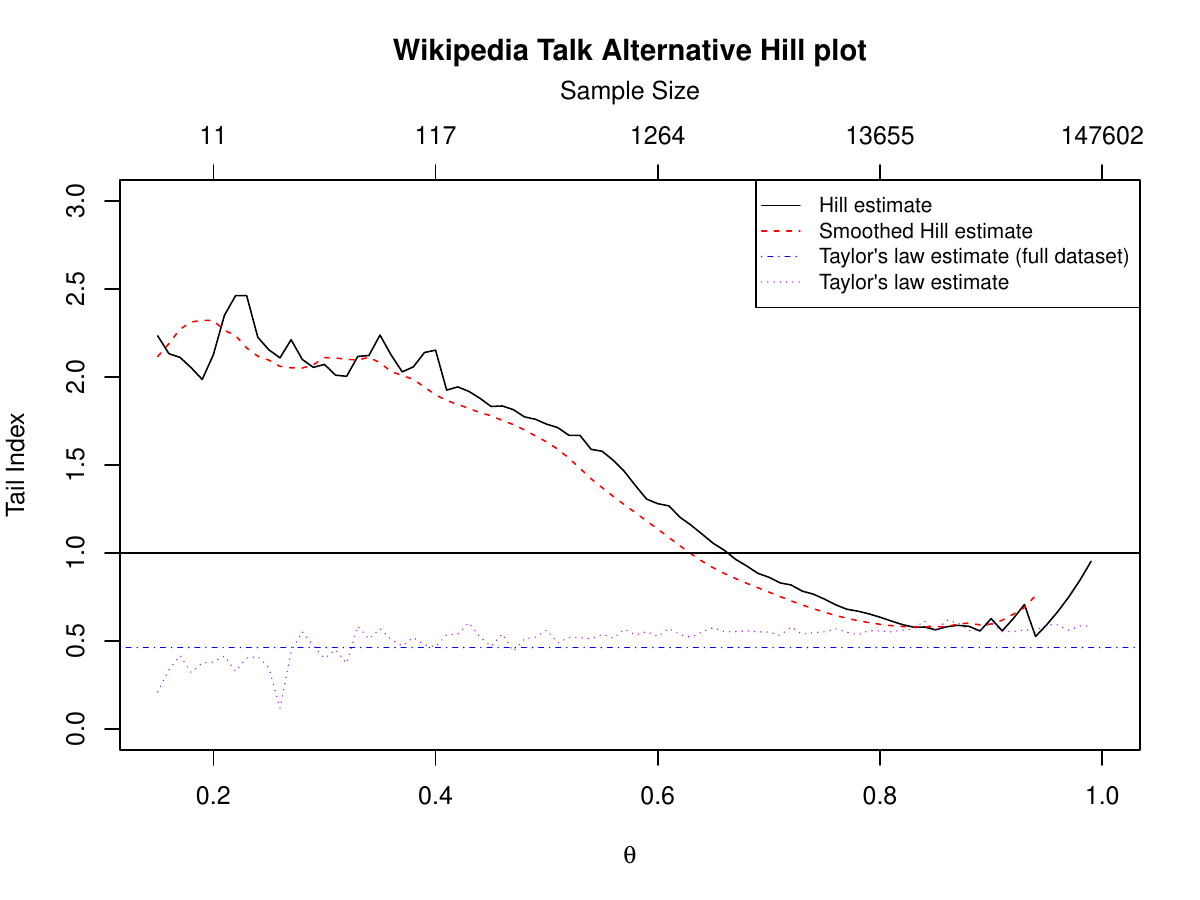}
\caption{Alternative Hill plot and tail index estimates implied by Taylor's law for Wikipedia talk dataset. The Hill and smoothed Hill estimates are plotted against the threshold parameter $\theta$, shown on the bottom $x$-axis. The Taylor’s law estimate based on subsamples of the full dataset has sample sizes indicated on the top $x$-axis. {
The Hill estimate at $\theta = 0.85$ is $0.563$ (99\% CI: 0.554--0.572).}}
\label{fig:wiki_hill}
\end{figure}

In Figure \ref{fig:wiki_hill}, the alternative Hill plot displays estimates of the tail index as a function of the threshold parameter $\theta \in [0,1]$, where the number of upper order statistics used in the Hill estimate is $\lceil n^{\theta} \rceil$. This rescaling facilitates interpretation and improves accuracy \citep{resnick1997smoothing}. 
In the figure, the black line represents the Hill estimate and the red line represents its smoothed version as proposed in \cite{resnick1997smoothing}, computed as the simple average of the subsequent $\lceil n^{\theta} \rceil$ Hill estimates.
The figure shows that the Hill estimate stabilizes for $\theta \in [0.8, 0.9]$.
The tail index inferred from this stabilization region closely aligns with that implied by Taylor's law using the full dataset (blue line), and is between $0$ and $1$.

Figure \ref{fig:wiki_hill} also plots the estimates from Taylor's law using increasingly large random subsamples of the full dataset (purple line), with the corresponding sample sizes indicated on the top $x$-axis. 
While the Hill plot is constructed using varying numbers of upper order statistics from the full dataset, the Taylor's law estimates are computed from fixed-size subsamples. We observe that the tail index estimates derived from Taylor's law converge toward and align with the estimate obtained from the Hill plot.

\begin{figure}[ht]
    \centering
    \includegraphics[width=0.5\linewidth]{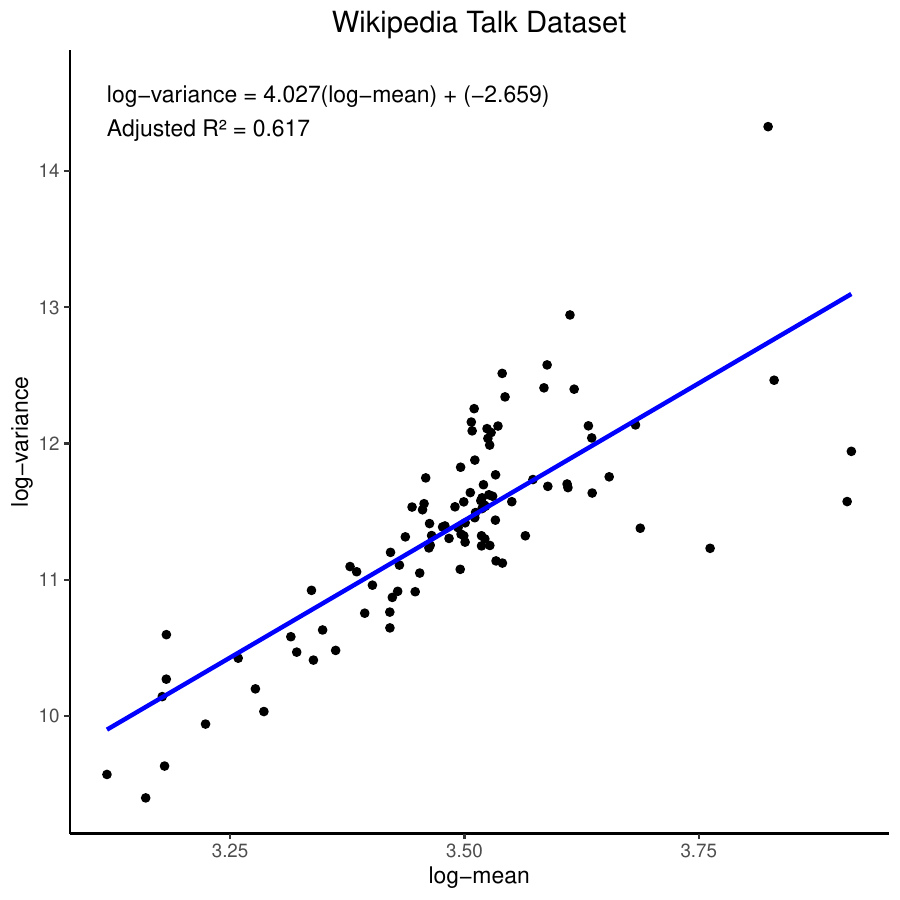}
    \caption{Log$_{10}$-variance versus log$_{10}$-mean across subsamples for Wikipedia talk dataset. 
    The analysis includes 100 pairs (log-mean, log-variance). The regression line was fitted using ordinary least squares.  The 95\% and 99\% confidence intervals for the intercept are $(-4.864 ,-0.455)$ and $(-5.577 ,0.258)$ respectively. For the slope, the 95\% confidence interval is $(3.396 , 4.658)$, while the 99\% confidence interval is $(3.191,  4.863)$. The adjusted $R^2$ value of the regression model is 0.617.}
    
    \label{fig:wiki_tl}
\end{figure}

To demonstrate Taylor's law, we generate a sequence of increasing sample sizes, evenly spaced on a logarithmic scale, drawn from the full dataset. The sequence is created by taking uniform increments from $\log_{10}(500)$ to $\log_{10}(N)$, where $N$ is the total number of samples. Each value is then transformed using the base-10 exponential function and rounded to the nearest integer. This approach results in denser sampling at smaller sizes and sparser sampling at larger sizes, minimizing redundancy where estimates stabilize. For each subsample, we compute the corresponding log-mean and log-variance. Figure \ref{fig:wiki_tl} plots the log-variances versus the log-means of these subsamples. 
A linear fit was approximately consistent with the points, in approximate agreement with Taylor's law.

\subsection{Epinions Dataset}
The Epinions dataset was obtained from the Network Data Repository with Interactive Graph Analytics and Visualization \citep{rossi2015network}\footnote{\url{https://networkrepository.com/rec_epinions_user_ratings.php}}. It is a bipartite rating network from Epinions, an online platform for product reviews.
Nodes represent users and products.
Edges connect users to the products they have rated; see \cite{richardson2003trust} and \cite{massa2005controversial} for more details. 
We use the number of edges connected to each product node (i.e., the number of reviews received by each product) to illustrate a network Taylor’s law. 
The dataset includes 635,268 user nodes and 120,492 product nodes. 
The maximum number of edges connected to a product node is 162,202, with an average of 113.4 edges per product, and a sample variance of 1,129,884. We apply the same methods as used for the Wikipedia talk dataset.

\begin{figure}[ht]
        \centering
        \includegraphics[width=0.7\linewidth]{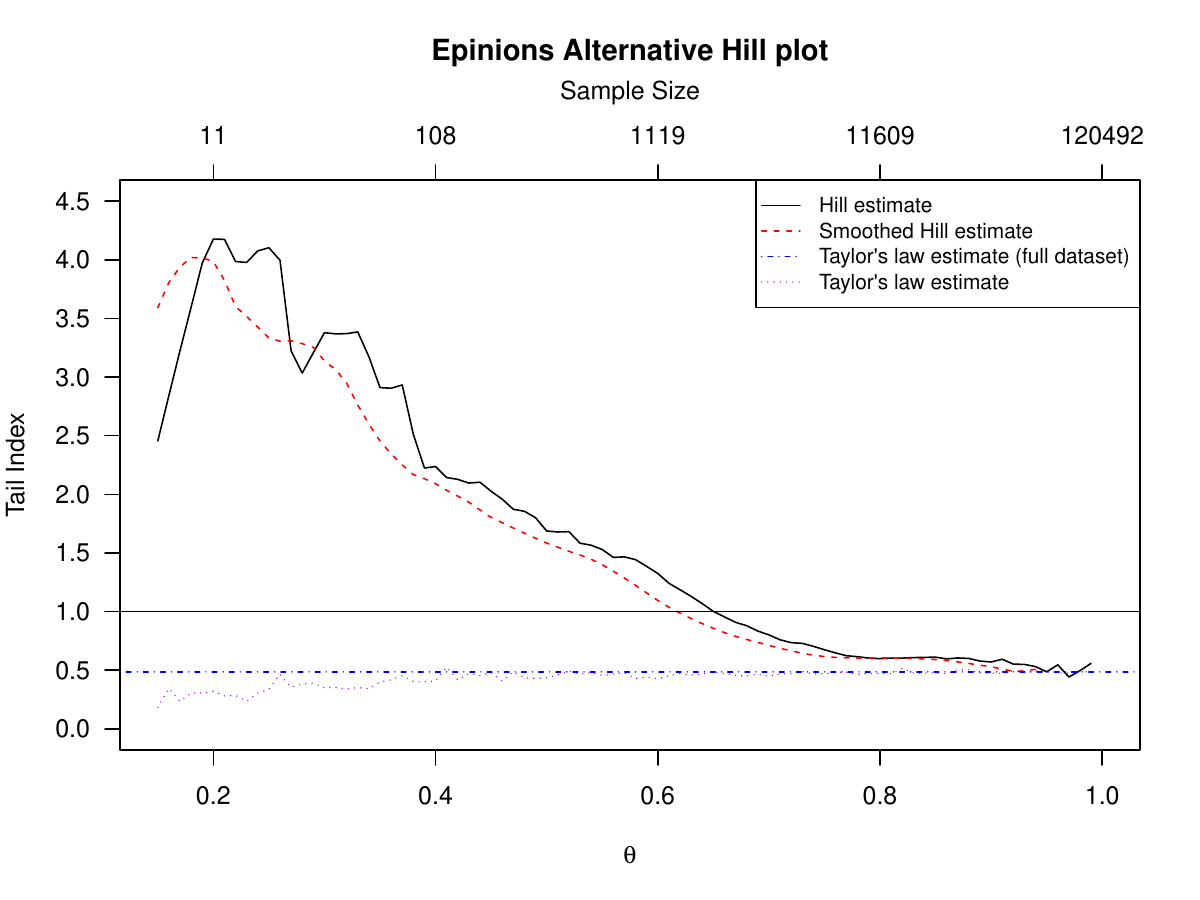}
        \caption{Alternative Hill plot and tail index estimates implied by Taylor's law for Epinions dataset. The Hill and smoothed Hill estimates are plotted against the threshold parameter $\theta$, shown on the bottom $x$-axis. The Taylor’s law estimate based on subsamples of the full dataset has sample sizes indicated on the top $x$-axis. {
        The Hill estimate at $\theta = 0.8$ is $0.539$ (99\% CI: 0.535--0.543).}}
        \label{fig:epinion hill}
    \end{figure}

\begin{figure}[ht]
        \centering
        \includegraphics[width=0.5\linewidth]{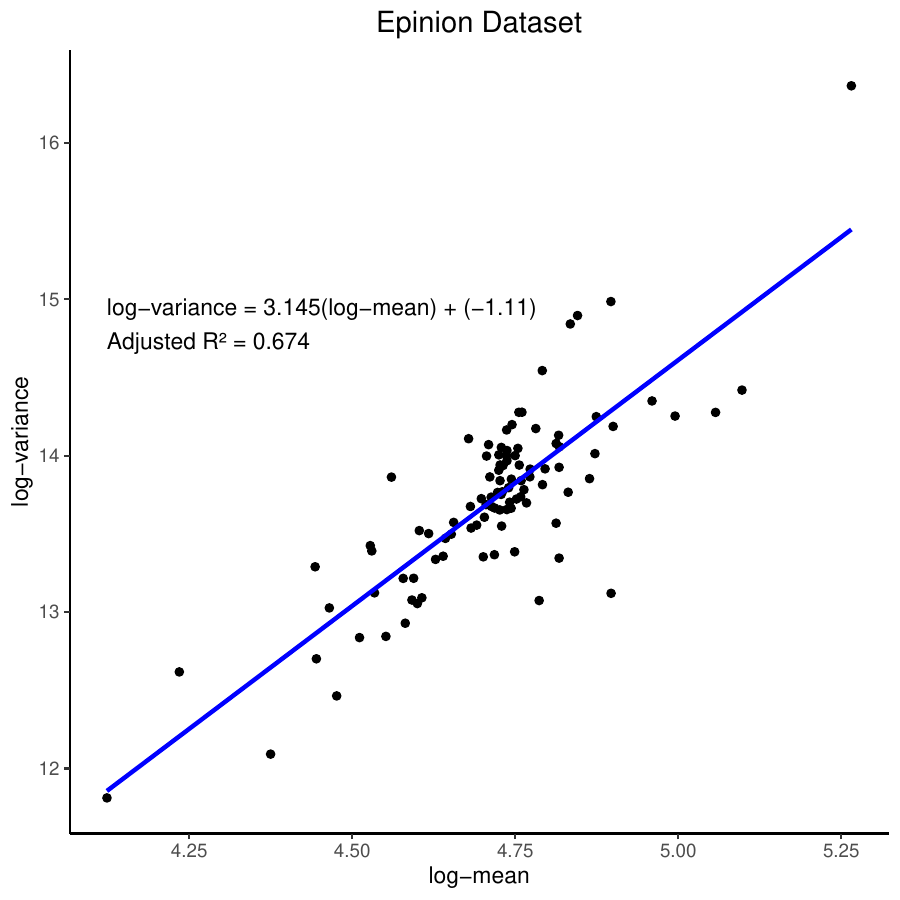}
        \caption{Log-variance versus log-mean across subsamples for Epinions dataset. A total of 100 pairs of (log-mean, log-variance) were included in the analysis. The regression line was fitted using ordinary least squares.  The 95\% and 99\% confidence intervals for the intercept are $(-3.166, 0.945)$ and $(-3.831 , 1.610)$ respectively. For the slope, the 95\% confidence interval is $( 2.709 , 3.580)$, while the 99\% confidence interval is $(2.568 , 3.721)$. The adjusted $R^2$ value of the regression model is 0.674.}
        \label{fig:epinion TL}
    \end{figure}

In Figure~\ref{fig:epinion hill}, the Hill estimate appears to stabilize around $\theta \approx 0.8$. Although there is a slight discrepancy between the tail index inferred from Taylor's law and that estimated from the Hill plot, both indicate that the tail index is less than $1$. Figure~\ref{fig:epinion TL} shows a clear linear relationship between the log-variances and log-means of the subsamples, consistent with Taylor's law.

\subsection{DBpedia Dataset}
The DBpedia dataset was obtained from Netzschleuder\footnote{\url{https://networks.skewed.de/net/dbpedia_country}}. It includes a bipartite network representing affiliations between notable individuals and countries worldwide, extracted from Wikipedia through the DBpedia project. The DBpedia project is a crowd-sourced community effort to extract structured content from information created in various Wikimedia projects.
The set of ``countries” in this dataset includes not only current nations but also former states, empires, kingdoms, and certain country-like entities; \cite{auer2007dbpedia} gives more details. 

We use the number of edges connected to each country node (i.e., the number of entities affiliated with each country) to illustrate a network Taylor’s law. 
The dataset contains 590,112 entity nodes and 2,302 country nodes. 
On average, each country node is connected to 276.8 edges, with a maximum of 111,132 edges and a sample variance of 9,700,633. We again apply the same analysis as before. 
The only difference is that, because this graph has fewer nodes, we generate the log-variance versus log-mean plot by progressively sampling the data starting from 10 nodes instead of 500 nodes in the logarithmic sequence.

\begin{figure}[ht]
    \centering
    \includegraphics[width=0.7\linewidth]{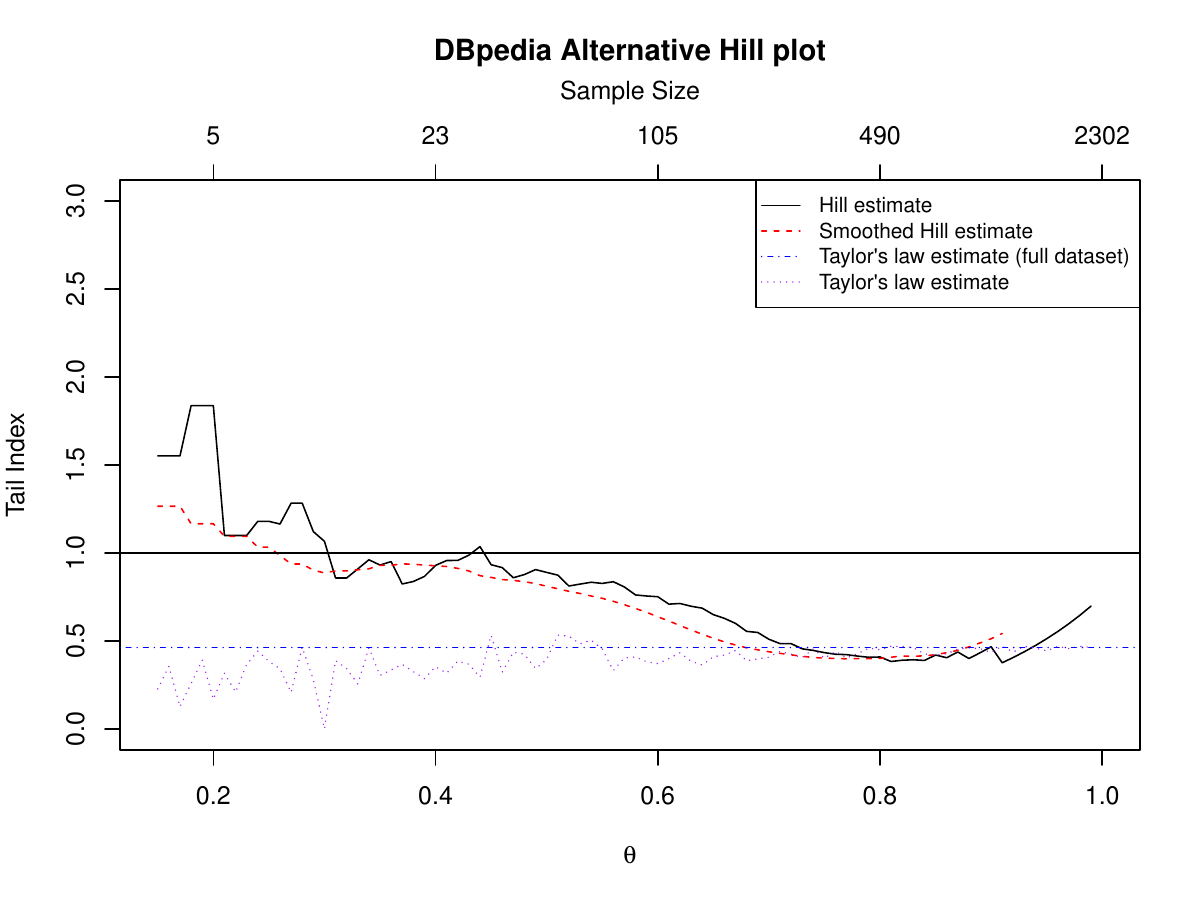}
        \caption{Alternative Hill plot and tail index estimates implied by Taylor's law for DBpedia dataset. The Hill and smoothed Hill estimates are plotted against the threshold parameter $\theta$, shown on the bottom $x$-axis. The Taylor’s law estimate based on subsamples of the full dataset has sample sizes indicated on the top $x$-axis. {
        The Hill estimate at $\theta = 0.8$ is $0.409$ (99\% CI: 0.366--0.463).}}
    \label{fig:DBpedia_hill}
\end{figure}

\begin{figure}[ht]
    \centering
    \includegraphics[width=0.5\linewidth]{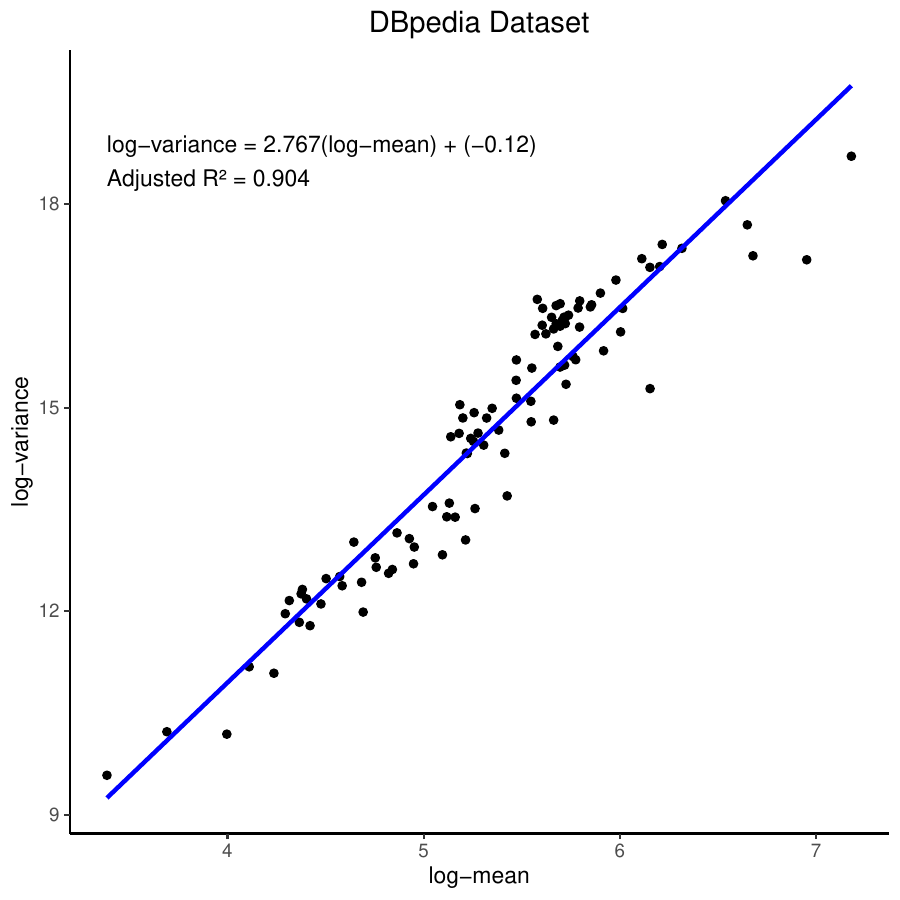}
    \caption{Log-variance versus log-mean across subsamples for DBpedia dataset. A total of 100 pairs of (log-mean, log-variance) were included in the analysis. The regression line was fitted using ordinary least squares.  The 95\% and 99\% confidence intervals for the intercept are $(-1.091, 0.850)$ and $(-1.405, 1.164)$ respectively. For the slope, the 95\% confidence interval is $(2.587, 2.946)$, while the 99\% confidence interval is $(2.528, 3.005)$. The adjusted $R^2$ value of the regression model is 0.904.}
    \label{fig:DBpedia_tl}
\end{figure}

In Figure~\ref{fig:DBpedia_hill}, the Hill estimate appears to stabilize around $\theta \approx 0.8$. The corresponding tail index aligns closely with the estimate inferred from Taylor's law; both suggest that the tail index is less than $1$. Figure~\ref{fig:DBpedia_tl} reveals a clear approximately linear relationship between the log-variances and log-means of the subsamples, approximately consistent with Taylor's law.


{
\subsection{Summary of empirical examples}
In all three examples, the lower bounds of the 99\% confidence intervals for the estimated slopes of plot of the log sample variance versus log sample mean exceed $2$ (see captions of Figures \ref{fig:wiki_tl}, \ref{fig:epinion TL}, and \ref{fig:DBpedia_tl}).
These estimated slopes are consistent with considerable heterogeneity in activity: most nodes have very few connections, while a few nodes have very many connections. 
That the data conform to Taylor’s law with such large slopes suggests that this variability follows a systematic scaling pattern reflecting structural properties of the underlying networks. 
In the Wikipedia dataset, this scaling implies that a small number of users carry out most of the communication. This may reflect the presence of highly active contributors or users with specific roles, such as moderators, who engage disproportionately in talk page interactions. 
In the Epinions dataset, the observed scaling is consistent with preferential attachment, where a few products accumulate the majority of attention as users are more likely to rate popular items. 
In the DBpedia dataset, the pattern suggests that a few countries or historical entities serve as central hubs with widespread affiliations. 
This unequal attention to countries likely reflects historical, cultural, or geopolitical prominence that contributes to their over-representation in structured knowledge sources.

The slopes exceeding $2$ also indicate that the observed counts in these examples are inconsistent with a Poisson distribution or with certain parameterizations of the negative binomial distribution, which are more commonly encountered in other applications of Taylor’s law.

To see directly whether the data are better approximated by an over-dispersed (variance > mean) light-tailed distribution or a heavy-tailed distribution, we fitted both the negative binomial distribution and the Pareto distribution to the node-associated values for the three data examples (Figure \ref{fig:three histograms}). 
We fitted the negative binomial distribution using \texttt{fitdist()} from the \texttt{fitdistrplus} package, which estimates parameters by maximum likelihood.
We fitted the Pareto distribution by least squares using \texttt{epareto()} from the \texttt{VGAM} package in \texttt{R}; see \cite{johnson1995continuous} for additional discussion of least-squares fitting for the Pareto distribution. 
To compare the data with the fitted models visually, we plotted the log of the estimated survival function based on the empirical distribution function and overlaid the corresponding functions from the fitted Pareto and negative binomial distributions.
We log-transformed the horizontal axis so that the Pareto distribution's survival function
would appear as a straight line.

\begin{figure}[ht]
    \centering
    \begin{subfigure}[b]{0.33\textwidth}
        \centering
        \includegraphics[width=\linewidth]{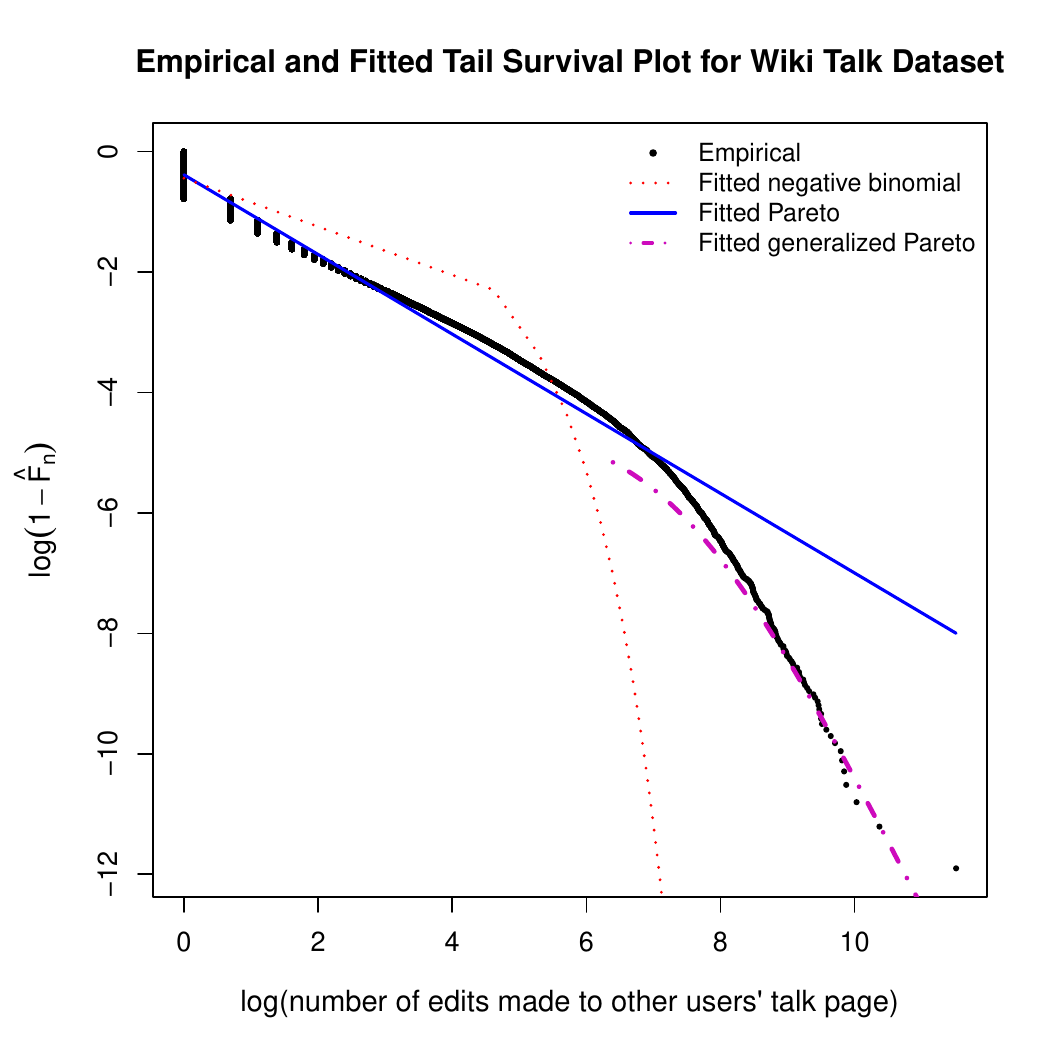}
        \label{fig:fig1}
    \end{subfigure}%
    \hfill
    \begin{subfigure}[b]{0.33\textwidth}
        \centering
        \includegraphics[width=\linewidth]{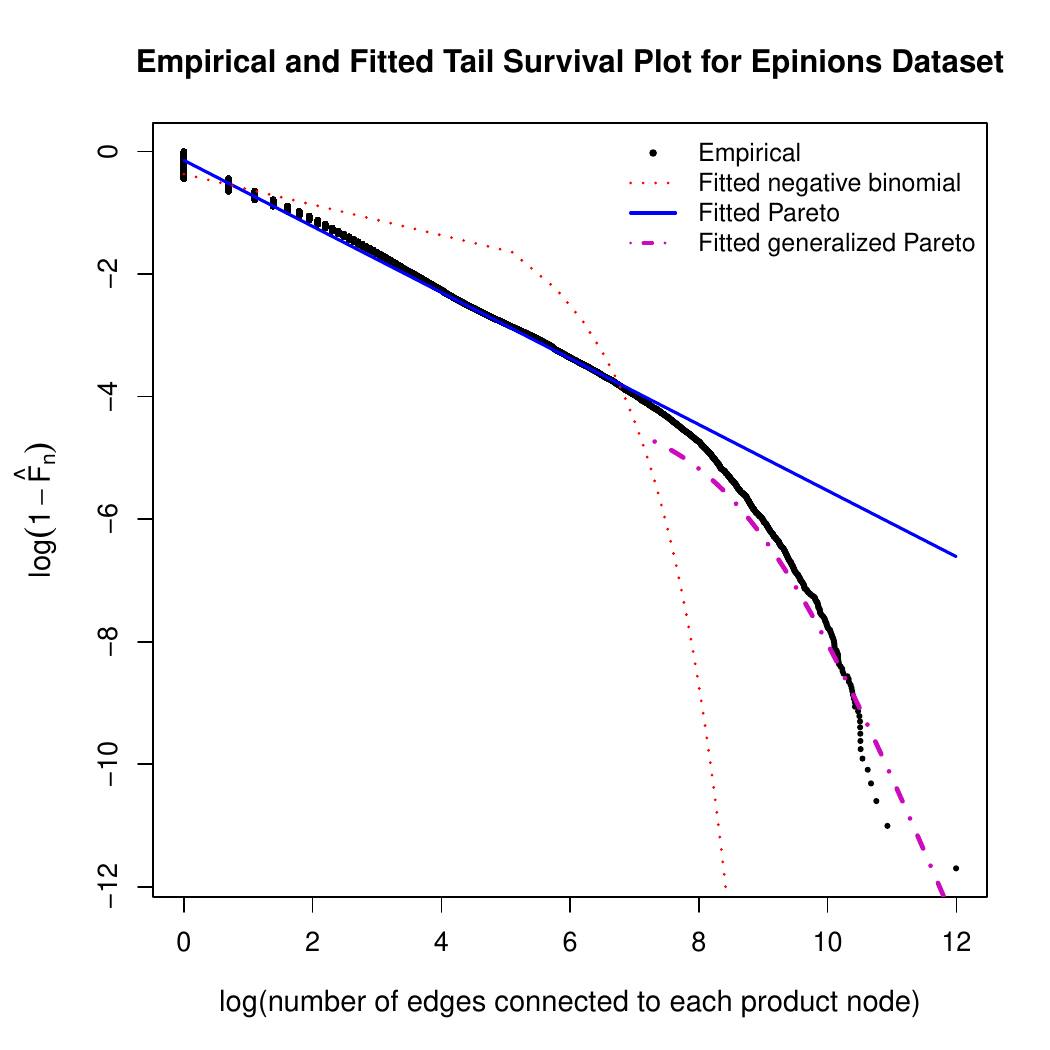}
        \label{fig:fig2}
    \end{subfigure}%
    \hfill
    \begin{subfigure}[b]{0.33\textwidth}
        \centering
        \includegraphics[width=\linewidth]{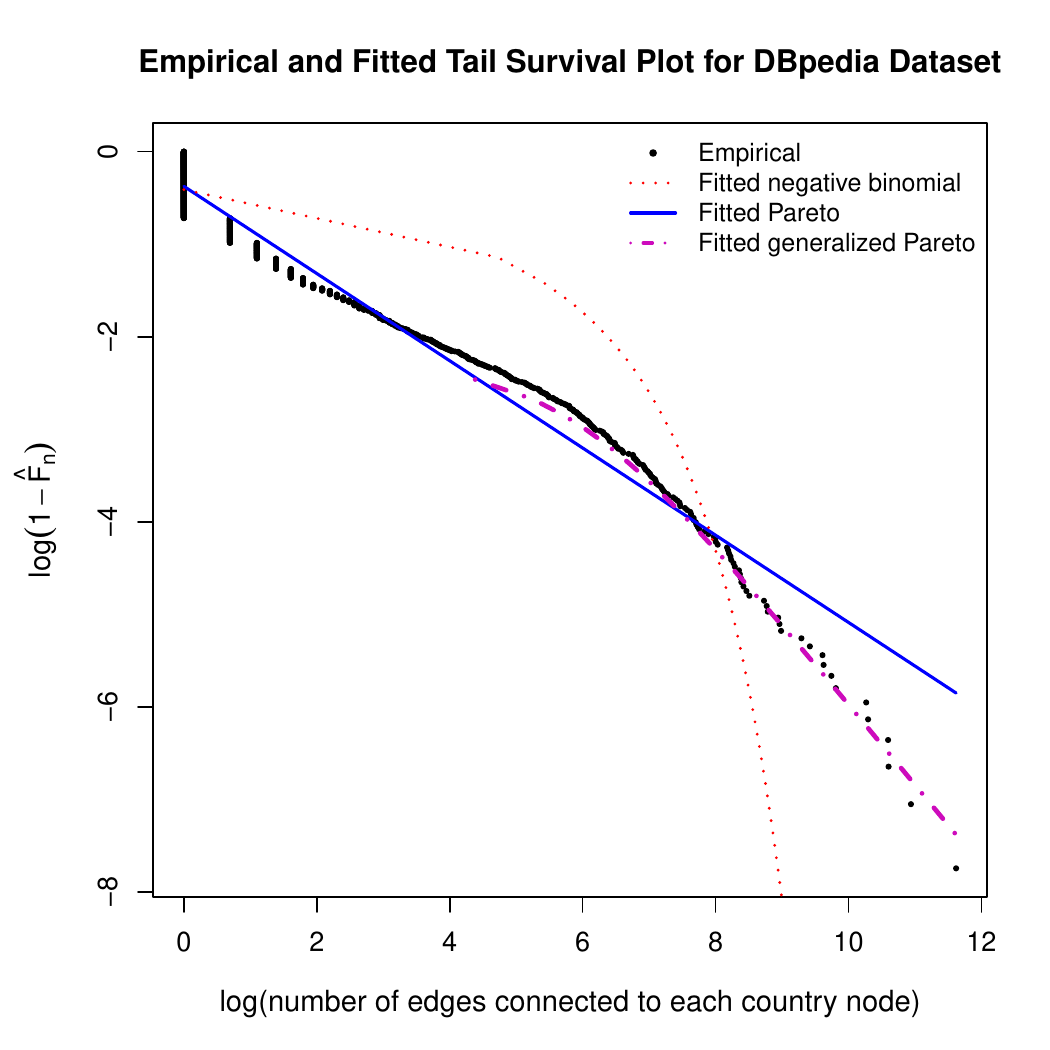}
        \label{fig:fig3}
    \end{subfigure}
    \caption{Plots of $\log(1 - \hat{F}_n)$, where $\hat{F}_n$ denotes the empirical distribution function or one of the fitted distributions (negative binomial, Pareto,  {
    or generalized Pareto}), for the Wiki Talk, Epinions, and DBpedia datasets. 
{
For the Wiki Talk, Epinions, and DBpedia datasets, the estimated tail indices (99\% CIs) based on the Pareto distribution are $0.661$ ($0.655$–$0.667$), $0.539$ ($0.535$–$0.543$), and $0.471$ ($0.443$–$0.512$), respectively, while the corresponding estimates based on the generalized Pareto distribution are $2.231$ ($1.748$–$2.675$), $2.584$ ($2.127$–$3.292$), and $0.786$ ($0.583$–$1.046$). The corresponding thresholds used in fitting the generalized Pareto distribution are $6.39$, $7.31$ and $4.38$ (in the log-scale).}
The confidence intervals of the estimated tail indices are obtained by 
applying a nonparametric bootstrap to the data.}

    \label{fig:three histograms}
\end{figure}

In all three datasets, the Pareto distribution fit the empirical data substantially better than the negative binomial distribution for moderate values of $x$,
but the data exhibit lighter tails than those implied by the fitted Pareto distribution.
The discrepancies between the empirical upper tails and the fitted Pareto upper tails,
though of uncertain origin,
are more likely to be due to heterogeneity and dependence in the network data than to the 
absence of a slowly varying factor from the Pareto distribution. 
In a regularly varying function, as we assume here, a slowly varying factor would be more influential for low to 
moderate values of the random variable and asymptotically negligible in affecting the slope for large values of the random variable.
Despite this uncertainty about the origin of the deviations from a Pareto distribution,
it is clear that the Pareto distribution approximates
the data better than the negative binomial distribution. 
{
In addition, the tail indices obtained from the Hill estimator and those from the Pareto distribution fitted by the least squares approach do not differ substantially.
The estimates are reported in the captions of Figures \ref{fig:wiki_hill}, \ref{fig:epinion hill}, \ref{fig:DBpedia_hill}, and \ref{fig:three histograms}.

To describe the right tail of the data better, we also applied the peaks-over-threshold approach and fitted a generalized Pareto distribution via maximum likelihood estimation using the \texttt{fitgpd()} function from the \texttt{POT} package in \texttt{R}. Threshold values were chosen approximately at the points where the fitted Pareto distribution no longer aligned well with the data and where the estimates appeared relatively stable across varying thresholds (figures not shown). The corresponding estimated survival functions for values above the threshold are plotted in Figure \ref{fig:three histograms}. We find that the generalized Pareto distribution provides a better fit to the right tail than the Pareto distribution.}
Because the right tail indices estimated by fitting the generalized Pareto
distribution to the Wiki Talk data and the Epinions data both exceed 2, the distributions
of nodal values are estimated to have finite population mean and finite population variance. 
Only for the DBpedia data does the tail index estimated by fitting the generalized Pareto
distribution fall below 1.

When the underlying data follow a power-law distribution with tail index $\alpha < 1$, 
the population mean and the population variance are infinite. 

Nonetheless, even in this case with tail index $\alpha < 1$, our theoretical results show that,
under our assumptions, the ratio of the log sample variance to the log sample mean converges to a finite limit that depends on $\alpha$, with the limit guaranteed to be at least 2. 
This ensures that the plot of log-variance versus log-mean remains meaningful even for infinite population moments, as shown first by \cite{brown2017taylor}. 
Moreover, as $\alpha \uparrow 1$, the limiting value of the slope increases, 
which explains the estimated slopes greater than 2 in our examples.

The observed scaling may provide a way to detect heterogeneity across subgraphs: if the log-variance versus log-mean relationship in a subset deviates notably from the overall trend, it may indicate a distinct structural pattern. Additionally, the scaling exponent may serve as a diagnostic tool for evaluating generative models of network data, as not all models reproduce the level of dispersion observed in these examples.
}

\section{Discussion}\label{sec:discussion}
We present a new framework for establishing the limits in probability 
of the ratios of moments that arise in various forms
of Taylor's law when the data can be heavy-tailed
with tail index $\alpha \in (0,1)$, 
dependent, and heterogeneous. 
We explore Taylor's law for higher (centered) moments and semivariances, along with its most commonly known form (\ref{eq:TL}). 
We introduce a covariance condition (Condition A($p$)) that leads to Taylor's law for dependent data, and we propose
mixing conditions that satisfy this covariance condition. 
We  extend our results to heterogeneous and correlated heavy-tailed dependent data and  
to an arbitrary index set,  illustrated by network-dependent data.
 To our knowledge, networks are a previously  unexplored setting for Taylor's law.
 The theoretical results are supported with three empirical examples and extensive simulations.

A future question  to consider is  whether the assumption that data are $i.i.d.$ 
when $\alpha \in (1,2)$, which is considered in \cite{cohen2022covid}, can 
be relaxed to allow weak dependence. 
In addition, more sophisticated network structures from real-world social network data
can be used to test the descriptive ability of Taylor's law in new domains.

\section{Disclosure Statement}\label{disclosure-statement}
No conflicts of interest exist.

\section{Data Availability Statement}\label{data-availability-statement}

Data used in this article are available at the following URLs: 
\begin{itemize}
    \item Stanford Network Analysis Project (\url{https://snap.stanford.edu/data/wiki-Talk.html})
    \item Network Data Repository with Interactive Graph Analytics and Visualization (\url{https://networkrepository.com/rec_epinions_user_ratings.php})
    
    \item Netzschleuder (\url{https://networks.skewed.de/net/dbpedia_country})
\end{itemize}

\setcounter{section}{0}

\renewcommand{\thesection}{\Alph{section}}          
\renewcommand{\thesubsection}{\thesection.\arabic{subsection}}

\numberwithin{equation}{section}                    
\renewcommand{\theequation}{\thesection.\arabic{equation}}

\section{Taylor's Law for Semivariances, Central Lower and Upper Moments of Heavy-tailed Data}\label{sec:semivariances}
In this section, we extend and discuss  Taylor's law for semivariances and central lower and upper moments of \emph{dependent} data; see \cite{brown2021taylor} for more discussion. 
The semivariance measures only the risk of unfavorable deviations from the mean, while the variance considers equally both upside and downside deviations from the mean.
Central lower and upper moments generalize the semivariance to powers other than $2$ of deviations from the mean.

Denote $a_+ := \max\{a,0\}$. Let $N^-_n := \sum^n_{i=1}\mathbbm{1}(X_i \leq M_{n,1})$ and $N^+_n := \sum^n_{i=1}\mathbbm{1}(X_i > M_{n,1})$.
For any $h > 0$, the sample $h$th central lower moments and sample $h$th central local lower moments are defined as
\begin{equation*}
	M^-_{n,h} := \frac{1}{n}\sum^n_{i=1}[(M_{n,1}- X_i)_+]^h \quad \text{ and } \quad M^{-*}_{n,h} := \frac{1}{N^-_n}\sum^n_{i=1}[(M_{n,1} - X_i)_+]^h.
\end{equation*}
Similarly, the sample $h$th central upper moments and sample $h$th central local upper moments are defined as
\begin{equation*}
	M^+_{n,h} := \frac{1}{n}\sum^n_{i=1}[(X_i - M_{n,1})_+]^h \quad \text{ and } \quad M^{+*}_{n,h} := \frac{1}{N^+_n}\sum^n_{i=1}[(X_i - M_{n,1})_+]^h.
\end{equation*}
When $h = 2$, $M^-_{n,2}, M^{-*}_{n,2}, M^+_{n,2}$, and $M^{+*}_{n,2}$ are also called the sample lower semivariance, sample local lower semivariance, sample upper semivariance and sample local upper semivariance, respectively.

\subsection{Central Lower and Upper Moments}\label{subsect:lower_upper}

We recall that $\{c_n\}$ is a positive sequence such that $nc_n \uparrow \infty$ and $\log c_n = o(\log n)$ and
that $\tilde{X} := X\mathbbm{1}(X < t_n)$.
Define $d_{n,1} := \mathbb{E}(\tilde{X}_1)$. For $\delta'>0$, define $b_n := {d_{n,1}}/{c_n^{2\delta'}}$ and $\tilde{b}_n := d_{n,1} c_n$.

\begin{theorem}[$h$th Central Lower Moments]\label{thm:hth_central_lower}
	Let $X_1,\ldots,X_n,\ldots$ be a sequence of 
	nonnegative random variables with a common survival function $\overline{F}(x) = x^{-\alpha}l(x)$, where $l$ is slowly varying and $\alpha \in (0, 1)$. Suppose that Condition A($p$) holds for $p = 1$.
	\begin{enumerate}[(i)]
		\item 
		Suppose that for some $a > 0$, \begin{equation}\label{eq:ass_lower_semi}
			\sum_{i\neq j} \Cov(\mathbbm{1}(X_i \leq a), \mathbbm{1}(X_j \leq a)) = o(n^2).
		\end{equation}
		Then, for any $h > 0$,
		\begin{equation*}
			\frac{\log M^-_{n,h}}{\log M_{n,1}} = h + O_p\left(\frac{1}{\log n}\right).
		\end{equation*}
		\item For some $\delta' > 2(1/\alpha-1) > 0$,
		suppose that
		\begin{equation}\label{eq:ass_lower_semi_seq}
			\sum_{i\neq j} \Cov(\mathbbm{1}(X_i \leq b_n), \mathbbm{1}(X_j \leq b_n)) = o(n^2).
		\end{equation}
		Then, for any $h > 0$, 
		\begin{equation*}
			\frac{M^-_{n,h}}{M^h_{n,1}} \stackrel{\mathbb{P}}{\rightarrow} 1.
		\end{equation*}
	\end{enumerate}
\end{theorem}

\begin{remark}
	In (\ref{eq:ass_lower_semi}),    \begin{equation*}
		\Cov(\mathbbm{1}(X_i \leq a), \mathbbm{1}(X_j \leq a)) = \mathbb{P}(X_i \leq a, X_j \leq a) - \mathbb{P}(X_i \leq a) \mathbb{P}(X_j \leq a).
	\end{equation*}
	This condition is ``weaker'' than strong mixing as we only require the condition to hold for one value of $a$ that is greater than the lower boundary of the support of $X$. 
	The condition  guarantees that $\frac{1}{n} \sum^n_{i=1}\mathbbm{1}(X_i \leq a)$ is bounded away from $0$ in probability.
\end{remark}

\begin{theorem}[$h$th Central Upper Moments]\label{thm:upper_central_moment}
	Let $X_1,\ldots,X_n,\ldots$ be a sequence of 
	nonnegative random variables with a common survival function $\overline{F}(x) = x^{-\alpha}l(x)$, where $l$ is slowly varying and $\alpha \in (0, 1)$. 
	\begin{enumerate}[(a)]
		\item 
		If Condition A($p$) holds for $p = h$ with  $\alpha < h < 1$, and if (\ref{eq:ass_lower_semi_seq}) also holds for some $\delta' > 2(1/\alpha-1) > 0$, then
		\begin{equation*}
			\frac{M_{n,h}^+}{M_{n,1}^h} = o_p(1).
		\end{equation*}
		\item Suppose that Condition A($p$) holds for $p = h$ with  $h >1$.
		\begin{enumerate}[(i)]
			\item Then
			\begin{equation*}
				R_{n,h,1}^{\text{cum}} :=   \frac{\log M_{n,h}^+}{\log M_{n,1}} - \imath(h, 1) =  O_p\left( \frac{\log c_n}{\log n}\right) + 
				O_p\left( \frac{|\log l(t_n)|}{\log n}\right).
			\end{equation*}

			\item 
			If $\lim_{x \rightarrow \infty}l(x) = \infty$, then
			$R_{n,h_1,h_2}^{\text{cum}} = O_p(s_{1n,l})$ and $(R_{n,h_1,h_2}^{\text{cum}})^{-1} = O_p(s_{2n,l}^{-1})$.
			\item
			If $\lim_{x \rightarrow \infty}l(x) = 0$, then $R_{n,h_1,h_2}^{\text{cum}} = O_p(s_{2n,l})$ and $(R_{n,h_1,h_2}^{\text{cum}})^{-1} = O_p(s_{1n,l}^{-1})$.
		\end{enumerate}
	\end{enumerate}
\end{theorem}

\subsection{Local Central Lower and Upper Moments}\label{subsect:local_lower_upper}

\begin{theorem}[$h$th Local Central Lower  Moments]\label{thm:local_lower_central_moment}
	Let $X_1,\ldots,X_n,\ldots$ be a sequence of nonnegative random variables with a common survival function $\overline{F}(x) = x^{-\alpha}l(x)$, where $l$ is slowly varying and $\alpha \in (0, 1)$.         
	If Condition A($p$) for $p = 1$ and (\ref{eq:ass_lower_semi}) hold, then for any $h > 0$,
	\begin{equation*}
		\frac{\log M^{-*}_{n,h}}{\log M_{n,1}} = h + O_p\left( \frac{1}{\log n}\right).
	\end{equation*}
\end{theorem}

To establish Taylor's law for the $h$th local central upper moment, we impose an additional condition for the truncated terms:
\begin{equation}\label{eq:cov_ass3}
	\sum_{i\neq j} \Cov(\mathbbm{1}(X_i > \tilde{b}_n), \mathbbm{1}(X_j > \tilde{b}_n)) = o\left(n^{2\alpha} \left(\frac{l(\tilde{b}_n)}{l(t_n) c_n}\right)^2 \right).
\end{equation}

\begin{theorem}[$h$th Local Upper Central Moments]\label{thm:local_upper_central_moment}
	Let $X_1,\ldots,X_n,\ldots$ be a sequence of nonnegative random variables with a common survival function $\overline{F}(x) = x^{-\alpha}l(x)$, where $l$ is slowly varying and $\alpha \in (0, 1)$.         Suppose that Condition A($p$) for $p = 1$ and  (\ref{eq:cov_ass3}) hold. 
	For $h > 1$, define  $\imath_+(h) := {(h-\alpha^2)}/{(1-\alpha)}$. Let
	\begin{equation*}
		L_n := \log c_n  + |\log l(t_n)| + |\log l(b_n)| + |\log l(\tilde{b}_n)|.
	\end{equation*}
	\begin{enumerate}[(a)]
		\item Then
		\begin{align*}
			R_{n,h,1}^{\text{lucm}} :=   \frac{\log M^{+*}_{n,h}}{\log M_{n,1}} - \imath_+(h) = O_p\left( \frac{L_n}{\log n}\right).
		\end{align*}
		
		\item Furthermore, suppose that $l(b_n)/l(t_n) \rightarrow c_1 $ and $l(\tilde{b}_n)/l(t_n) \rightarrow c_2$ for some constants $c_1, c_2 > 0$. 
		
		\begin{enumerate}[(i)]
			\item  If $\lim_{x \rightarrow \infty}l(x) = \infty$, then $R_{n,h_1,h_2}^{\text{lucm}} = O_p(s_{1n,l})$ and $(R_{n,h_1,h_2}^{\text{lucm}})^{-1} = O_p(s_{2n,l}^{-1})$.
			
			\item If $\lim_{x \rightarrow \infty}l(x) = 0$, then  $R_{n,h_1,h_2}^{\text{lucm}} = O_p(s_{2n,l})$ and $(R_{n,h_1,h_2}^{\text{lucm}})^{-1} = O_p(s_{1n,l}^{-1})$.
		\end{enumerate}
		
	\end{enumerate}
\end{theorem}

\begin{remark}
	For some choices of $l$, we have $l(b_n)/l(t_n) \rightarrow c_1 > 0$ and $l(\tilde{b}_n)/l(t_n) \rightarrow c_2 > 0$. 
	For example, consider $l(t) = \log t$.  
	As $n \rightarrow \infty$,
	\begin{align*}
		d_{n,1}c_n \sim C (t_n^\alpha)^{1/\alpha - 1} l(t_n) = C n^{1/\alpha - 1} c_n^{1/\alpha} l(t_n)^{1/\alpha}.
	\end{align*}
	Therefore, as $n \rightarrow \infty$,
	\begin{align*}
		\frac{l(\tilde{b}_n)}{l(t_n)} = \frac{\log \tilde{b}_n}{\frac{1}{\alpha} \log t_n^\alpha} = \frac{o(1) + \log C + (\frac{1}{\alpha}-1)\log n+ \frac{1}{\alpha} \log c_n + \frac{1}{\alpha} \log l(t_n) }{\frac{1}{\alpha}(\log n + \log l(t_n) + \log c_n)} \rightarrow \frac{\frac{1}{\alpha}-1}{\frac{1}{\alpha}} = 1-\alpha.
	\end{align*}
	Similarly, as $n \rightarrow\infty$,
	\begin{equation*}
		\frac{l(b_n)}{l(t_n)} \rightarrow 1-\alpha.
	\end{equation*}
\end{remark}

\section{Simulation Studies}\label{sec:simulation}
The following simulations illustrate the results discussed in the main text. We use the \texttt{R} packages 
\texttt{EnvStats}, \texttt{stabledist} to simulate heavy-tailed distributions and \texttt{ggplot2} to visualize the results.

In each of the following figures, we simulate  samples of size $n = 10^p$ for $p=1,\ldots,7$. 
In panels (a)--(d) of each figure, we  plot the locations of (a) $(\log(M_{n,1}),\log(V_n))$ for variance, (b) $(\log(M_{n,1}),\log(M^{+}_{n,2}))$ for upper semivariance, (c) $(\log(M_{n,1}),\log(M^c_{n,3}))$ for third central moment, and (d) $(\log(M_{n,1}),\log(M^{+*}_{n,2}))$ for local upper semivariance.
In each plot,  three occurrences of $1,\ 2,\ \ldots,\ 7$ represent the location of the ratios with sample size $10^1,\ 10^2,\ \ldots,\ 10^7$ in three different samples, respectively. 

The dashed straight line in each plot has a slope equal to  the theoretical limit and an intercept estimated from the samples when the slope is fixed to be the theoretical limit. 
From the logarithmic form of Taylor's law in (\ref{eq:TL}),
it is evident that  the slopes represent the asymptotic ratio
of the moments being plotted. 
For example, in each panel (a), where $\frac{\log(V_n))}{\log(M_{n,1})} \stackrel{\mathbb{P}}{\rightarrow} \frac{2-\alpha}{1-\alpha}$,  the slope is set to be $b=\frac{2-\alpha}{1-\alpha}$. Theoretically, the points in each graph should lie close to the dashed straight line.

\subsection{Taylor's Law for $i.i.d.$ Data}
Simulations of  $i.i.d.$ heavy-tailed random variables, generated by three different distributions with three different choices of $\alpha$ each,
illustrate  the asymptotic results 
of Theorem \ref{thm:higher_central_moments} with $h_1=2$ and $h_2=1$, Theorem \ref{thm:upper_central_moment}, Theorem \ref{thm:higher_central_moments} with $h_1=3$ and $h_2=1$, and Theorem \ref{thm:local_upper_central_moment},
in the following nine figures.

We observe that the points lie closer to the dashed line  when $\alpha=0.2$
than when $\alpha=0.5$ and 0.7. 
Moreover, the points with ``7'' are also generally closer to the line than the other points, as they are based on the largest sample size.

\begin{enumerate}
	
	\item 
	Suppose that $X_1,\ldots,X_n$ are $i.i.d.$ Pareto random variables with location parameter $x_m = 0.5$ and scale parameter $\alpha$, denoted by Pareto($0.5,\alpha$), with survival function $\mathbb{P}(X_i>x)=x^{-\alpha}$ for $i=1,\ldots,n$, $0<\alpha<1$,  and $x\geq x_m$. 
	{
		Figures \ref{fig:pareto_iid_1}, \ref{fig:pareto_iid_2}, and \ref{fig:pareto_iid_3} show the results 
		when $\alpha=0.2, 0.5,$ and $0.8$, respectively. We employ the function \texttt{rpareto} in the \textbf{R} package \textbf{EnvStats} to simulate Pareto random variables.}

	\begin{figure}[h]
		\centering
		\subfloat[Taylor's law]{\includegraphics[width=0.3\textwidth]{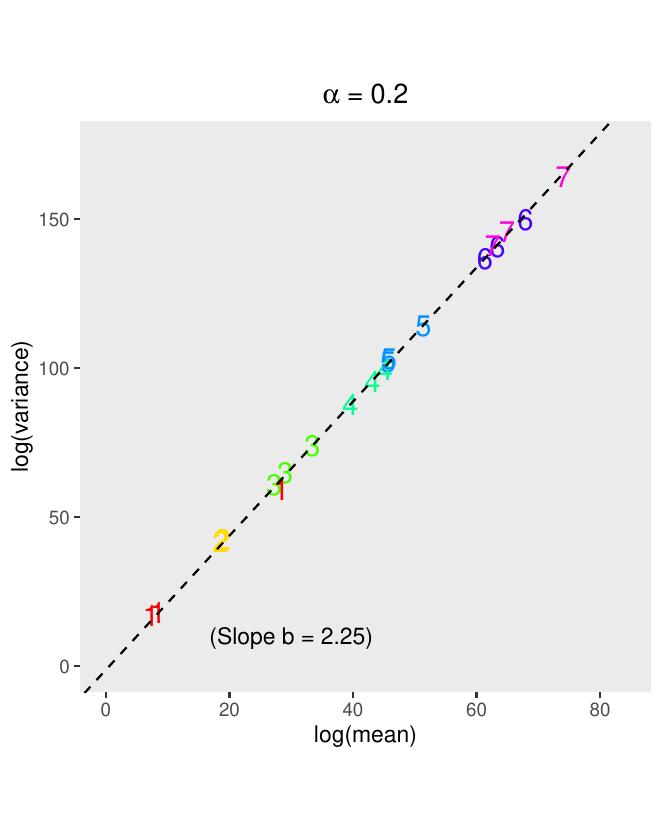}}
		\qquad
		\subfloat[Taylor's law for upper semivariance]{\includegraphics[width=0.3\textwidth]{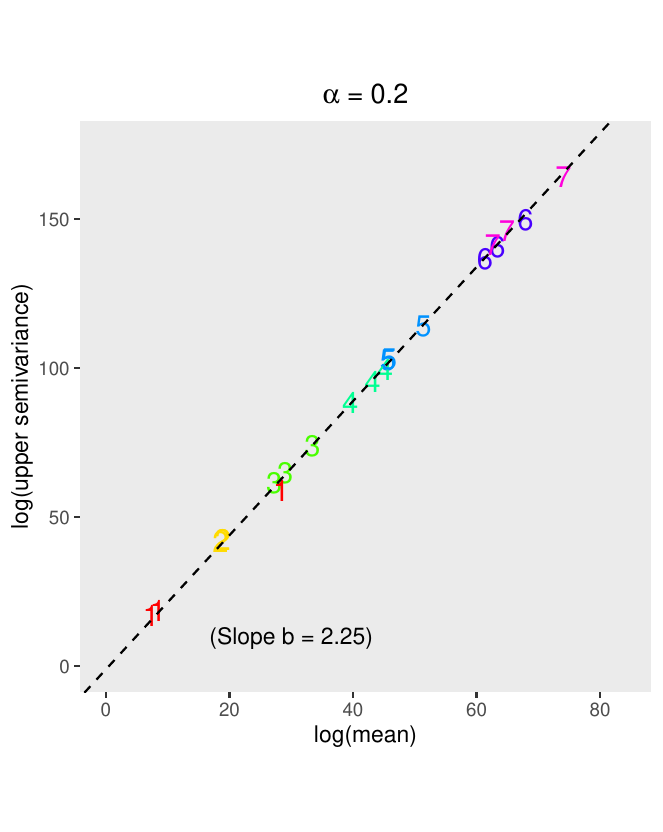}}
		\\
		\subfloat[Taylor's law for third central moment]{\includegraphics[width=0.3\textwidth]{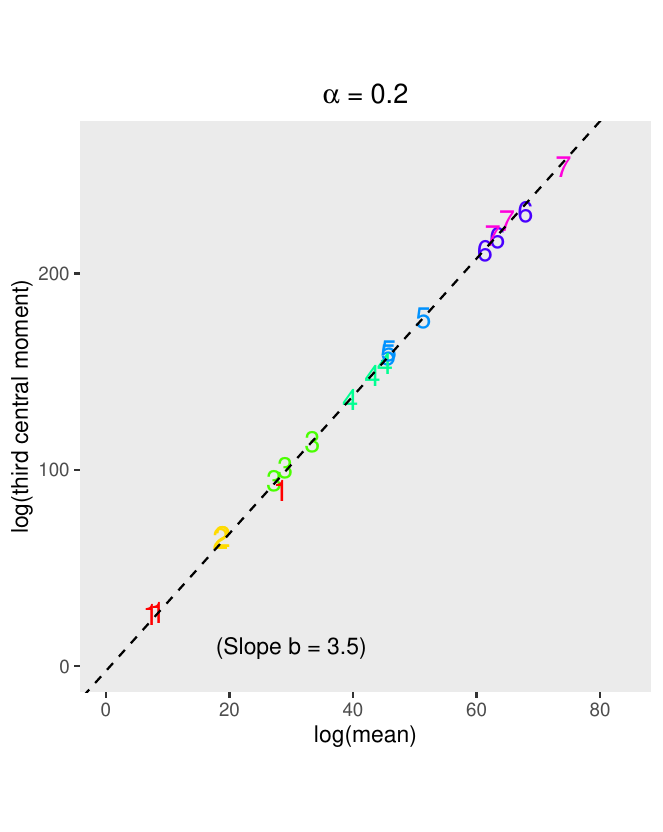}}
		\qquad
		\subfloat[Taylor's law for local upper semivariance]{\includegraphics[width=0.3\textwidth]{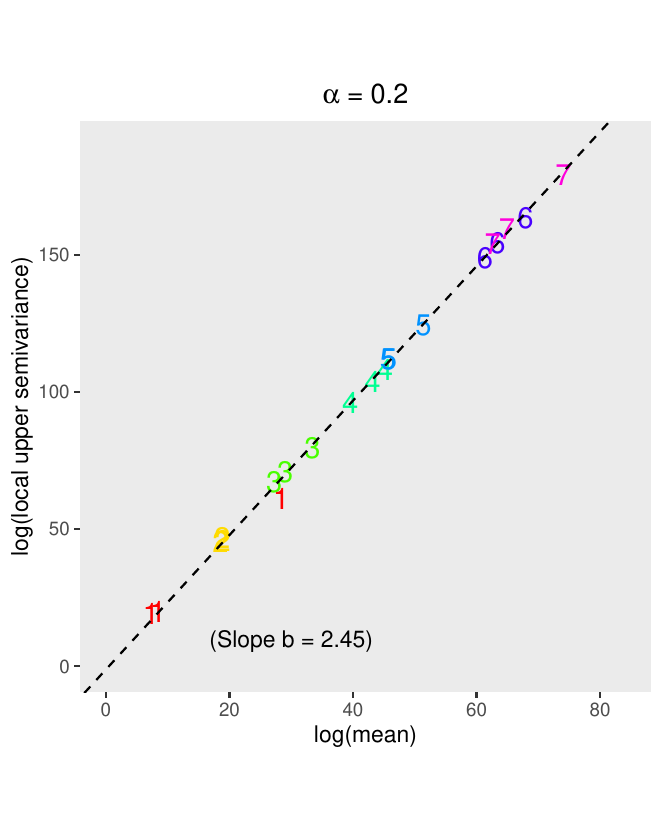}}
		
		\caption{Scatterplots of (a) $(\log(M_{n,1}),\log(V_n))$, (b) $(\log(M_{n,1}),\log(M^{+}_{n,2}))$, (c) $(\log(M_{n,1}),\log(M^c_{n,3}))$, and (d) $(\log(M_{n,1}),\log(M^{+*}_{n,2}))$, corresponding to Theorem \ref{thm:higher_central_moments} with $h_1=2$ and $h_2=1$, Theorem \ref{thm:upper_central_moment}, Theorem \ref{thm:higher_central_moments} with $h_1=3$ and $h_2=1$, and Theorem \ref{thm:local_upper_central_moment}, respectively, when $X_i \overset{\mathrm{i.i.d.}}{\sim}$   Pareto$(0.5,0.2)$ with sample size $10^p$ for $p=1,2,\ldots,7$.}
		\label{fig:pareto_iid_1}
	\end{figure}
	
	\begin{figure}[h]
		\centering
		\subfloat[Taylor's law]{\includegraphics[width=0.3\textwidth]{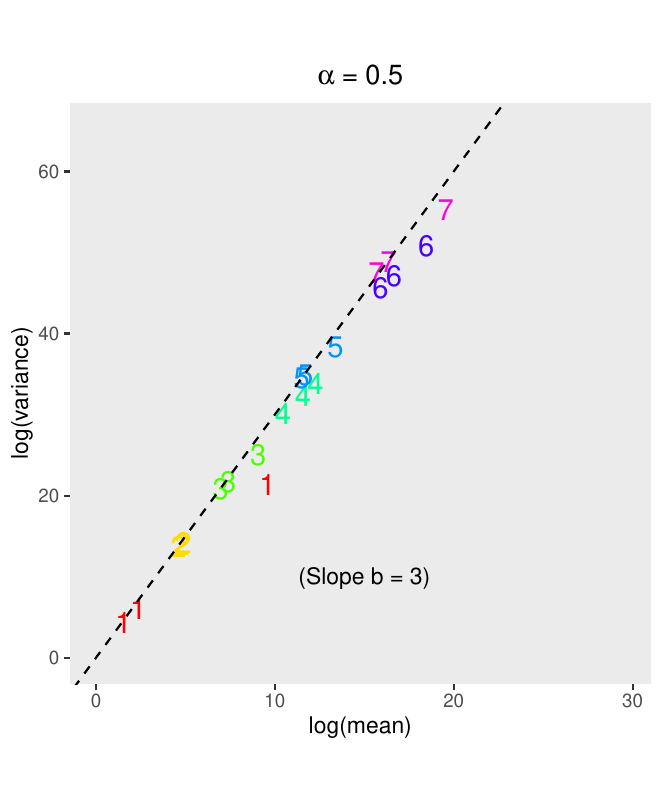}}
		\qquad
		\subfloat[Taylor's law for upper semivariance]{\includegraphics[width=0.3\textwidth]{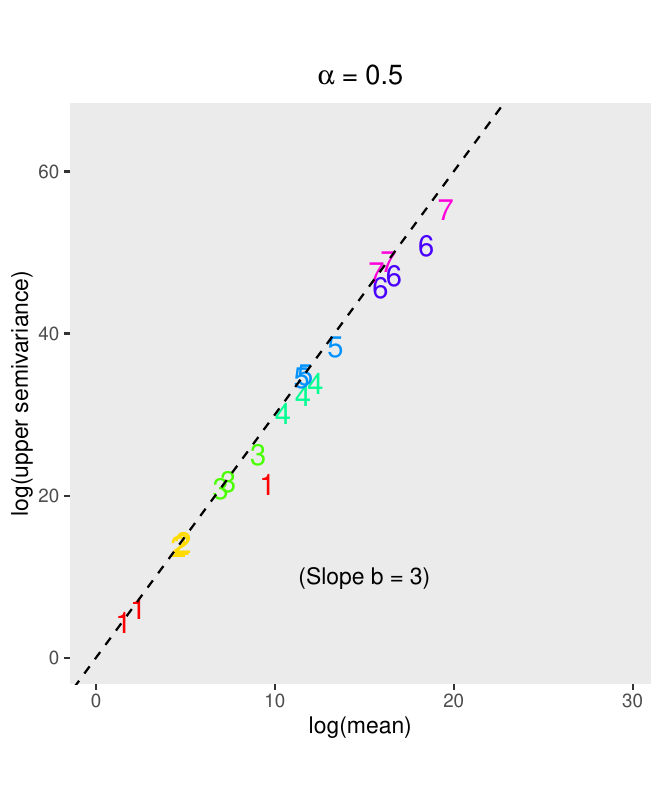}}
		\\
		\subfloat[Taylor's law for third central moment]{\includegraphics[width=0.3\textwidth]{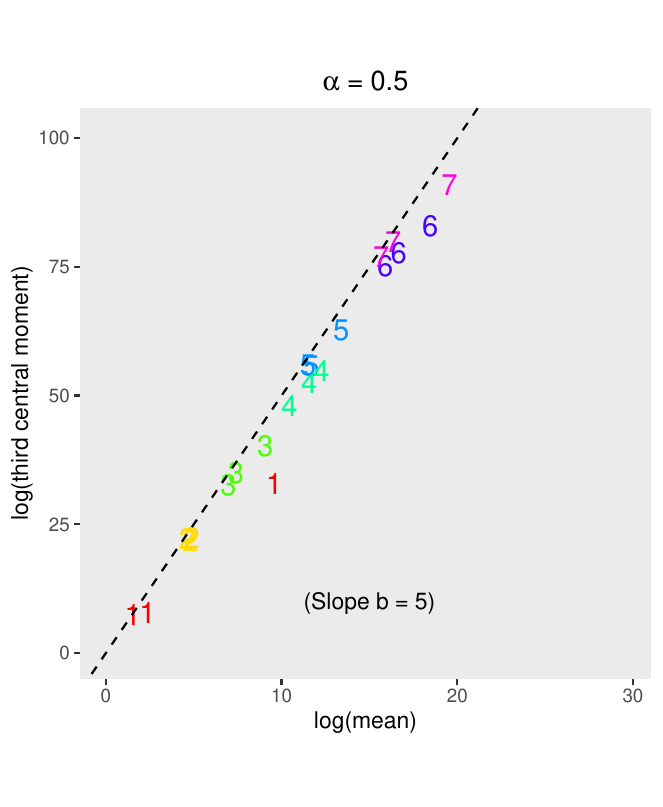}}
		\qquad
		\subfloat[Taylor's law for local upper semivariance]{\includegraphics[width=0.3\textwidth]{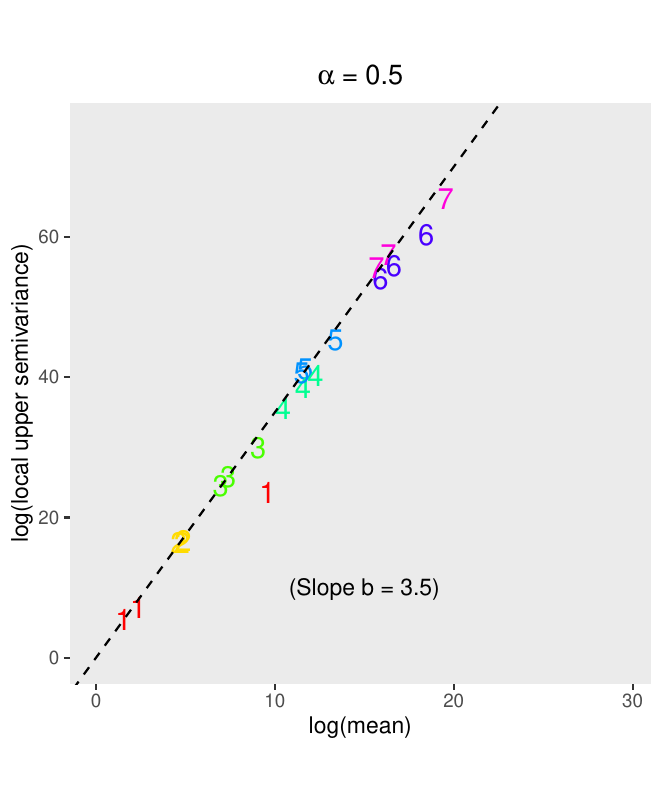}}

		\caption{Scatterplots of (a) $(\log(M_{n,1}),\log(V_n))$, (b) $(\log(M_{n,1}),\log(M^{+}_{n,2}))$, (c) $(\log(M_{n,1}),\log(M^c_{n,3}))$, and (d) $(\log(M_{n,1}),\log(M^{+*}_{n,2}))$, corresponding to Theorem \ref{thm:higher_central_moments} with $h_1=2$ and $h_2=1$, Theorem \ref{thm:upper_central_moment}, Theorem \ref{thm:higher_central_moments} with $h_1=3$ and $h_2=1$, and Theorem \ref{thm:local_upper_central_moment}, respectively, when $X_i \overset{\mathrm{i.i.d.}}{\sim}$ Pareto$(0.5,0.5)$ with sample size $10^p$ for $p=1,2,\ldots,7$.}
		\label{fig:pareto_iid_2}
	\end{figure}

	\begin{figure}[h]
		\centering
		\subfloat[Taylor's law]{\includegraphics[width=0.3\textwidth]{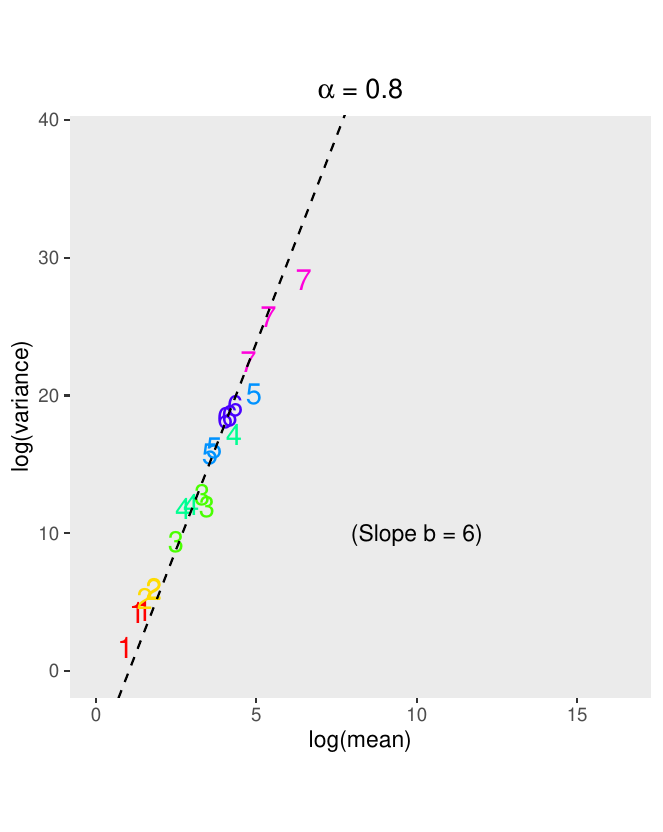}}
		\qquad
		\subfloat[Taylor's law for upper semivariance]{\includegraphics[width=0.3\textwidth]{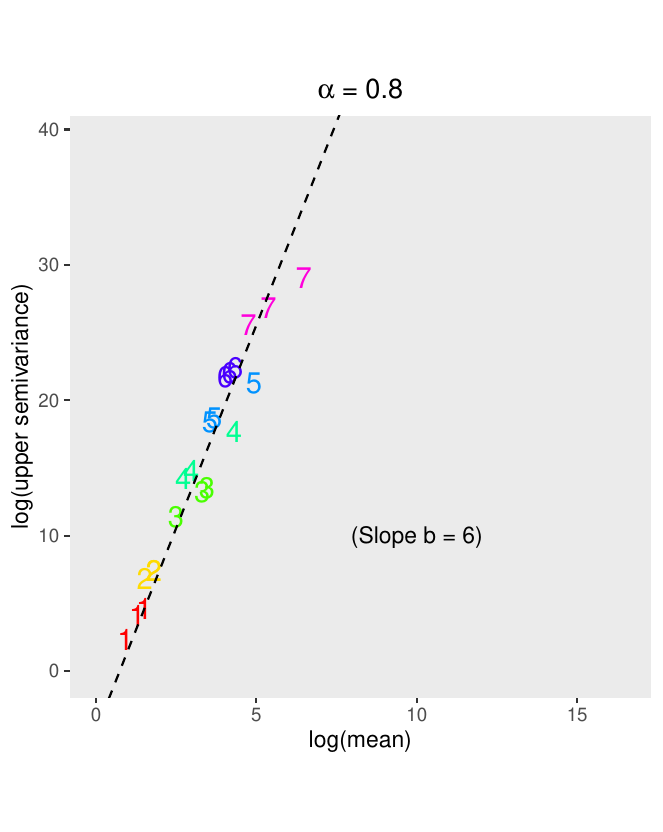}}
		\\
		\subfloat[Taylor's law for third central moment]{\includegraphics[width=0.3\textwidth]{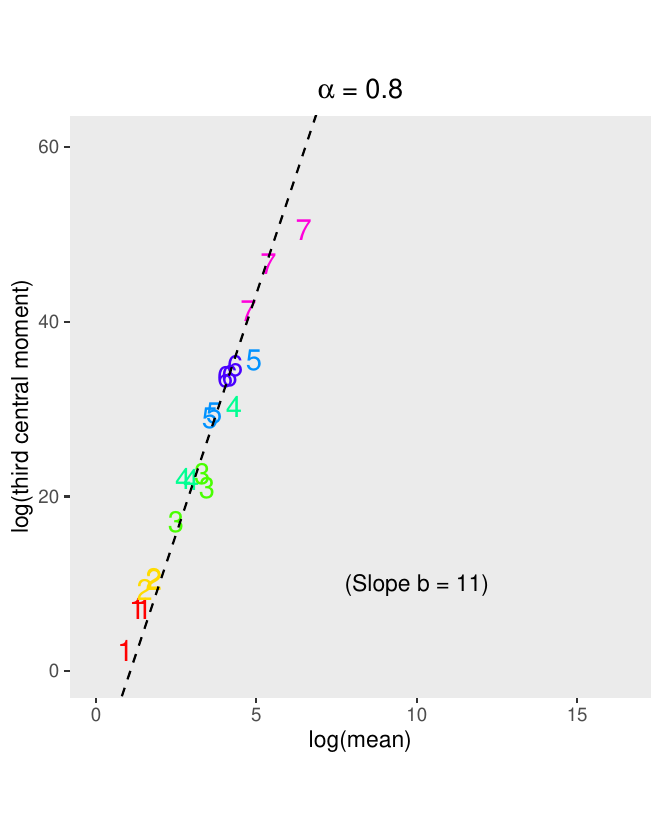}}
		\qquad
		\subfloat[Taylor's law for local upper semivariance]{\includegraphics[width=0.3\textwidth]{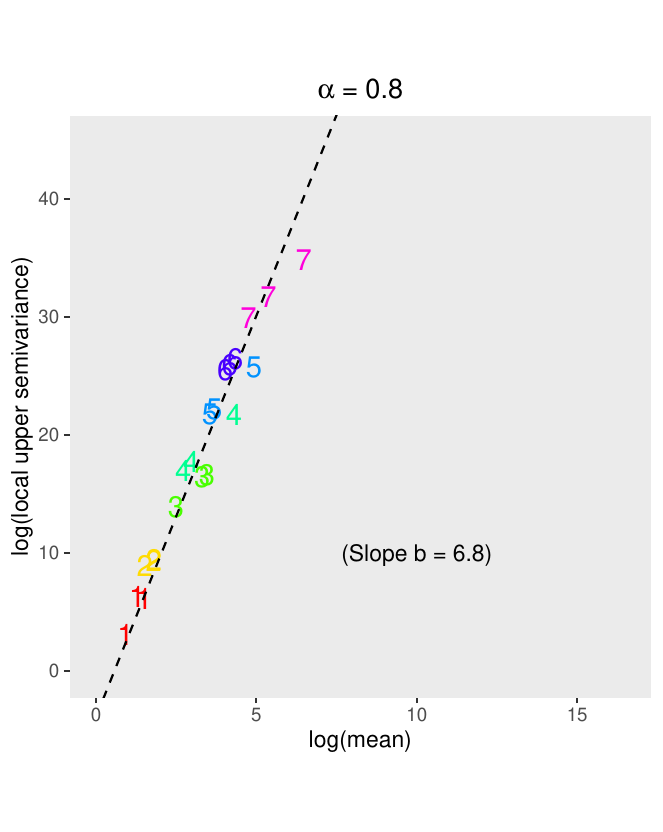}}
		
		\caption{Scatterplots of (a) $(\log(M_{n,1}),\log(V_n))$, (b) $(\log(M_{n,1}),\log(M^{+}_{n,2}))$, (c) $(\log(M_{n,1}),\log(M^c_{n,3}))$, and (d) $(\log(M_{n,1}),\log(M^{+*}_{n,2}))$, corresponding to Theorem \ref{thm:higher_central_moments} with $h_1=2$ and $h_2=1$, Theorem \ref{thm:upper_central_moment}, Theorem \ref{thm:higher_central_moments} with $h_1=3$ and $h_2=1$, and Theorem \ref{thm:local_upper_central_moment}, respectively, when $X_i \overset{\mathrm{i.i.d.}}{\sim}$ Pareto$(0.5,0.8)$ with sample size $10^p$ for $p=1,2,\ldots,7$.}
		\label{fig:pareto_iid_3}
	\end{figure}
	
	\item{
		As a special case of the survival function $\overline{F}(x)=x^{-\alpha}l(x)$,
		where $l$ is slowly varying and $\alpha>0$,
		we define the cumulative distribution function
		\begin{align}\label{eq:F1}
			F_{1,\alpha}(x)=1-x^{-\alpha}\log(x) e\alpha, \quad x \geq \exp(\alpha^{-1}).
		\end{align}
		Let $x_{m,\alpha} := \exp(\alpha^{-1})$.
		We used the acceptance-rejection method \citep{ross2022simulation} to simulate random variables with the distribution function $F_{1,\alpha}$. Specifically, a Pareto distribution with the same location parameter $x_{m,\alpha}$ and a shape parameter of $\alpha -0.05$ (i.e., Pareto($x_{m,\alpha}$, $\alpha - 0.05$)) was used as the proposal distribution. This ensures that the two distributions have the same support, while the proposal distribution has a heavier tail than $F_{1,\alpha}$.}

	Figures \ref{fig:inverse_iid_1}, \ref{fig:inverse_iid_2}, and \ref{fig:inverse_iid_3} show the results
	when $\alpha=0.2, 0.5,$ and $0.8$, respectively.
	\begin{figure}[h]
		\centering
		\subfloat[Taylor's law]{\includegraphics[width=0.3\textwidth]{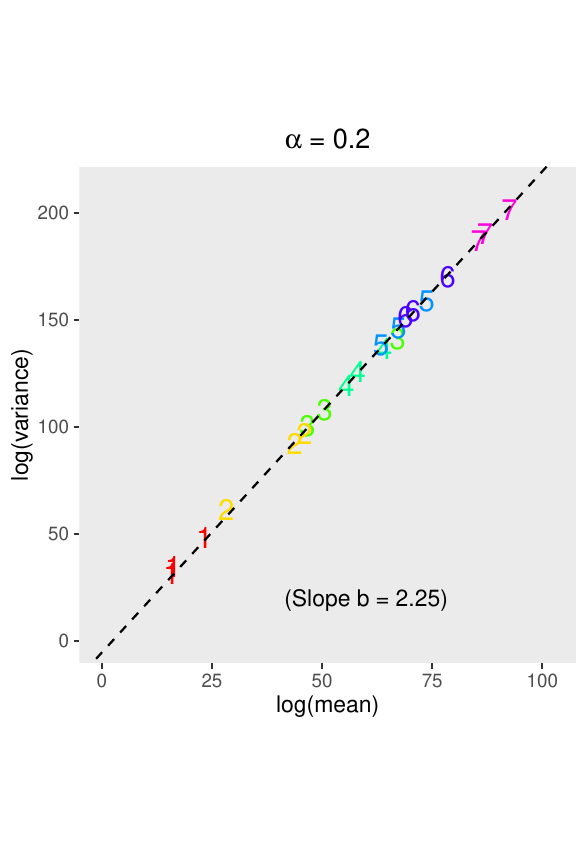}}
		\qquad
		\subfloat[Taylor's law for upper semivariance]{\includegraphics[width=0.3\textwidth]{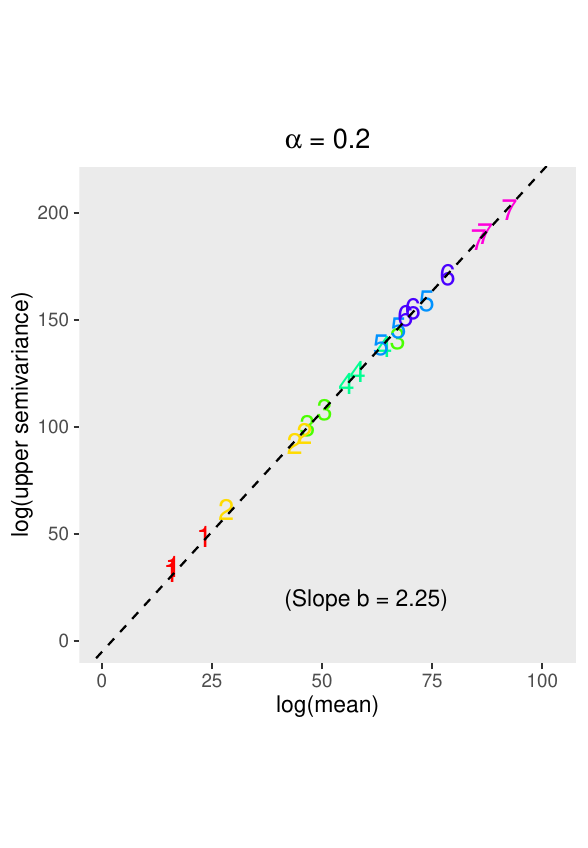}}
		\\
		\subfloat[Taylor's law for third central moment]{\includegraphics[width=0.3\textwidth]{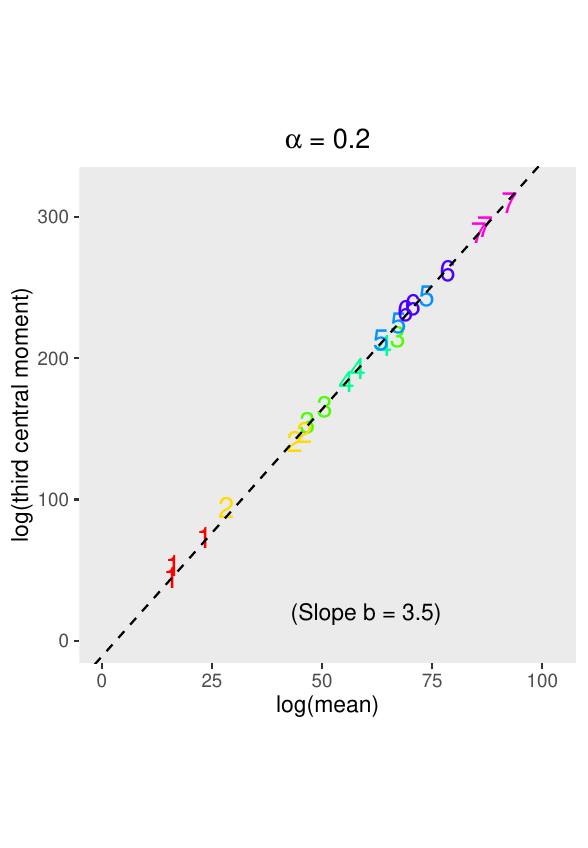}}
		\qquad
		\subfloat[Taylor's law for local upper semivariance]{\includegraphics[width=0.3\textwidth]{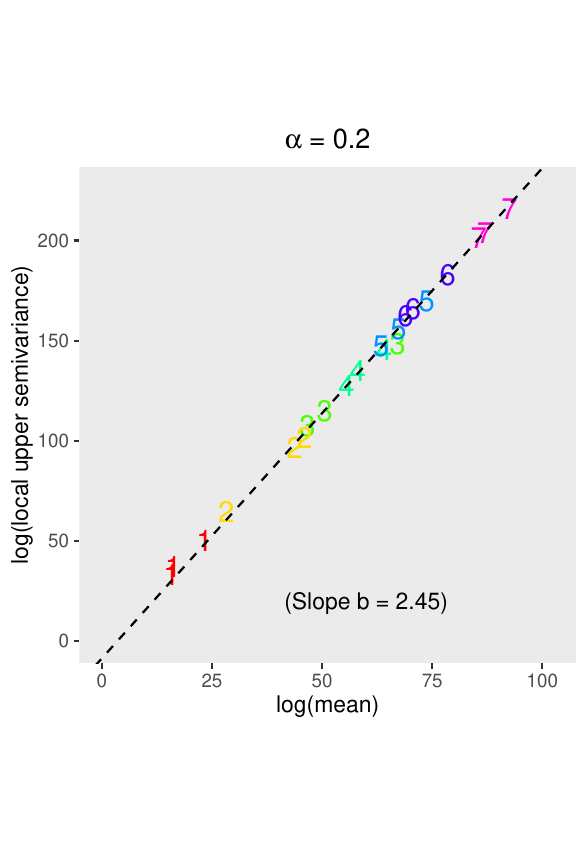}}
		
		\caption{Scatterplots of (a) $(\log(M_{n,1}),\log(V_n))$, (b) $(\log(M_{n,1}),\log(M^{+}_{n,2}))$, (c) $(\log(M_{n,1}),\log(M^c_{n,3}))$, and (d) $(\log(M_{n,1}),\log(M^{+*}_{n,2}))$, corresponding to Theorem \ref{thm:higher_central_moments} with $h_1=2$ and $h_2=1$, Theorem \ref{thm:upper_central_moment}, Theorem \ref{thm:higher_central_moments} with $h_1=3$ and $h_2=1$, and Theorem \ref{thm:local_upper_central_moment}, respectively, when $F_{1,0.2}$ 
			in \eqref{eq:F1} is the distribution function of $i.i.d.$ $X_i$ with sample size $10^p$ for $p=1,2,\ldots,7$.}
		\label{fig:inverse_iid_1}
	\end{figure}
	
	\begin{figure}[h]
		\centering
		\subfloat[Taylor's law]{\includegraphics[width=0.3\textwidth]{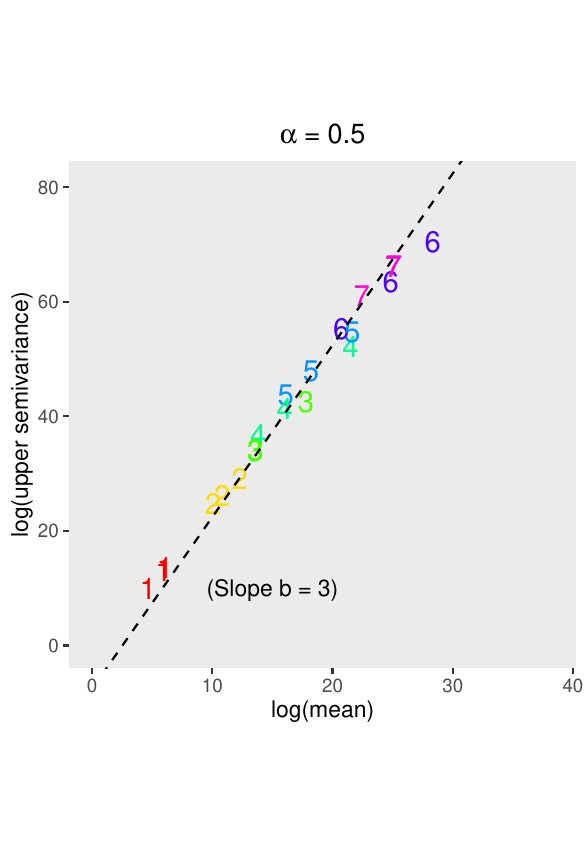}}
		\qquad
		\subfloat[Taylor's law for upper semivariance]{\includegraphics[width=0.3\textwidth]{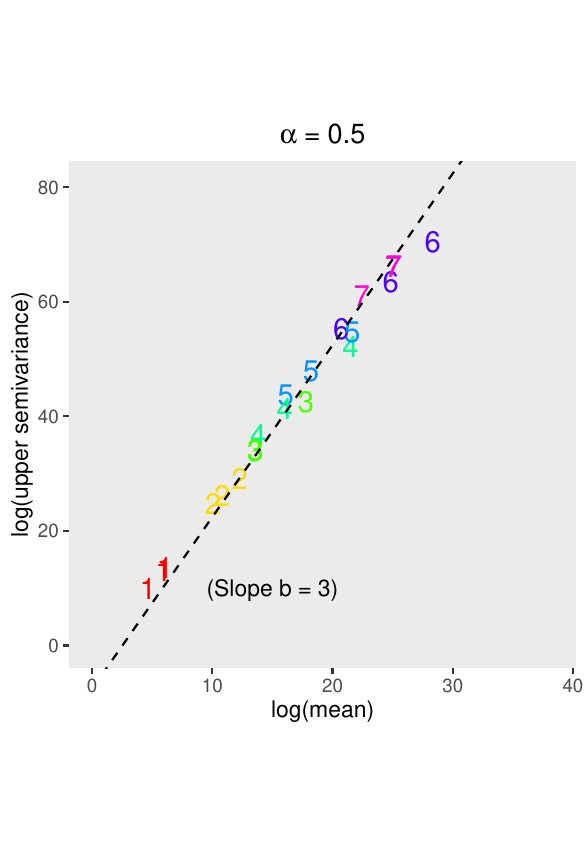}}
		\\
		\subfloat[Taylor's law for third central moment]{\includegraphics[width=0.3\textwidth]{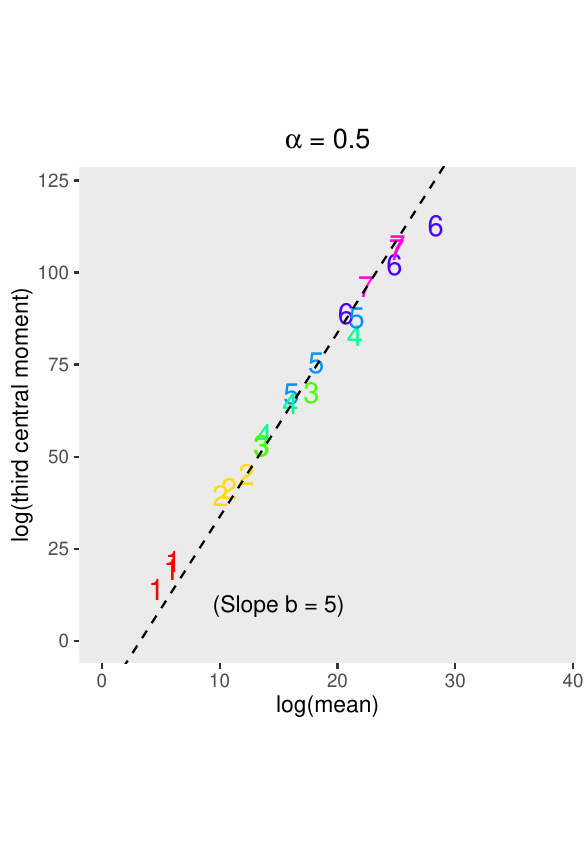}}
		\qquad
		\subfloat[Taylor's law for local upper semivariance]{\includegraphics[width=0.3\textwidth]{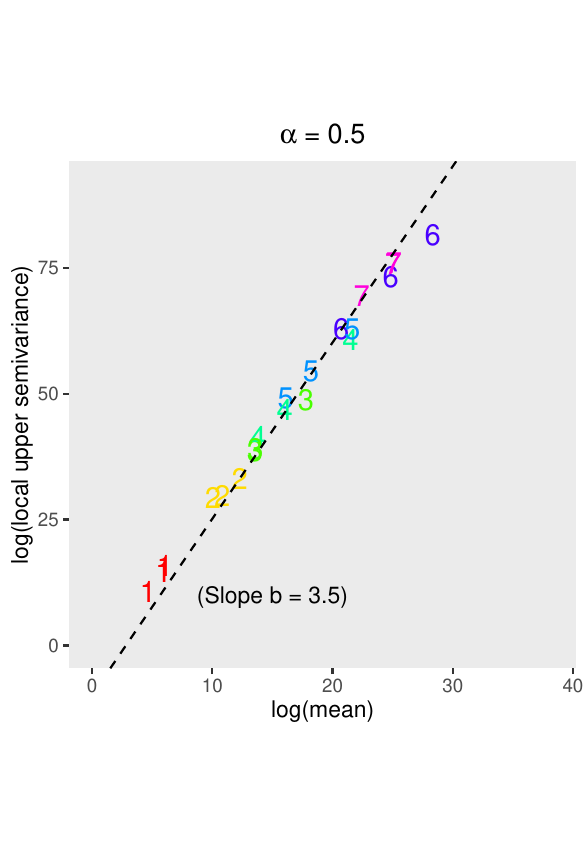}}
		\caption{Scatterplots of (a) $(\log(M_{n,1}),\log(V_n))$, (b) $(\log(M_{n,1}),\log(M^{+}_{n,2}))$, (c) $(\log(M_{n,1}),\log(M^c_{n,3}))$, and (d) $(\log(M_{n,1}),\log(M^{+*}_{n,2}))$, corresponding to Theorem \ref{thm:higher_central_moments} with $h_1=2$ and $h_2=1$, Theorem \ref{thm:upper_central_moment}, Theorem \ref{thm:higher_central_moments} with $h_1=3$ and $h_2=1$, and Theorem \ref{thm:local_upper_central_moment}, respectively, when $F_{1,0.5}$ 
			in \eqref{eq:F1} is the distribution function of $i.i.d.$ $X_i$ with sample size $10^p$ for $p=1,2,\ldots,7$.}
		\label{fig:inverse_iid_2}
	\end{figure}
	
	\begin{figure}[h]
		\centering
		\subfloat[Taylor's law]{\includegraphics[width=0.3\textwidth]{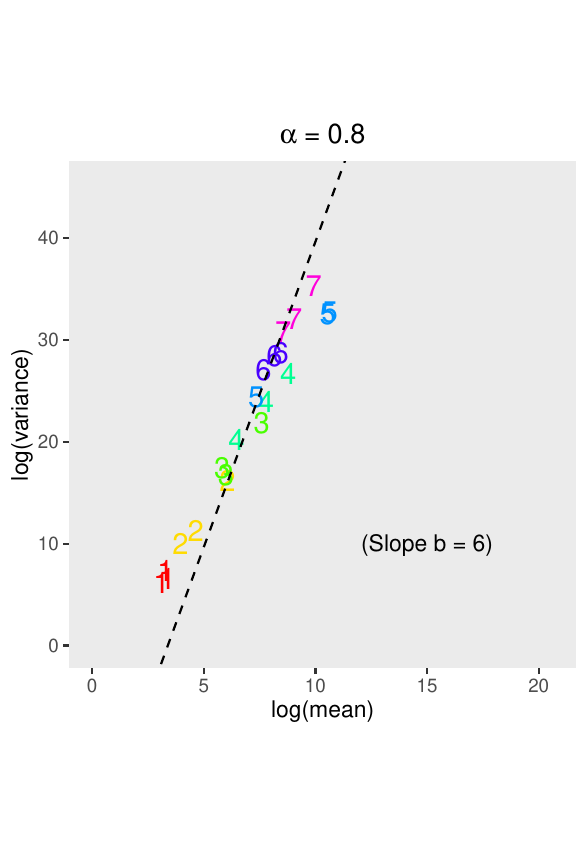}}
		\qquad
		\subfloat[Taylor's law for upper semivariance]{\includegraphics[width=0.3\textwidth]{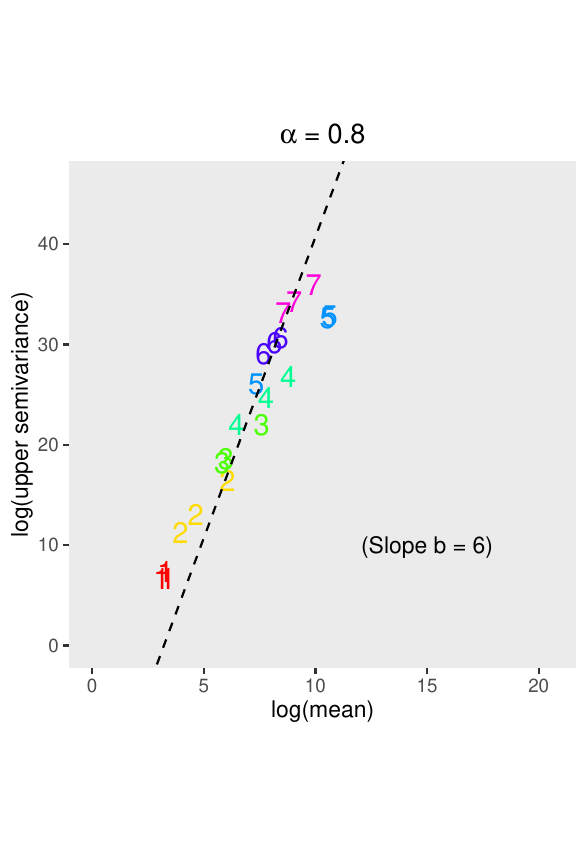}}
		\\
		\subfloat[Taylor's law for third central moment]{\includegraphics[width=0.3\textwidth]{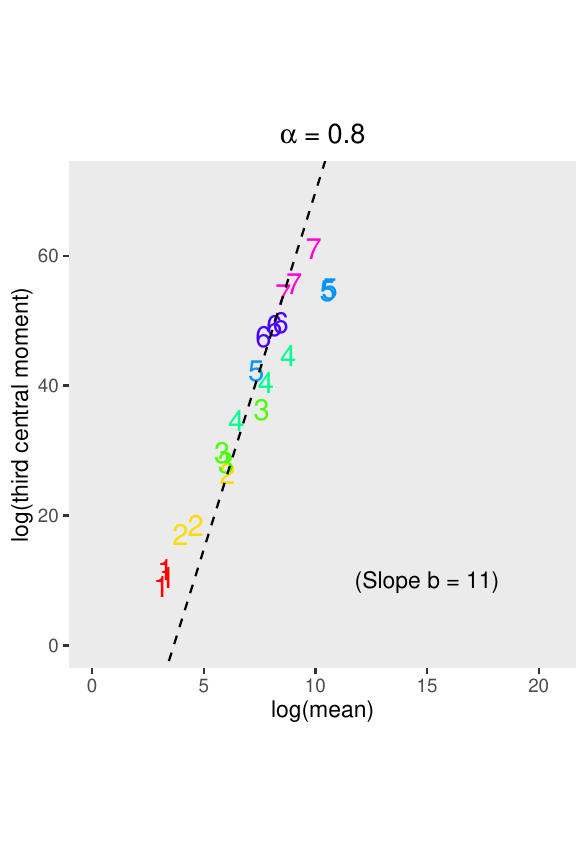}}
		\qquad
		\subfloat[Taylor's law for local upper semivariance]{\includegraphics[width=0.3\textwidth]{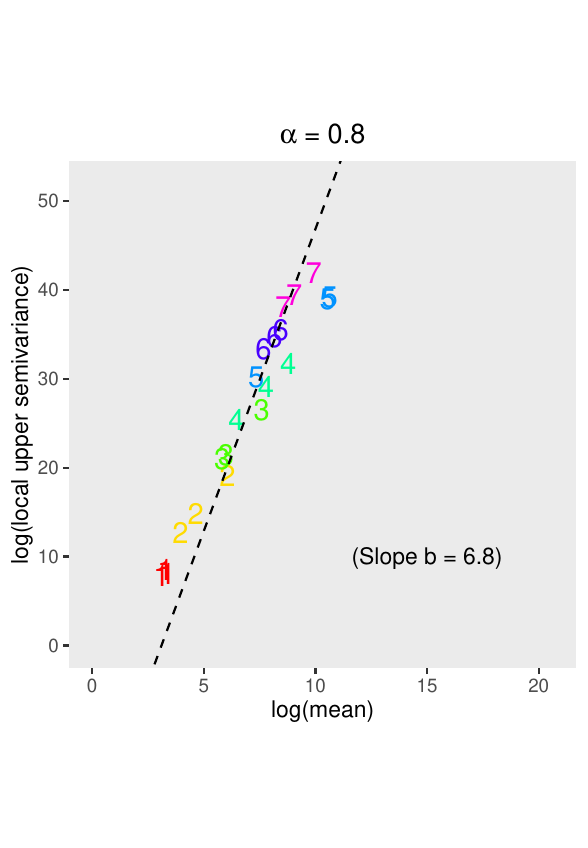}}
		\caption{Scatterplots of (a) $(\log(M_{n,1}),\log(V_n))$, (b) $(\log(M_{n,1}),\log(M^{+}_{n,2}))$, (c) $(\log(M_{n,1}),\log(M^c_{n,3}))$, and (d) $(\log(M_{n,1}),\log(M^{+*}_{n,2}))$, corresponding to Theorem \ref{thm:higher_central_moments} with $h_1=2$ and $h_2=1$, Theorem \ref{thm:upper_central_moment}, Theorem \ref{thm:higher_central_moments} with $h_1=3$ and $h_2=1$, and Theorem \ref{thm:local_upper_central_moment}, respectively, when $F_{1,0.8}$ 
			in \eqref{eq:F1} is the distribution function of $i.i.d.$ $X_i$ with sample size $10^p$ for $p=1,2,\ldots,7$.}
		\label{fig:inverse_iid_3}
	\end{figure}

	\item Consider the stable distribution $F(c, \alpha)${
		on the nonnegative half  line with} 
	Laplace transform
	\begin{align*}
		\mathcal{L}(s) = \mathbb{E}(e^{-(sX)})= e^{-(cs)^\alpha}, \quad s\geq0, \quad c>0, \quad 0<\alpha<1.
	\end{align*}
	According to the tail approximation of the stable distribution \citep{nolan2020univariate}, the survival 
	function of $X \sim F(c, \alpha)$ satisfies $\mathbb{P}(X > x) \sim C x^{-\alpha}$ for some constant $C > 0$. Thus, it has a tail index of $\alpha$.
	Figures \ref{fig:stable_iid_1}, \ref{fig:stable_iid_2}, and \ref{fig:stable_iid_3} present the results of simulations for $\alpha = 0.2$, $0.5$, and $0.8$. We used the \texttt{rstable} function from the \textbf{R} package \textbf{stabledist} to generate random variables following the stable distribution.

	\begin{figure}[h]
		\centering
		\subfloat[Taylor's law]{\includegraphics[width=0.3\textwidth]{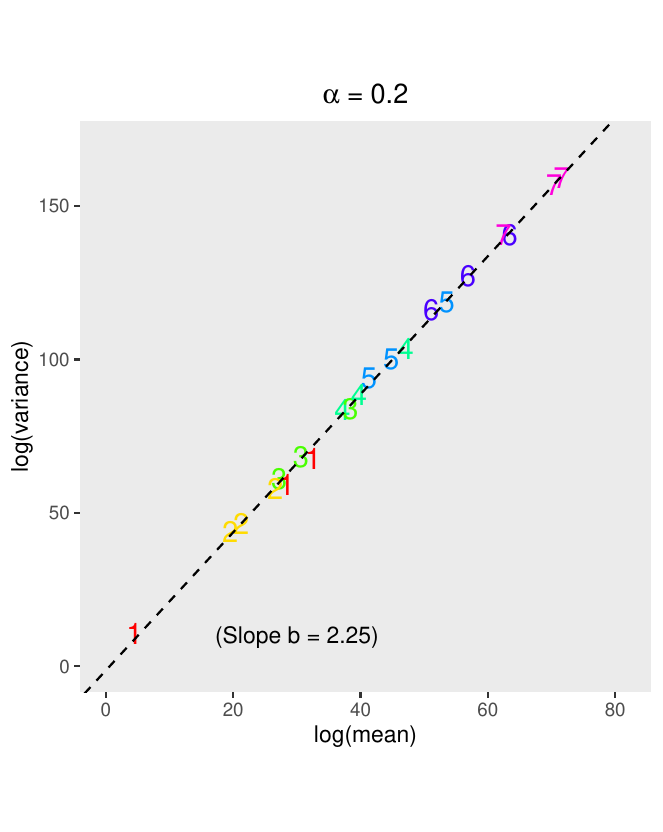}}
		\qquad
		\subfloat[Taylor's law for upper semivariance]{\includegraphics[width=0.3\textwidth]{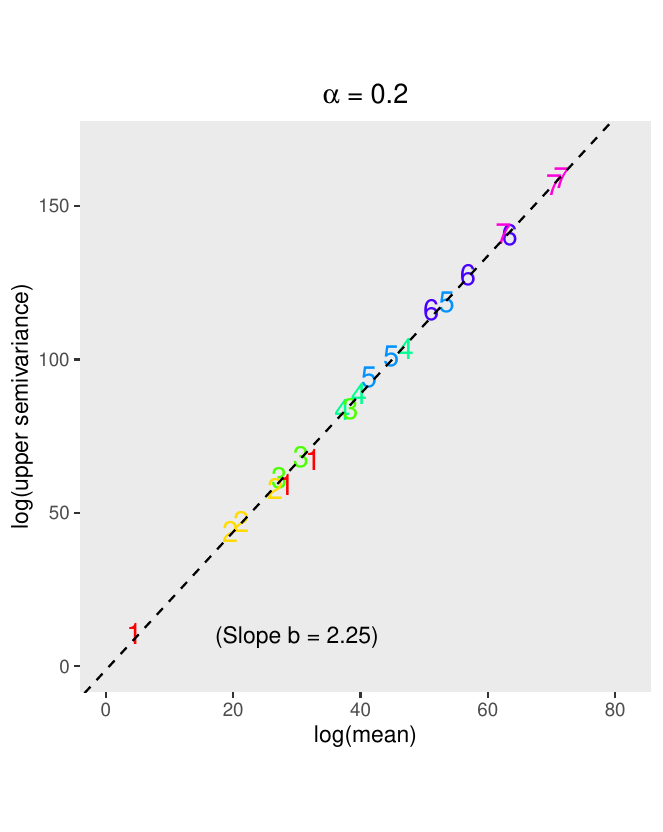}}
		\\
		\subfloat[Taylor's law for third central moment]{\includegraphics[width=0.3\textwidth]{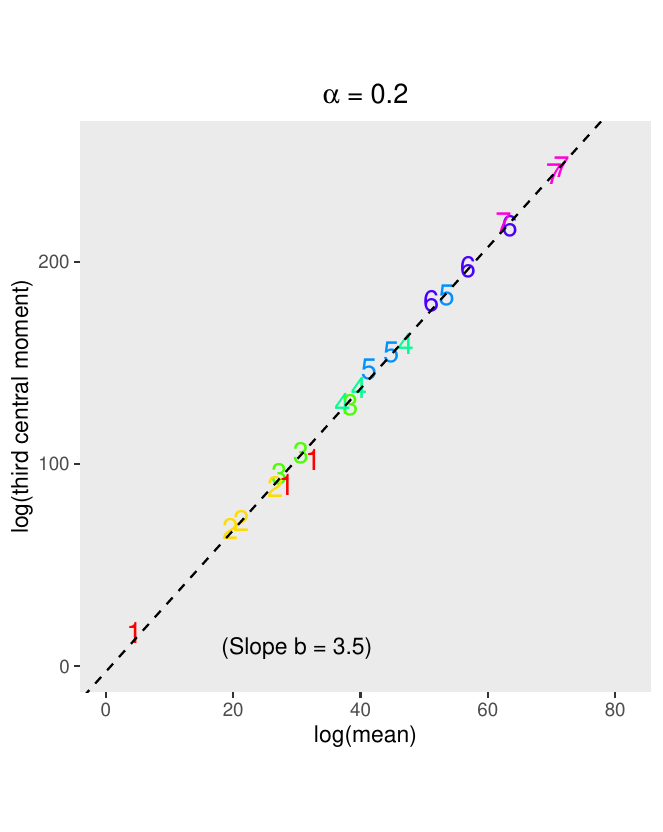}}
		\qquad
		\subfloat[Taylor's law for local upper semivariance]{\includegraphics[width=0.3\textwidth]{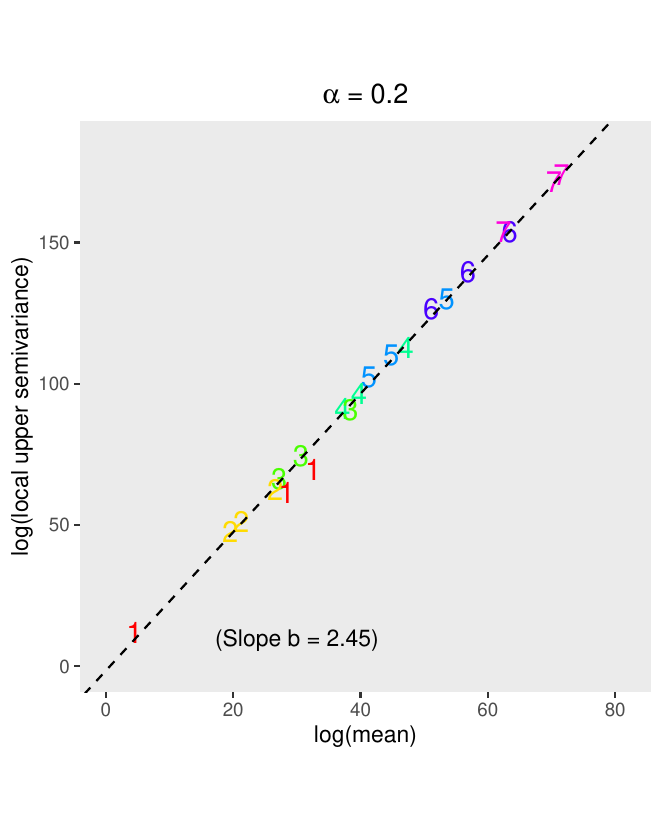}}
		\caption{Scatterplots of (a) $(\log(M_{n,1}),\log(V_n))$, (b) $(\log(M_{n,1}),\log(M^{+}_{n,2}))$, (c) $(\log(M_{n,1}),\log(M^c_{n,3}))$, and (d) $(\log(M_{n,1}),\log(M^{+*}_{n,2}))$, corresponding to Theorem \ref{thm:higher_central_moments} with $h_1=2$ and $h_2=1$, Theorem \ref{thm:upper_central_moment}, Theorem \ref{thm:higher_central_moments} with $h_1=3$ and $h_2=1$, and Theorem \ref{thm:local_upper_central_moment}, respectively, when $X_i \overset{\mathrm{i.i.d.}}{\sim} F(1,0.2)$ with sample size $10^p$ for {$p=1,2,\ldots,7$}.}
		\label{fig:stable_iid_1}
	\end{figure}
	
	\begin{figure}[h]
		\centering
		\subfloat[Taylor's law]{\includegraphics[width=0.3\textwidth]{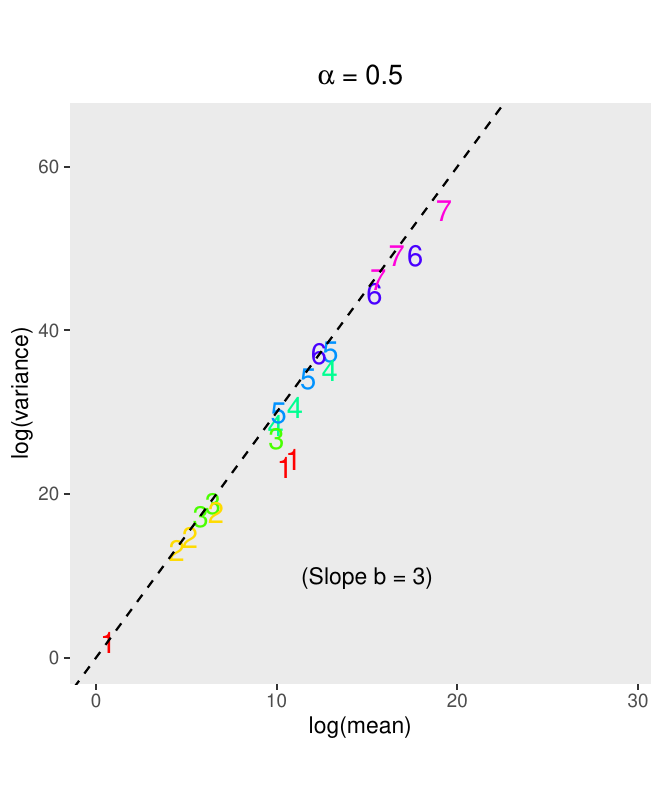}}
		\qquad
		\subfloat[Taylor's law for upper semivariance]{\includegraphics[width=0.3\textwidth]{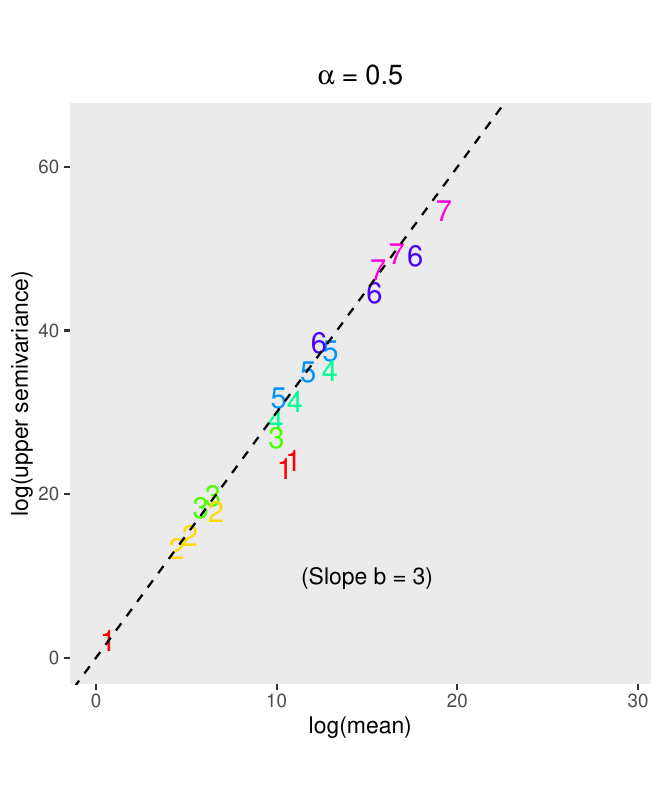}}
		\\
		\subfloat[Taylor's law for third central moment]{\includegraphics[width=0.3\textwidth]{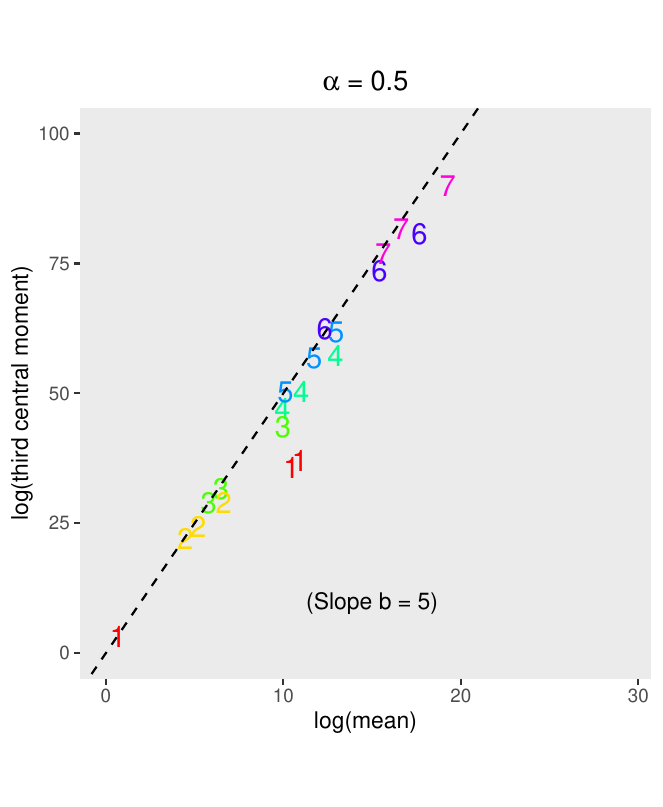}}
		\qquad
		\subfloat[Taylor's law for local upper semivariance]{\includegraphics[width=0.3\textwidth]{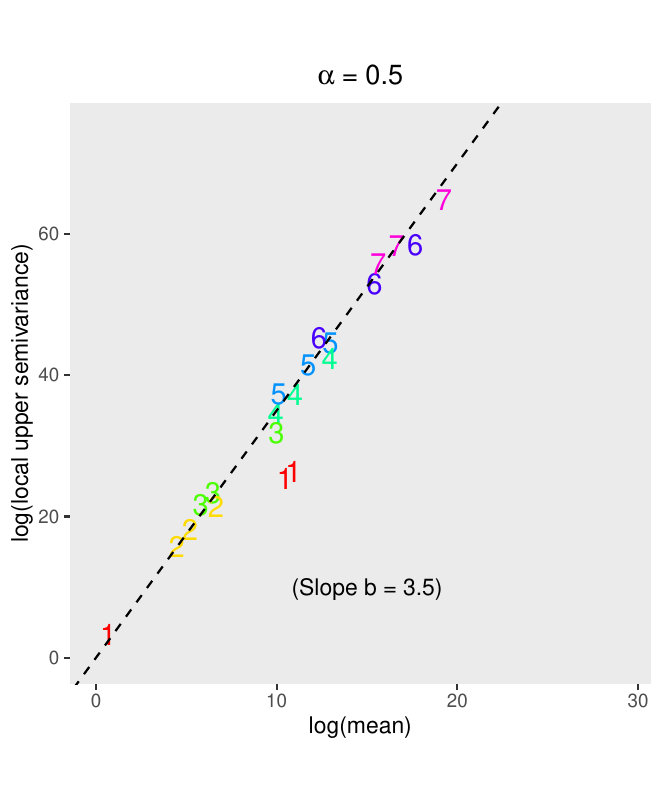}}
		\caption{Scatterplots of (a) $(\log(M_{n,1}),\log(V_n))$, (b) $(\log(M_{n,1}),\log(M^{+}_{n,2}))$, (c) $(\log(M_{n,1}),\log(M^c_{n,3}))$, and (d) $(\log(M_{n,1}),\log(M^{+*}_{n,2}))$, corresponding to Theorem \ref{thm:higher_central_moments} with $h_1=2$ and $h_2=1$, Theorem \ref{thm:upper_central_moment}, Theorem \ref{thm:higher_central_moments} with $h_1=3$ and $h_2=1$, and Theorem \ref{thm:local_upper_central_moment}, respectively, when $X_i \overset{\mathrm{i.i.d.}}{\sim} F(1,0.5)$ with sample size $10^p$ for {$p=1,2,\ldots,7$}.}
		\label{fig:stable_iid_2}
	\end{figure}
	
	\begin{figure}[h]
		\centering
		\subfloat[Taylor's law]{\includegraphics[width=0.3\textwidth]{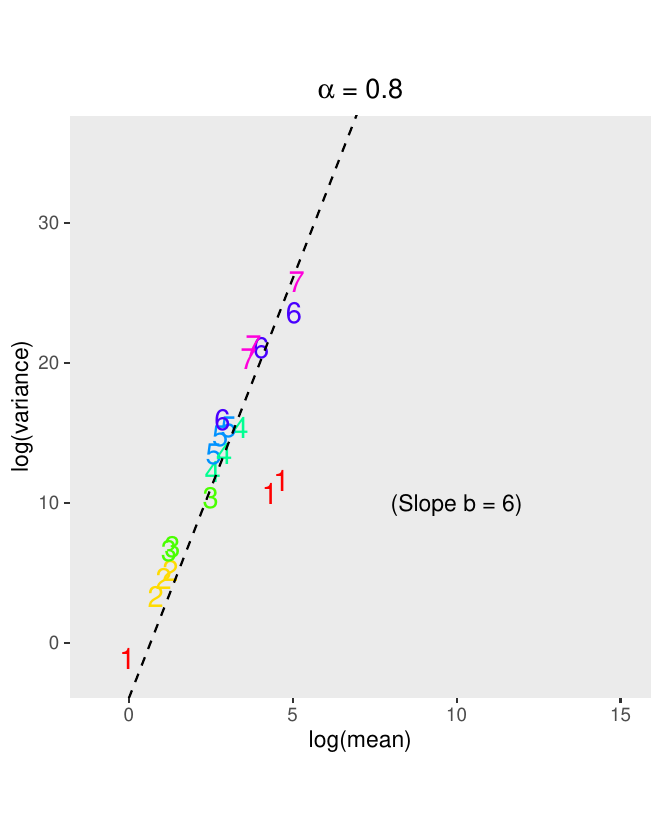}}
		\qquad
		\subfloat[Taylor's law for upper semivariance]{\includegraphics[width=0.3\textwidth]{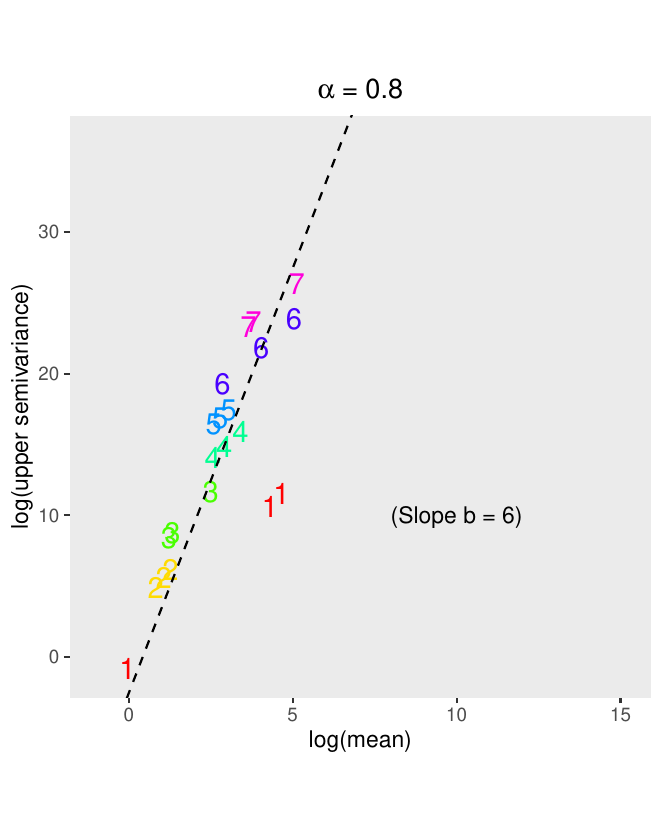}}
		\\
		\subfloat[Taylor's law for third central moment]{\includegraphics[width=0.3\textwidth]{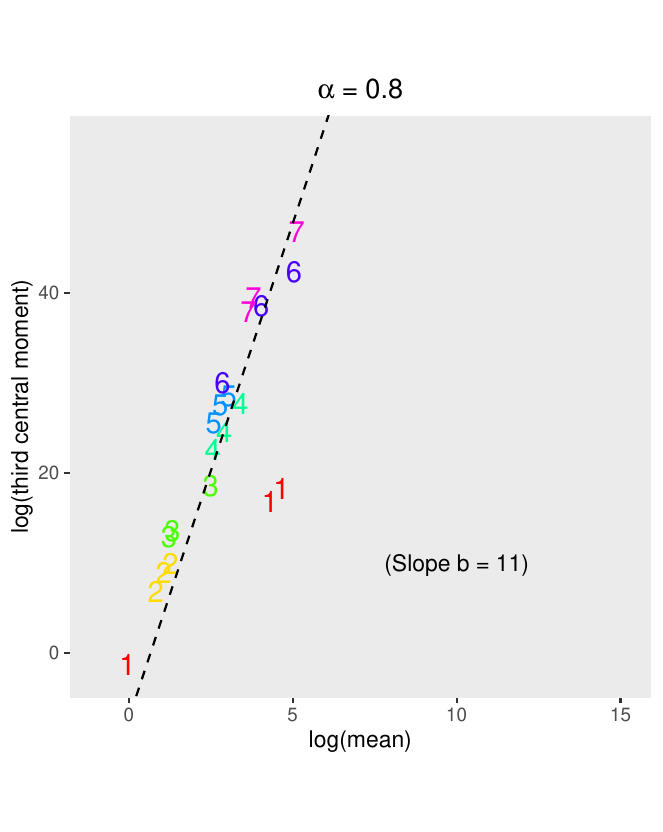}}
		\qquad
		\subfloat[Taylor's law for local upper semivariance]{\includegraphics[width=0.3\textwidth]{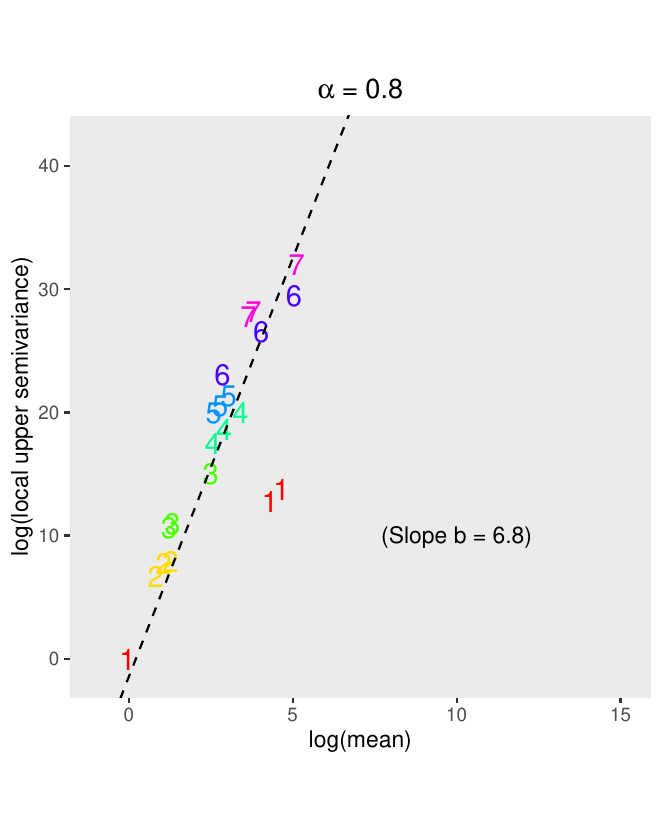}}
		\caption{Scatterplots of (a) $(\log(M_{n,1}),\log(V_n))$, (b) $(\log(M_{n,1}),\log(M^{+}_{n,2}))$, (c) $(\log(M_{n,1}),\log(M^c_{n,3}))$, and (d) $(\log(M_{n,1}),\log(M^{+*}_{n,2}))$, corresponding to Theorem \ref{thm:higher_central_moments} with $h_1=2$ and $h_2=1$, Theorem \ref{thm:upper_central_moment}, Theorem \ref{thm:higher_central_moments} with $h_1=3$ and $h_2=1$, and Theorem \ref{thm:local_upper_central_moment}, respectively, when $X_i \overset{\mathrm{i.i.d.}}{\sim} F(1,0.8)$ with sample size $10^p$ for {$p=1,2,\ldots,7$}.}
		\label{fig:stable_iid_3}
	\end{figure}
	
\end{enumerate}

\subsection{Taylor's Law for an AR(1) Model}
\begin{figure}[h]
	\centering
	\subfloat[Taylor's law]{\includegraphics[width=0.3\textwidth]{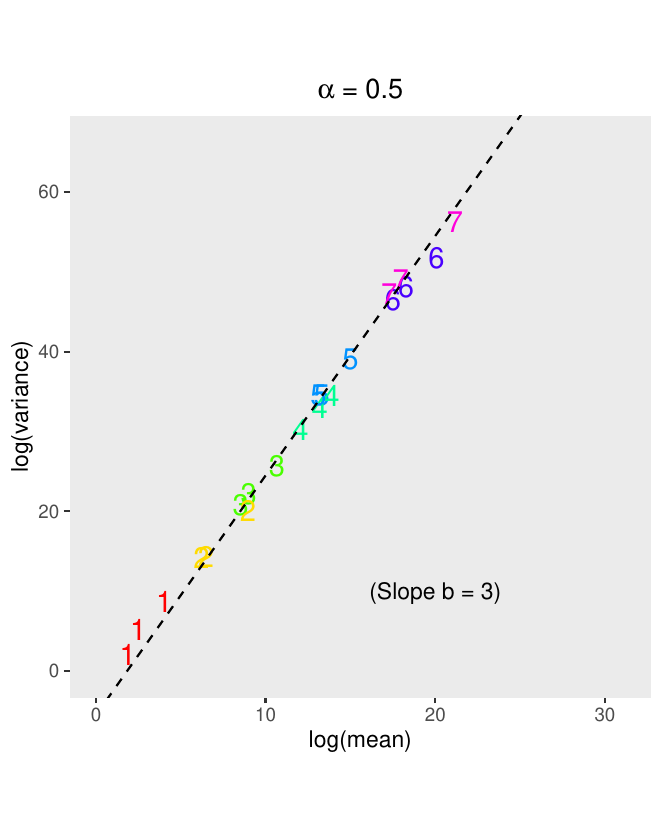}}
	\qquad
	\subfloat[Taylor's law for upper semivariance]{\includegraphics[width=0.3\textwidth]{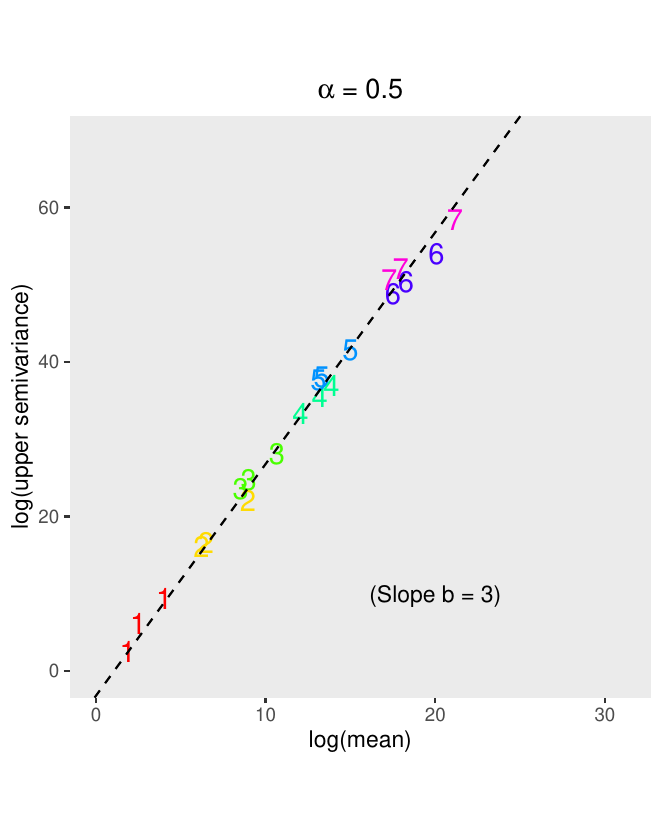}}
	\\
	\subfloat[Taylor's law for third central moment]{\includegraphics[width=0.3\textwidth]{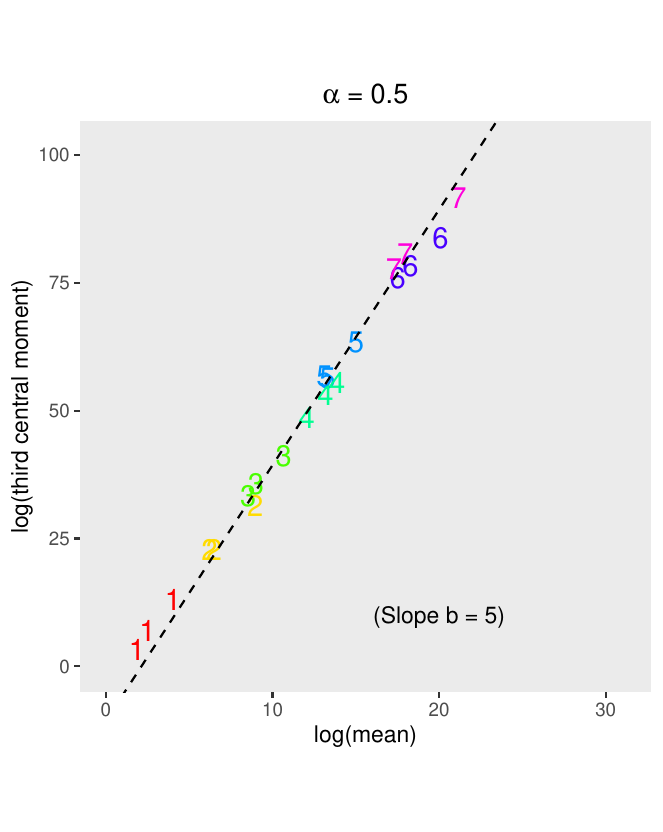}}
	\qquad
	\subfloat[Taylor's law for local upper semivariance]{\includegraphics[width=0.3\textwidth]{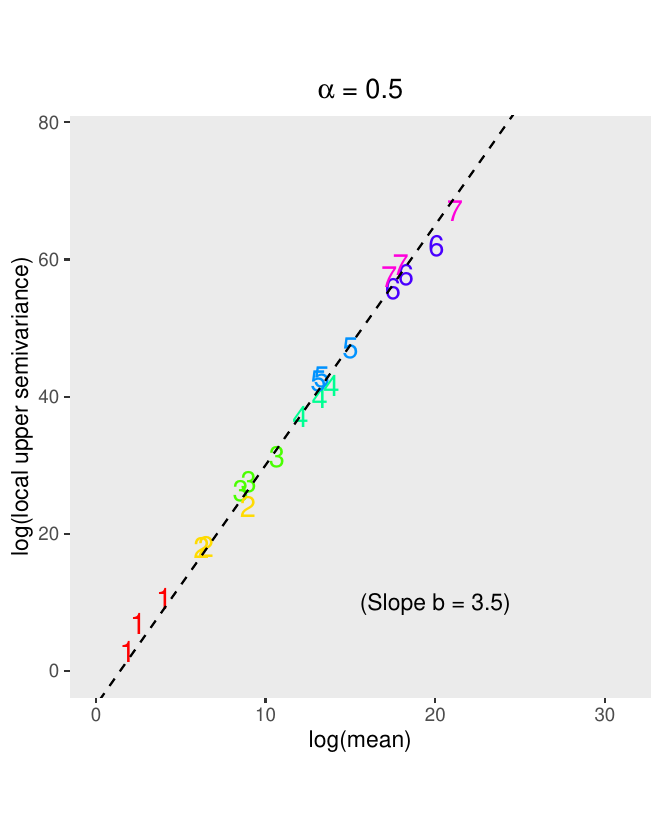}}
	\caption{Scatterplots of (a) $(\log(M_{n,1}),\log(V_n))$, (b) $(\log(M_{n,1}),\log(M^{+}_{n,2}))$, (c) $(\log(M_{n,1}),\log(M^c_{n,3}))$, and (d) $(\log(M_{n,1}),\log(M^{+*}_{n,2}))$, corresponding to Theorem \ref{thm:higher_central_moments} with $h_1=2$ and $h_2=1$, Theorem \ref{thm:upper_central_moment}, Theorem \ref{thm:higher_central_moments} with $h_1=3$ and $h_2=1$, and Theorem \ref{thm:local_upper_central_moment}, respectively, when $X_i$ follows an AR(1) model with $\beta_1=0.8$, and $\epsilon_i \overset{\mathrm{i.i.d.}}{\sim} $Pareto$(1,0.5)$ with sample size $10^p$ for $p=1,2,\ldots,7$.}
	\label{fig:AR_pareto_1}
\end{figure}

\begin{figure}[h]
	\centering
	\subfloat[Taylor's law]{\includegraphics[width=0.3\textwidth]{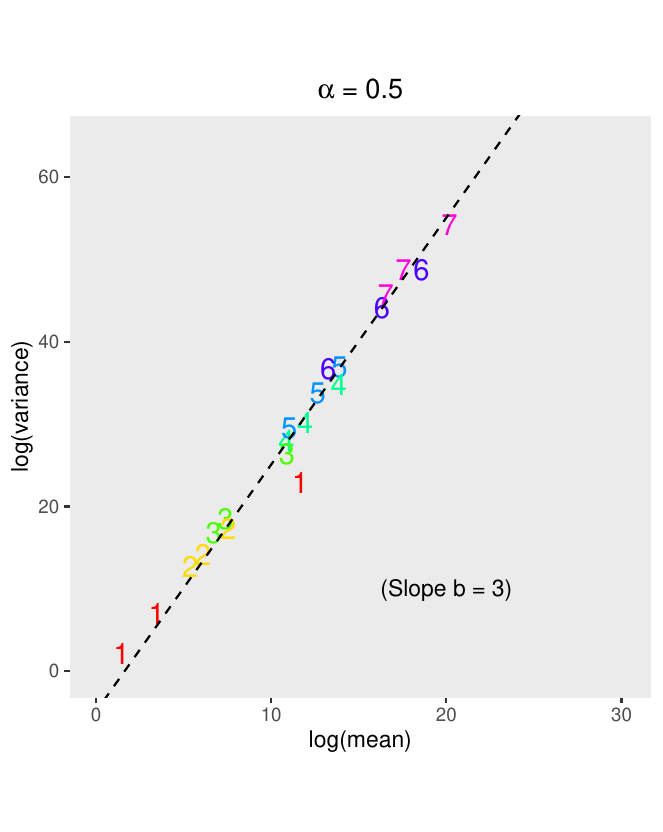}}
	\qquad
	\subfloat[Taylor's law for upper semivariance]{\includegraphics[width=0.3\textwidth]{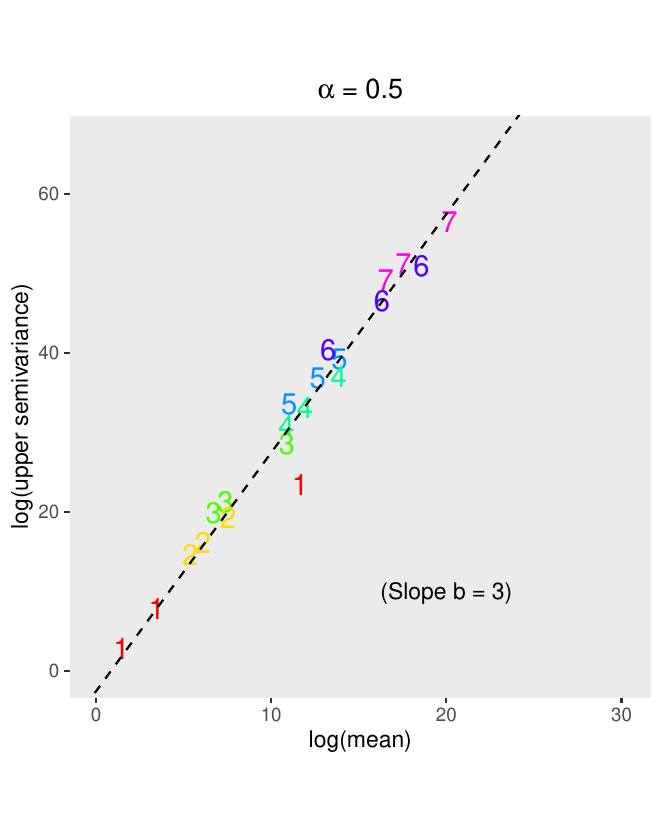}}
	\\
	\subfloat[Taylor's law for third central moment]{\includegraphics[width=0.3\textwidth]{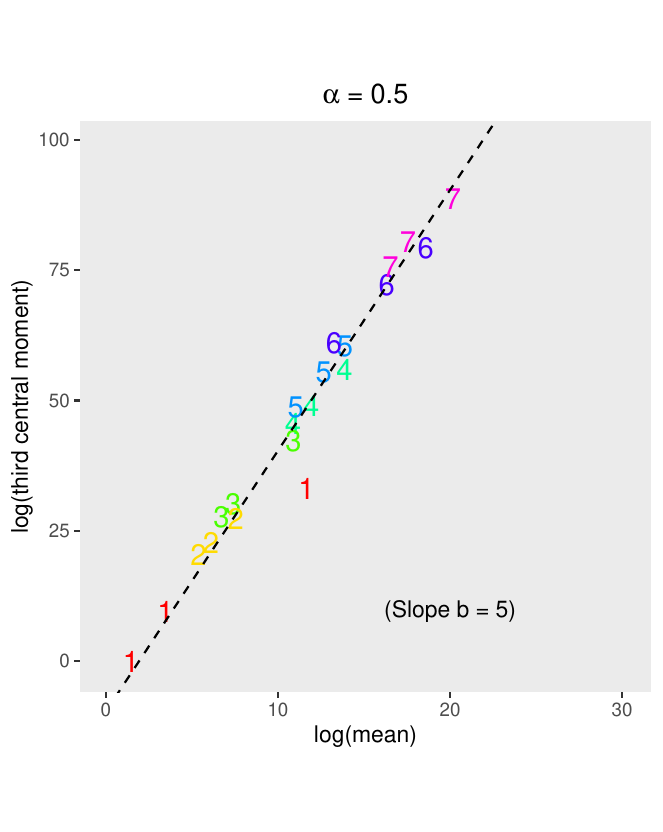}}
	\qquad
	\subfloat[Taylor's law for local upper semivariance]{\includegraphics[width=0.3\textwidth]{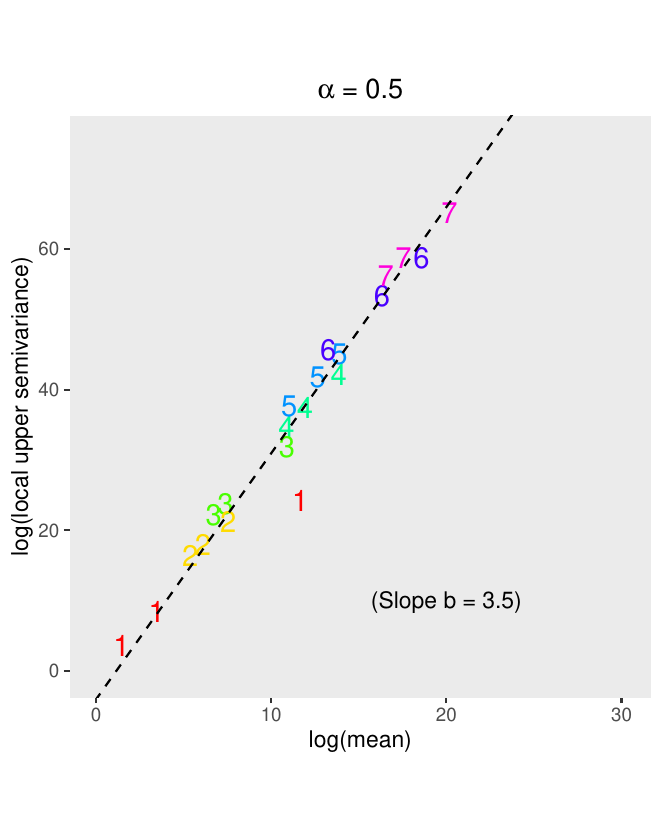}}
	\caption{Scatterplots of (a) $(\log(M_{n,1}),\log(V_n))$, (b) $(\log(M_{n,1}),\log(M^{+}_{n,2}))$, (c) $(\log(M_{n,1}),\log(M^c_{n,3}))$, and (d) $(\log(M_{n,1}),\log(M^{+*}_{n,2}))$, corresponding to Theorem \ref{thm:higher_central_moments} with $h_1=2$ and $h_2=1$, Theorem \ref{thm:upper_central_moment}, Theorem \ref{thm:higher_central_moments} with $h_1=3$ and $h_2=1$, and Theorem \ref{thm:local_upper_central_moment}, respectively, when $X_i$ follows an AR(1) model with $\beta_1=0.8$, and $\epsilon_i \overset{\mathrm{i.i.d.}}{\sim} $$F(1,0.5)$ with sample size $10^p$ for $p=1,2,\ldots,7$.}
	
	\label{fig:AR_stable_1}
\end{figure}

{We simulated two cases of the first-order autoregressive model discussed in Section \ref{subsec:cov condition}, Theorem \ref{thm:AR_cov}:
	\begin{align*}
		X_t=\beta_1 X_{t-1}+\epsilon_t
	\end{align*}
	for $t=1,\ldots,n$, $\beta_1 \in (0,1)$. 
	In the first case, $\epsilon_t \sim Pareto(1,\alpha)$. 
	In the second case, $\epsilon_t \sim F(1,\alpha)$ with $\alpha=0.5$. 
	Both distributions satisfy the conditions given in  Theorem \ref{thm:AR_cov}. 
	Indeed, the Pareto distribution has a decreasing density and any stable distribution is unimodal \citep{nolan2020univariate}, and both densities are also bounded.
	
	Figures \ref{fig:AR_pareto_1} and \ref{fig:AR_stable_1} show the 
	results of simulations of
	AR(1) models when the noise follows Pareto and stable distributions, respectively, with $\alpha=0.5$, and $\beta_1 = 0.8$. 
	
	The sample size was $n=10^p$ for $p=1,2,\ldots,7$. }

When there is autoregressive dependence in the simulated time series, 
the points in the figures  lie further from the straight dashed line than 
when the random variables are  $i.i.d.$, 
suggesting that, even though our results hold for weak dependence
asymptotically, the asymptotic results may not be very accurate for 
finite samples.

\subsection{Taylor's Law with Heterogeneous Data}
\begin{figure}[h]
	\centering
	\subfloat[Taylor's law]{\includegraphics[width=0.3\textwidth]{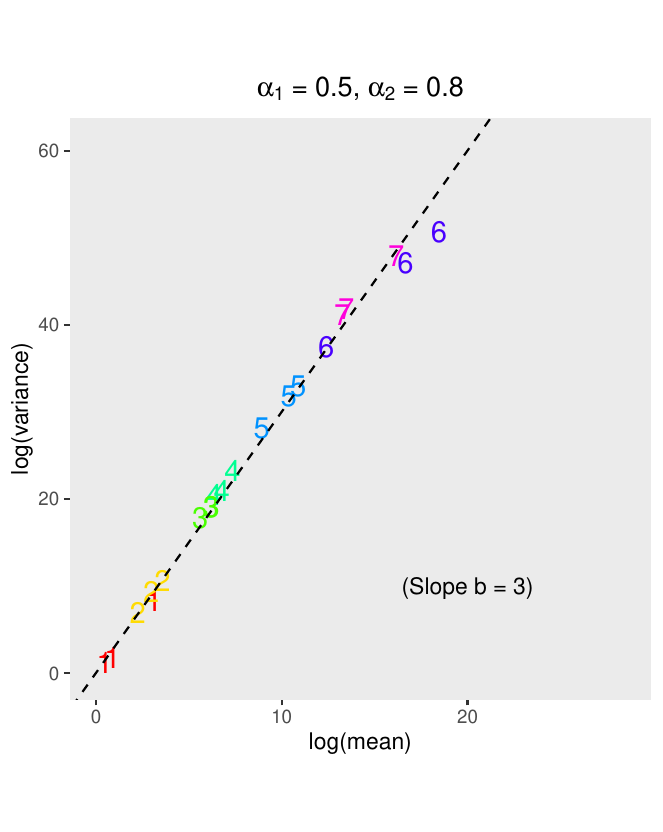}}
	\qquad
	\subfloat[Taylor's law for upper semivariance]{\includegraphics[width=0.3\textwidth]{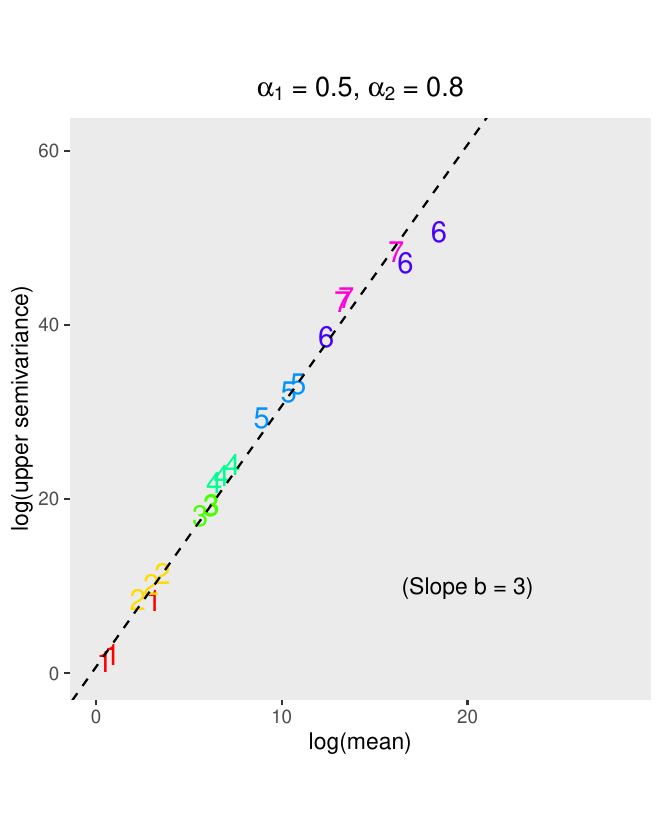}}
	\\
	\subfloat[Taylor's law for third central moment]{\includegraphics[width=0.3\textwidth]{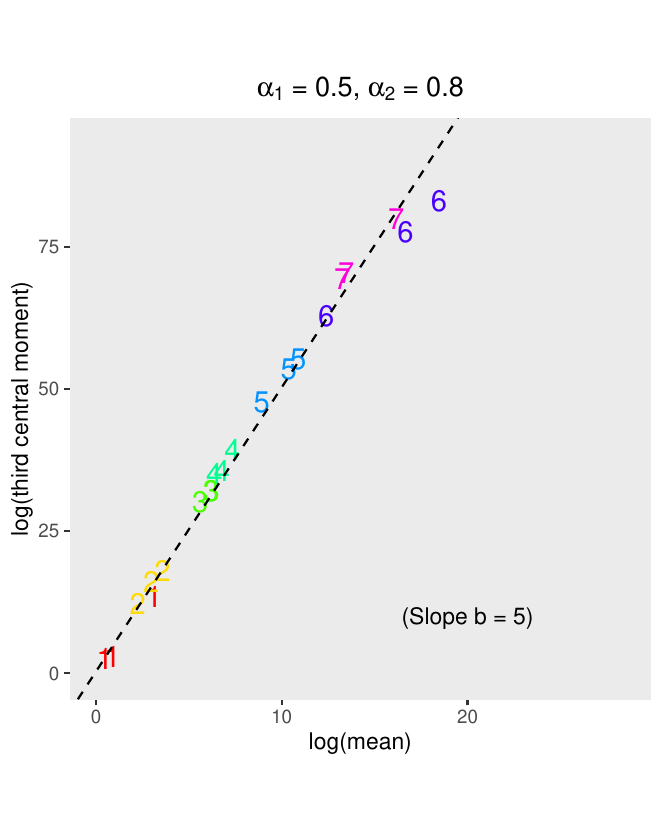}}
	\qquad
	\subfloat[Taylor's law for local upper semivariance]{\includegraphics[width=0.3\textwidth]{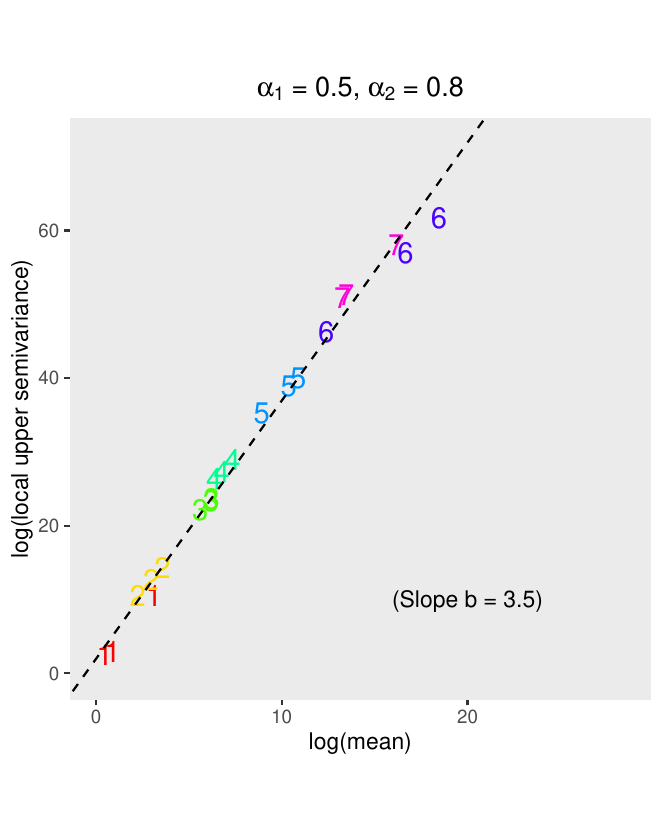}}
	\caption{Scatterplots of (a) $(\log(M_{n,1}),\log(V_n))$, 
		(b) $(\log(M_{n,1}),\log(M^{+}_{n,2}))$, (c) $(\log(M_{n,1}),\log(M^c_{n,3}))$, and (d) $(\log(M_{n,1}),\log(M^{+*}_{n,2}))$, corresponding to Theorem \ref{thm:higher_central_moments} with $h_1=2$ and $h_2=1$, Theorem \ref{thm:upper_central_moment}, Theorem \ref{thm:higher_central_moments} with $h_1=3$ and $h_2=1$, and Theorem \ref{thm:local_upper_central_moment}, respectively, when $X_{i,U} \overset{\mathrm{i.i.d.}}{\sim}$Pareto$(0.5,0.5)$, $X_{i,V} \overset{\mathrm{i.i.d.}}{\sim}$ Pareto$(0.8,0.8)$, and $p^*=0.6$, with sample size $10^p$ for $p=1,2,\ldots,7$.}
	\label{fig:heter_pareto_1}
\end{figure}

\begin{figure}[h]
	\centering
	\subfloat[Taylor's law]{\includegraphics[width=0.3\textwidth]{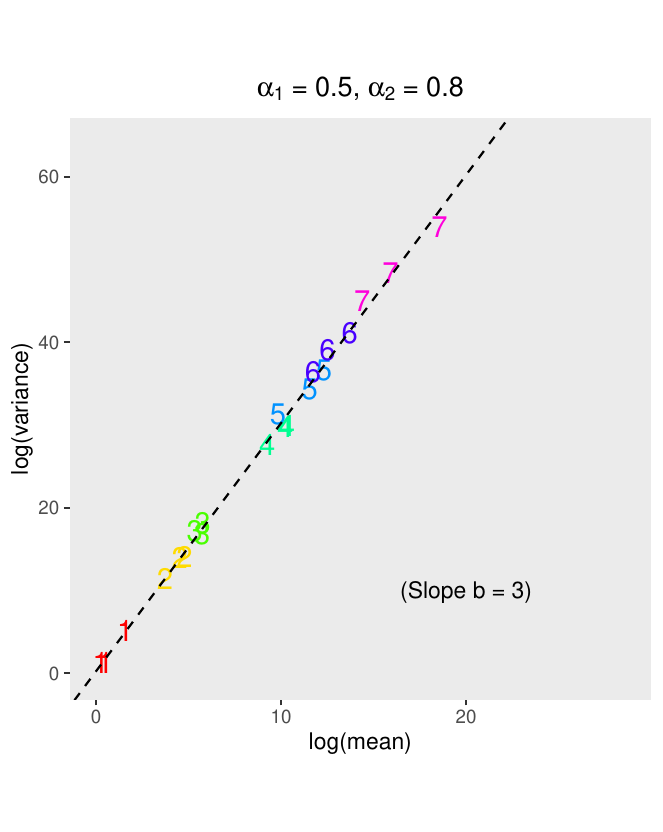}}
	\qquad
	\subfloat[Taylor's law for upper semivariance]{\includegraphics[width=0.3\textwidth]{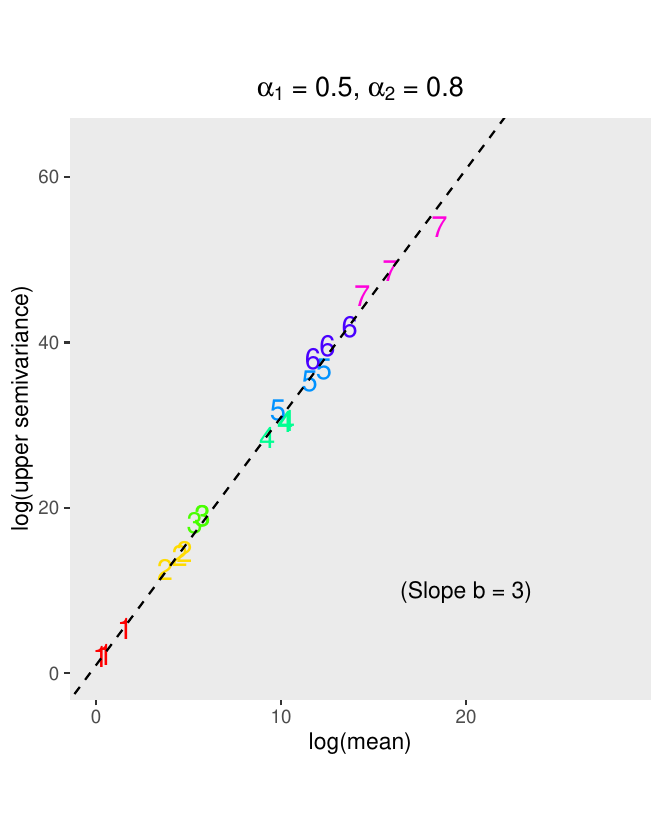}}
	\\
	\subfloat[Taylor's law for third central moment]{\includegraphics[width=0.3\textwidth]{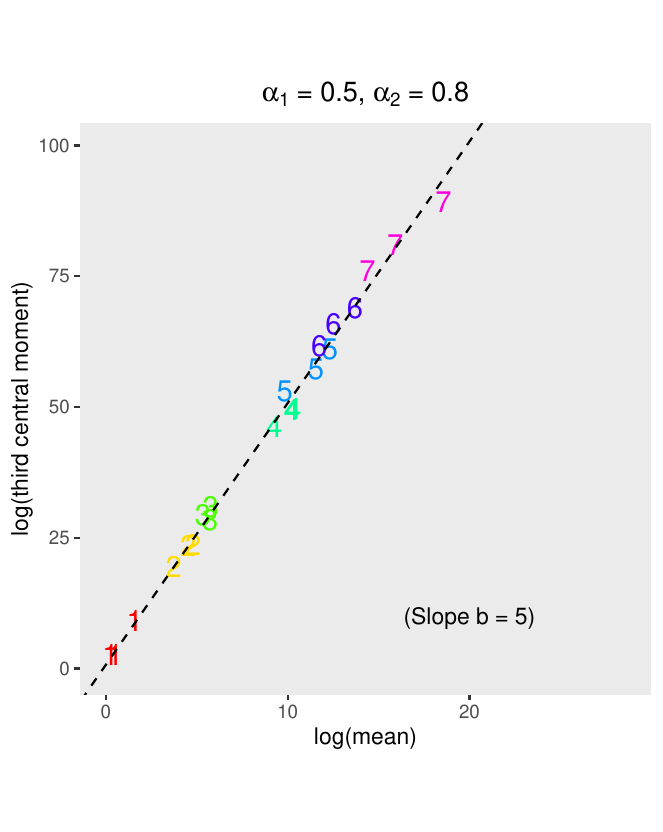}}
	\qquad
	\subfloat[Taylor's law for local upper semivariance]{\includegraphics[width=0.3\textwidth]{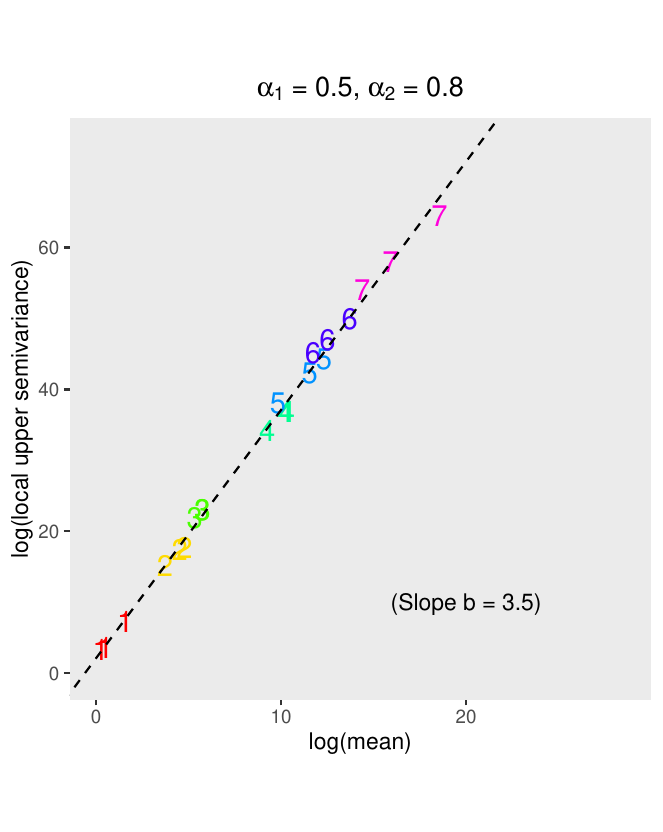}}
	\caption{Scatterplots of (a) $(\log(M_{n,1}),\log(V_n))$, (b) $(\log(M_{n,1}),\log(M^{+}_{n,2}))$, (c) $(\log(M_{n,1}),\log(M^c_{n,3}))$, and (d) $(\log(M_{n,1}),\log(M^{+*}_{n,2}))$, corresponding to Theorem \ref{thm:higher_central_moments} with $h_1=2$ and $h_2=1$, Theorem \ref{thm:upper_central_moment}, Theorem \ref{thm:higher_central_moments} with $h_1=3$ and $h_2=1$, and Theorem \ref{thm:local_upper_central_moment}, respectively, when $X_{i,U} \overset{\mathrm{i.i.d.}}{\sim} F(1,0.5)$, $X_{i,V} \overset{\mathrm{i.i.d.}}{\sim} F(1,0.8)$, and $p^*=0.3$, with sample size $10^p$ for $p=1,2,\ldots,7$.}
	\label{fig:heter_stable_1}
\end{figure}

\begin{figure}[h]
	\centering
	\subfloat[Taylor's law]{\includegraphics[width=0.3\textwidth]{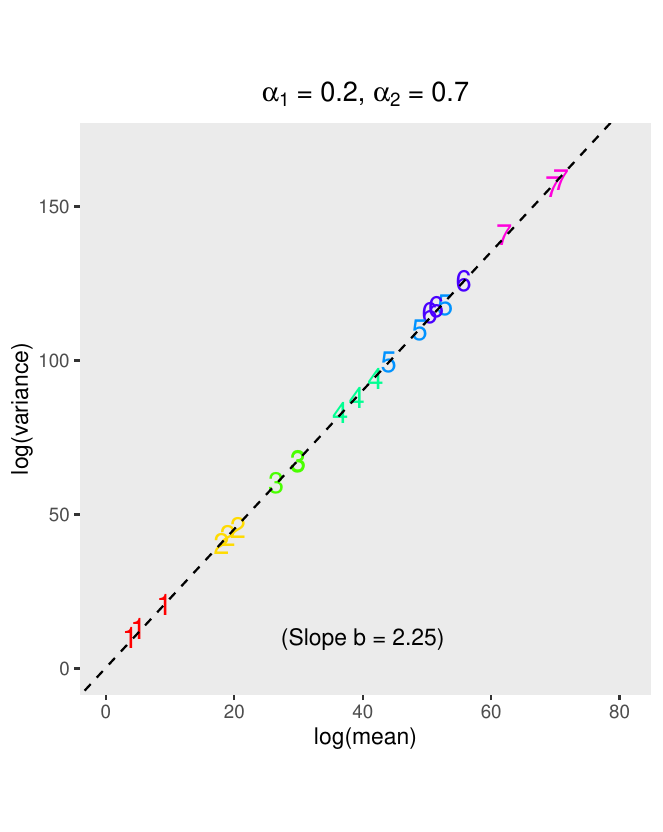}}
	\qquad
	\subfloat[Taylor's law for upper semivariance]{\includegraphics[width=0.3\textwidth]{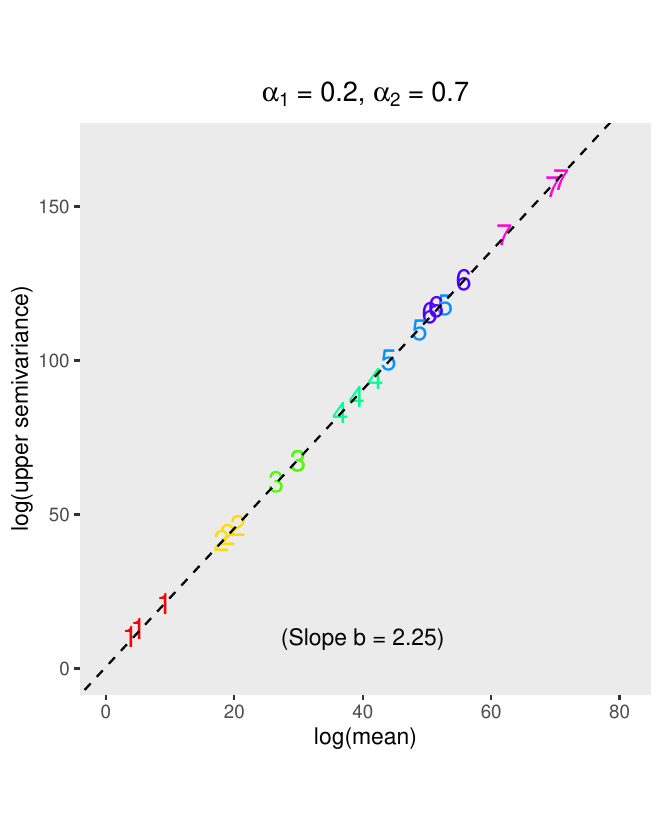}}
	\\
	\subfloat[Taylor's law for third central moment]{\includegraphics[width=0.3\textwidth]{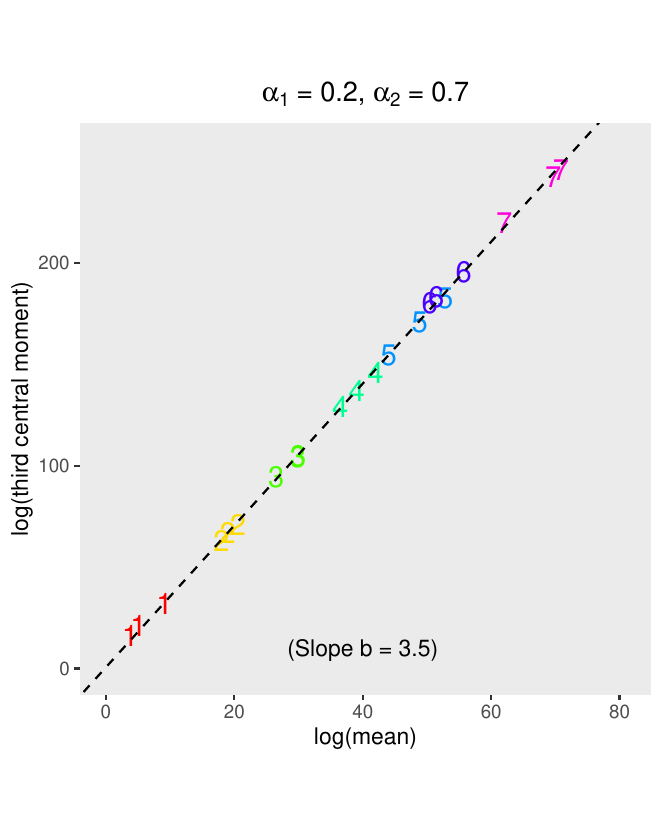}}
	\qquad
	\subfloat[Taylor's law for local upper semivariance]{\includegraphics[width=0.3\textwidth]{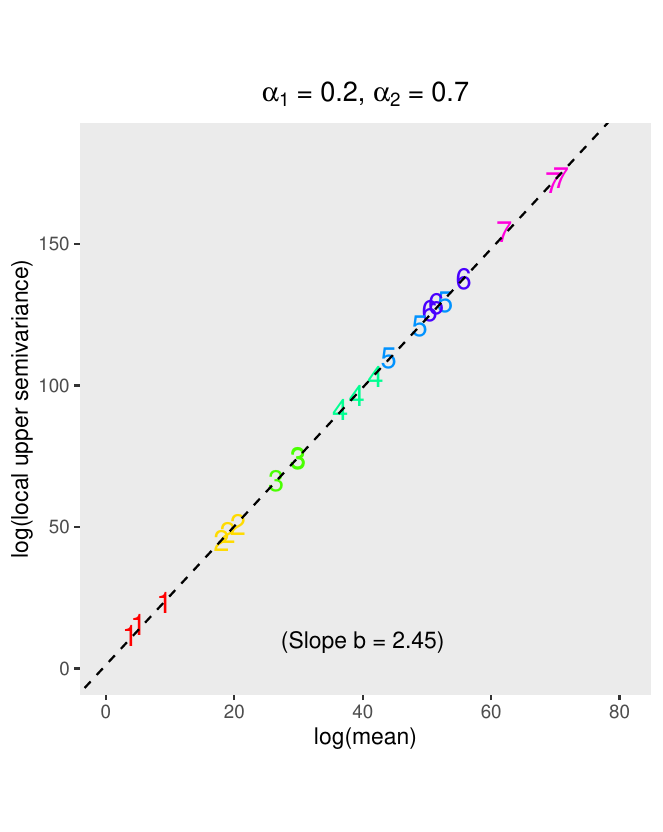}}
	\caption{Scatterplots of (a) $(\log(M_{n,1}),\log(V_n))$, (b) $(\log(M_{n,1}),\log(M^{+}_{n,2}))$, (c) $(\log(M_{n,1}),\log(M^c_{n,3}))$, and (d) $(\log(M_{n,1}),\log(M^{+*}_{n,2}))$, corresponding to Theorem \ref{thm:higher_central_moments} with $h_1=2$ and $h_2=1$, Theorem \ref{thm:upper_central_moment}, Theorem \ref{thm:higher_central_moments} with $h_1=3$ and $h_2=1$, and Theorem \ref{thm:local_upper_central_moment}, respectively, when $X_{i,U} \overset{\mathrm{i.i.d.}}{\sim} F(1,0.2)$, $X_{i,V} \overset{\mathrm{i.i.d.}}{\sim} F(1,0.7)$, and $p^*=0.3$, with sample size $10^p$ for $p=1,2,\ldots,7$.}
	\label{fig:heter_stable_2}
\end{figure}

\begin{figure}[h]
	\centering
	\subfloat[Taylor's law]{\includegraphics[width=0.3\textwidth]{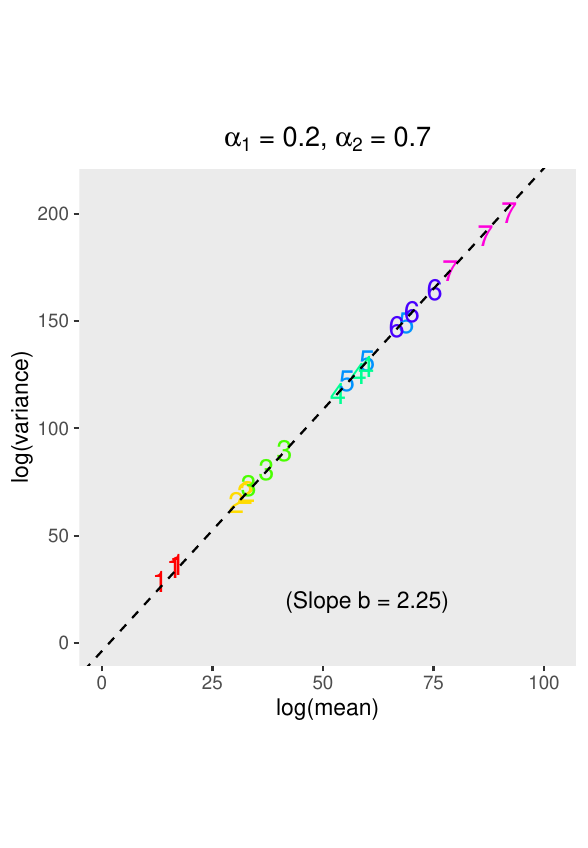}}
	\qquad
	\subfloat[Taylor's law for upper semivariance]{\includegraphics[width=0.3\textwidth]{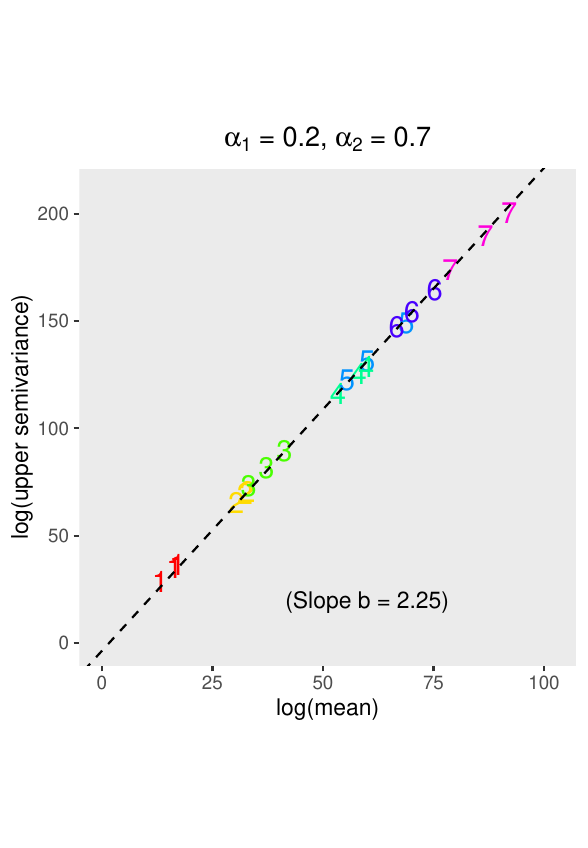}}
	\\
	\subfloat[Taylor's law for third central moment]{\includegraphics[width=0.3\textwidth]{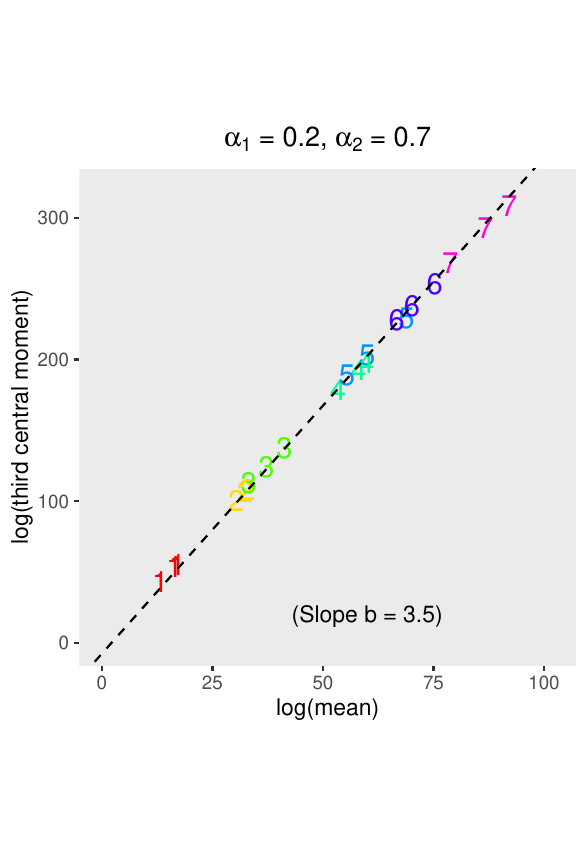}}
	\qquad
	\subfloat[Taylor's law for local upper semivariance]{\includegraphics[width=0.3\textwidth]{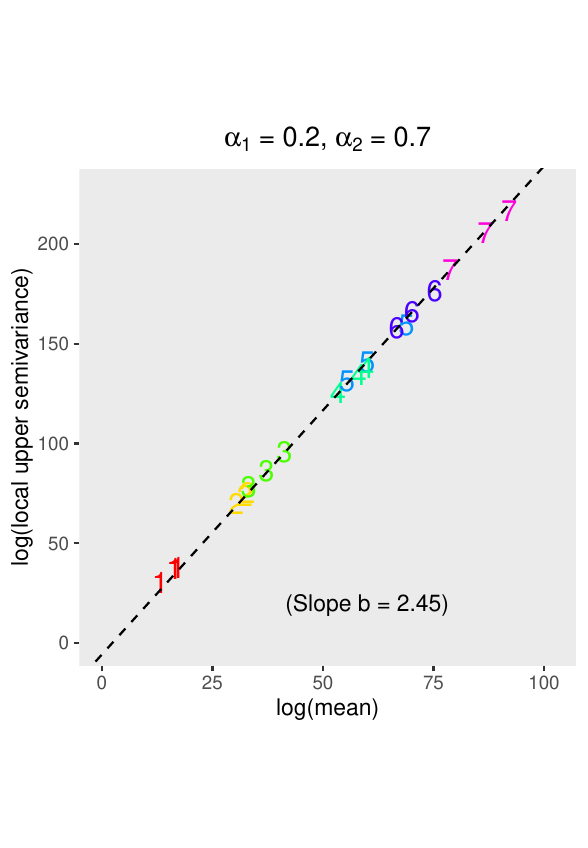}}
	\caption{Scatterplots of (a) $(\log(M_{n,1}),\log(V_n))$, (b) $(\log(M_{n,1}),\log(M^{+}_{n,2}))$, (c) $(\log(M_{n,1}),\log(M^c_{n,3}))$, and (d) $(\log(M_{n,1}),\log(M^{+*}_{n,2}))$, corresponding to Theorem \ref{thm:higher_central_moments} with $h_1=2$ and $h_2=1$, Theorem \ref{thm:upper_central_moment}, Theorem \ref{thm:higher_central_moments} with $h_1=3$ and $h_2=1$, and Theorem \ref{thm:local_upper_central_moment}, respectively, when $X_{i,U} \overset{\mathrm{i.i.d.}}{\sim} F_{1,0.2}$, $X_{i,V} \overset{\mathrm{i.i.d.}}{\sim} F_{1, 0.7}$, and $p^*=0.3$, with sample size $10^p$ for $p=1,2,\ldots,7$.}
	\label{fig:heter_inverse_1}
\end{figure}

Here we simulate the setting described in Section \ref{sec:heterogeneous}. 
Given a value of $p^* \in (0, 1)$, we first generate $u_n$ from a binomial distribution with $n$ trials and a 
probability of success $p^*$ on each trial. 
Next, we generate $u_n$ random variables $X_{i,U}$ for $i = 1, \ldots, u_n$, each with a survival function $\overline{F}_U(x) = x^{-\alpha_1}l_1(x)$, where $l_1$ is slowly varying and $\alpha_1 \in (0, 1)$. Similarly, we generate $n - u_n$ random variables $X_{i,V}$ for $i = 1, \ldots, n - u_n$, each with a survival function $\overline{F}_V(x) = x^{-\alpha_2}l_2(x)$, where $l_2$ is slowly varying and $\alpha_2 > \alpha_1$. Consequently, 
\begin{equation*}
	\lim_{x \rightarrow \infty} \frac{\mathbb{P}(X_{1,V} > x)}{\mathbb{P}(X_{1,U} > x)} = 0.
\end{equation*}
For each simulation setting, we specify the value of $p^*$, the distributions, and the corresponding values of 
$\alpha_1$ and $\alpha_2$ in the captions of Figures \ref{fig:heter_pareto_1}--\ref{fig:heter_inverse_1}.
The straight dashed line in each panel represents the corresponding limit of the smaller value among $\alpha_1$ and $\alpha_2$, in agreement with Theorem \ref{thm:moment_hetero}.
In general,  the points lie close to the line.
Figures \ref{fig:heter_pareto_1} and \ref{fig:heter_stable_1} show the 
simulation results
from heterogeneous heavy-tailed data generated by Pareto and stable distributions, respectively, with $\alpha_1=0.5$ and $\alpha_2=0.8$.
Figures \ref{fig:heter_stable_2} and \ref{fig:heter_inverse_1} show the 
simulation results from heterogeneous heavy-tailed data generated from a stable distribution and  
the distribution function $F_{1,\alpha}$ defined in \eqref{eq:F1}, with $\alpha_1=0.2$ and $\alpha_2=0.7$. 
The asymptotic results depend only on $\alpha_1$, the smaller $\alpha$, in the heterogeneous case. 
Hence there is only one straight dashed  line in each panel of these figures. 
In  all these figures, the points lie close to or  on the line in almost all settings, suggesting that 
the asymptotic theory usefully approximates the behavior of our finite samples.

\subsection{Taylor's Law with Correlated Heavy-tailed Data}
\begin{figure}[h]
	\centering
	\subfloat[Taylor's law]{\includegraphics[width=0.35\textwidth]{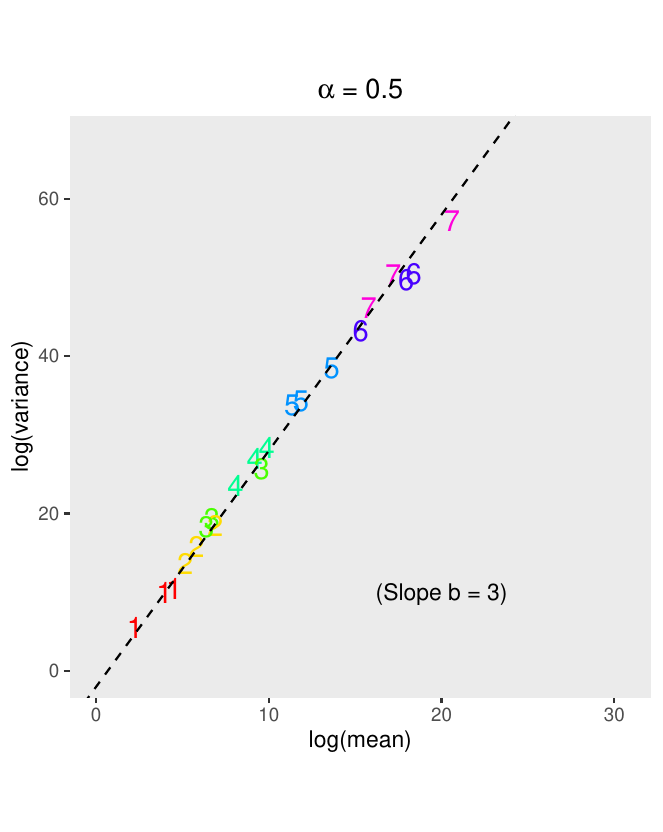}}
	\qquad
	\subfloat[Taylor's law for upper semivariance]{\includegraphics[width=0.35\textwidth]{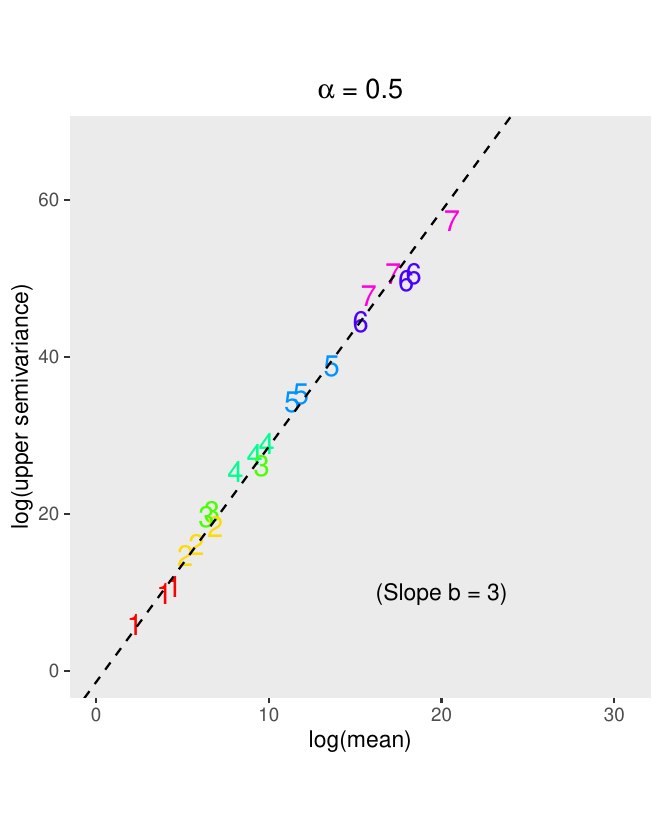}}
	\\
	\subfloat[Taylor's law for third central moment]{\includegraphics[width=0.35\textwidth]{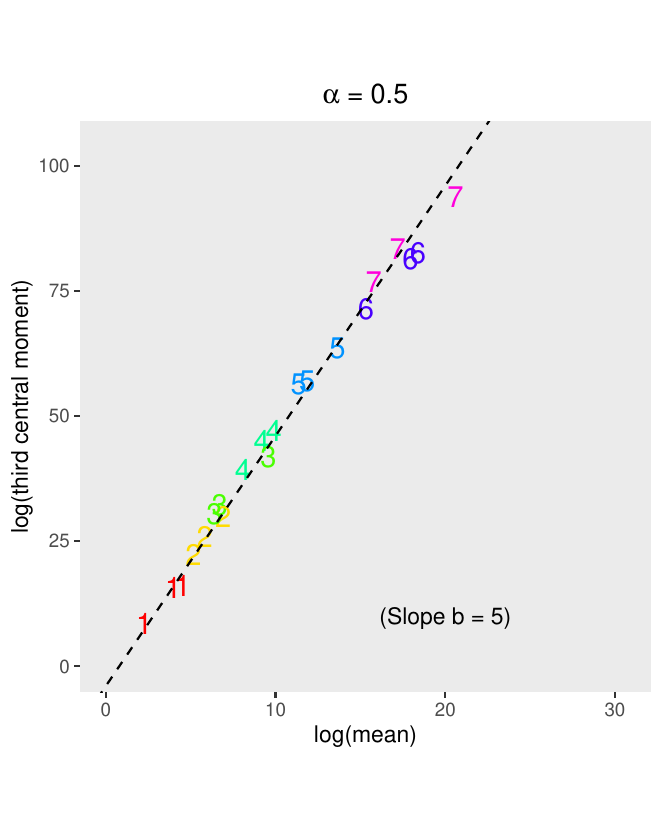}}
	\qquad
	\subfloat[Taylor's law for local upper semivariance]{\includegraphics[width=0.35\textwidth]{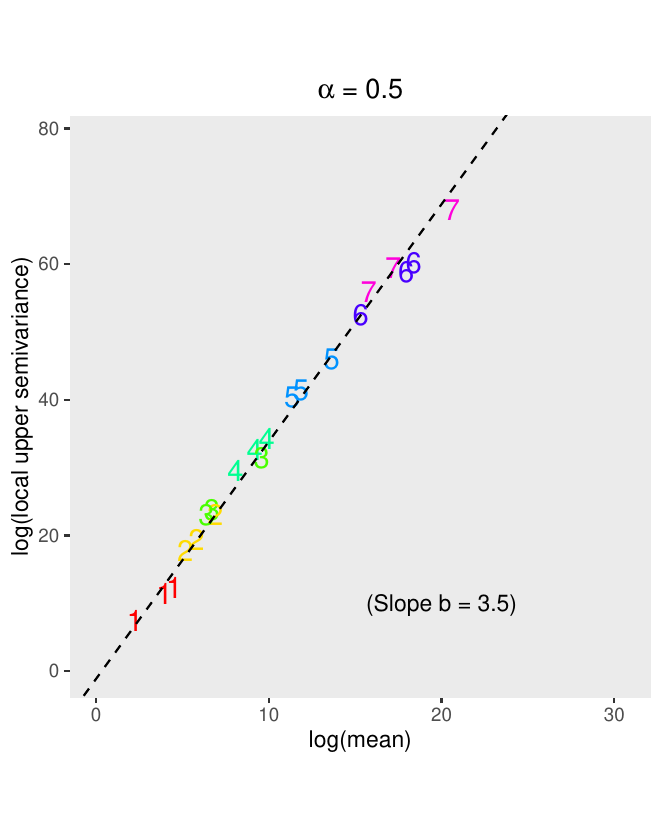}}
	\caption{Scatterplots of (a) $(\log(M_{n,1}),\log(V_n))$, (b) $(\log(M_{n,1}),\log(M^{+}_{n,2}))$, (c) $(\log(M_{n,1}),\log(M^c_{n,3}))$, and (d) $(\log(M_{n,1}),\log(M^{+*}_{n,2}))$, corresponding to Theorem \ref{thm:higher_central_moments} with $h_1=2$ and $h_2=1$, Theorem \ref{thm:upper_central_moment}, Theorem \ref{thm:higher_central_moments} with $h_1=3$ and $h_2=1$, and Theorem \ref{thm:local_upper_central_moment}, respectively. 
		The simulation follows Example \ref{example:correlated_2020} with $\alpha=0.5$ and $\rho=0.1$ and sample size $10^p$ for $p=1,2,\ldots,7$.}
	\label{fig:cor1}
\end{figure}

\begin{figure}[h]
	\centering
	\subfloat[Taylor's law]{\includegraphics[width=0.3\textwidth]{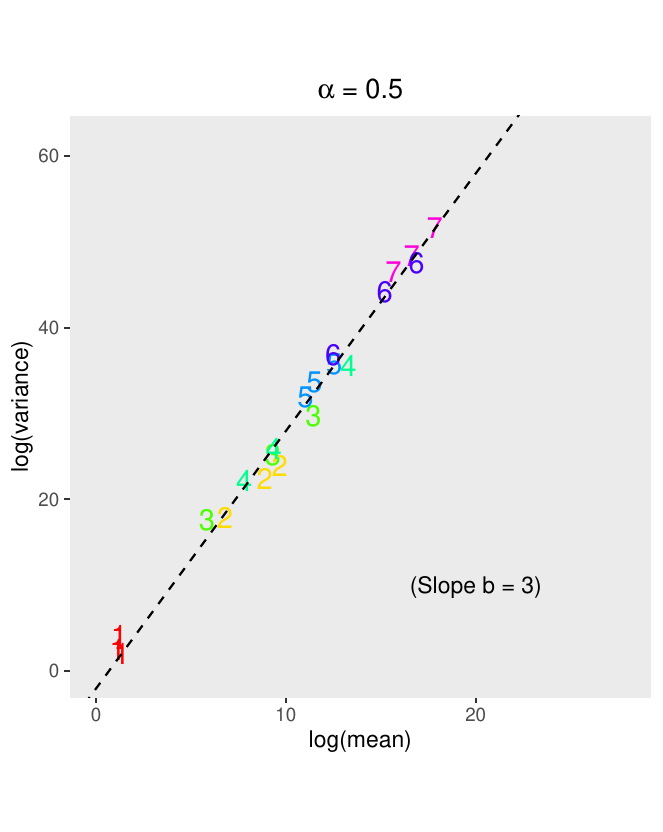}}
	\qquad
	\subfloat[Taylor's law for upper semivariance]{\includegraphics[width=0.3\textwidth]{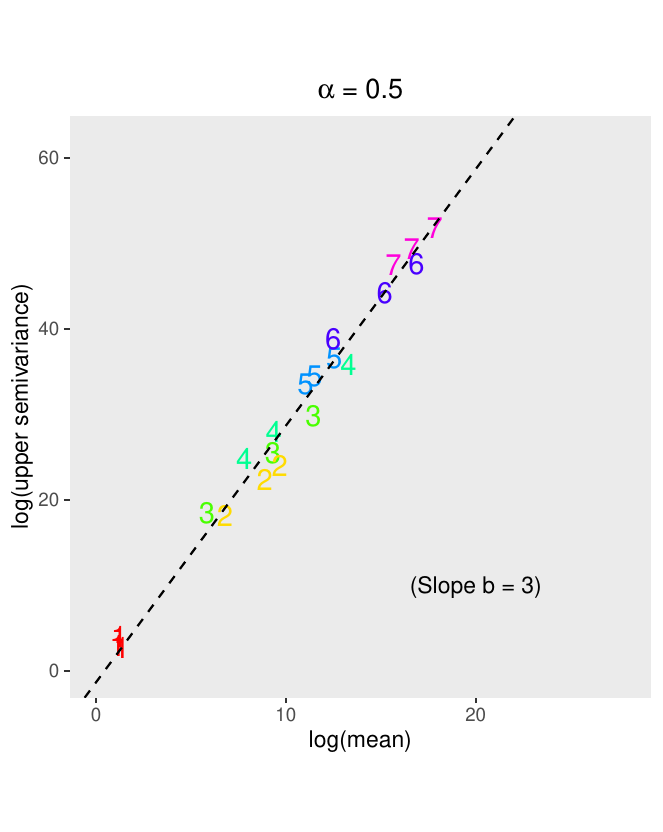}}
	\\
	\subfloat[Taylor's law for third central moment]{\includegraphics[width=0.3\textwidth]{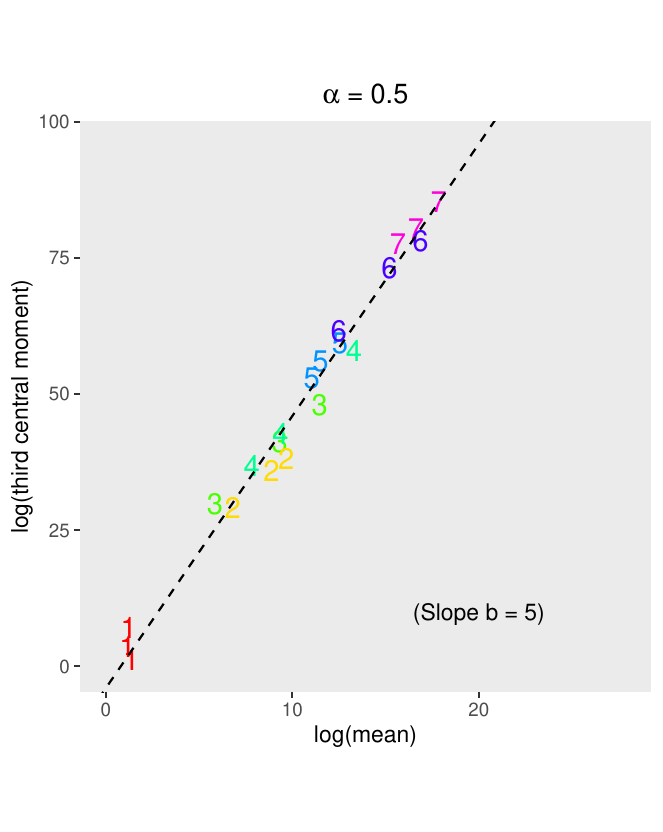}}
	\qquad
	\subfloat[Taylor's law for local upper semivariance]{\includegraphics[width=0.3\textwidth]{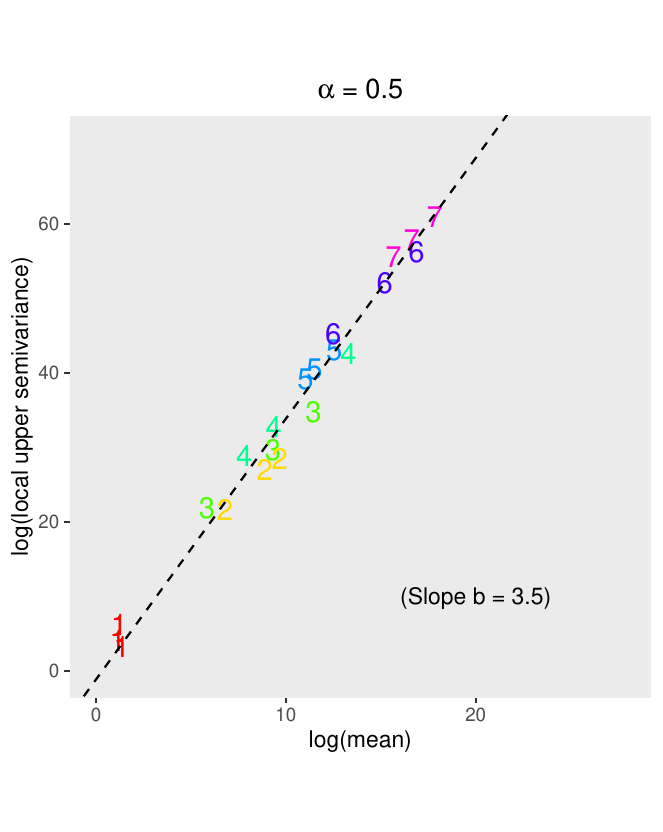}}
	\caption{Scatterplots of (a) $(\log(M_{n,1}),\log(V_n))$, (b) $(\log(M_{n,1}),\log(M^{+}_{n,2}))$, (c) $(\log(M_{n,1}),\log(M^c_{n,3}))$, and (d) $(\log(M_{n,1}),\log(M^{+*}_{n,2}))$, corresponding to Theorem \ref{thm:higher_central_moments} with $h_1=2$ and $h_2=1$, Theorem \ref{thm:upper_central_moment}, Theorem \ref{thm:higher_central_moments} with $h_1=3$ and $h_2=1$, and Theorem \ref{thm:local_upper_central_moment}, respectively. 
		The simulation follows Example \ref{example:correlated_2022} with $\alpha=0.5$ and  $\gamma(X_i,X_j)=\exp(-|i - j| / 100)$} and sample size $10^p$ for $p=1,2,\ldots,7$. 
	\label{fig:cor2}
\end{figure}

Simulations of correlated heavy-tailed data  specified 
in Examples \ref{example:correlated_2020} and \ref{example:correlated_2022}
yielded  Figures \ref{fig:cor1} to \ref{fig:cor2}.
The asymptotic ratio of moments describes the finite-sample simulations well.

\subsection{Taylor's Law for Heavy-tailed Network Data}\label{subsec:networksimulation}
\begin{figure}[h]
	\centering
	\subfloat[Example 1]{\includegraphics[width=0.45\textwidth]{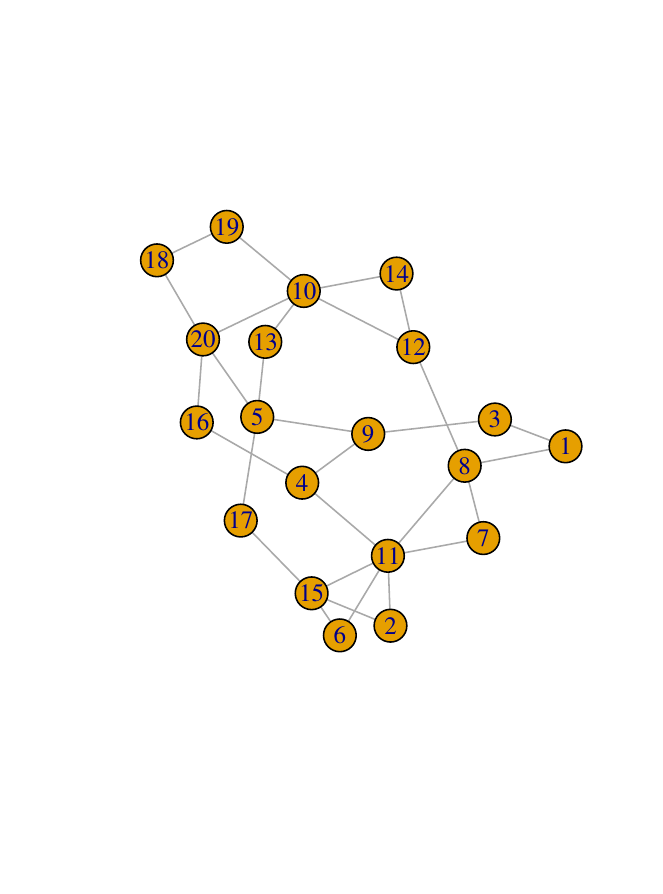}}
	\qquad
	\subfloat[Example 2]{\includegraphics[width=0.45\textwidth]{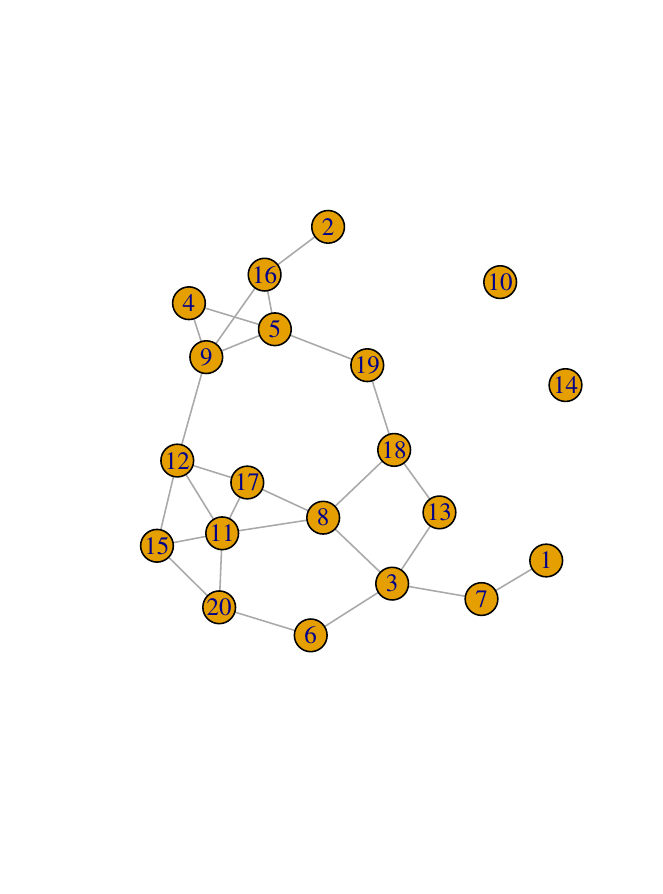}}
	\caption{Two subgraphs each containing 20 nodes from the network generated by the random graph}
	\label{figure:network_example}
\end{figure}

{Here we illustrate Taylor's law using network data. We  generate random graphs based 
	on the Erd\H{o}s–R\'enyi–Gilbert model \citep{gilbert1959random}. 
	In this model, each edge has a fixed probability of being present or absent, independently of the other edges. These random graphs are generated with the \texttt{erdos.renyi.game} function from the R package \textbf{igraph}. The edges for each node are generated through independent Bernoulli trials with a success probability $p$, that is, each pair of nodes has a probability $p$ of being joined by an edge. We  illustrate these random graphs with  two examples in Figure \ref{figure:network_example}, each containing 20 nodes with $p = 0.2$.

	Define $\Xi_{j,1}^n$ to be the set of nodes in the graph that are at distance 1 from node $j$, i.e., $\Xi_{j,1}^n = \{k : d(j,k) = 1, 1 \leq k \leq n\}$, where the distance between two nodes is defined as the length of the shortest path connecting them. We associate each node with a random variable defined via:
	\begin{align}\label{eq:network_sim}
		X_j := Z_j + \frac{\sum_{k \in \Xi_{j,1}^n} Z_k}{|\Xi_{j,1}^n|},
	\end{align}
	where $Z_j$ are the i.i.d. heavy-tailed random variables. We consider $n = 100, 500, 1000, 5000$, $10000$ with $p = \frac{10}{n-1}$ for $n = 100$ and $500$, and $p = \frac{100}{n-1}$ for $n = 1000, 5000$, and $10000$. 
	We also limit the maximum number of edges to $20$ when $n = 100$ or $500$ and to  $200$ when $n = 1000, 5000$, or $10000$. 
	These bounds are larger than twice the expected number of edges. 
	The actual number of edges is unlikely to exceed these bounds for large $n$. 
	If the number of simulated edges exceeds these bounds, we regenerate another random graph until the constraint is satisfied.

	Under this setting, if $d(X_i, X_j) \geq 3$, then $\Cov(X_i, X_j) = 0$, since $X_i = Z_i + \frac{\sum_{k \in \Xi^n_{i,1}} Z_k}{|\Xi^n_{i,1}|}$ and $X_j = Z_j + \frac{\sum_{k \in \Xi^n_{j,1}} Z_k}{|\Xi^n_{j,1}|}$ have no overlapping $Z_k$ for $k \in \mathbb{N}$. On the other hand, we restrict $|\Xi^n_{i,1}| \leq K$ for a positive constant $K$ across different cases. Thus, this setting satisfies the two conditions in Example \ref{Example:network}, and hence the covariance condition (\ref{eq:network_cov_assumption}) holds.
	
	If each edge $(j, k)$ were associated with a positive weight, $w(j, k) > 0$, then a natural generalization would be to define
	\begin{align*}
		X_j := Z_j + \frac{\sum_{k \in \Xi_{j,1}^n} w(j, k) Z_k}{\sum_{k \in \Xi_{j,1}^n} w(j, k)}.
	\end{align*}
	
}

\begin{figure}[h]
	\centering
	\subfloat[Taylor's law]{\includegraphics[width=0.3\textwidth]{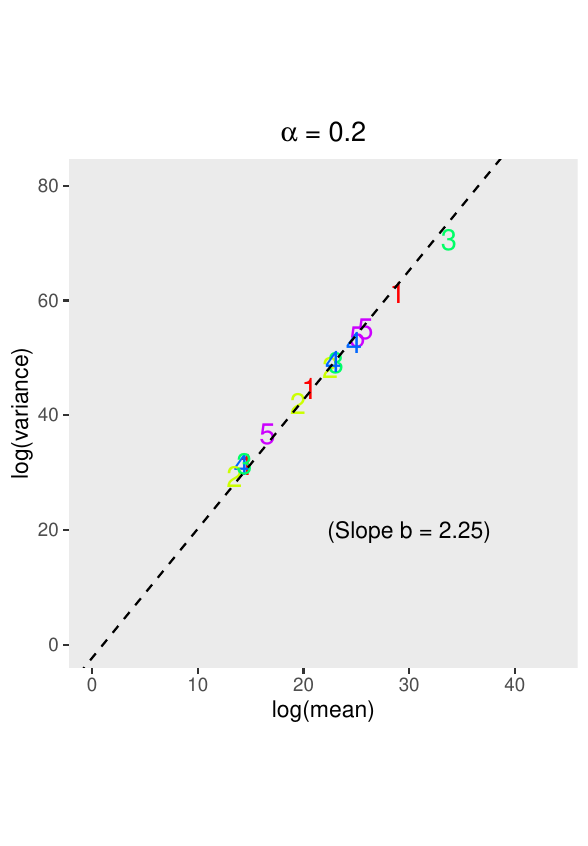}}
	\qquad
	\subfloat[Taylor's law for upper semivariance]{\includegraphics[width=0.3\textwidth]{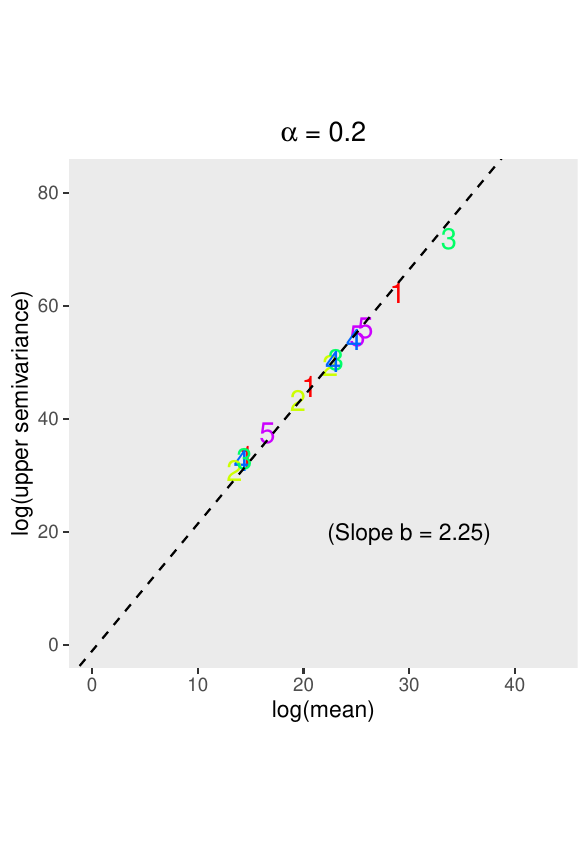}}
	\\
	\subfloat[Taylor's law for third central moment]{\includegraphics[width=0.3\textwidth]{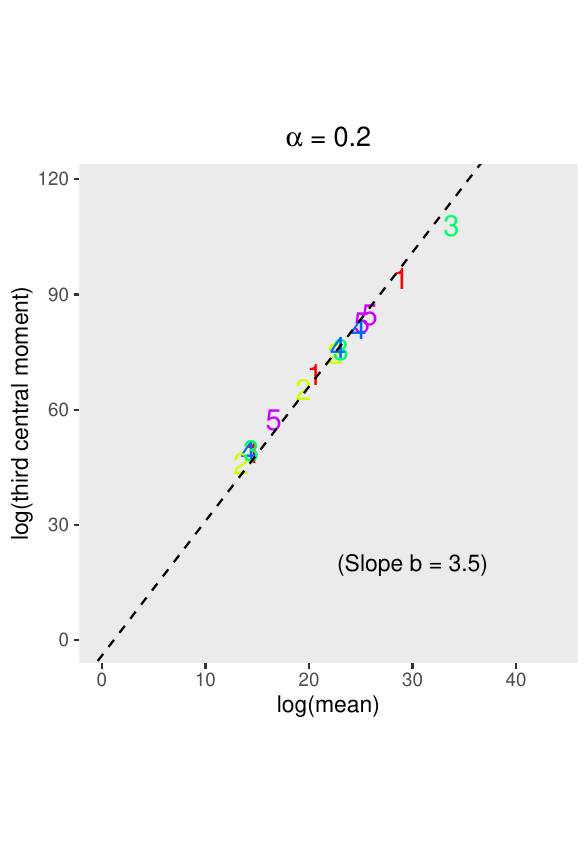}}
	\qquad
	\subfloat[Taylor's law for local upper semivariance]{\includegraphics[width=0.3\textwidth]{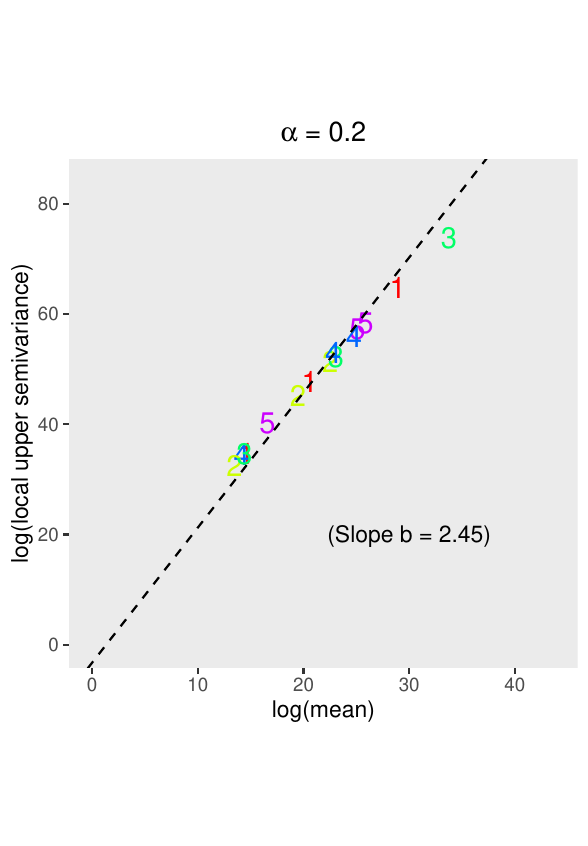}}
	\caption{Scatterplots of (a) $(\log(M_{n,1}),\log(V_n))$, 
		(b) $(\log(M_{n,1}),\log(M^{+}_{n,2}))$, (c) $(\log(M_{n,1}),\log(M^c_{n,3}))$, and (d) $(\log(M_{n,1}),\log(M^{+*}_{n,2}))$, corresponding to Theorem \ref{thm:higher_central_moments} with $h_1=2$ and $h_2=1$, Theorem \ref{thm:upper_central_moment}, Theorem \ref{thm:higher_central_moments} with $h_1=3$ and $h_2=1$, and Theorem \ref{thm:local_upper_central_moment}, respectively, when $Z_i \overset{\mathrm{i.i.d.}}{\sim}$ Pareto$(0.5,0.2)$ with a sample size of $100, 500, 1000, 5000, 10000$, which are labeled as ``1'', $\ldots$, ``5'' in the figures, respectively.}
	\label{fig:network_1}
\end{figure}

\begin{figure}[h]
	\centering
	\subfloat[Taylor's law]{\includegraphics[width=0.3\textwidth]{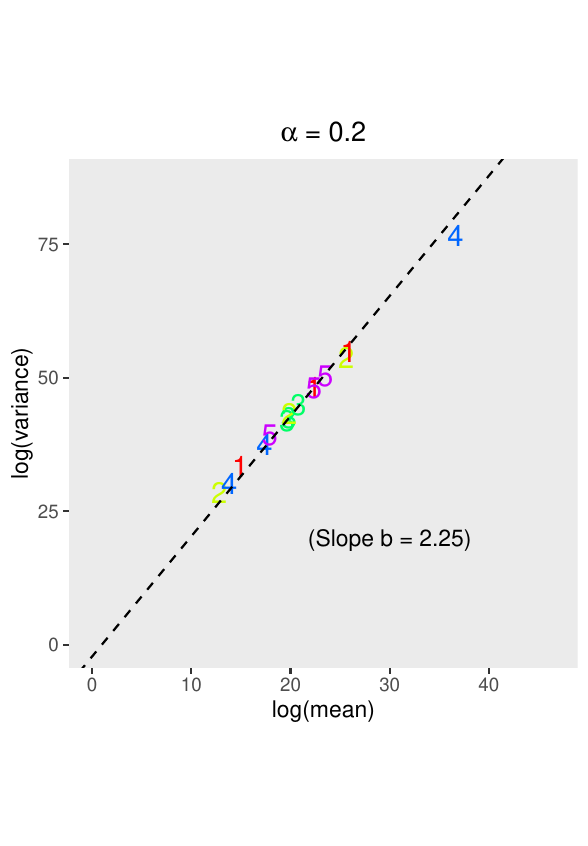}}
	\qquad
	\subfloat[Taylor's law for upper semivariance]{\includegraphics[width=0.3\textwidth]{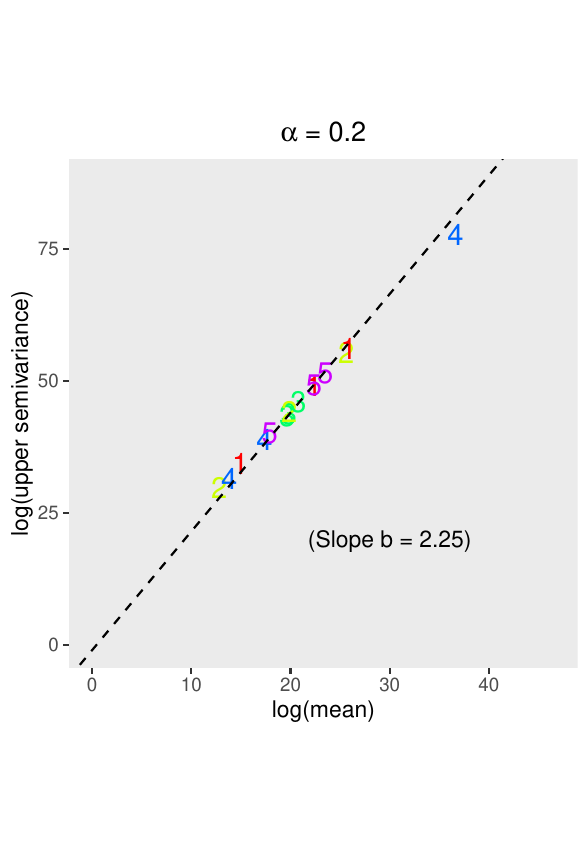}}
	\\
	\subfloat[Taylor's law for third central moment]{\includegraphics[width=0.3\textwidth]{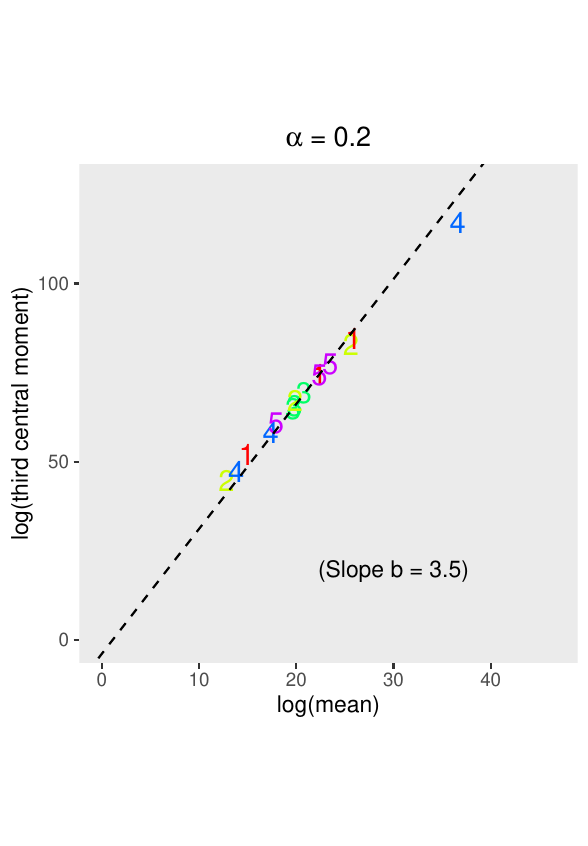}}
	\qquad
	\subfloat[Taylor's law for local upper semivariance]{\includegraphics[width=0.3\textwidth]{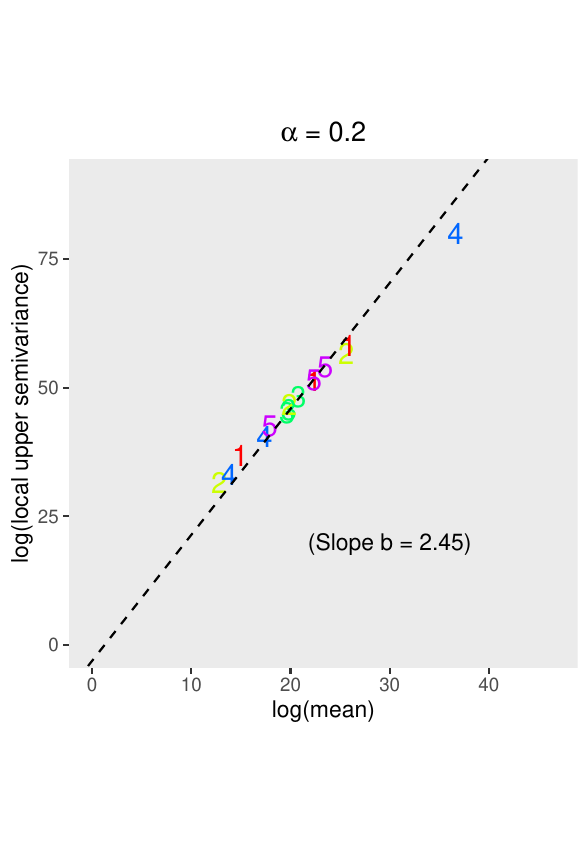}}
	\caption{Scatterplots of (a) $(\log(M_{n,1}),\log(V_n))$, (b) $(\log(M_{n,1}),\log(M^{+}_{n,2}))$, (c) $(\log(M_{n,1}),\log(M^c_{n,3}))$, and (d) $(\log(M_{n,1}),\log(M^{+*}_{n,2}))$, corresponding to Theorem \ref{thm:higher_central_moments} with $h_1=2$ and $h_2=1$, Theorem \ref{thm:upper_central_moment}, Theorem \ref{thm:higher_central_moments} with $h_1=3$ and $h_2=1$, and Theorem \ref{thm:local_upper_central_moment}, respectively, when $Z_i \overset{\mathrm{i.i.d.}}{\sim}F(1,0.2)$ with a sample size of $100, 500, 1000, 5000, 10000$, which are labeled as ``1'', $\ldots$, ``5'' in the figures, respectively.}
	\label{fig:network_2}
\end{figure}

\begin{figure}[h]
	\centering
	\subfloat[Taylor's law]{\includegraphics[width=0.3\textwidth]{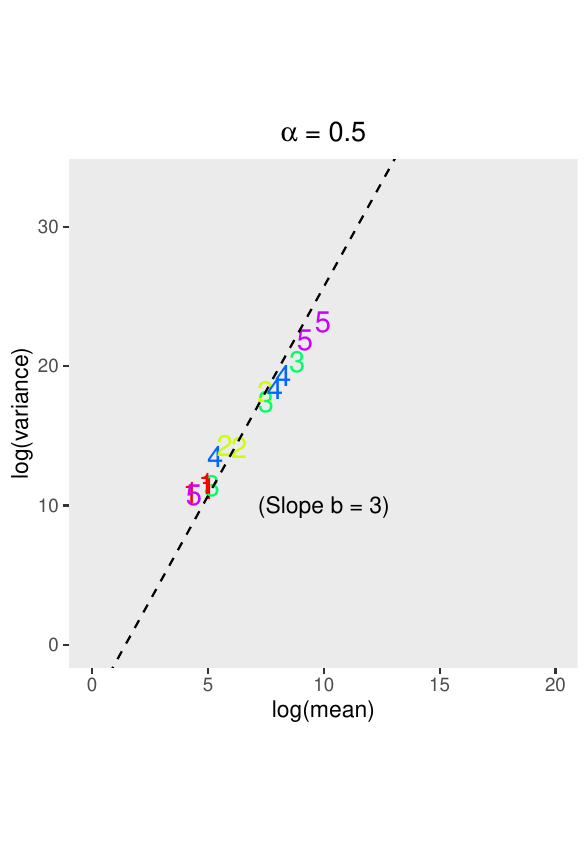}}
	\qquad
	\subfloat[Taylor's law for upper semivariance]{\includegraphics[width=0.3\textwidth]{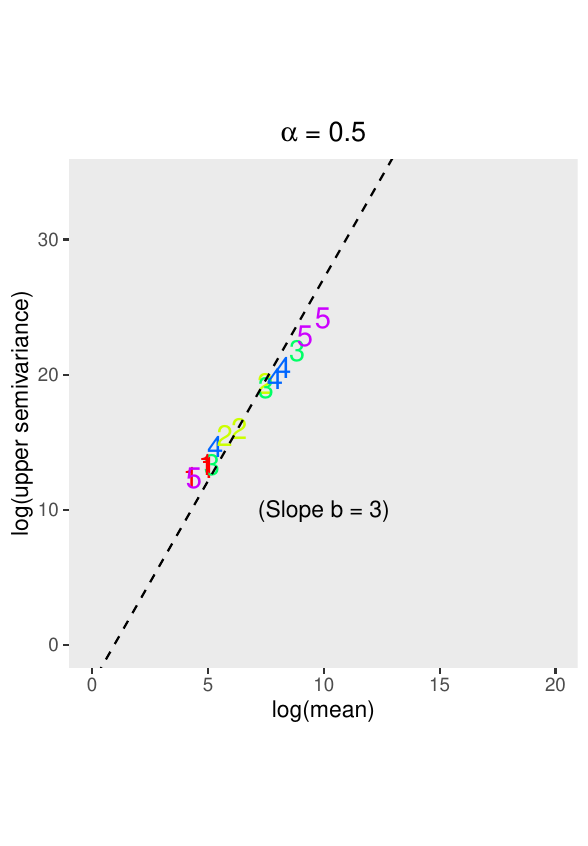}}
	\\
	\subfloat[Taylor's law for third central moment]{\includegraphics[width=0.3\textwidth]{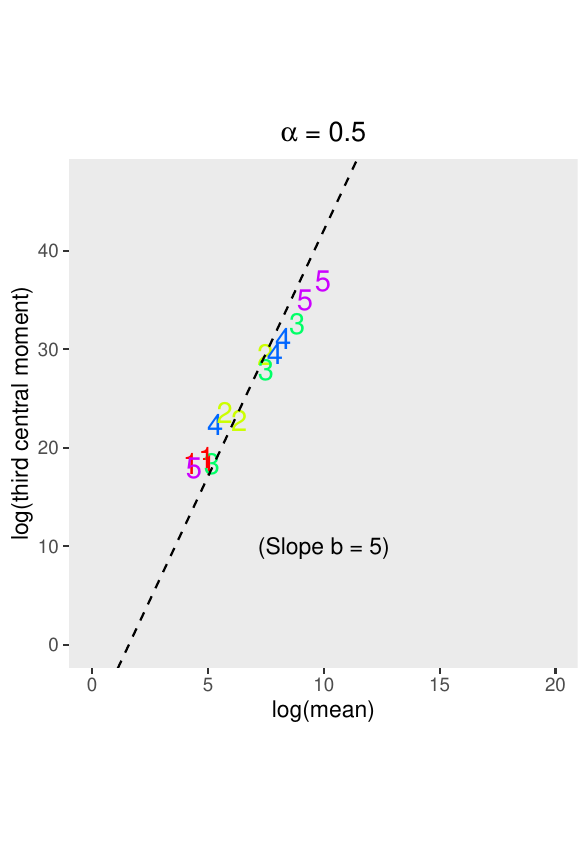}}
	\qquad
	\subfloat[Taylor's law for local upper semivariance]{\includegraphics[width=0.3\textwidth]{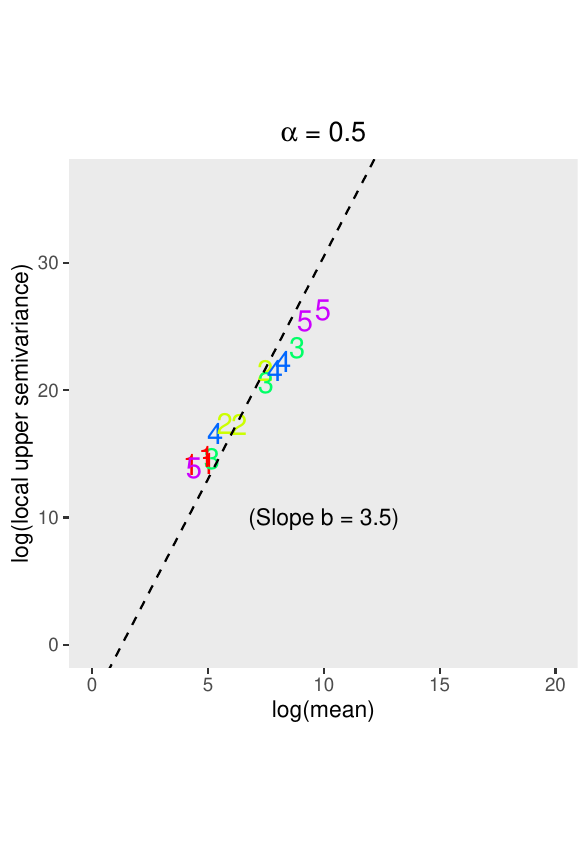}}
	\caption{Scatterplots of (a) $(\log(M_{n,1}),\log(V_n))$, (b) $(\log(M_{n,1}),\log(M^{+}_{n,2}))$, (c) $(\log(M_{n,1}),\log(M^c_{n,3}))$, and (d) $(\log(M_{n,1}),\log(M^{+*}_{n,2}))$, corresponding to Theorem \ref{thm:higher_central_moments} with $h_1=2$ and $h_2=1$, Theorem \ref{thm:upper_central_moment}, Theorem \ref{thm:higher_central_moments} with $h_1=3$ and $h_2=1$, and Theorem \ref{thm:local_upper_central_moment}, respectively, when $Z_i \overset{\mathrm{i.i.d.}}{\sim}$ Pareto$(0.5,0.5)$ with a sample size of $100, 500, 1000, 5000, 10000$, which are labeled as ``1'', $\ldots$, ``5'' in the figures, respectively.}
	\label{fig:network_3}
\end{figure}

\begin{figure}[h]
	\centering
	\subfloat[Taylor's law]{\includegraphics[width=0.35\textwidth]{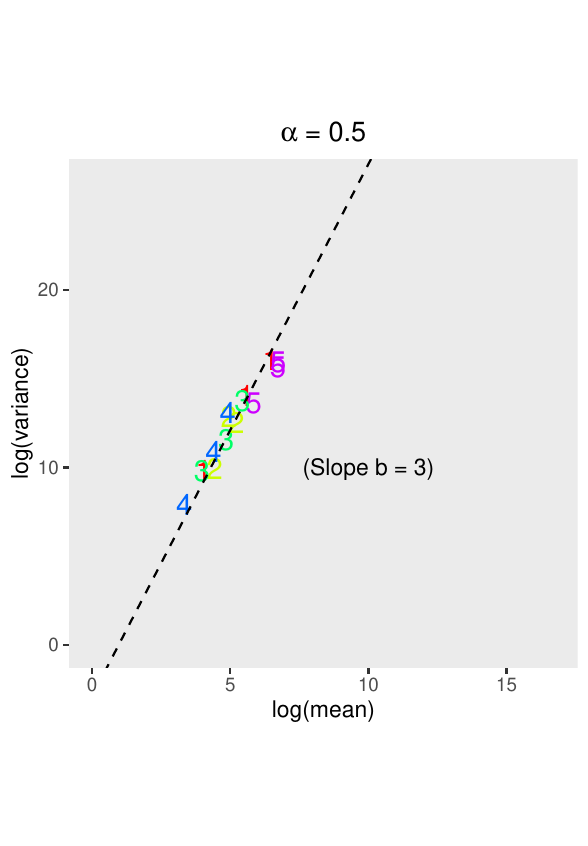}}
	\qquad
	\subfloat[Taylor's law for upper semivariance]{\includegraphics[width=0.35\textwidth]{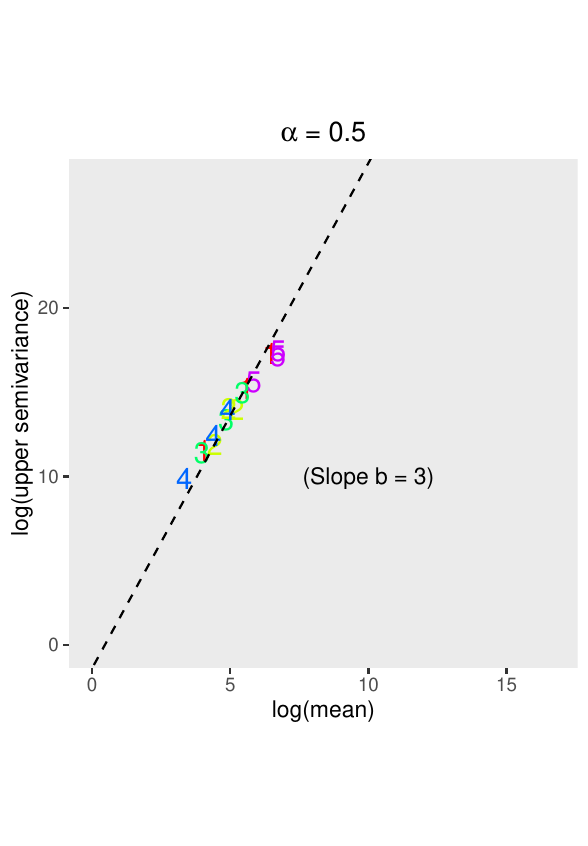}}
	\\
	\subfloat[Taylor's law for third central moment]{\includegraphics[width=0.35\textwidth]{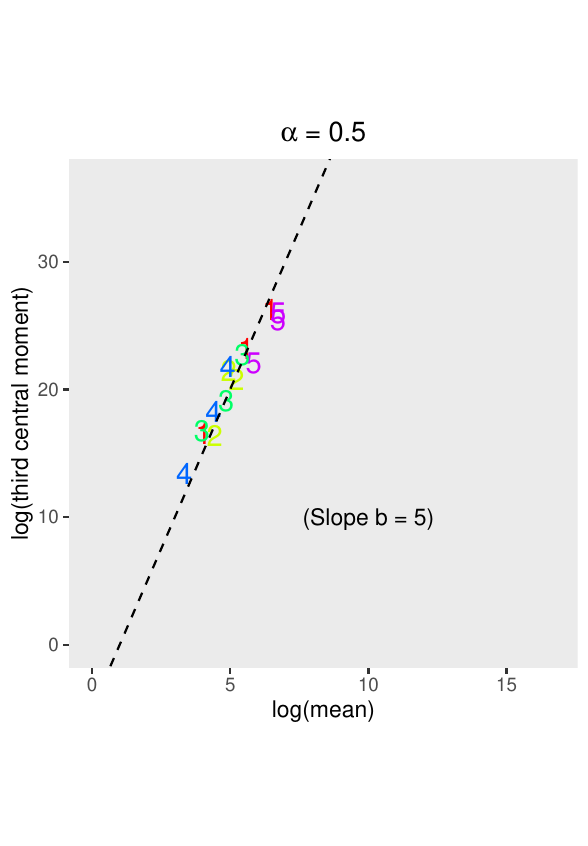}}
	\qquad
	\subfloat[Taylor's law for local upper semivariance]{\includegraphics[width=0.35\textwidth]{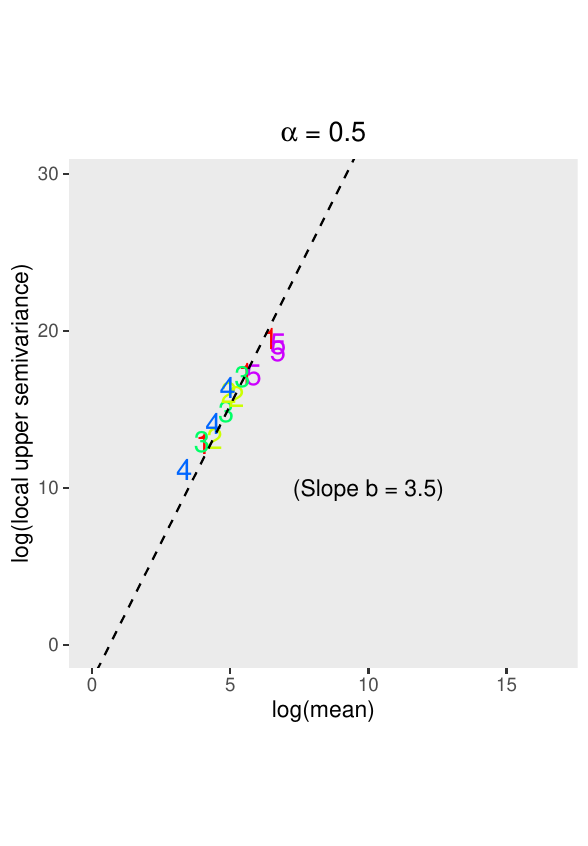}}
	\captionsetup{font=small}
	\caption{Scatterplots of (a) $(\log(M_{n,1}),\log(V_n))$, (b) $(\log(M_{n,1}),\log(M^{+}_{n,2}))$, (c) $(\log(M_{n,1}),\log(M^c_{n,3}))$, and (d) $(\log(M_{n,1}),\log(M^{+*}_{n,2}))$, corresponding to Theorem \ref{thm:higher_central_moments} with $h_1=2$ and $h_2=1$, Theorem \ref{thm:upper_central_moment}, Theorem \ref{thm:higher_central_moments} with $h_1=3$ and $h_2=1$, and Theorem \ref{thm:local_upper_central_moment}, respectively, when $Z_i \overset{\mathrm{i.i.d.}}{\sim} F(1,0.5)$ with a sample size of $100, 500, 1000, 5000, 10000$, which are labeled as ``1'', $\ldots$, ``5'' in the figures, respectively.}
	\label{fig:network_4}
\end{figure}

\begin{figure}[h]
	\centering
	\subfloat[Taylor's law]{\includegraphics[width=0.3\textwidth]{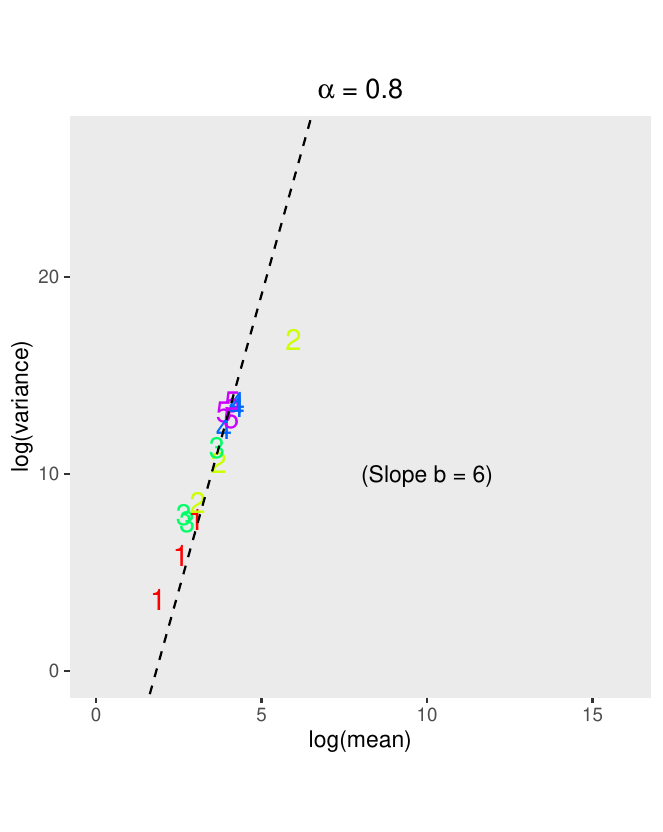}}
	\qquad
	\subfloat[Taylor's law for upper semivariance]{\includegraphics[width=0.3\textwidth]{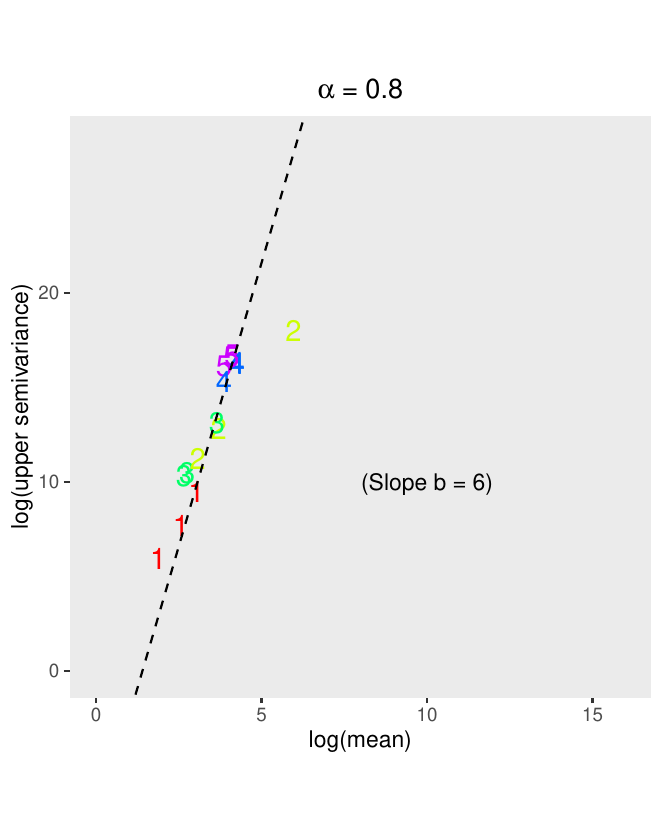}}
	\\
	\subfloat[Taylor's law for third central moment]{\includegraphics[width=0.3\textwidth]{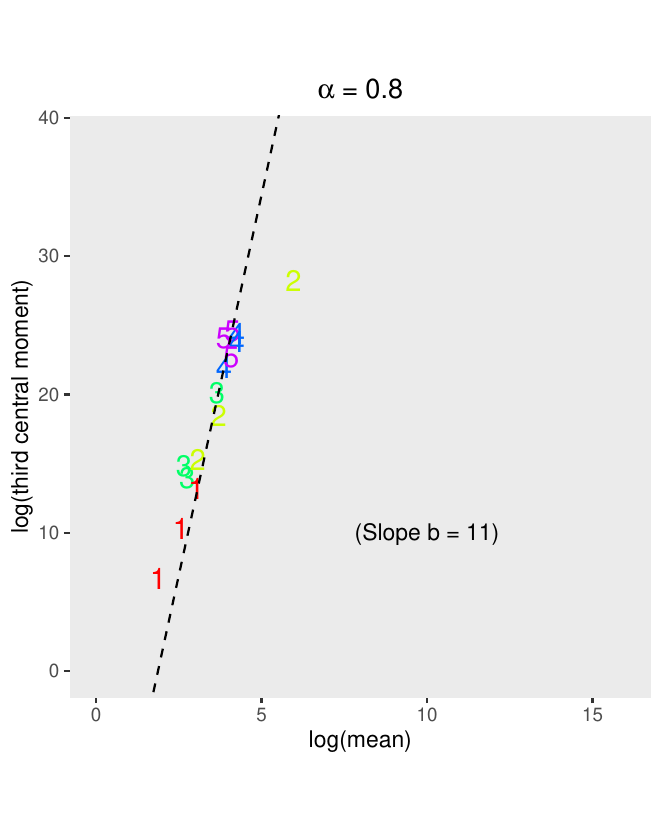}}
	\qquad
	\subfloat[Taylor's law for local upper semivariance]{\includegraphics[width=0.3\textwidth]{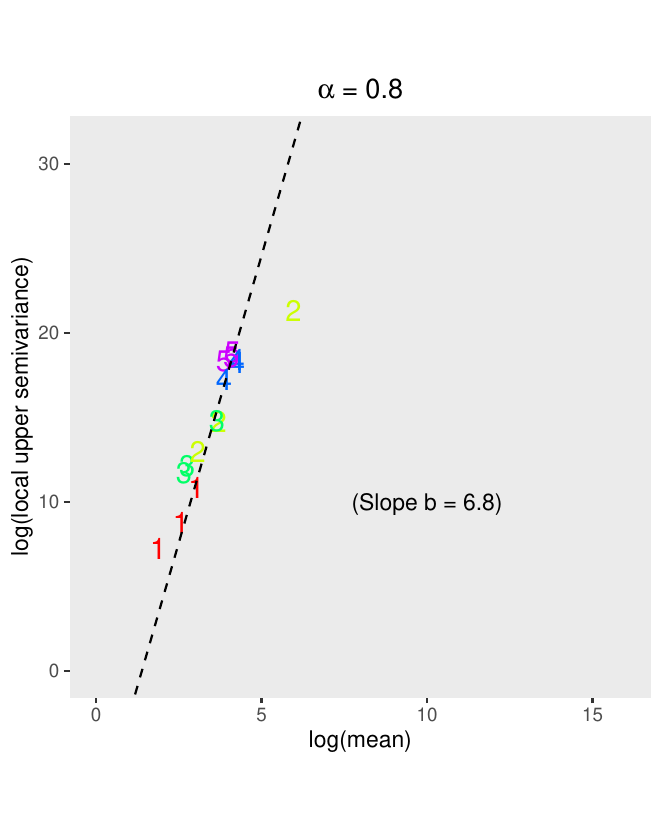}}

	\caption{Scatterplots of (a) $(\log(M_{n,1}),\log(V_n))$, (b) $(\log(M_{n,1}),\log(M^{+}_{n,2}))$, (c) $(\log(M_{n,1}),\log(M^c_{n,3}))$, and (d) $(\log(M_{n,1}),\log(M^{+*}_{n,2}))$, corresponding to Theorem \ref{thm:higher_central_moments} with $h_1=2$ and $h_2=1$, Theorem \ref{thm:upper_central_moment}, Theorem \ref{thm:higher_central_moments} with $h_1=3$ and $h_2=1$, and Theorem \ref{thm:local_upper_central_moment}, respectively, when $Z_i \overset{\mathrm{i.i.d.}}{\sim}$Pareto$(0.5,0.8)$ with a sample size of $100, 500, 1000, 5000, 10000$ which are labeled as ``1'', $\ldots$, ``5'' in the figures, respectively.}
	\label{fig:network_5}
\end{figure}

\begin{figure}[h]
	\centering
	\subfloat[Taylor's law]{\includegraphics[width=0.3\textwidth]{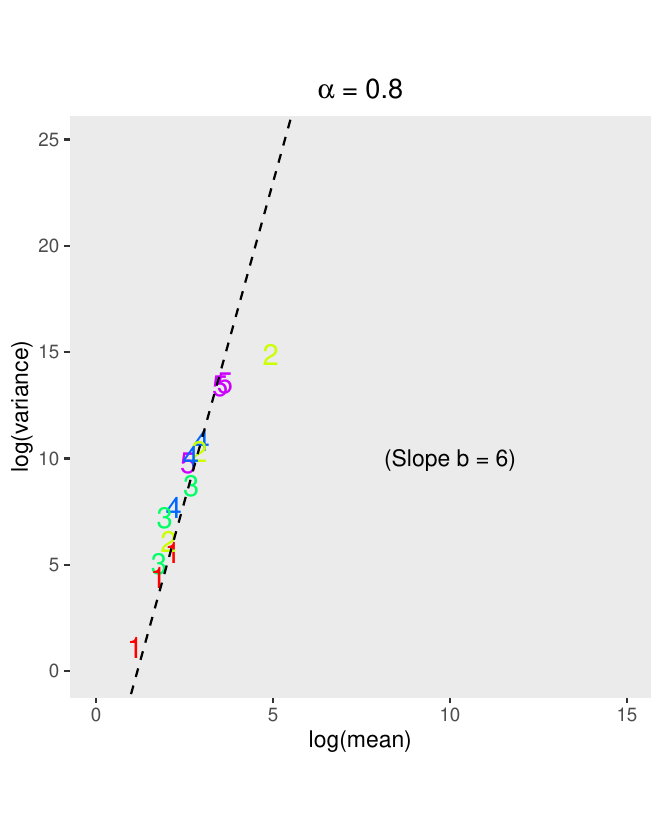}}
	\qquad
	\subfloat[Taylor's law for upper semivariance]{\includegraphics[width=0.3\textwidth]{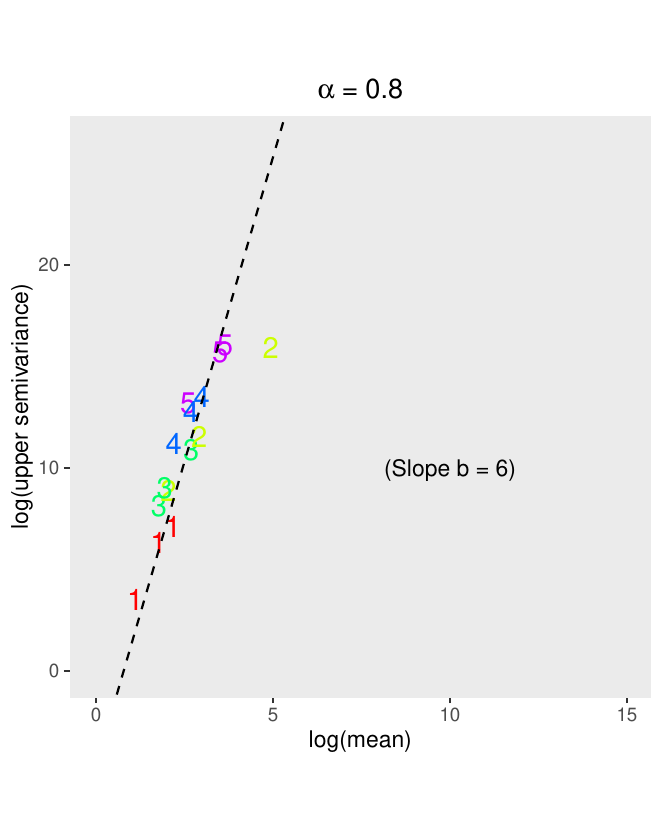}}
	\\
	\subfloat[Taylor's law for third central moment]{\includegraphics[width=0.3\textwidth]{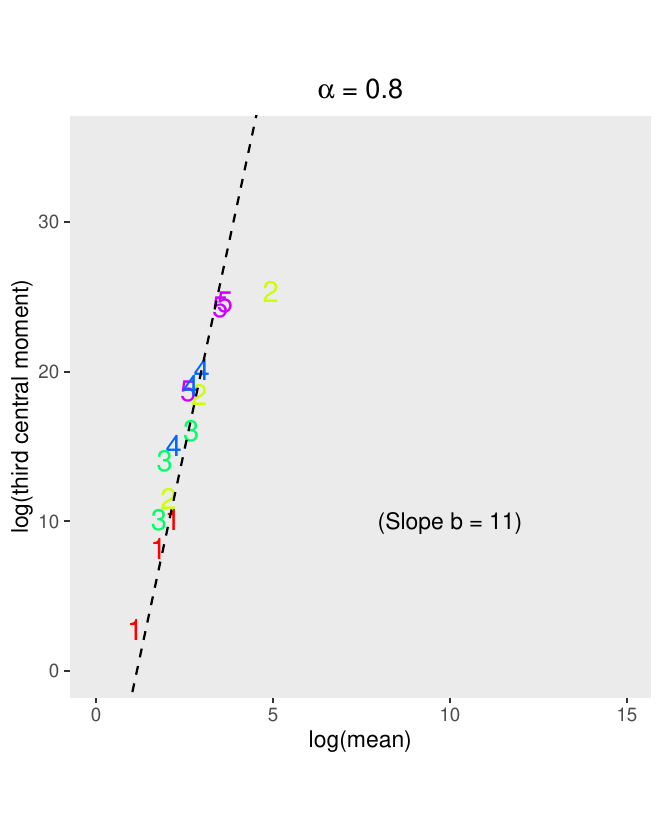}}
	\qquad
	\subfloat[Taylor's law for local upper semivariance]{\includegraphics[width=0.3\textwidth]{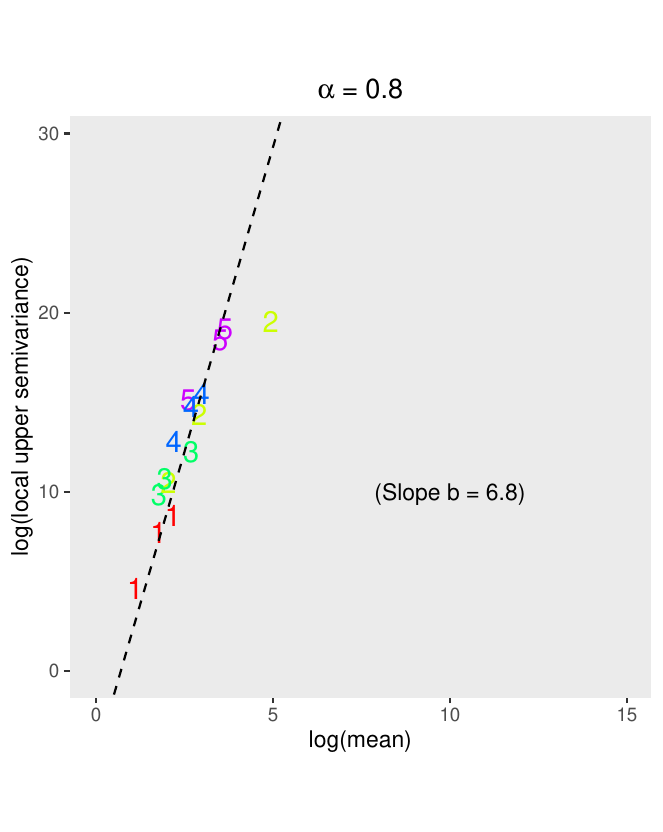}}
	
	\caption{Scatterplots of (a) $(\log(M_{n,1}),\log(V_n))$, (b) $(\log(M_{n,1}),\log(M^{+}_{n,2}))$, (c) $(\log(M_{n,1}),\log(M^c_{n,3}))$, and (d) $(\log(M_{n,1}),\log(M^{+*}_{n,2}))$, corresponding to Theorem \ref{thm:higher_central_moments} with $h_1=2$ and $h_2=1$, Theorem \ref{thm:upper_central_moment}, Theorem \ref{thm:higher_central_moments} with $h_1=3$ and $h_2=1$, and Theorem \ref{thm:local_upper_central_moment}, respectively, when $Z_i \overset{\mathrm{i.i.d.}}{\sim} F(1,0.8)$ with a sample size of $100, 500, 1000, 5000, 10000$, which are labeled as ``1'', $\ldots$, ``5'' in the figures, respectively.}
	\label{fig:network_6}
\end{figure}

Figures \ref{fig:network_1}--\ref{fig:network_6} illustrate Taylor's laws when the heavy-tailed data are connected by a network graph, and $Z_i$'s are generated from Pareto or stable distributions, with $\alpha = 0.2, 0.5,$ and $0.8$. In general, under these network settings, the theoretical asymptotic results of Taylor's laws fit the finite-sample heavy-tailed network data quite well.

\FloatBarrier

\section{Proofs for Section \ref{sec:main results}}

\subsection{Proofs for Section \ref{subsect:auxiliary}}

We first state two well-known results in probability theory.
\begin{lemma}\label{lemma:expectation of y^p}
	Let $Y \geq 0$ be a nonnegative random variable and $p > 0$. 
	Then $\mathbb{E}(Y^p)=\int_0^\infty py^{p-1}\mathbb{P}(Y>y)dy$.
\end{lemma}
For a proof of Lemma \ref{lemma:expectation of y^p}, see \cite{chakraborti2019higher} or \cite{liu2020general}.

\begin{proposition}[Karamata's Theorem, e.g., Proposition 1.5.8 in \cite{bingham1989regular}]\label{prop:Karamata}
	If the function $l$ is slowly varying and $l$ is locally bounded on $[c_0, \infty)$ and if $\alpha > -1$, then
	\begin{equation*}
		\int^x_{c_0} t^\alpha l(t)dt \sim \frac{x^{\alpha+1}}{\alpha+1}l(x).
	\end{equation*}
\end{proposition}
\begin{proof}[Proof of Lemma \ref{lemma:truncate_moment_asy}]
	\begin{enumerate}[(a)]
		\item The survival function of $\tilde{X}$ for any $x \in \mathbb{R}$ is
		\begin{equation*}
			\mathbb{P}(\tilde{X} >x) = \left[ 1 - \frac{F(x)}{F(t_n)} \right] \mathbbm{1}(x < t_n) = \left[ \frac{\overline{F}(x) - \overline{F}(t_n)}{F(t_n)}\right] \mathbbm{1}(x < t_n).
		\end{equation*}
		Thus,  by Lemma \ref{lemma:expectation of y^p},
		\begin{align*}
			\mathbb{E}[\tilde{X}^p] &=  p \int^\infty_0 y^{p-1}\mathbb{P}(\tilde{X} > y) dy = \frac{p}{F(t_n)} \int^{t_n}_0 (y^{p-1} \overline{F}(y) - y^{p-1}\overline{F}(t_n))dy \\
			&= \frac{p}{F(t_n)} \left[ \int^{t_n}_0 y^{p-1}\overline{F}(y)dy - \frac{t_n^p}{p} \overline{F}(t_n)\right].
		\end{align*}
		Since $\overline{F}(y) = y^{-\alpha}l(y)$,
		\begin{align*}
			\mathbb{E}[\tilde{X}^p] = 
			\frac{p}{F(t_n)} \left[ \int^{t_n}_0 y^{p-1-\alpha}l(y) dy - \frac{t_n^{p-\alpha}}{p}l(t_n) \right].
		\end{align*}
		By Proposition \ref{prop:Karamata}, for $p - 1- \alpha > -1$, or equivalently, $p > \alpha$,
		\begin{equation*}
			\int^{t_n}_0 y^{p-1-\alpha}l(y) dy \sim \frac{t_n^{p-\alpha}l(t_n)}{p-\alpha}, \quad \text{as } n \rightarrow \infty.
		\end{equation*}
		Thus, as $n \rightarrow \infty$, 
		\begin{equation*}
			\frac{1}{t_n^{p-\alpha}l(t_n)}\left[   \int^{t_n}_0 y^{p-1-\alpha}l(y) dy - \frac{t_n^{p-\alpha}}{p}l(t_n)  \right] \rightarrow \frac{1}{p-\alpha} - \frac{1}{p} = \frac{\alpha}{p(p-\alpha)}.
		\end{equation*}
		As $\overline{F}(t_n) \rightarrow 0$, we have
		\begin{equation*}
			\mathbb{E}[\tilde{X}^p]  \sim \frac{\alpha }{p-\alpha} t_n^{p-\alpha}l(t_n), \quad \text{as } n \rightarrow \infty.
		\end{equation*}
		
		\item Part (a) implies that
		\begin{equation*}
			\log \mathbb{E}[\tilde{X}^p] - \log \left[ \frac{\alpha }{p-\alpha} t_n^{p-\alpha}l(t_n)\right] \rightarrow 0, \quad \text{as } n \rightarrow \infty.
		\end{equation*}
		Thus, as $n \rightarrow \infty$,
		\begin{equation*}
			\frac{\log \mathbb{E}[\tilde{X}^p]}{\log t_n} \rightarrow p - \alpha,
		\end{equation*}
		as  clearly $\frac{\log l(t_n)}{\log t_n} \rightarrow 0$; also see Proposition 1.3.6 (i) in \cite{bingham1989regular}.
	\end{enumerate}

\end{proof}

\begin{proof}[Proof of Lemma \ref{lemma:existence_of_tn}]
	Let
	\begin{equation*}
		g_n(x) := \frac{n c_n l(x)}{x^\alpha} - 1.
	\end{equation*}
	Recall that $n c_n \uparrow \infty$. Consider a sufficiently large $n$ such that $nc_n>\frac{2}{\overline{F}_0(0+)}$. For small enough $x$, since $\overline{F}_0(0+)>0$, we have $g_n(x) > 0$. On the other hand, by Proposition 1.3.6 (v) in \cite{bingham1989regular}, $l(x)/x^\alpha \rightarrow 0$ as $x \rightarrow \infty$. Thus,
	\begin{equation*}
		\lim_{x \rightarrow \infty}g_n(x) = -1.
	\end{equation*}
	Since $g_n$ is continuous, there exists a root $t_n$ such that $g_n(t_n) = 0$.
\end{proof}

\subsection{Proofs for Section \ref{subsec:cov condition}}
\begin{proof}[Proof of Theorem \ref{thm:m-dependence}]
	Fix $p > \alpha$. First, for any $i, j$,
	\begin{align*}
		\left|\Cov(\breve{X}^p_i,\breve{X}^p_j)\right| \leq \sqrt{\Var(\breve{X}^p_i)\Var(\breve{X}^p_j)}=\Var(\breve{X}^p_1) \leq  \mathbb{E}(\breve{X}^{2p}_1).
	\end{align*}
	Then, as $\{X_i,i\geq1\}$ is a $m$-dependent sequence,
	\begin{align*}
		\frac{\sum_{i\neq j} \Cov(\breve{X}^p_i,\breve{X}^p_j)}{v^{2p}_n c_n^2}
		\leq \frac{2 n m  \mathbb{E}(\breve{X}^{2p}_1)}{v^{2p}_n c_n^2} \sim \frac{2nm \frac{\alpha}{2p-\alpha} v_n^{2p-\alpha} l(v_n)}{v_n^{2p}c_n^2} =  \frac{2 m\alpha}{2p-\alpha} \frac{n l(v_n)}{v_n^\alpha c_n^2} \rightarrow 0,
	\end{align*}
	where the asymptotic equivalence follows from Lemma \ref{lemma:truncate_moment_asy} and the last convergence follows from (\ref{eq:v_n_choice1}) and $c_n \rightarrow \infty$.
	Therefore, Condition A($p$) is satisfied.
\end{proof}

\begin{proof}[Proof of Theorem \ref{thm:strong_mixing}]
	Fix $p > 0$.    As $\breve{X}_i \leq v_n$, by Lemma 1.2 in \cite{Ibragimov1962} we have for any $i \neq j$,
	\begin{align*}
		|\Cov(\breve{X}^p_i,\breve{X}^p_j)|=     |\mathbb{E}(\breve{X}_i^p\breve{X}_j^p) -\mathbb{E}(\breve{X}_i^p)\mathbb{E}(\breve{X}_j^p)| \leq 4v_n^{2p}\alpha(|j-i|).
	\end{align*}
	Then
	\begin{align*}
		\frac{\sum_{i\neq j}\Cov(\breve{X}^p_i, \breve{X}^p_j) }{v^{2p}_n c_n^2} \leq \frac{\sum_{i\neq j} 4 v_n^{2p}\alpha(|j-i|)}{v_n^{2p} c_n^2} = \frac{8 \sum^n_{k=1} \alpha(k)}{c_n^2} 
		\rightarrow 0.
	\end{align*}
\end{proof}

\begin{proof}[Proof of Corollary \ref{corollary:mixing}]
	According to \cite{Bradley2005basic},
	\begin{enumerate}[(i)]
		\item $\psi$-mixing implies $\phi$-mixing;
		\item $\phi$-mixing implies $\rho$-mixing and $\beta$-mixing;
		\item $\rho$-mixing and $\beta$-mixing each imply strong-mixing.
	\end{enumerate}
	By transitivity, $\{X_i,i\geq 1\}$ must be a strong mixing sequence. 
	For any two $\sigma$-fields $\mathcal{A}$ and $\mathcal{B} \subset \mathcal{F}$,  we also have
	\begin{align*}
		2\alpha(\mathcal{A},\mathcal{B}) \leq \beta(\mathcal{A},\mathcal{B}) \leq \phi(\mathcal{A},\mathcal{B}) \leq \frac{1}{2}\psi(\mathcal{A},\mathcal{B})
	\end{align*}
	and
	\begin{align*}
		4\alpha(\mathcal{A},\mathcal{B}) \leq \rho(\mathcal{A},\mathcal{B}) \leq \psi(\mathcal{A},\mathcal{B}).
	\end{align*}
	The inequalities immediately imply that
	\begin{align*}
		2\alpha(n) \leq \beta(n) \leq \phi(n) \leq \frac{1}{2}\psi(n)
	\end{align*}
	and
	\begin{align*}
		4\alpha(n) \leq \rho(n) \leq \psi(n).
	\end{align*}
	Then $\lim_{n \to \infty} \frac{1}{c_n^2}\sum^\infty_{n=1} \psi(n) = 0$ implies that $\lim_{n \to \infty} \frac{1}{c_n^2}\sum^\infty_{n=1} \alpha(n) = 0$, since $ \alpha(\mathcal{A},\mathcal{B}) \leq \frac{1}{4}\psi(\mathcal{A},\mathcal{B})$. 
	The same result holds for other types of mixing. We omit the details. 
	As a result, Theorem \ref{thm:strong_mixing}  implies that Condition A($p$) holds.
\end{proof}

\begin{proof}[Proof of Theorem \ref{thm:AR_cov}]
	To show that Condition A($p$) holds for any $p > 0$, by Theorem \ref{thm:strong_mixing}, it suffices to show that as $n \rightarrow \infty$, 
	\begin{align*}
		\frac{\sum_{k=1}^n \alpha(k)}{c_n^2} \to 0.
	\end{align*}
	We first obtain an upper bound for $\alpha(k)$. Under the conditions for $\epsilon_t$, by Theorems 1 and 2 and the remark following Theorem 2 in \cite{andrews1983first}, $\{X_t\}$ is a strong mixing process, and there exists a positive integer $s_0$ such that the strong mixing coefficient satisfies
	\begin{align}\label{eq:ar_mix_bound1}
		\alpha(k) \leq C \mathbb{E}(X_1^q) (\beta_1^q)^k \quad \text{for}\quad k \geq s_0
	\end{align}
	and
	\begin{align}\label{eq:ar_mix_bound2}
		\alpha(k) \leq 1 \quad \text{for} \quad 1 \leq k < s_0,
	\end{align}
	where $C > 0$ is a constant. 
	We have  $\mathbb{E}(X_1^q) < \infty$ because $\mathbb{E}(\epsilon^q_1) < \infty$. 
	Furthermore, since $\beta_1 \in (0, 1)$, $\sum_{k=s_0}^n (\beta_1^q)^k \leq \frac{\beta_1^q}{1-\beta_1^q} < \infty$.
	Since $c_n \rightarrow \infty$ as $n \rightarrow \infty$,
	\begin{align*}
		\frac{\sum_{k=1}^n \alpha(k)}{c_n^2} \leq \frac{ (s_0-1) + C \mathbb{E}(X_1^q) \sum_{k=s_0}^n (\beta_1^q)^k}{c_n^2} \rightarrow 0.
	\end{align*}
\end{proof}
\subsection{Proof of Lemma \ref{lemma:log_M_np/log_n} and Theorem \ref{thm:moment}}
We introduce the following Lemmas \ref{lemma:log_t_n/log_n_limit}--\ref{lemma:log_Mnp_dnp} before proving Lemma \ref{lemma:log_M_np/log_n} and Theorem \ref{thm:moment}. 
Recall that $t_n$ satisfies (\ref{eq:choice_of_tn}) and that
$c_n$ satisfies $nc_n \rightarrow \infty$ and $\log c_n = o(\log n)$.

\begin{lemma}\label{lemma:log_t_n/log_n_limit}
	If $l_0>0$ is continuous and slowly varying with 
	$\overline{F}_0(0+):=\lim_{x \downarrow 0} x^{-\alpha}l_0(x)>0$, then
	\begin{equation*}
		\lim_{n \rightarrow \infty} \frac{\log t_n}{\log n} = \frac{1}{\alpha}.
	\end{equation*}
\end{lemma}

\begin{proof}[Proof of Lemma \ref{lemma:log_t_n/log_n_limit}]
	We have
	\begin{align*}
		\frac{\log t_n}{\log n} &=\frac{\log t_n^\alpha}{\alpha \log n} = \frac{1}{\alpha} \left( \frac{\log n + \log c_n + \log l_0(t_n)}{\log n}\right) = \frac{1}{\alpha} \left(1 + \frac{\log c_n}{\log n} + \frac{\log l_0(t_n)}{\log n}\right).
	\end{align*}
	Since $\mathbb{P}(X > t_n) = 1/(nc_n)$ and $n c_n \uparrow\infty$, $t_n$ must also go to $\infty$ as $n \rightarrow \infty$.  Then, by Proposition 1.3.6 in \cite{bingham1989regular} and the choice of $c_n$ that $\log c_n = o(\log n)$, we have 
	\begin{equation*}
		\lim_{n \rightarrow \infty} \frac{\log l_0(t_n)}{\log t_n} = 0 \text{ and } \lim_{n \rightarrow \infty} \frac{\log c_n}{\log n}  = 0,
	\end{equation*}
	and the claim in the lemma follows.
\end{proof}

For $i=1,\ldots,n$, let $\tilde{X}_i:= X_i \mathbbm{1}(X_i < t_n)$. For $p > \alpha$, denote the $p$th truncated sample moment $\tilde{M}_{n,p} := \frac{1}{n}\sum^n_{i=1}\tilde{X}_i^p$ and its corresponding population moment $d_{n,p} := \mathbb{E}(\tilde{X}_1^{p})$. 
The following lemma provides an upper bound for the ratio ${M_{n,p}}/{d_{n,p}}$ by truncating $X_i$ at $t_n$. 
\begin{lemma}[Upper bound for ${M_{n,p}}/{d_{n,p}}$]\label{lemma:UB}
	Let $X_1,\ldots,X_n,\ldots$ be a sequence of 
	nonnegative random variables with a common survival function
	$\overline{F}(x) = x^{-\alpha}l(x)$, where $l$ is slowly varying and $\alpha > 0$ (not necessarily in $(0, 1)$).
	For $p > \alpha$, we have
	\begin{equation*}
		\frac{M_{n,p}}{d_{n,p}} = O_p(1).
	\end{equation*}
\end{lemma}

\begin{proof}[Proof of Lemma \ref{lemma:UB}]
	Let $\varepsilon > 0$. Choose a constant $C>0$ such that $1/C < \varepsilon /2$. By (\ref{eq:t_n_to_0}), there exists $N$ such that for all $n \geq N$, $nl(t_n)/t^\alpha_n < \varepsilon/2$. 
	Then
	\begin{align}
		\mathbb{P} \left( \frac{M_{n,p}}{d_{n,p}}> C \right) 
		& = \mathbb{P} \left( \frac{M_{n,p}}{d_{n,p}}> C , \tilde{M}_{n,p} = M_{n,p} \right) 
		+
		\mathbb{P} \left( \frac{M_{n,p}}{d_{n,p}}> C , \tilde{M}_{n,p} \neq M_{n,p} \right)  \nonumber
		\\
		&\leq  \mathbb{P} \left( \frac{\tilde{M}_{n,p}}{d_{n,p}} > C \right)  + \mathbb{P}(\tilde{M}_{n,p} \neq M_{n,p}).  \label{eq:upper_bound_1}
	\end{align}
	For the first term in (\ref{eq:upper_bound_1}), since ${\tilde{M}_{n,p}}/{d_{n,p}}$ is non-negative, by Markov's inequality,
	\begin{equation*}
		\mathbb{P} \left( \frac{\tilde{M}_{n,p}}{d_{n,p}}> C \right) \leq \frac{\mathbb{E}(\tilde{M}_{n,p}/d_{n,p})}{C} =  \frac{1}{C} < \frac{\varepsilon}{2}.
	\end{equation*}
	For the second term in (\ref{eq:upper_bound_1}),
	\begin{equation*}
		\mathbb{P}(\tilde{M}_{n,p} \neq M_{n,p})
		\leq \mathbb{P} \left( \cup^n_{i=1} \{X_i > t_n\} \right) \leq n \mathbb{P}(X > t_n) = \frac{n l(t_n)}{t_n^\alpha} < \frac{\varepsilon}{2}.
	\end{equation*}
\end{proof}

Next Lemma \ref{lemma:LB} provides a lower bound for ${M_{n,p}}/{d_{n,p}}$. Consider the truncated term $\breve{X}_i := X_i \mathbbm{1}(X_i < v_n)$. Recall that $\{v_n\}$ is a sequence satisfying \eqref{eq:v_n_choice1}.
Define $\breve{M}_{n,p} := \frac{1}{n}\sum^n_{i=1}\breve{X}_i^p$  and $\breve{d}_{n,p} := \mathbb{E}(\breve{X}^p_1)$. 

\begin{lemma}[Lower Bound for ${M_{n,p}}/{d_{n,p}}$]\label{lemma:LB}
	Let $X_1,\ldots,X_n,\ldots$ be a sequence of 
	nonnegative random variables with a common survival function $\overline{F}(x) = x^{-\alpha}l(x)$, where $l$ is slowly varying and $\alpha > 0$ (not necessarily in $(0, 1)$). 
	If  Condition A($p$) holds for $p > \alpha$, then
	\begin{equation*}
		\frac{M_{n,p}}{d_{n,p}} \geq \frac{c_{p,\alpha,\delta'}(1+o_p(1))}{ c_n^{\delta'}},
	\end{equation*}
	for any $\delta' > 2(p/\alpha  -1)>0$ and $c_{p,\alpha, \delta'} > 0$ is some constant depending only on $p,\alpha,\delta'$.
\end{lemma}

\begin{proof}[Proof of Lemma \ref{lemma:LB}]
	Write
	\begin{equation}\label{eq:lower_bound1}
		n M_{n, p} \geq (n\breve{M}_{n,p} - n \breve{d}_{n,p}) + n \breve{d}_{n,p}.
	\end{equation}
	For any $p > \alpha$, by Lemma \ref{lemma:truncate_moment_asy},
	\begin{equation}\label{eq:LB2_new}
		\frac{n \breve{d}_{n,p}}{v_n^p} \sim \frac{ \frac{\alpha}{p - \alpha} n v_n^{p-\alpha}l(v_n)}{v_n^p} \sim c_n \frac{\alpha}{p-\alpha}.
	\end{equation}    
	For any $x > 0$, by Markov's inequality,
	\begin{align}\label{eq:lower_bound3}
		\mathbb{P} \left(  \frac{|n\breve{M}_{n,p} - n \breve{d}_{n,p}|}{v_n^p c_n} > x \right)  
		&\leq \frac{\Var(n\breve{M}_{n,p} - n \breve{d}_{n,p})}{v^{2p}_n c_n^2 x^2} \nonumber \\
		&\leq \frac{n^2\Var(\breve{M}_{n,p})}{v^{2p}_n c_n^2 x^2} \nonumber \\
		&\leq \frac{n \breve{d}_{n,2p} + \sum_{i\neq j}\Cov(\breve{X}^p_i, \breve{X}^p_j)}{v^{2p}_n c_n^2 x^2} \rightarrow 0,
	\end{align}
	where the last convergence follows from (\ref{eq:LB2_new}) since 
	${n \breve{d}_{n,2p}}/{v_n^{2p}x^2} \to 0$, the fact that $c_n \rightarrow \infty$ and  (\ref{eq:cov_ass2}).
	In view of (\ref{eq:lower_bound1})--(\ref{eq:lower_bound3}),
	\begin{equation*}
		\frac{n M_{n,p}}{v_n^p c_n} \geq o_p\left(\frac{1}{c_n}\right) + c_n\frac{\alpha}{p-\alpha}.
	\end{equation*}
	Then
	\begin{align*}
		\frac{ M_{n,p}}{d_{n,p}} = \frac{n M_{n,p}}{v_n^p c_n} \cdot \frac{v_n^p c_n}{n d_{n,p}} \geq \left(o_p\left(\frac{1}{c_n}\right) + \frac{\alpha}{p-\alpha} \right)\frac{v_n^p c_n}{n d_{n,p}}.
	\end{align*}
	By Lemma \ref{lemma:truncate_moment_asy},
	\begin{align*}
		\frac{v_n^p c_n}{n d_{n,p}} \sim \frac{v_n^p c_n}{n t_n^{p-\alpha} l(t_n)\frac{\alpha}{p-\alpha} } = \frac{t_n^\alpha}{n l(t_n)} c_n \left( \frac{v_n}{t_n}\right)^p
		= c_n^2 \left( \frac{v_n}{t_n}\right)^p.
	\end{align*}
	By the definitions of $v_n$ and $t_n$,
	\begin{align}\label{eq:LB_b_1}
		\left( \frac{v_n}{t_n}\right)^p 
		&= \left( \frac{nl_0(v_n)}{c_n} \cdot \frac{1}{n c_n l_0(t_n)} \right)^{p/\alpha} 
		= \left( \frac{1}{c_n^2} \frac{l_0(v_n)}{l_0(t_n)} \right)^{p/\alpha}.
	\end{align}
	By Potter's Theorem (see Theorem 1.5.6 in \cite{bingham1989regular}), for any $0 < \delta < \alpha$, there exists $N$ such that for all $n \geq N$,
	\begin{equation*}
		\frac{l_0(t_n)}{l_0(v_n)} \leq 2 \max \left\{ \left( \frac{t_n}{v_n}\right)^\delta, \left( \frac{t_n}{v_n}\right)^{-\delta}\right\}.
	\end{equation*}    
	Since $v_n \leq t_n$ (otherwise, $M_{n,p}/d_{n,p} \geq (o_p(1) + \frac{\alpha}{p-\alpha})(1+o(1))c_n^2 \rightarrow \infty$ could contradict the fact that $M_{n,p}/{d_{n,p}} = O_p(1)$ established in Lemma \ref{lemma:UB}), we have
	\begin{equation}\label{eq:LB_b_2}
		\frac{l_0(t_n)}{l_0(v_n)} \leq 2 \left( \frac{t_n}{v_n}\right)^\delta.
	\end{equation}
	Then (\ref{eq:LB_b_1}) and (\ref{eq:LB_b_2}) give
	\begin{equation*}
		\left( \frac{v_n}{t_n}\right)^{p} \geq \left( \frac{1}{c_n^2}\right)^{p/\alpha} \cdot
		\frac{1}{2^{p/\alpha}}\left( \frac{v_n}{t_n}\right)^{\delta p /\alpha}.
	\end{equation*}
	Therefore
	\begin{equation*}
		\left( \frac{v_n}{t_n}\right)^{p(1-\delta/\alpha)} \geq \left( \frac{1}{2 c_n^2}\right)^{p/\alpha} .
	\end{equation*}
	Thus
	\begin{align*}
		c_n^2 \left( \frac{v_n}{t_n}\right)^p \geq c_n^2 
		\left( \frac{1}{2 c_n^2}\right)^{\frac{p}{\alpha - \delta}}  = \frac{1}{2^{\frac{p}{\alpha-\delta}}} \left( \frac{1}{c_n^2} \right)^{ \frac{p}{\alpha-\delta}-1}.
	\end{align*}
	The claim follows by letting $\delta' := 2 \left( \frac{p}{\alpha -\delta} - 1\right) > 2(p/\alpha  -1)$ and $c_{p, \alpha, \delta'} := \frac{\alpha}{p-\alpha} \frac{1}{2^{p/(\alpha-\delta)}}$.
\end{proof}

\begin{lemma}\label{lemma:log_Mnp_dnp}
	Under the conditions in Lemma \ref{lemma:LB}, 
	\begin{equation*}
		\log \frac{M_{n,p}}{d_{n,p}} = O_p(\log c_n).
	\end{equation*}
\end{lemma}

\begin{proof}[Proof of Lemma \ref{lemma:log_Mnp_dnp}]
	Fix $\varepsilon \in (0, 1)$. By Lemmas \ref{lemma:UB} and \ref{lemma:LB}, with probability greater than $1 - \varepsilon$, for all large $n$, there exists some $C > 0$ such that
	\begin{equation*}
		\frac{c_{p, \alpha, \delta'}}{2 c_n^{\delta'}} \leq \frac{M_{n,p}}{d_{n,p}} \leq  C.
	\end{equation*}
	Therefore, for all sufficiently large $n$,
	\begin{equation*}
		\left|\log \frac{M_{n,p}}{d_{n,p}} \right| \leq  |\log C| + |\log(c_{p, \alpha, \delta'}/2)|  + \delta' \log c_n.
	\end{equation*}
	This completes the proof.
\end{proof}

\begin{proof}[Proof of Lemma \ref{lemma:log_M_np/log_n}]
	We write:
	\begin{equation*}
		\log \frac{M_{n,p}}{n^{\frac{p-\alpha}{\alpha}} } = 
		\log \frac{M_{n,p}}{d_{n,p}}
		+ 
		\log \frac{d_{n,p}}{n^{\frac{p-\alpha}{\alpha}} }.
	\end{equation*}
	By Lemma \ref{lemma:truncate_moment_asy},
	\begin{equation*}
		\frac{d_{n,p}}{n^{\frac{p-\alpha}{\alpha}} } \sim \frac{\frac{\alpha}{p-\alpha}t_n^{p-\alpha}l(t_n)  }{n^{\frac{p-\alpha}{\alpha}}} = \frac{ \frac{\alpha}{p-\alpha} (n l(t_n) c_n)^{\frac{p-\alpha}{\alpha}} l(t_n)  }{n^{\frac{p-\alpha}{\alpha}} } = \frac{\alpha}{p-\alpha} c_n^{\frac{p-\alpha}{\alpha} }l(t_n)^{\frac{p}{\alpha}}.
	\end{equation*}
	Thus, by Lemma \ref{lemma:log_Mnp_dnp},
	\begin{equation*}
		\frac{\log M_{n,p}}{\log n} - \frac{p-\alpha}{\alpha}  =
		\frac{ \log \frac{M_{n,p}}{d_{n,p}} }{\log n} + \frac{\log \frac{d_{n,p}}{n^{\frac{p-\alpha}{\alpha}} }}{\log n} = 
		O_p\left(\frac{\log c_n}{\log n}\right) + O\left( \frac{|\log l(t_n)|}{\log n}\right).
	\end{equation*}
	Finally, (\ref{eq:estimator_moment_n}) follows from $\log l(t_n) / \log t_n \rightarrow 0$ and Lemma \ref{lemma:log_t_n/log_n_limit}.
\end{proof}
Now we are ready to prove Theorem \ref{thm:moment}.
\begin{proof}[Proof of Theorem \ref{thm:moment}]
	\begin{enumerate}[(a)]
		\item We  write
		\begin{align}\label{eq:thm_moment_decomp}
			\log \frac{M_{n,h_1}}{M_{n,h_2}^{\imath(h_1,h_2)}} =
			\log \frac{M_{n,h_1}}{d_{n,h_1}} +
			\log \frac{d_{n,h_1}}{d_{n,h_2}^{\imath(h_1,h_2)}} +
			\imath(h_1,h_2) \log \frac{d_{n,h_2}}{M_{n,h_2}}.
		\end{align}
		By Lemma \ref{lemma:log_Mnp_dnp}, 
		\begin{equation}\label{eq:thm_moment_2}
			\log \frac{M_{n, h_1}}{d_{n,h_1}} = O_p(\log c_n) \text{ and }  \log \frac{d_{n,h_2}}{M_{n,h_2}} = O_p(\log c_n).
		\end{equation}
		Also
		\begin{align*}
			\frac{d_{n,h_1}}{d_{n,h_2}^{\imath(h_1,h_2)}} \sim \frac{\frac{\alpha}{h_1-\alpha} t_n^{h_1-\alpha}l(t_n)}{ [\frac{\alpha}{h_2-\alpha}t_n^{h_2-\alpha}l(t_n) ]^{\imath(h_1,h_2)} }  
			= \left( \frac{\alpha}{h_1-\alpha} \right) \left( \frac{\alpha}{h_2-\alpha}\right)^{-\imath(h_1,h_2)} l(t_n)^{(h_2-h_1)/(h_2-\alpha)}.
		\end{align*}
		Thus
		\begin{equation}\label{eq:thm_moment_3}
			\log \frac{d_{n,h_1}}{d_{n,h_2}^{\imath(h_1,h_2)}} = o(1) +  C_1 + C_2 \log l(t_n),
		\end{equation}
		where $C_1, C_2$ are some constants depending only on $\alpha,h_1,h_2$.
		Combining (\ref{eq:thm_moment_decomp})--(\ref{eq:thm_moment_3}), we obtain
		\begin{equation}\label{eq:thm_moment_4}
			\log \frac{M_{n,h_1}}{M_{n,h_2}^{\imath(h_1,h_2)}} =    \log M_{n,h_1} - \imath(h_1,h_2) \log M_{n, h_2} = O_p(\log c_n) + O(|\log l(t_n)|).
		\end{equation}
		Dividing both sides of  (\ref{eq:thm_moment_4}) by $\log M_{n,h_2}$ gives
		\begin{equation}\label{eq:log_Mnk_Mn1}
			\frac{\log M_{n,h_1}}{\log M_{n,h_2}} - \imath(h_1,h_2) = O_p\left( \frac{\log c_n}{\log M_{n,h_2}} \right) + O \left( \frac{|\log l(t_n)|}{\log M_{n,h_2}} \right).
		\end{equation}
		The claim in \eqref{eq:pareto rate} then follows since $\log M_{n,h_2}/ \log  n \stackrel{\mathbb{P}}{\rightarrow} \frac{h_2-\alpha}{\alpha}$ as shown in Lemma \ref{lemma:log_M_np/log_n}.
		\item \begin{enumerate}[(i)]
			\item 
			As $l(x) \rightarrow \infty$ and $t_n \rightarrow \infty$, we have $ l(t_n) > 1$ and $l(x) > 1$ for all sufficiently large $n$ and $x$. Let  $c_n \rightarrow \infty$ be such that $1 \leq c_n \leq l(n)$ and
			\begin{equation}\label{eq:cn_choice_rate}
				\frac{\log (\inf_{ x \in I_{n,l} } l(x))}{\log c_n} \rightarrow \infty.
			\end{equation} 
			From $n l(t_n)/t_n^\alpha \sim 1/c_n$, we have $t_n^\alpha \leq 2 n c_n l(t_n)$ for all sufficiently large $n$. As $l(x)/x^{\delta/2} \rightarrow 0$ and $c_n \leq l(n)$, we have for all sufficiently large $n$,
			\begin{equation*}
				t_n^{\alpha-\delta} \leq 2 n c_n \frac{l(t_n)}{t_n^\delta} \leq 2 n l(n) \frac{l(t_n)}{t_n^\delta} \leq 2n^{1+\delta/2} \frac{l(t_n)}{t_n^\delta} \leq n^{1+\delta}.
			\end{equation*}
			Thus, $t_n \leq n^{(1+\delta)/(\alpha-\delta)}$ for all sufficiently large $n$.    
			On the other hand, as $c_n \geq 1$, 
			\begin{align*}
				\frac{nl(t_n)}{t_n^\alpha} \sim \frac{1}{c_n} \leq 1,
			\end{align*}
			and so
			\begin{align}\label{eq:t_n choice}
				t_n^{\alpha +\delta} \geq \frac{1}{2} n t_n^{\delta} l(t_n) \geq n t_n^{\delta/2} l(t_n).
			\end{align}
			As $t_n^{\delta/2}l(t_n) \to \infty$ by Proposition 1.3.6 (v) in \cite{bingham1989regular}, we have $t_n^{\alpha+\delta} \geq  n$ and so $t_n \geq n^{1/(\alpha+\delta)}$ for all large $n$. From $n^{1/(\alpha+\delta)} \leq t_n \leq n^{(1+\delta)/(\alpha-\delta)}$ and $t_n^\alpha \sim c_n n l(t_n)$, we have
			\begin{equation}\label{eq:tn_alpha bound}
				\frac{1}{2} n \inf_{n^{1/(\alpha+\delta)} \leq x \leq n^{(1+\delta)/(\alpha-\delta)}} l(x) \leq t_n^\alpha \leq 2 l(n)n \sup_{n^{1/(\alpha+\delta)} \leq x \leq n^{(1+\delta)/(\alpha-\delta)}} l(x).
			\end{equation}
			Thus $t_n \in I_{n,l}$. Then, in view of (\ref{eq:cn_choice_rate}), 
			for all large $n$, 
			\begin{equation*}
				\log c_n + |\log l(t_n)| \leq 2 \log \left( \sup_{x \in I_{n,l}}l(x)\right).
			\end{equation*}
			Then, using Theorem \ref{thm:moment} (a), we obtain $R_{n,h_1,h_2}^{\text{hm}} = O_p(s_{1n,l})$.
			
			Next, we shall show that $(R_{n,h_1,h_2}^{\text{hm}})^{-1} = O_p(s_{2n,l}^{-1})$.
			By (\ref{eq:thm_moment_decomp}) and the triangle inequality, we have
			\begin{align*}
				\left| \log \frac{M_{n,h_1}}{M_{n,h_2}^{\imath(h_1,h_2)}}\right| & \geq 
				\left|\log \frac{d_{n,h_1}}{d_{n,h_2}^{\imath(h_1,h_2)}}\right| 
				- 
				\left|\log \frac{M_{n,h_1}}{d_{n,h_1}}\right| -
				\left|\imath(h_1,h_2) \log \frac{d_{n,h_2}}{M_{n,h_2}}\right|.
			\end{align*}
			For any $\varepsilon > 0$, by (\ref{eq:thm_moment_3}) and Lemma \ref{lemma:log_Mnp_dnp}, there exists some $C > 0$ independent of $n$ such that  for all large $n$, with probability at least $1 - \varepsilon$,
			\begin{align*}
				\left| \log \frac{M_{n,h_1}}{M_{n,h_1}^{\imath(h_1,h_2)}}\right| 
				& \geq 
				\left|o(1) + C_1 + C_2 \log l(t_n)\right|  - C \log c_n\\
				& \geq 
				\left|o(1) + C_1 + C_2 \log\left(\inf_{x \in I_{n,l}\textbf{}} l(x)\right) \right|  - C \log c_n\\
				&= \left|o(1) + C_1 + C_2 \log\left(\inf_{x \in I_{n,l}} l(x)\right) \right|  - C o\left(\log\left(\inf_{x \in I_{n,l}} l(x)\right)\right),
			\end{align*}
			where the last equality follows from  (\ref{eq:cn_choice_rate}). Thus,
			\begin{align*}
				\left| \log \frac{M_{n,h_1}}{M_{n,h_1}^{\imath(h_1,h_2)}}\right| & \geq 
				\frac{1}{2}\left|C_2 \log\left(\inf_{x \in I_{n,l}} l(x)\right)\right|.
			\end{align*}
			Finally, dividing both sides of the above inequality by $|\log M_{n,h_1}|$, with  probability at least $1 - \varepsilon$,  we have
			\begin{align*}
				\left|  \frac{\log M_{n,h_1}}{\log M_{n,h_1}} - \imath(h_1,h_2) \right| & \geq 
				\frac{|C_2|}{2} \frac{|\log(\inf_{x \in I_{n,l}} l(x))|}{\log M_{n, h_1}}.
			\end{align*}
			Thus, $(R_{n,h_1,h_2}^{\text{hm}})^{-1} = O_p(s_{2n,l}^{-1})$  in view of Lemma \ref{lemma:log_M_np/log_n}.\\
			
			\item
			As $l(x) \rightarrow 0$ and $t_n \rightarrow \infty$, we have $l(t_n) < 1$ and $l(x) < 1$ for all large $n$ and $x$.  Let  $c_n \rightarrow \infty$ be such that $1 \leq c_n \leq l(n)^{-1}$ and
			\begin{equation}\label{eq:cn_choice_rate2}
				\frac{|\log (\sup_{ x \in I_{n,l} } l(x))|}{\log c_n} \rightarrow \infty.
			\end{equation} 
			Using the same argument as in the proof of (b)(i), we have $t_n \in I_{n,l}$ for all sufficiently large $n$.
			Thus, by (\ref{eq:cn_choice_rate2}),  for all sufficiently large $n$, 
			\begin{equation*}
				\log c_n + |\log l(t_n)| \leq 2 \left| \log \left( \inf_{x \in I_{n,l}}l(x)\right)\right|.
			\end{equation*}
			In view of Theorem \ref{thm:moment} (a), we obtain $R_{n,h_1,h_2}^{\text{hm}} = O_p(s_{2n,l})$.
			
			Next, we shall show that $(R_{n,h_1,h_2}^{\text{hm}})^{-1} = O_p(s_{1n,l}^{-1})$. 
			By (\ref{eq:thm_moment_decomp}) and the triangle inequality, we have
			\begin{align*}
				\left| \log \frac{M_{n,h_1}}{M_{n,h_2}^{\imath(h_1,h_2)}}\right| & \geq 
				\left|\log \frac{d_{n,h_1}}{d_{n,h_2}^{\imath(h_1,h_2)}}\right| 
				- 
				\left|\log \frac{M_{n,h_1}}{d_{n,h_1}}\right| -
				\left|\imath(h_1,h_2) \log \frac{d_{n,h_2}}{M_{n,h_2}}\right|.
			\end{align*}
			As $t_n \geq n$ for all large $n$, we have
			\begin{equation*}
				0 < l(t_n) \leq \sup_{x \in I_{n,l}} l(x) < 1,
			\end{equation*}
			and so
			\begin{equation*}
				|    \log l(t_n) | \geq \left|\log \left(\sup_{x \in I_{n,l}}l(x)\right)\right|.
			\end{equation*}
			
			For any $\varepsilon > 0$, by (\ref{eq:thm_moment_3}) and Lemma \ref{lemma:log_Mnp_dnp}, there exists some $C > 0$ independent of $n$ such that  for all large $n$, with probability at least $1 - \varepsilon$,
			\begin{align*}
				\left| \log \frac{M_{n,h_1}}{M_{n,h_1}^{\imath(h_1,h_2)}}\right| 
				& \geq 
				\left|o(1) + C_1 + C_2 \log l(t_n)\right|  - C \log c_n\\
				& \geq 
				\left|o(1) + C_1 + C_2 \log\left(\sup_{x \in I_{n,l}} l(x)\right) \right|  - C \log c_n\\
				& = \left|o(1) + C_1 + C_2 \log\left(\sup_{x \in I_{n,l}} l(x)\right) \right|  + o\left( \log\left(\sup_{x \in I_{n,l}} l(x)\right)\right),
			\end{align*}
			where the last equality follows from  (\ref{eq:cn_choice_rate2}). Thus
			\begin{align*}
				\left| \log \frac{M_{n,h_1}}{M_{n,h_1}^{\imath(h_1,h_2)}}\right| & \geq 
				\frac{1}{2}\left|C_2 \log\left(\sup_{x \in I_{n,l}} l(x)\right)\right|.
			\end{align*}
			Dividing both sides of the above inequality by $|\log M_{n,h_1}|$, we have, 
			with  probability at least $1 - \varepsilon$,  
			\begin{align*}
				\left|  \frac{\log M_{n,h_1}}{\log M_{n,h_1}} - \imath(h_1,h_2) \right| & \geq 
				\frac{|C_2|}{2} \frac{|\log(\sup_{x \in I_{n,l}} l(x))|}{\log M_{n, h_1}}.
			\end{align*}
			Thus, $(R_{n,h_1,h_2}^{\text{hm}})^{-1} = O_p(s_{1n,l}^{-1})$  follows in view of Lemma \ref{lemma:log_M_np/log_n}.    
		\end{enumerate}
	\end{enumerate}
\end{proof}

\subsection{Proof of Theorem \ref{thm:higher_central_moments}}
To  prove Theorem \ref{thm:higher_central_moments}, we  introduce  Lemma \ref{lemma:higher_central_moments}, which is 
analogous to  Lemma \ref{lemma:log_Mnp_dnp}.
\begin{lemma}\label{lemma:higher_central_moments}
	Let $X_1,\ldots,X_n,\ldots$ be a sequence of 
	nonnegative random variables with a common survival function $\overline{F}(x) = x^{-\alpha}l(x)$, where $l$ is slowly varying and $\alpha > 0$ (not necessarily in $(0, 1)$).  If Condition A($p$) holds for $k > \alpha$, where $k \in \mathbb{N}$, then
	\begin{equation*}
		\log \frac{|M^c_{n,k}|}{d_{n,k}} = O_p\left( \log c_n\right).
	\end{equation*}
\end{lemma}

\begin{proof}[Proof of Lemma \ref{lemma:higher_central_moments}]
	Let $M_{n,0} := 1$. As
	\begin{equation*}
		M^c_{n,k} = \sum^k_{j=0} C^k_j(-1)^{k-j} M_{n,j} M_{n,1}^{k-j},
	\end{equation*}
	we have an upper bound
	\begin{equation}\label{eq:higher_central_3}
		|M^c_{n,k}| \leq M_{n,k} + \sum^{k-1}_{j=0} C^k_j M_{n,j} M_{n,1}^{k-j}.
	\end{equation}
	By Lemma \ref{lemma:truncate_moment_asy}, for some constant $C$, 
	\begin{align}
		\frac{d_{n,j} d^{k-j}_{n,1}}{d_{n,k}} &\sim \frac{Ct_n^{j-\alpha} t_n^{(1-\alpha)(k-j)} l(t_n) l(t_n)^{k-j}}{t_n^{k-\alpha}l(t_n) } = \frac{C l(t_n)^{k-j}}{t_n^{\alpha(k-j)}} \nonumber \\
		&= \frac{C l(t_n)^{k-j}}{(n l(t_n)c_n)^{k-j}} = \frac{C}{(nc_n)^{k-j}}. \label{eq:higher_central_2}
	\end{align}
	By  (\ref{eq:higher_central_2}) and Lemma 
	\ref{lemma:UB},
	\begin{align}
		\frac{1}{d_{n,k}} \sum^{k-1}_{j=0} C^k_j M_{n,j} M_{n,1}^{k-j}  &= \sum^{k-1}_{j=0} C^k_j \frac{M_{n,j}}{d_{n,j}} \cdot \frac{M_{n,1}^{k-j}}{d_{n,1}^{k-j}} \cdot \frac{d_{n,j} d^{k-j}_{n,1}}{d_{n,k}}\nonumber \\
		&= \sum^{k-1}_{j=0} O_p(c_n) O_p((c_n)^{k-j}) O\left( \frac{1}{(n c_n)^{(k-j)}}\right) \nonumber \\
		&= O_p(\frac{\frac{c_n}{n^k}(n^k-1)}{n-1}) \nonumber \\
		&= O_p(\frac{c_n}{n-1})-\frac{c_n}{n^k(n-1)}\nonumber \\
		&= O_p \left( \frac{c_n}{n} \right). \label{eq:higher_central_4}
	\end{align}
	By (\ref{eq:higher_central_3}), (\ref{eq:higher_central_4}) and Lemma \ref{lemma:UB},
	\begin{equation}\label{eq:higher_central_UB}
		\frac{|M^c_{n,k}|}{d_{n,k}} \leq O_p( c_n) + O_p \left( \frac{c_n}{n} \right) = O_p(c_n).
	\end{equation}
	On the other hand, we have a lower bound
	\begin{equation}\label{eq:higher_central_LB1}
		|M^c_{n,k}| \geq M_{n,k} - \sum^{k-1}_{j=0} C^k_j M_{n,j} M_{n,1}^{k-j}.
	\end{equation}
	By Lemma \ref{lemma:LB} and (\ref{eq:higher_central_LB1}),
	\begin{equation}\label{eq:higher_central_LB}
		\frac{|M^c_{n,k}|}{d_{n,k}} \geq \frac{M_{n,k}}{d_{n,k} }- \frac{1}{d_{n,k}} \sum^{k-1}_{j=0} C^k_j M_{n,j} M_{n,1}^{k-j} \geq 
		\frac{c_{p,\alpha,\delta'}(1+o_p(1))}{ c_n^{\delta'}} - O_p\left( \frac{c_n}{n}\right).
	\end{equation}
	The claim in the lemma then follows from (\ref{eq:higher_central_UB}) and (\ref{eq:higher_central_LB}).
\end{proof}

\begin{proof}[Proof of Theorem \ref{thm:higher_central_moments}]
	\begin{enumerate}[(a)]
		\item Write
		\begin{equation*}
			\log \frac{|M^c_{n,k}|}{M_{n,1}^{\imath(k, 1)}}  = \log \frac{|M^c_{n,k}|}{d_{n,k}} + \log \frac{d_{n,k}}{d_{n,1}^{\imath(k, 1)}} + \imath(k, 1) \log \frac{d_{n,1}}{M_{n,1}}
		\end{equation*}
		and
		\begin{equation*}
			\log\frac{ |M^c_{n,h_1}|}{ |M^c_{n,h_2}|} = \log \frac{|M^c_{n,h_1}|}{d_{n,h_1}} + \log \frac{d_{n,h_1}}{d_{n,h_2}^{\imath(h_1,h_2)}} + \imath(h_1, h_2) \log \frac{d_{n, h_2}}{|M^c_{n,h_2}|}.
		\end{equation*}
		With Lemma \ref{lemma:higher_central_moments}, the rest of the proof follows the same argument as in the proof of Theorem \ref{thm:moment} (a).
		\item It follows the same argument as in the proof of Theorem \ref{thm:moment} (b).
		
	\end{enumerate}
	
\end{proof}

\section{Proofs for Section \ref{sec:semivariances}}
\subsection{Proofs for Section \ref{subsect:lower_upper}}
\begin{proof}[Proof of Theorem \ref{thm:hth_central_lower}]
	\begin{enumerate}[(i)]
		\item We have an upper bound:
		\begin{equation}\label{eq:hth_central_lower1}
			\frac{M^-_{n,h}}{M^h_{n,1}}  \leq \frac{\frac{1}{n} \sum_{i : X_i \leq M_{n,1}}M_{n,1}^h}{ M_{n,1}^h } \leq 1.
		\end{equation}
		Fix $a > 0$. We also have a lower bound:
		\begin{align}
			\frac{M^-_{n,h}}{M_{n,1}^h} &= \frac{1}{n}\sum_{i : X_i \leq M_{n,1}} 
			\left(1 - \frac{X_i}{M_{n,1}}\right)^h 
			\geq \frac{1}{n}\sum_{i : X_i \leq M_{n,1}} 
			\left(1 - \frac{X_i}{M_{n,1}}\right)^h \mathbbm{1}(M_{n,1} > a)\nonumber \\
			&\geq \frac{1}{n}\sum_{i : X_i \leq a} 
			\left(1 - \frac{X_i}{M_{n,1}}\right)^h \mathbbm{1}(M_{n,1} > a)
			\geq \frac{1}{n} \sum_{i: X_i \leq a } \left(1 - \frac{a}{M_{n,1}}\right)^h \mathbbm{1}(M_{n,1} > a) \nonumber \\
			&= \left[\frac{1}{n} \sum^n_{i=1}\mathbbm{1}(X_i \leq a) \right] \left(1 - \frac{a}{M_{n,1}}\right)^h \mathbbm{1}(M_{n,1} > a). \label{eq:hth_central_lower3}
		\end{align}
		Since Condition A($p$) holds for $p=1$, from the lower bound for $M_{n,1}/d_{n,1}$ in Lemma \ref{lemma:LB}, we know that $M_{n,1}$ diverges to $\infty$ in probability. Thus, as $a > 0$ is fixed, 
		\begin{equation}\label{eq:hth_central_lower4}
			\left( 1-\frac{a}{M_{n,1}}\right)^h \mathbbm{1}(M_{n,1} > a) \stackrel{\mathbb{P}}{\rightarrow} 1.
		\end{equation}
		For any $\varepsilon > 0$, by Markov's inequality,
		\begin{align*}
			\mathbb{P} \left( \left| \frac{1}{n}\sum^n_{i=1} \mathbbm{1}(X_i \leq a) - F(a) \right| > \varepsilon \right) 
			&\leq \frac{\Var\left(\frac{1}{n}\sum^n_{i=1} \mathbbm{1}(X_i \leq a)\right)}{\epsilon^2} \\
			&= \frac{ F(a)(1-F(a))}{n \varepsilon^2} + \frac{\sum_{i\neq j} \Cov(\mathbbm{1}(X_i \leq a), \mathbbm{1}(X_j \leq a))}{n^2 \varepsilon^2}.
		\end{align*}
		
		The last expression goes to $0$ by (\ref{eq:ass_lower_semi}). This together with (\ref{eq:hth_central_lower3}) and (\ref{eq:hth_central_lower4}) give
		\begin{equation}\label{eq:hth_central_lower2}
			\frac{M^-_{n,h}}{M^h_{n,1}} \geq (F(a) + o_P(1))(1+o_p(1)) = F(a) + o_p(1).
		\end{equation}
		In view of (\ref{eq:hth_central_lower1}) and (\ref{eq:hth_central_lower2}), we have
		\begin{equation}\label{eq:hth_central_lower}
			\log \frac{M^-_{n,h}}{M^h_{n,1}}  = O_p(1),
		\end{equation}
		which implies that
		\begin{equation*}
			\log M_{n,h}^- - h \log M_{n,1} = O_p(1).
		\end{equation*}
		Dividing both sides of the last equation by $\log M_{n,1}$, we obtain
		\begin{equation*}
			\frac{\log M_{n,h}^-}{\log M_{n,1}} = h + O_p\left( \frac{1}{\log M_{n,1}} \right).
		\end{equation*}
		The claim in the theorem then follows since $\log M_{n,1}/ \log  n \stackrel{\mathbb{P}}{\rightarrow} \frac{1-\alpha}{\alpha}$ as shown in Lemma \ref{lemma:log_M_np/log_n}.
		
		\item As in (i), we have
		\begin{equation}\label{eq:hth_1}
			\left[\frac{1}{n} \sum^n_{i=1}\mathbbm{1}(X_i \leq b_n) \right] \left(1 - \frac{b_n}{M_{n,1}}\right)^h \mathbbm{1}(M_{n,1} > b_n) \leq \frac{M^-_{n,h}}{M_{n,1}^h} \leq 1.
		\end{equation}
		By Lemma \ref{lemma:LB},
		\begin{equation*}
			\frac{c_{\alpha,1,\delta'}(1+o_p(1))}{c_n^{\delta'}} \leq \frac{M_{n,1}}{d_{n,1}}.
		\end{equation*}
		Thus, by the definition of $b_n$,
		\begin{equation*}
			\frac{b_n}{M_{n,1}} \leq     \frac{1}{c_n^{\delta'} c_{\alpha,1,\delta'}(1+o_p(1))}\rightarrow 0. 
		\end{equation*}
		As a result,
		\begin{equation}\label{eq:hth_2}
			\mathbbm{1}(M_{n,1} > b_n) \stackrel{\mathbb{P}}{\rightarrow} 1 \text{ and } \left(1 - \frac{b_n}{M_{n,1}}\right)^h  \stackrel{\mathbb{P}}{\rightarrow} 1.
		\end{equation}
		By the same argument as in (i),
		\begin{equation}\label{eq:hth_3}
			\frac{1}{n} \sum^n_{i=1}\mathbbm{1}(X_i \leq b_n) = o_p(1) + F(b_n) \rightarrow 1,
		\end{equation}
		as $b_n \rightarrow \infty$. In view of (\ref{eq:hth_1})--(\ref{eq:hth_3}),
		\begin{equation*}
			\frac{M^-_{n,h}}{M_{n,1}^h} \stackrel{\mathbb{P}}{\rightarrow} 1.
		\end{equation*}
	\end{enumerate}
	
\end{proof}

To prove Theorem \ref{thm:upper_central_moment}, we  establish a  lemma  analogous to Lemma \ref{lemma:log_Mnp_dnp}.
\begin{lemma}\label{lemma:central_upper_moment} 
	Let $X_1,\ldots,X_n,\ldots$ be a sequence of 
	nonnegative random variables with a common survival function $\overline{F}(x) = x^{-\alpha}l(x)$, where $l$ is slowly varying and $\alpha \in (0, 1)$. Suppose that Condition A($p$) holds for $p = h > 1$. Then
	\begin{equation*}
		\log \frac{M^+_{n,h}}{d_{n,h}} = O_p\left( \log c_n \right).
	\end{equation*}   
\end{lemma}

\begin{proof}[Proof of Lemma \ref{lemma:central_upper_moment}]
	To establish the claimed result, we shall obtain upper and lower bounds of $M_{n,h}^+ / d_{n,h}$.
	
	For an upper bound for $M_{n,h}^+ / d_{n,h}$, we first note that:
	\begin{align}
		M^+_{n,h} &= \frac{1}{n}\sum^n_{i=1}[(X_i-M_{n,1})_+]^h \leq \frac{1}{n}\sum^n_{i=1}|X_i - M_{n,1}|^h \nonumber \\
		& \leq 2^{h-1}\left( \frac{1}{n} \sum^n_{i=1}X_i^h +  M_{n,1}^h\right) = 2^{h-1} (M_{n,h}  + M_{n,1}^h). \label{eq:M_nh+_UB}
	\end{align}
	By Lemma \ref{lemma:UB},
	\begin{equation}\label{eq:M_{nh}_d_{nh}O_p}
		\frac{M_{n,h}}{d_{n,h}} = O_p(1).
	\end{equation}
	As
	\begin{equation}\label{eq:dhn1dnh}
		\frac{d_{n,1}^h}{d_{n,h}} \sim \frac{t_n^{h-h\alpha}l(t_n)^h}{t_n^{h-\alpha}l(t_n)} = \frac{l(t_n)^{h-1}}{t_n^{(h-1)\alpha}} = \frac{l(t_n)^{h-1}}{[n l(t_n) c_n]^{h-1}} = \frac{1}{(nc_n)^{h-1}},
	\end{equation}
	we have
	\begin{equation}\label{eq:M_{nh}_d_{nh}O_p2}
		\frac{M^h_{n,1}}{d_{n,h}} = \left( \frac{M_{n,1}}{d_{n,1}}\right)^h \frac{d_{n,1}^h}{d_{n,h}} = O_p\left(\frac{1}{(nc_n)^{h-1}}\right).
	\end{equation}
	Thus, by \eqref{eq:M_nh+_UB}, \eqref{eq:M_{nh}_d_{nh}O_p}, and \eqref{eq:M_{nh}_d_{nh}O_p2}, we have
	\begin{equation}\label{eq:central_upper_moment1}
		\frac{M^+_{n,h}}{d_{n,h}} \leq 2^{h-1} \left( \frac{M_{n,h}}{d_{n,h}} + \frac{M^h_{n,1}}{d_{n,h}} \right) = 2^{h-1} \left( O_p(1) + O_p\left(\frac{1}{(nc_n)^{h-1}}\right) \right) =  O_p(1),
	\end{equation}
	where the last equality follows as $h > 1$ and $c_n \uparrow \infty$.
	
	We now establish a lower bound for $M_{n,h}^+ / d_{n,h}$. 
	As $|a-b| = (a-b)_+ + (b-a)_+$ for any $a,b \in \mathbb{R}$, we have
	\begin{equation}\label{eq:central_upper_moment2}
		M^+_{n,h} = \frac{1}{n}\sum^n_{i=1}|X_i - M_{n,1}|^h - \frac{1}{n}\sum^n_{i=1}(M_{n,1}- X_i)_+^h. 
	\end{equation}
	Using the inequality $(a+b)^p \leq 2^{p-1} (a^p + b^p)$ for $a,b \geq 0 $ and $p \geq 1$, we have
	\begin{align*}
		\frac{1}{n}\sum_{i=1}^n X_i^h &= \frac{1}{n}\sum_{i=1}^n |X_i - M_{n,1} + M_{n,1}|^h 
		\leq 2^{h-1} \left( \frac{1}{n}\sum_{i=1}^n|X_i - M_{n,1}|^h + M_{n,1}^h\right).
	\end{align*}
	Thus,
	\begin{equation}\label{eq:central_upper_moment3}
		\frac{1}{n}\sum^n_{i=1}|X_i - M_{n,1}|^h  \geq \frac{1}{n 2^{h-1}}\sum^n_{i=1}X_i^h - M_{n,1}^h = \frac{M_{n,h}}{2^{h-1}}  - M_{n,1}^h.
	\end{equation}
	As $(M_{n,1}-X_i)_+ \leq M_{n,1}$, we have
	\begin{equation}\label{eq:central_upper_moment4}
		\frac{1}{n}\sum^n_{i=1}(M_{n,1}- X_i)_+^h \leq M_{n,1}^h.
	\end{equation}
	In view of (\ref{eq:central_upper_moment2})--(\ref{eq:central_upper_moment4}),
	\begin{equation*}
		M^+_{n,h} \geq \frac{M_{n,h}}{2^{h-1}} - 2 M_{n,1}^h.
	\end{equation*}
	Then, by Lemma \ref{lemma:LB} and \eqref{eq:M_{nh}_d_{nh}O_p2},
	\begin{equation}\label{eq:central_upper_moment5}
		\frac{M^+_{n,h}}{d_{n,h}} \geq \frac{c_{h,\alpha, \delta'}(1+o_p(1))}{ 2^{h-1}c_n^{\delta'}} - O_p\left( \frac{1}{(nc_n)^{h-1}}\right),
	\end{equation}
	for some constant $\delta' > 2(h/\alpha - 1) > 0$ and $c_{h,\alpha,\delta'} > 0$.
	
	Finally, in view of (\ref{eq:central_upper_moment1}) and (\ref{eq:central_upper_moment5}), as $h > 1$, for any $\varepsilon \in (0, 1)$, with probability greater than $1 - \varepsilon$, for all sufficiently large $n$, there exists some $C > 0$ such that
	\begin{equation*}
		\frac{1}{2} \cdot \frac{c_{h, \alpha, \delta'}}{2^{h-1} c_n^{\delta'}} \leq \frac{M_{n,h}^+}{d_{n,h}} \leq  C.
	\end{equation*}
	Therefore, for all sufficiently large $n$,
	\begin{equation*}
		\left|\log \frac{M_{n,h}^+}{d_{n,h}} \right| \leq  |\log C| + |\log(c_{h, \alpha, \delta'}/2)|  + h \log 2+ \delta' \log c_n.
	\end{equation*}
	This completes the proof.
\end{proof}

\begin{proof}[Proof of Theorem \ref{thm:upper_central_moment}]
	\begin{enumerate}[(a)]
		\item 
		We first show that 
		\begin{equation*}
			\frac{ \frac{1}{n}\sum^n_{i=1}|X_i-M_{n,1}|^h}{M_{n,1}^h} - 1 = o_p(1).
		\end{equation*}
		Following the proof of Theorem B.1 of \cite{brown2021taylor}, we have
		\begin{align*}
			\left| \frac{ \frac{1}{n}\sum^n_{i=1}|X_i-M_{n,1}|^h}{M_{n,1}^h} - 1 \right|
			\leq \frac{M_{n,h}}{M^h_{n,1}}.
		\end{align*}
		By Lemmas \ref{lemma:UB} and \ref{lemma:LB}, and (\ref{eq:dhn1dnh}),
		\begin{equation*}
			\frac{M_{n,h}}{M^h_{n,1}} =  \frac{M_{n,h}}{d_{n,h}} \cdot \frac{d_{n,h}}{d_{n,1}^h} \cdot \left(\frac{d_{n,1}}{M_{n,1}}\right)^h = O_p(1) O_p(c_n^{h \delta'})\frac{d_{n,h}}{d_{n,1}^h} \sim O_p(c_n^{h \delta'})(nc_n)^{h-1} = o_p(1),
		\end{equation*}
		as $h < 1$. By Theorem \ref{thm:hth_central_lower} (ii),
		\begin{equation*}
			\frac{M^+_{n,h}}{M_{n,1}^h} = \frac{ \frac{1}{n}\sum^n_{i=1}|X_i-M_{n,1}|^h}{M_{n,1}^h} - \frac{M^-_{n,h}}{M_{n,1}^h} =
			(1+o_p(1)) - (1 +o_p(1)) = o_p(1).
		\end{equation*}

		\item With Lemma \ref{lemma:central_upper_moment}, the proof of the result is similar to that of Theorem \ref{thm:moment} and is therefore omitted.
	\end{enumerate}
\end{proof}

\subsection{Proofs for Section \ref{subsect:local_lower_upper}}
Let $M^*_{n,1} := \frac{1}{n}\sum^n_{i=1}X^*_i$, where $X^*_i$'s have the same joint distribution as $X_i$'s.  
Recall that $d_{n,1} = \mathbb{E}(X)\mathbbm{1}(X < t_n)$, $b_n = d_{n,1}/c_n^{2\delta'}$, and $\tilde{b}_n = d_{n,1}c_n$ for  $\delta' > 0$. 
To  prove Theorems \ref{thm:local_lower_central_moment} and \ref{thm:local_upper_central_moment}, 
we establish Lemmas \ref{lemma:indep_n_alpha_order}--\ref{lemma:N-n_n converge to 1}.

\begin{lemma}\label{lemma:indep_n_alpha_order}
	Let $X_1,\ldots,X_n,\ldots$ be a sequence of 
	nonnegative random variables with a common survival function $\overline{F}(x) = x^{-\alpha}l(x)$, where $l$ is slowly varying and $\alpha \in (0, 1)$.
	\begin{enumerate}[(i)]
		\item 
		If  Condition A($p$) holds for $p = 1$,  then
		\begin{equation*}
			\frac{\sum^n_{i=1}\mathbbm{1}(X_i > M_{n,1}^*) }{n^\alpha} = O_p\left(\frac{ l(b_n) c_n^{2\delta' \alpha + \alpha - 1}}{l(t_n) }\right),
		\end{equation*}
		where $\delta' > 2(1/\alpha -1)$.
		\item 
		If (\ref{eq:cov_ass3}) holds, then
		\begin{equation*}
			\left(   \frac{\sum^n_{i=1}\mathbbm{1}(X_i > M^*_{n,1}) }{n^\alpha} \right)^{-1} = O_p\left(\frac{l(t_n) c_n}{l(\tilde{b}_n)}\right).
		\end{equation*}
		
	\end{enumerate}
	
\end{lemma}

\begin{proof}[Proof of Lemma \ref{lemma:indep_n_alpha_order}]
	\begin{enumerate}[(i)]
		\item Let $C > 0$ and 
		\begin{equation*}
			g(n) := \frac{n \mathbb{P}(X_1 > b_n)}{n^\alpha}.
		\end{equation*}
		We have
		\begin{align}
			\mathbb{P} \left( \frac{\sum^n_{i=1}\mathbbm{1}(X_i > M_{n,1}^*) }{n^\alpha} > C g(n) \right)\nonumber 
			&= \int^{b_n}_0 \mathbb{P} \left( \frac{\sum^n_{i=1}\mathbbm{1}(X_i > t) }{n^\alpha} > C g(n) \right) dF_{M_{n,1}^*}(t) \nonumber \\
			& \quad + \int^\infty_{b_n} \mathbb{P} \left( \frac{\sum^n_{i=1}\mathbbm{1}(X_i > t) }{n^\alpha} > C g(n) \right) dF_{M_{n,1}^*}(t), \label{eq:n+bound1}
		\end{align}
		where $F_{M_{n,1}^*}(\cdot)$ is the cumulative distribution function of $M_{n,1}^*$. 
		For the first term on the right side of (\ref{eq:n+bound1}), since any probability is bounded by $1$, we have
		\begin{align}\label{eq:n+bound2}
			\int^{b_n}_0 \mathbb{P} \left( \frac{\sum^n_{i=1}\mathbbm{1}(X_i > t)}{n^\alpha} > C g(n) \right) dF_{M_{n,1}^*}(t) \leq \mathbb{P}\left( M_{n,1}^* \leq b_n\right) = 1 - \mathbb{P}\left( M_{n,1} > b_n\right),
		\end{align}
		as $M^*_{n,1} \stackrel{d}{=} M_{n,1}$. By the definition of $b_n$ and Lemma \ref{lemma:LB}, for some constant $C'$ independent of $n$,
		\begin{equation*}
			\frac{M_{n,1}}{b_n} = c_n^{2\delta'} \frac{M_{n,1}}{d_{n,1}} \geq C' c_n^{\delta'} (1 + o_p(1)).
		\end{equation*}
		As $c_n \rightarrow \infty$, we have 
		\begin{equation}\label{eq:n+bound3}
			\mathbb{P}(M_{n,1} > b_n) \rightarrow 1.     
		\end{equation}
		Since for any $t > b_n$,
		\begin{align*}
			\mathbb{P} \left( \frac{\sum^n_{i=1}\mathbbm{1}(X_i > t) }{n^\alpha} > C g(n) \right) 
			\leq  \mathbb{P} \left( \frac{\sum^n_{i=1}\mathbbm{1}(X_i > b_n) }{n^\alpha} > C g(n) \right).
		\end{align*}
		For the second term on the right  side of (\ref{eq:n+bound1}), we have
		\begin{align}\label{eq:n+bound4}
			\int^\infty_{b_n} \mathbb{P} \left( \frac{\sum^n_{i=1}\mathbbm{1}(X_i > t) }{n^\alpha} > C g(n) \right) dF_{M^*_{n,1}}(t)
			\leq  \mathbb{P} \left( \frac{\sum^n_{i=1}\mathbbm{1}(X_i > b_n) }{n^\alpha} > C g(n) \right).
		\end{align}
		By Markov's inequality and the definition of $g(n)$,
		\begin{align}\label{eq:n+bound5}
			\mathbb{P} \left( \frac{\sum^n_{i=1}\mathbbm{1}(X_i > b_n) }{n^\alpha} > C g(n) \right) &\leq \frac{n \mathbb{P}(X_1 > b_n) }{ C n^\alpha g(n)} =\frac{1}{C}.
		\end{align}
		In view of (\ref{eq:n+bound1})--(\ref{eq:n+bound5}), we have shown that
		\begin{equation*}
			\frac{\sum^n_{i=1}\mathbbm{1}(X_i > M^*_{n,1})}{n^\alpha} = O_p(g(n)).
		\end{equation*}
		Finally, by Lemma \ref{lemma:truncate_moment_asy}, for some constant $C'$,
		\begin{equation*}
			g(n)   = \frac{n^{1-\alpha}l(b_n) c_n^{2\delta' \alpha}}{d_n^\alpha}
			\sim C'\frac{n^{1-\alpha}l(b_n) c_n^{2\delta' \alpha} }{ (t_n^{1-\alpha}l(t_n))^\alpha} 
			= C'\frac{n^{1-\alpha}l(b_n) c_n^{2\delta' \alpha}}{ (nl(t_n)c_n)^{1-\alpha} l(t_n)^\alpha } 
			= C' \frac{ l(b_n) c_n^{2\delta' \alpha + \alpha - 1}}{l(t_n) }.
		\end{equation*}
		Thus, the result in part (i) of the lemma follows.

		\item Let $C > 0$ and 
		\begin{equation*}
			g_2(n) := \frac{l(\tilde{b}_n)}{l(t_n) c_n}.
		\end{equation*}
		We have
		\begin{align}
			&\mathbb{P}\left(  \left( \frac{\sum^n_{i=1}\mathbbm{1}(X_i > M^*_{n,1}) }{n^\alpha} \right)^{-1} > \frac{C}{g_2(n)} \right)\nonumber \\
			&= \int^{\tilde{b}_n}_0
			\mathbb{P}\left(  \frac{\sum^n_{i=1}\mathbbm{1}(X_i > t) }{n^\alpha} < \frac{g_2(n)}{C} \right)d F_{M^*_{n,1}}(t) \nonumber\\
			&\quad + \int^\infty_{\tilde{b}_n}
			\mathbb{P}\left(  \frac{\sum^n_{i=1}\mathbbm{1}(X_i > t) }{n^\alpha} < \frac{g_2(n)}{C} \right)d F_{M^*_{n,1}}(t). \label{eq:local_term_2}
		\end{align}
		For the second term on the right hand side of (\ref{eq:local_term_2}),
		\begin{equation*}
			\int^\infty_{\tilde{b}_n}
			\mathbb{P}\left(  \frac{\sum^n_{i=1}\mathbbm{1}(X_i > t) }{n^\alpha} < \frac{g_2(n)}{C} \right)d F_{M^*_{n,1}}(t) \leq 1 - F_{M^*_{n,1}}(\tilde{b}_n).
		\end{equation*}
		We have
		\begin{align*}
			1 -   F_{M^*_{n,1}}(\tilde{b}_n) = \mathbb{P}(M_{n,1}^*(t) > \tilde{b}_n) 
			= \mathbb{P} \left( \frac{M_{n,1}}{d_{n,1}c_n}  > 1\right) \rightarrow 0,
		\end{align*}
		where the last convergence follows as $\frac{M_{n,1}}{d_{n,1}} = O_p(1)$ (by Lemma \ref{lemma:UB}) so that $\frac{M_{n,1}}{d_{n,1}c_n} = o_p(1)$.
		For the first term on the right hand side of (\ref{eq:local_term_2}), when $0 \leq t \leq \tilde{b}_n$,
		\begin{equation*}
			\mathbb{P}\left(  \frac{\sum^n_{i=1}\mathbbm{1}(X_i > t) }{n^\alpha} < \frac{g_2(n)}{C} \right)
			\leq \mathbb{P}\left(  \frac{\sum^n_{i=1}\mathbbm{1}(X_i > \tilde{b}_n) }{n^\alpha} < \frac{g_2(n)}{C} \right)
		\end{equation*}
		and so
		\begin{equation*}
			\int^{\tilde{b}_n}_0
			\mathbb{P}\left(  \frac{\sum^n_{i=1}\mathbbm{1}(X_i > t) }{n^\alpha} < \frac{g_2(n)}{C} \right)d F_{M^*_{n,1}}(t) \leq 
			\mathbb{P}\left(  \frac{\sum^n_{i=1}\mathbbm{1}(X_i > \tilde{b}_n) }{n^\alpha} < \frac{g_2(n)}{C} \right),
		\end{equation*}
		as $F_{M^*_{n,1}}(\tilde{b}_n) \leq 1$.
		Write
		\begin{equation*}
			\frac{\sum^n_{i=1}\mathbbm{1}(X_i > \tilde{b}_n) }{n^\alpha} =  \frac{\sum^n_{i=1}(\mathbbm{1}(X_i > \tilde{b}_n) - \overline{F}(\tilde{b}_n)) }{n^\alpha} + \frac{n \overline{F}(\tilde{b}_n)}{n^\alpha}.
		\end{equation*}
		For some constant $C' > 0$ independent of $n$, by Lemma \ref{lemma:truncate_moment_asy},
		\begin{equation*}
			\frac{n \overline{F}(\tilde{b}_n)}{n^\alpha} 
			= \frac{n^{1-\alpha}l(\tilde{b}_n) }{d_{n,1}^\alpha c_n^\alpha} \sim C' \frac{n^{1-\alpha}l(\tilde{b}_n)}{t_n^{(1-\alpha)\alpha} l(t_n)^{\alpha} c_n^\alpha }
			= C' \frac{l(\tilde{b}_n)}{l(t_n) c_n} = C' g_2(n).
		\end{equation*}
		For any $\varepsilon > 0$, by Markov's inequality,
		\begin{align}
			& \mathbb{P} \left( \frac{\sum^n_{i=1}(\mathbbm{1}(X_i > \tilde{b}_n) - \overline{F}(\tilde{b}_n)) }{n^\alpha} > \varepsilon g_2(n) \right)\nonumber \\
			& \leq \frac{n \overline{F}(\tilde{b}_n)(1-\overline{F}(\tilde{b}_n))}{n^{2\alpha } \varepsilon^2 g^2_2(n)} + \frac{\sum_{i\neq j} \Cov(\mathbbm{1}(X_i > \tilde{b}_n), \mathbbm{1}(X_j > \tilde{b}_n)) }{n^{2\alpha} \varepsilon^2 g^2_2(n) } \label{eq:n+N_2}
		\end{align}
		The first term on the right side of  (\ref{eq:n+N_2}) obeys 
		\begin{equation*}
			\frac{n \overline{F}(\tilde{b}_n)(1-\overline{F}(\tilde{b}_n))}{n^{2\alpha } \varepsilon^2 g^2_2(n)} \sim \frac{1}{n^\alpha \varepsilon^2 g_2(n)} \rightarrow 0.
		\end{equation*}
		The second term on the right side of (\ref{eq:n+N_2})  goes to $0$ by  (\ref{eq:cov_ass3}). Thus we have shown that
		\begin{equation*}
			\frac{\sum^n_{i=1}\mathbbm{1}(X_i > \tilde{b}_n) }{n^\alpha} = (1+o_p(1)) g_2(n) + (C' + o_p(1))g_2(n).
		\end{equation*}
		Finally,
		\begin{equation*}
			\mathbb{P}\left(  \frac{\sum^n_{i=1}\mathbbm{1}(X_i > \tilde{b}_n) }{n^\alpha} <\frac{g_2(n)}{C} \right) 
			= \mathbb{P}\left(  (1+o_p(1)) + (C'+o_p(1))< \frac{1}{C} \right).
		\end{equation*}
		The last probability can be made arbitrarily small by picking $C$ large enough.
	\end{enumerate}
\end{proof}

\begin{lemma}\label{lemma:N+_n_alpha_order}
	Let $X_1,\ldots,X_n,\ldots$ be a sequence of 
	nonnegative random variables with a common survival function $\overline{F}(x) = x^{-\alpha}l(x)$, where $l$ is slowly varying and $\alpha \in (0, 1)$.
	\begin{enumerate}[(i)]
		\item 
		If Condition A($p$) holds for $p = 1$, then
		\begin{equation*}
			\frac{N^+_n}{n^\alpha} = O_p\left(\frac{ l(b_n) c_n^{2\delta' \alpha + \alpha - 1}}{l(t_n) }\right),
		\end{equation*}
		where $\delta' > 2(1/\alpha -1)$.
		
		\item     If (\ref{eq:cov_ass3}) holds, then 
		\begin{equation*}
			\left( \frac{N^+_n}{n^\alpha}\right)^{-1} = O_p\left(\frac{l(t_n) c_n}{l(\tilde{b}_n)}\right).
		\end{equation*}        
	\end{enumerate}
\end{lemma}

\begin{proof}[Proof of Lemma \ref{lemma:N+_n_alpha_order}]
	\begin{enumerate}[(i)]
		\item 
		With Lemma \ref{lemma:indep_n_alpha_order} (i), the proof follows using the argument in the proof of Lemma A.3 in \cite{brown2021taylor} and is therefore omitted.
		
		\item

		Consider $k > 1$ independent samples, each of size $n$, having the same joint distribution as $X_1,\ldots,X_n$:
		\begin{equation*}
			\{X^{(1)}_1,\ldots,X^{(1)}_n\}, \ldots, \{X^{(k)}_1,\ldots,X^{(k)}_n\}.
		\end{equation*}
		Denote $M^{(j)}_{n,1} := \frac{1}{n}\sum^n_{i=1}X^{(j)}_i$. Then
		\begin{align}
			N^+_n &= \sum^n_{i=1}\mathbbm{1}(X_i > M_{n,1}) 
			\geq \sum^n_{i=1}  \mathbbm{1} \left(X_i > M_{n,1} > \max_{1 \leq j \leq k} M^{(j)}_{n,1}\right)  \nonumber \\
			&= \sum^n_{i=1} \left[ \mathbbm{1}\left(X_i > \max_{1 \leq j \leq k} M^{(j)}_{n,1}\right) - \mathbbm{1}\left(\max_{1 \leq j \leq k}M^{(j)}_{n,1} < X_i \leq M_{n,1}\right) \right]. \label{eq:Nn+_decomp}
		\end{align}
		We  claim that
		\begin{equation}\label{eq:max_min1}
			\sum^n_{i=1}  \mathbbm{1}\left(X_i > \max_{1 \leq j \leq k} M^{(j)}_{n,1}\right) = \min_{1 \leq j \leq k} \sum^n_{i=1} \mathbbm{1}\left(X_i > M^{(j)}_{n,1}\right).
		\end{equation}
		Indeed, for any $j = 1,\ldots,k$,
		\begin{equation*}
			\sum^n_{i=1}  \mathbbm{1}\left(X_i > \max_{1 \leq j \leq k} M^{(j)}_{n,1}\right) \leq \sum^n_{i=1} \mathbbm{1} \left(X_i > M^{(j)}_{n,1}\right).
		\end{equation*}
		Thus
		\begin{equation}\label{eq:N+nalpha_claim1}
			\sum^n_{i=1}  \mathbbm{1}\left(X_i > \max_{1 \leq j \leq k} M^{(j)}_{n,1}\right) \leq \min_{1 \leq j \leq k} \sum^n_{i=1} \mathbbm{1}\left(X_i > M^{(j)}_{n,1}\right).
		\end{equation}
		On the other hand, for any $j = 1,\ldots,k$,
		\begin{equation*}
			\min_{1 \leq j \leq k} \sum^n_{i=1} \mathbbm{1}\left(X_i > M^{(j)}_{n,1} \right) \leq \sum^n_{i=1} \mathbbm{1}\left(X_i > M^{(j)}_{n,1}
			\right).
		\end{equation*}
		Since $\max_{1 \leq j \leq k} M^{(j)}_{n,1}$ is equal to one of $M^{(j)}_{n,1}$, we still have
		\begin{equation}\label{eq:N+nalpha_claim2}
			\min_{1 \leq j \leq k} \sum^n_{i=1} \mathbbm{1}\left(X_i > M^{(j)}_{n,1}\right) \leq \sum^n_{i=1} \mathbbm{1}\left(X_i > \max_{1 \leq j \leq k} M^{(j)}_{n,1}\right).
		\end{equation}
		In view of (\ref{eq:N+nalpha_claim1}) and (\ref{eq:N+nalpha_claim2}), we have (\ref{eq:max_min1}). From (\ref{eq:Nn+_decomp}) and (\ref{eq:max_min1}), we obtain
		\begin{align}
			N^+_n & \geq \min_{1 \leq j \leq k} \sum^n_{i=1} \mathbbm{1}\left(X_i > M^{(j)}_{n,1}\right) - \sum^n_{i=1}\mathbbm{1}\left(\max_{1 \leq j \leq k}M^{(j)}_{n,1} < X_i \leq M_{n,1}\right) \nonumber \\
			& \geq 
			\min_{1 \leq j \leq k} \sum^n_{i=1} \mathbbm{1}\left(X_i > M^{(j)}_{n,1}\right) - \sum^n_{i=1} \mathbbm{1}\left(\max_{1 \leq j \leq k}M^{(j)}_{n,1} < M_{n,1}\right) \nonumber \\
			& = \min_{1 \leq j \leq k} \sum^n_{i=1} \mathbbm{1}\left(X_i > M^{(j)}_{n,1}\right) - n \mathbbm{1}\left(\max_{1 \leq j \leq k}M^{(j)}_{n,1} < M_{n,1}\right), \label{eq:N^+_n_LB}
		\end{align}
		where the second inequality follows as 
		\begin{equation*}
			\mathbbm{1}\left(\max_{1 \leq j \leq k}M^{(j)}_{n,1} < X_i \leq M_{n,1}\right) \leq \mathbbm{1}\left(\max_{1 \leq j \leq k}M^{(j)}_{n,1} < M_{n,1}\right).
		\end{equation*}
		Denote
		\begin{equation*}
			g_2(n) := \frac{l(\tilde{b}_n)}{l(t_n) c_n}.
		\end{equation*}
		Let $\varepsilon > 0$. For any $C > 0$, by (\ref{eq:N^+_n_LB}),
		\begin{align}
			&
			\mathbb{P}\left( \left(\frac{N^+_n}{n^\alpha}\right)^{-1} > \frac{C}{g_2(n)} \right) = 
			\mathbb{P}\left( \frac{N^+_n}{n^\alpha} < \frac{g_2(n)}{C} \right)\nonumber \\
			&\leq \mathbb{P} \left( \frac{1}{n^\alpha} \min_{1 \leq j \leq k} \sum^n_{i=1}\mathbbm{1}\left(X_i > M_{n,i}^{(j)}\right) - n^{1-\alpha} \mathbbm{1}\left(\max_{1 \leq j \leq k}M^{(j)}_{n,1} < M_{n,1} \right) < \frac{g_2(n)}{C}\right)\nonumber   \\
			&= \mathbb{P} \left( \frac{1}{n^\alpha} \min_{1 \leq j \leq k}\sum^n_{i=1}\mathbbm{1}\left(X_i > M_{n,i}^{(j)}\right) - n^{1-\alpha} \mathbbm{1}\left(\max_{1 \leq j \leq k} M^{(j)}_{n,1} < M_{n,1} \right) < \frac{g_2(n)}{C}, \max_{1 \leq j \leq k} M^{(j)}_{n,1} \geq M_{n,1}  \right) \nonumber  \\
			& \quad + \mathbb{P} \left( \frac{1}{n^\alpha} \min_{1 \leq j \leq k} \sum^n_{i=1}\mathbbm{1}\left(X_i > M_{n,i}^{(j)}\right) - n^{1-\alpha} \mathbbm{1}\left(\max_{1 \leq j \leq k} M^{(j)}_{n,1} < M_{n,1} \right) < \frac{g_2(n)}{C}, \max_{1 \leq j \leq k} M^{(j)}_{n,1} < M_{n,1} \right)\nonumber  \\
			&\leq \mathbb{P} \left( \frac{1}{n^\alpha} \min_{1 \leq j \leq k}\sum^n_{i=1}\mathbbm{1}\left(X_i > M_{n,i}^{(j)}\right) - 0 < \frac{g_2(n)}{C}\right) \nonumber  \\
			& \quad + \mathbb{P} \left( \frac{1}{n^\alpha} \min_{1 \leq j \leq k} \sum^n_{i=1}\mathbbm{1}\left(X_i > M_{n,i}^{(j)}\right) - n^{1-\alpha} \mathbbm{1}\left(\max_{1 \leq j \leq k} M^{(j)}_{n,1} < M_{n,1} \right) < \frac{g_2(n)}{C}, \max_{1 \leq j \leq k} M^{(j)}_{n,1} < M_{n,1} \right)\nonumber  \\
			& \leq 
			\mathbb{P} \left( \frac{1}{n^\alpha} \min_{1 \leq j \leq k} \sum^n_{i=1}\mathbbm{1}\left(X_i > M_{n,i}^{(j)}\right)  < \frac{g_2(n)}{C}\right)  + \mathbb{P}\left( \max_{1 \leq j \leq k} M^{(j)}_{n,1} < M_{n,1}\right), \label{eq:N+_nalpha22}
		\end{align}
		where the second inequality holds because $\mathbbm{1}\left(\max_{1 \leq j \leq k} M^{(j)}_{n,1} < M_{n,1} \right) = 0$ as \\
		$\max_{1 \leq j \leq k} M^{(j)}_{n,1} \geq M_{n,1}$. 
		Since $M_{n,1}^{(1)},\ldots,M_{n,1}^{(k)},M_{n,1}$ are independent and have a common distribution, by symmetry,
		\begin{equation}\label{eq:N+_nalpha33}
			\mathbb{P}\left( \max_{1 \leq j \leq k} M^{(j)}_{n,1} < M_{n,1}\right)  = \frac{1}{k+1} < \frac{\varepsilon}{2},
		\end{equation}
		where the last inequality follows by choosing a large enough $k$.
		For any random variables $Y_1,\ldots,Y_k$, $\mathbb{P}( \min_{j=1,\ldots,k}Y_j < c) \leq \mathbb{P}(\cup_j \{Y_j < c\}) \leq \sum^k_{j=1}\mathbb{P}(Y_j < c)$. Thus
		\begin{align*}
			\mathbb{P} \left( \frac{1}{n^\alpha} \min_{1 \leq j \leq k} \sum^n_{i=1}\mathbbm{1}(X_i > M_{n,i}^{(j)})  < \frac{g_2(n)}{C}\right) &\leq k \mathbb{P} \left( \frac{1}{n^\alpha}  \sum^n_{i=1}\mathbbm{1}(X_i > M_{n,i}^{(1)})  < \frac{g_2(n)}{C}\right)\\
			&= k \mathbb{P} \left( \left( \frac{1}{n^\alpha}  \sum^n_{i=1}\mathbbm{1}(X_i > M_{n,i}^{(1)}) \right)^{-1}  < \frac{C}{g_2(n)}\right).
		\end{align*}
		By Lemma \ref{lemma:indep_n_alpha_order} (ii), there exists a $C> 0$ such that for all large $n$, the last probability is smaller than $\frac{\varepsilon}{2k}$. 
		This together with (\ref{eq:N+_nalpha22}) and (\ref{eq:N+_nalpha33}) show that there exists $C > 0$ such that for all large $n$,
		\begin{equation*}
			\mathbb{P}\left( \left(\frac{N^+_n}{n^\alpha}\right)^{-1} > \frac{C}{g_2(n)} \right) < \varepsilon,
		\end{equation*}
		and the proof is completed.
	\end{enumerate}
\end{proof}

\begin{lemma}\label{lemma:N-n_n converge to 1}
	Let $X_1,\ldots,X_n,\ldots$ be a sequence of 
	nonnegative random variables with a common survival function $\overline{F}(x) = x^{-\alpha}l(x)$, where $l$ is slowly varying and $\alpha \in (0, 1)$.         If  Condition A($p$) holds for $p = 1$, then 
	\begin{equation*}
		\frac{N^-_n}{n} \stackrel{\mathbb{P}}{\rightarrow} 1.
	\end{equation*}
\end{lemma}

\begin{proof}[Proof of Lemma \ref{lemma:N-n_n converge to 1}]
	By Lemma \ref{lemma:N+_n_alpha_order} (i), we have
	\begin{equation*}
		\frac{N^+_n}{n^\alpha} = O_p\left(\frac{ l(b_n) c_n^{2\delta' \alpha + \alpha - 1}}{l(t_n) }\right).
	\end{equation*}
	Thus, $\frac{N^+_n}{n} = o_p(1)$ as $\frac{ l(b_n) c_n^{2\delta' \alpha + \alpha - 1}}{l(t_n) n^{1-\alpha}} \rightarrow 0$.
	The proof is then completed by noting that
	\begin{equation*}
		1 - \frac{N^+_n}{n} = \frac{N^-_n}{n} \leq 1.
	\end{equation*}
\end{proof}

\begin{proof}[Proof of Theorem \ref{thm:local_lower_central_moment}]
	By Lemma \ref{lemma:N-n_n converge to 1}, $     \frac{N^-_n}{n} \stackrel{\mathbb{P}}{\rightarrow} 1$ and so $\log (     \frac{N^-_n}{n}) \stackrel{\mathbb{P}}{\rightarrow} 0$. Since \eqref{eq:ass_lower_semi} holds, we obtain (\ref{eq:hth_central_lower}) and thus
	\begin{equation*}
		\log \frac{M^{-*}_{n,h}}{M_{n,1}^h} = \log \frac{M^-_{n,h}}{M_{n,1}^h} - \log \frac{N^-_n}{n} = O_p(1) + o_p(1) = O_p(1),
	\end{equation*}
	which implies the required result as in the proof of Theorem \ref{thm:moment}.
\end{proof}

\begin{proof}[Proof of Theorem \ref{thm:local_upper_central_moment}]

	\begin{enumerate}[(a)]
		\item 
		Write
		\begin{equation}
			\log \frac{M^{+*}_{n,h}}{M_{n,1}^{\imath_+(h)}} = \log \frac{M^+_{n,h}}{d_{n,h}} + \log \frac{d_{n,h}}{d_{n,1}^{\imath_+(h)}} + \imath_+(h) \log \frac{d_{n,1}}{M_{n,1}} + \log \frac{n}{N^+_n}.
		\end{equation}
		By Lemmas \ref{lemma:central_upper_moment} and \ref{lemma:log_Mnp_dnp},
		\begin{equation*}
			\log \frac{M^+_{n,h}}{d_{n,h}} + \imath_+(h) \log \frac{d_{n,1}}{M_{n,1}} = O_p(\log c_n).
		\end{equation*}
		By Lemma \ref{lemma:truncate_moment_asy}, for some constant $C' > 0$,
		\begin{equation*}
			\frac{d_{n,h}}{d_{n,1}^{\imath_+(h)}} \sim C' 
			\frac{t_n^{h-\alpha}l(t_n)}{t_n^{\imath+(h)(1-\alpha)} l(t_n)^{\imath_+(h)}} = C' t_n^{\alpha^2-\alpha} l(t_n)^{1-\imath_+(h)}
		\end{equation*}
		and thus
		\begin{align*}
			\log \frac{d_{n,h}}{d_{n,1}^{\imath_+(h)}} &= o(1) + \log C' + (\alpha-1) \log t^\alpha_n + (1-\imath_+(h))\log l(t_n)\\
			&= o(1) + \log C' + (\alpha-1)  \log n + (2-\alpha-\imath_+(h))\log l(t_n) +(1-\alpha)\log c_n.
		\end{align*}
		By Lemma \ref{lemma:N+_n_alpha_order},
		\begin{align*}
			& \log \frac{d_{n,h}}{d_{n,1}^{\imath_+(h)}}+ \log \frac{n}{N^+_n} \\
			&= o(1) + \log C' + \log(n^{(\alpha-1)}) + \log \frac{n}{N^+_n}+ (2-\alpha-\imath_+(h))\log l(t_n) +(1-\alpha)\log c_n\\
			&= o(1) + \log C' + \log \frac{n^\alpha}{N^+_n}+ (2-\alpha-\imath_+(h))\log l(t_n) +(1-\alpha)\log c_n\\
			&= O(|\log l(t_n)|) + O_p(\log c_n) + O_p(|\log l(b_n)|) + O(|\log l(\tilde{b}_n)|).
		\end{align*}
		Therefore
		\begin{align*}
			\frac{\log M^{+*}_{n,h}}{\log M_{n,1}} - \imath_+(h) = O_p\left( \frac{L_n}{\log M_{n,1}}\right) = O_p\left( \frac{L_n}{\log n}\right).
		\end{align*}
		\item 
		In this case, we have
		\begin{equation*}
			\frac{N^+_n}{n^\alpha} = O_p\left( c_n^{2\delta' \alpha + \alpha - 1}\right) \text{ and } \left( \frac{N^+_n}{n^\alpha}\right)^{-1} = O_p\left(c_n\right).
		\end{equation*}
		The rest of the proof is similar to that of Theorem \ref{thm:moment} and  is therefore omitted.
		
	\end{enumerate}
\end{proof}

\section{Proofs for Section \ref{sec:heterogeneous}}

To prove Theorem \ref{thm:moment_hetero}, we  establish  Lemmas  \ref{lemma:heter_UV_mean}--\ref{lemma:LB_heter}.

\begin{lemma}\label{lemma:heter_UV_mean}
	Let $X_U $ be a random variable with survival function $\overline{F}_U(x) = x^{-\alpha}l(x)$, where $l$ is slowly varying and $\alpha \in (0, 1)$. Let $X_V$ be a random variable such that
	\begin{equation}\label{eq:heter_heavy_tail2}
		\lim_{x \rightarrow \infty} \frac{\mathbb{P}(X_V > x)}{\mathbb{P}(X_U > x)}= 0.
	\end{equation}
	Let $d_{n, p}^V := \mathbb{E}(X_U \mathbbm{1}(X_U < t_n))$ and $d_{n, p}^U := \mathbb{E}(X_V \mathbbm{1}(X_V < t_n))$.
	Then for all large enough $n$,
	\begin{equation*}
		d^V_{n,p} \leq C_0  d^U_{n,p},
	\end{equation*}
	where $C_0> 0$  is a constant that depends only on $\alpha, p$, but not on $n$.
\end{lemma}

\begin{proof}[Proof of Lemma \ref{lemma:heter_UV_mean}]
	In view of (\ref{eq:heter_heavy_tail2}), there exists $x_0 > 0$ such that for all $x \geq x_0$, $\overline{F}_V(x) \leq \overline{F}_U(x)$. Consider large $n$ such that $t_n \geq x_0$. Then, by Lemma \ref{lemma:expectation of y^p},
	\begin{align*}
		d^V_{n,p} &= \mathbb{E}(\tilde{X}^p_{1,V}) \\
		& = p \int^\infty_0 y^{p-1}\mathbb{P}(\tilde{X}_{1,V} > y) dy \\
		&= \frac{p}{F_V(t_n)} \left[ \int^{t_n}_0 y^{p-1} \overline{F}_V(y) - y^{p-1}\overline{F}_V(t_n) dy \right]\\
		& \leq \frac{p}{F_V(t_n)}\int^{t_n}_0 y^{p-1} \overline{F}_V(y)dy \\
		&= \frac{p}{F_V(t_n)}\left[ \int^{x_0}_0 y^{p-1} \overline{F}_V(y)dy  + 
		\int^{t_n}_{x_0} y^{p-1} \overline{F}_V(y)dy \right] \\
		& \leq \frac{px_0^{p}}{F_V(t_n)} + \frac{p}{F_V(t_n)}\int^{t_n}_{x_0} y^{p-1} \overline{F}_U(y)dy.
	\end{align*}
	By Proposition \ref{prop:Karamata},
	\begin{equation*}
		\int^{t_n}_{x_0} y^{p-1} \overline{F}_U(y)dy \sim \frac{t_n^{p-\alpha}l(t_n)}{p-\alpha}.
	\end{equation*}
	Recall that
	\begin{equation*}
		d^U_{n,p} \sim \frac{\alpha}{p-\alpha} t_n^{p-\alpha}l(t_n).
	\end{equation*}
	Thus, as $F_V(t_n) \rightarrow 1$ and $d^U_{n,p} \rightarrow \infty$,
	\begin{align*}
		\frac{d^V_{n,p}}{d^U_{n,p}} \leq o_p(1) + (1+o_p(1))\frac{p}{p-\alpha} \cdot \frac{p-\alpha}{\alpha} = o_p(1) + \frac{p}{\alpha}.
	\end{align*}
	Therefore, for all large $n$,
	\begin{equation*}
		d^V_{n,p} \leq \left( 1 + \frac{p}{\alpha}\right) d^U_{n,p}.
	\end{equation*}
	The claim in the lemma holds with  $C_0 = 1 + \frac{p}{\alpha}$.
\end{proof}

In the remainder of this section, we consider the setting described in Section \ref{sec:heterogeneous}. 
We define some additional notation. 
Let $\tilde{X}_{i,U} := X_{i,U}\mathbbm{1}(X_{i,U} < t_n)$ and $\tilde{X}_{i,V} := X_{i,V} \mathbbm{1}(X_{i,V}< t_n)$. 
{While $\breve{X}_{i,U}:= X_{i,U} \mathbbm{1}(X_{i,U} < v_n)$ and $\breve{X}_{i,V}:= X_{i,V} \mathbbm{1}(X_{i,V} < v_n)$ have similar definitions, they differ
	because the thresholds in the indicator functions are $v_n$ instead of $t_n$.}

Define $\tilde{M}_{n,p} := \frac{1}{n}( \sum^{u_n}_{i=1}\tilde{X}_{i,U}^p + \sum^{n-u_n}_{i=1}\tilde{X}_{i,V}^p)$ and $d_{n,p} = \mathbb{E}(\tilde{M}_{n,p})$.

The setting in Lemma \ref{lemma:UB_heterogenous}
is the same as that in Theorem \ref{thm:moment_hetero} except that 
Lemma \ref{lemma:UB_heterogenous} does not assume that \eqref{eq:cov_cond_hetero} holds.

\begin{lemma}\label{lemma:UB_heterogenous}
	For any positive integer $n>1$, assume that each of 
	$u_n < n$ random variables $X_{i,U}$, $i=1,\ldots,u_n$, follows the same  survival function $\overline{F}_U(x) = x^{-\alpha}l(x)$, where $l$ is slowly varying for some $\alpha \in (0, 1)$; 
	and all of $n -u_n$ random variables $X_{i,V}$, $i=1,\ldots,n-u_n$,  have the 
	same distribution such  that (\ref{eq:heter_heavy_tail}) holds. For $p > \alpha$, we have
	\begin{equation*}
		\frac{M_{n,p}}{d_{n,p}} = O_p(1).
	\end{equation*}
\end{lemma}

\begin{proof}[Proof of Lemma \ref{lemma:UB_heterogenous}]
	The proof is similar to that of Lemma \ref{lemma:UB}.  Let $\varepsilon > 0$. We shall show that $\mathbb{P} \left( \frac{M_{n,p}}{d_{n,p}}> C \right) < \varepsilon$ for some $C$ and all large $n$.
	We have
	\begin{align*}
		\mathbb{P}(\tilde{M}_{n,p} \neq M_{n,p})
		&\leq \mathbb{P} \left( ( \cup^{u_n}_{i=1} \{X_{i,U} > t_n\}) \cup  (\cup^{n-u_n}_{i=1}\{X_{i,V} > t_n\}) \right) \\
		&\leq u_n\mathbb{P}(X_{1,U} > t_n) + (n-u_n) \mathbb{P}(X_{1,V} > t_n).
	\end{align*}
	Thus
	\begin{align}
		\mathbb{P} \left( \frac{M_{n,p}}{d_{n,p}}> C \right) 
		& \leq  \mathbb{P} \left( \frac{\tilde{M}_{n,p}}{d_{n,p}} > C \right)  + \mathbb{P}(\tilde{M}_{n,p} \neq M_{n,p}) \nonumber \\
		&\leq  \mathbb{P} \left( \frac{\tilde{M}_{n,p}}{d_{n,p}}> C \right) + \bigg( u_n\mathbb{P}(X_{1,U} > t_n) + (n-u_n) \mathbb{P}(X_{1,V} > t_n) \bigg). \label{eq:heter_1}
	\end{align}
	For the first term in (\ref{eq:heter_1}), by Markov's inequality and Lemma \ref{lemma:heter_UV_mean}, we have for all large $n$,
	\begin{align*}
		\mathbb{P} \left( \frac{\tilde{M}_{n,p}}{d_{n,p}}> C \right) & \leq \frac{u_n d^U_{n,p} + (n-u_n)d^V_{n,p}}{C n d^U_{n,p}} \leq \frac{1 + C_0}{C}.
	\end{align*}
	Choose $C$ such that $\frac{1+C_0}{C} \leq \varepsilon / 2$. For the second term in (\ref{eq:heter_1}), 
	\begin{align*}
		u_n\mathbb{P}(X_{1,U} > t_n) + (n-u_n) \mathbb{P}(X_{1,V} > t_n) & \leq \frac{n l(t_n)}{t_n^\alpha} + n \mathbb{P}(X_{1,U} > t_n) \cdot  \frac{\mathbb{P}(X_{1,V} > t_n)}{\mathbb{P}(X_{1,U} > t_n)},
	\end{align*}
	which goes to $0$ by the choice of $t_n$ and (\ref{eq:heter_heavy_tail}). 
	Thus, for all large $n$, it is smaller than $\varepsilon / 2$. The proof is then completed.
\end{proof}

\begin{lemma}\label{lemma:LB_heter}
	Under the conditions in Theorem \ref{thm:moment_hetero}, we have
	\begin{equation*}
		\frac{M_{n,p}}{d_{n,p}} \geq \frac{c_{p,\alpha,\delta'}(1+o_p(1))}{ c_n^{\delta'}},
	\end{equation*}
	for any $\delta' \in (0, \frac{p-\alpha}{\alpha})$, where $c_{p,\alpha, \delta'} > 0$ is some constant depending on $p,\alpha,\delta'$.
\end{lemma}

\begin{proof}[Proof of Lemma \ref{lemma:LB_heter}]
	Let $\breve{M}_{n,p}:=\frac{1}{n} (\sum^{u_n}_{i=1}\breve{X}^p_{i,U} + \sum^{n-u_n}_{i=1} \breve{X}^p_{i,V})$. Define $\breve{d}^U_{n,p} := \mathbb{E}(\breve{X}^p_{1,U})$ and $\breve{d}^V_{n,p} := \mathbb{E}(\breve{X}^p_{1,V})$.     Write
	\begin{equation*}
		n M_{n, p} \geq n \breve{M}_{n, p} \geq (n\breve{M}_{n,p} - u_n \breve{d}^U_{n,p} - (n-u_v)\breve{d}^V_{n,p}) + u_n\breve{d}^U_{n,p}.
	\end{equation*}
	For any $p > \alpha$,
	\begin{equation}\label{eq:LB2_hetero}
		\frac{n \breve{d}^U_{n,p}}{v_n^p} \sim \frac{ \frac{\alpha}{p - \alpha} n v_n^{p-\alpha}l(t_n)}{v_n^p} \sim c_n \frac{\alpha}{p-\alpha}.
	\end{equation}
	For any $x > 0$, by Markov's inequality,
	\begin{align*}
		& \mathbb{P} \left(  \frac{|n\breve{M}_{n,p} - u_n \breve{d}^U_{n,p} - (n-u_v)\breve{d}^V_{n,p})|}{v_n^p c_n} > x \right) \\
		& \qquad \leq 
		\frac{u_n \breve{d}_{n,2p}^U  + (n-u_n) \breve{d}_{n,2p}^V +\sum_{i\neq j}\Cov(\breve{Y}^p_i, \breve{X}^Y_j)}{v^{2p}_n c_n^2 x^2} \rightarrow 0,
	\end{align*}
	where the last convergence follows from (\ref{eq:LB2_hetero}), Lemma \ref{lemma:heter_UV_mean}, the fact that $c_n \rightarrow \infty$ and Condition (\ref{eq:cov_cond_hetero}).
	Thus,
	\begin{equation*}
		\frac{n M_{n,p}}{v_n^p c_n} \geq o_p(1) + \frac{\alpha}{p-\alpha}.
	\end{equation*}
	The rest of the proof is similar to that of Lemma \ref{lemma:LB} and is therefore omitted.
\end{proof}

\begin{proof}[Proof of Theorem \ref{thm:moment_hetero}]
	The  proof  is similar to that of Theorem \ref{thm:moment} with the use of Lemmas \ref{lemma:UB_heterogenous} and \ref{lemma:LB_heter}. Therefore we omit the proof.
\end{proof}

\section{Proof for Section \ref{sec:corr}}
\begin{proof}[Proof of Theorem \ref{thm:correlated_case}]
	Following the proof of Theorem \ref{thm:moment}, it suffices to show that 
	\begin{equation*}
		\log \frac{M_{n,p}}{d_{n,p}} = O_p(\log c_n),
	\end{equation*}
	which will be satisfied if we can show 
	\begin{equation*}
		\frac{M_{n,p}}{d_{n,p}} = O_p(1) \text{ and } \frac{M_{n,p}}{d_{n,p}} \geq \frac{c_{p,\alpha,\delta'}(1+o_p(1)}{c_n^{\delta'}}
	\end{equation*}
	for some $\delta' > 0$ and constant $c_{p,\alpha,\delta'}$ as in Lemmas \ref{lemma:UB} and \ref{lemma:LB}. The same proof as in Lemma \ref{lemma:UB} can show the former result is true. For the latter result, we only need to modify the proof of (\ref{eq:lower_bound3}):
	\begin{equation*}
		\mathbb{P} \left(  \frac{|n\breve{M}_{n,p} - n \breve{d}_{n,p}|}{v_n^p c_n} > x \right)  \rightarrow 0.
	\end{equation*}
	The rest of the proof then follows as in the proof of Lemma \ref{lemma:LB}. 
	To this end, by the conditional Markov's inequality,
	\begin{align*}
		\mathbb{P} \left(  \frac{|n\breve{M}_{n,p} - n \breve{d}_{n,p}|}{v_n^p c_n} > x \right) 
		&= \mathbb{E} \left( \mathbb{P} \left(  \frac{|n\breve{M}_{n,p} - n \breve{d}_{n,p}|}{v_n^p c_n} > x \bigg| \mathcal{G} \right) \right)\\
		&\leq 
		\mathbb{E}\left( \frac{\mathbb{E}(|n\breve{M}_{n,p} - n \breve{d}_{n,p}|^2|\mathcal{G}) + \sum_{i\neq j}\Cov(\breve{X}^p_i, \breve{X}^p_j|\mathcal{G})}{v^{2p}_n c_n^2 x^2} \right)\\
		&\leq \frac{n \breve{d}_{n,2p}}{v^{2p}_n c_n^2 x^2} + o(1) \rightarrow 0,
	\end{align*}
	where the last convergence follows from (\ref{eq:LB2_new}).
	
\end{proof}

\section{Appendix for Section \ref{sec:network TL}}\label{app:network}
In this section, we provide the details of Example \ref{Example:network}.   Fix $p > \alpha$ and $k \in [n]$, by Conditions (1) and (2), we have
\begin{align*}
	\sum_{j: j \neq k}^n \Cov(\breve{X}_k^p,\breve{X}_j^p) 
	&=  \sum_{m=1}^{M} \sum_{j:(k,j) \in U_{n,m}} \Cov(\breve{X}_k^p,\breve{X}_j^p) 
	+ \sum_{m=M+1}^{n-1} \sum_{j:(k,j) \in U_{n,m}} \Cov(\breve{X}_k^p,\breve{X}_j^p)\\
	&\leq \sum_{m=1}^M K^m \Var(\breve{X}_k^p) \\
	&\leq C_1 \mathbb{E}(\breve{X}_k^{2p})
\end{align*}
where $C_1 := \sum_{m=1}^M K^m > 0$.
The first equality holds because we can separate the summation based on the distance from  $v_1$.
The first inequality follows as $\sum_{j :(k,j)\in U_{n,m}} 1 \leq K^m$ in view of Condition (2).  Hence, it follows that
\begin{align*}
	\sum_{i\neq j}\Cov(\breve{X}^p_i, \breve{X}^p_j) 
	\leq nC_1 \mathbb{E}(\breve{X}_1^{2p}).
\end{align*}
As a result,
\begin{align*}
	\frac{\sum_{i\neq j} \Cov(\breve{X}^p_i,\breve{X}^p_j)}{v^{2p}_n c_n^2}
	\leq \frac{C_1 n  \mathbb{E}(\breve{X}^{2p}_1)}{v^{2p}_n c_n^2} \sim \frac{C_1n \frac{\alpha}{2p-\alpha} v_n^{2p-\alpha} l(v_n)}{v_n^{2p}c_n^2} =  \frac{C_1\alpha}{2p-\alpha} \frac{n l(v_n)}{v_n^\alpha c_n^2} \rightarrow 0,
\end{align*}
where the asymptotic equivalence follows from Lemma \ref{lemma:truncate_moment_asy} and the last convergence follows from (\ref{eq:v_n_choice1}) and $c_n \rightarrow \infty$.
Therefore, Condition \eqref{eq:network_cov_assumption} is satisfied.






  \bibliography{references20241008.bib}

\end{document}